\documentclass[a4paper,9pt]{article}
\usepackage[latin1]{inputenc}
\usepackage[T1]{fontenc}
\usepackage{amsmath,amssymb,array,amsthm,graphics}
\usepackage{dsfont}\let\mathbb\mathds
\usepackage[english,frenchb]{babel}
\usepackage[all,dvips]{xy}
\newenvironment{preuve}{\vskip 2mm\noindent {\it Preuve: }}{\flushright $\Box$ \vskip 2mm}%\flushright
\newcommand{\calo}{\mathcal{O}}
\newcommand{\calv}{\mathcal{V}}
\newcommand{\calm}{\mathcal{M}}

\newcommand{\calf}{\mathcal{F}}
\newcommand{\calp}{\mathcal{P}}
\newcommand{\calu}{\mathcal{U}}
\newcommand{\cale}{\mathcal{E}}
\newcommand{\calg}{\mathcal{G}}
\newcommand{\cald}{\mathcal{D}}
\newcommand{\cals}{\mathcal{S}}
\newcommand{\caly}{\mathcal{Y}}
\newcommand{\calx}{\mathcal{X}}
\newcommand{\calh}{\mathcal{H}}
\newcommand{\call}{\mathcal{L}}
\newcommand{\calc}{\mathcal{C}}
\newcommand{\calk}{\mathcal{K}}
\newcommand{\cali}{\mathcal{I}}
\newcommand{\calw}{\mathcal{W}}
\newcommand{\calt}{\mathcal{T}}
\newcommand{\calr}{\mathcal{R}}
\newtheorem{theoreme}{Théorème} %ecrit les théorème avec le numéro de la section devant
\newtheorem{lemme}{Lemme} 

\newtheorem{corollaire}{Corollaire}
\newtheorem{definition}{Définition}
\newtheorem{proposition}{Proposition}

\theoremstyle{plain}

\newtheorem*{theoremese}{ Th\'eor\`eme} % Ne met pas de numéro de théorème
\newtheorem*{corollaire*}{Corollaire}

\newtheorem{propdef}{Proposition-Definition}

\newcommand{\crzmhs}{{\mathcal Cat}_{{\bf CRZ}-{\mathcal MHS}}}
\newcommand{\crmhs}{{\mathcal Cat}_{{\bf CR}-{\mathcal MHS}}}
\newcommand{\cmhs}{{\mathcal Cat}_{{\bf C}-{\mathcal MHS}}}
\newcommand{\trif}{{\mathcal Cat}_{Trif}}

\begin{document}
\thispagestyle{empty}
\begin{center}
{\large{\bf STRUCTURES DE HODGE MIXTES ET\\
FIBRES SUR LE PLAN PROJECTIF\\
COMPLEXE}} 
\end{center}
\begin{center}
Olivier PENACCHIO
\end{center}
${}$\\
${}$\\
\begin{center}
{\bf Abstract}
\end{center}
 
The purpose of this work is to geometrize the notion of mixed Hodge structure. Therefore, we associate equivariant vector bundles on the projective plane to trifiltered vector spaces. Making this Rees construction with filtrations arising from mixed Hodge structures gives an equivalence of categories between the category of mixed Hodge structures and a category of equivariant vector bundles on ${\bf P}^{2}_{\bf C}$ that are semistable for a semistability condition with respect to a divisor. Using this dictionary we can recover geometrically some properties of mixed Hodge structures like the abelian caracter of the category of mixed Hodge structures or the fact that higher extensions vanish.

Next, looking at the second Chern class of the Rees bundles, we define the ${\bf R}$-split level of the mixed Hodge structures that generalize the notion of ${\bf R}$-split mixed Hodge structure. This invariant can be expressed in terms of the Hodge numbers $h^{p,q}$ and some Hodge-type integers $s^{p,q}$ that are given by the relative position of the Hodge and conjugate filtrations. When the Hodge numbers are constant for a family of mixed Hodge structures, these numbers vary semicountinously. This allows to define stratifications of the bases of variations of mixed Hodge structures. \\ 
${}$\\
{\bf Keywords}: Structures de Hodge mixtes, Fibr\'es semistables, Fibr\'es \'equivariants.
${}$\\
{\bf MSC-class:} 14C30-32G20.

\tableofcontents
\newpage

\section{Introduction}
Le développement de la théorie de Hodge a permis, à travers les structures additionnelles qu'elle apporte, une compréhension plus profonde de nombreux invariants topologiques, groupes d'homotopie, groupes de cohomologie, dans le cadre de la géométrie algébrique complexe. Dans ce travail nous nous proposons de géométriser la notion de structure de Hodge mixte. Nous associons aux structures de Hodge mixtes, qui sont essentiellement la donnée d'un espace vectoriel muni de trois filtrations vérifiant certaines relations d'incidence, des objets géométriques que sont les fibrés vectoriels. Cette construction, "de Rees", est largement inspirée d'une part par la construction des structures twistorielles mixtes sur la droite projective complexe de C.Simpson \cite{sim2}, d'autre part par la construction de C.Sabbah \cite{sab}, toujours sur la droite projective complexe, de fibrés vectoriels pour étudier les structures de Frobenius sur les variétés. Dans \cite{sab}, on construit des fibrés sur ${\bf P}^1$ à partir de la filtration de Hodge et de la filtration par le poids d'une structure de Hodge. Les structures twistorielles mixtes de \cite{sim2}, fibrés vectoriels sur ${\bf P}^1$, sont construites, lorsqu'elles viennent de la théorie de Hodge à partir de la filtration de Hodge et de sa conjuguée ; les positions relatives des filtrations se traduisent en termes de semistabilité des fibrés. Les objets géométriques que nous associons aux espaces vectoriels trifiltrés sont des fibrés sur le plan projectif complexe ${\bf P}^2$ qui sont équivariants pour l'action de tores et ont certaines propriétés de semistabilité. Cette construction rassemble en quelque sorte celle de \cite{sab} et de \cite{sim2} puisque lorsque les objets proviennent de structures de Hodge les fibrés de Rees sont construits à partir de la filtration de Hodge, sa conjuguée et la filtration par le poids. Cette construction est similaire à celle de Klyachko qui a montré dans \cite{kly} que les fibrés vectoriels équivariants sur les variétés toriques sont équivalents à certains types d'ensembles de filtrations d'espaces vectoriels et de sa généralisation par \cite{per}. 

Pour que la correspondance soit pertinente, il ne suffit pas de décrire les objets associés aux structures de Hodge mixtes mais aussi les morphismes entre ces objets. La puissance de la théorie de Hodge telle qu'elle a été développée par P.Deligne dans \cite{del2} et \cite{del3} tient au fait que ses constructions sont fonctorielles et au fait que les morphismes entre structures de Hodge qui apparaissent naturellement ont une grande rigidité, on dit, et nous l'expliciterons plus bas, qu'ils sont strictement compatibles. Nous décrivons ici les morphismes équivariants entre fibrés qui correspondent à ces morphismes strictement compatibles d'espaces vectoriels filtrés.

La motivation, outre la curiosité naturelle, pour s'appliquer à géométriser ces structures et écrire un "dictionnaire
" théorie de Hodge-faisceaux équivariants est de pouvoir traduire les problèmes étudiés en théorie de Hodge, variations de structures de Hodge, limites de structures de Hodge, modules de structures de Hodge en termes de problèmes "classiques" de géométrie algébrique, familles de faisceaux équivariants, limites dans les espaces de modules ou les champs de faisceaux cohérents équivariants. Le problème de donner un sens à une limite de structures de Hodge mixtes par exemple se traduirait en la question suivante : le champs des objets équivariants associés est-il compact ? Malheureusement, comme nous le décrivons plus bas, bien que nous puissions \'ecrire un dictionnaire entre structures de Hodge mixtes et fibr\'es \'equivariants, nous n'arrivons pas ici réellement à définir d'objet géométrique équivalent à une variation de structures de Hodge mixtes.    

\begin{center}
{\bf *}
\end{center}

Ce texte se décompose en quatre parties et une annexe. La première et la deuxième partie sont vouées à l'étude ponctuelle et respectivement en famille des espaces vectoriels filtrés. La troisième et la quatrième sont les applications de la première et de la deuxième partie à la théorie de Hodge dans le sens où les espaces filtrés étudiés proviennent de structures de Hodge. On verra que la traduction des propriétés ponctuelles est aisée, de la partie 1 à partie 3 alors que l'anti-holomorphicité est un obstacle à décrire les familles de structures de Hodge à l'aide de familles algébriques ou holomorphes de fibrés équivariants.

\begin{center}
{\bf **}
\end{center}

Décrivons plus précisement le contenu des quatre parties.\\

${\bullet}$ La première aboutit au résultat clé de ce travail puisqu'il décrit une équivalence de catégories entre la catégorie des espaces vectoriels munis de trois filtrations et la catégorie des fibrés équivariants pour l'action d'un tore munis d'une certaine condition de semistabilité. 

Elle débute par de nombreuses sections expositoires des notions utilisées par la suite. On se place sur un corps $k$ algébriquement clos de caractéristique nulle. On rappelle différentes définitions sur les filtrations suivant \cite{del2}, la principale étant la suivante : un morphisme entre espaces vectoriels filtrés $f:(V,F^{\bullet })\rightarrow (V',F^{\bullet }{}')$ est strictement compatible aux filtrations s'il est non seulement compatible aux filtrations, i.e. pour tout entier $p$, $f(F^{p}) \subset F^{p}{}'$, mais plus, si tout élément de l'image par $f$ de $V$ qui est dans $F^{p}{}'$ provient d'un élément de $F^{p}{}$, i.e.  $$f(F^{p})=\text{Im}(f) \cap F^{p}{}'.$$

La notion cruciale définie ensuite est celle de filtrations opposées. Trois filtrations finies et décroissantes d'un même espace vectoriel $V$, $W^{\bullet }$, $F^{\bullet }$ et $G^{\bullet }$ sont dites opposées si $$ \forall p,q,n, \,\,\, \, , p+q+n\neq 0 \Rightarrow \, Gr^{p}_{F^{\bullet }}Gr^{q}_{G^{\bullet }}Gr^{n}_{W^{\bullet }}V =0.$$
On prendra garde, comme noté dans \cite{del2}, et nous le verrons dans les propriétés des fibrés associés à des triplets de filtrations opposées, que l'ordre des filtrations importe dans la définition. Par le lemme de Zassenhaus $F^{\bullet }$ et $G^{\bullet }$ jouent ici un rôle symétrique. On pourra lui préférer la définition équivalente,
$$ (W^{\bullet }, F^{\bullet }, G^{\bullet }) \text{ opposées } \Leftrightarrow \forall \,\,\, n,\,\, ( F^{\bullet }, G^{\bullet }) \text{ sont } -n\text{-opposées sur } Gr_{W^{\bullet}}^{n}V,$$
dans laquelle la dissymétrie est plus explicite.

On aborde ensuite les constuctions de Rees qui forment la base de la géométrisation des espaces vectoriels filtrés. La philosophie générale de ces construction est la suivante : à un vectoriel filtré on associe un fibré équivariant pour l'action d'un tore, l'action permet r\'eciproquement de retrouver la filtration sur une des fibres du fibré. Plus précisement, considérons un espace vectoriel muni d'une filtration décroissante et exhaustive $(V,F^{\bullet })$ on associe un faisceau sur la droite
${\bf A}^{1}=\text{Spec}\,k[u]$, sous-faisceau de $j^{*}(\calo_{{\bf G}_m}\otimes_{k} V)$ (où $j: {\bf G}_{m} \rightarrow { \bf A}^{1}$) engendré par les éléments de la forme $u^{-p}.w$ pour $w \in F^{p}$. La filtration fournit des "pôles" en $0$. Le faisceau obtenu est le faisceau de Rees $\xi_{{\bf A}^{1}}(V,F^{\bullet })$, il est localement libre et équivariant pour l'action du tore ${\bf G}_{m}$ par translation sur la base. Réciproquement, si l'on part d'un fibré équivariant sur la droite affine, on peut reconstruire une filtration sur l'espace vectoriel fibre en $1$ du fibré en regardant l'ordre des p\^oles des sections équivariantes. La construction de Rees et sa construction inverse donnent une correspondance (cf.\cite{sim2}),

$$
\xymatrix{
\{ \text{Espaces vectoriels filtrés} \} \ar@<2pt>[r]^-{\Phi_{R}} &  \{  \text{Fibrés }{\bf G}_{m}\text{ équivariants sur }{\bf A}^{1}  \} \ar@<2pt>[l]^-{\Phi_{I}}
}
.$$

On peut poursuivre cette correspondance aux espaces vectoriels munis de deux filtrations. A gauche les morphismes sont les morphismes d'espaces vectoriels filtrés strictement compatibles aux filtrations, à droite les morphismes sont les morphismes équivariants de fibrés vectoriels dont le conoyau est singulier au plus en codimension 2, le lieu de singularité étant contenu dans l'origine du plan affine,
$$
\xymatrix{
\{ \text{Espaces vectoriels bifiltrés} \} \ar@<2pt>[r]^-{\Phi_{R}} &  \{  \text{Fibrés }{\bf G}_{m}^{2}\text{ équivariants sur }{\bf A}^{2}  \} \ar@<2pt>[l]^-{\Phi_{I}}
}
.$$

Si l'on part d'un espace vectoriel trifiltré  $(V,F_{0}^{\bullet },F_{1}^{\bullet },F_{2}^{\bullet })$, on peut alors définir un fibré sur le plan projectif ${\xi}_{{\bf P}^2}(V,F_{0}^{\bullet },F_{1}^{\bullet },F_{2}^{\bullet })$ en recollant les fibrés sur les affines standards obtenus à partir des trois couples de filtrations.\ C'est le fibré de Rees associé à l'espace vectoriel trifiltré $(V,F_{0}^{\bullet },F_{1}^{\bullet },F_{2}^{\bullet })$. Les actions des ${\bf G}_{m}^{2}$ se recollent pour former une action du tore quotienté par la diagonale, ${\bf T}={\bf G}_{m}^{3}/\Delta({\bf G}_{m}^{3})$. On peut alors exhiber une équivalence de catégories entre la catégorie des espaces vectoriels trifiltrés munie des morphismes strictement compatibles entre les filtrations et la catégorie des fibrés ${\bf T}$-équivariants sur le plan projectif munie des morphismes équivariants dont le conoyau est singulier au plus en codimension 2, donc sans torsion. Le conoyau d'un morphisme dans cette catégorie est défini comme étant le faisceau réflexif associé au conoyau faisceautique, et est donc localement libre puisque l'on est sur une surface. L'équivalence de catégories découle donc du fait que le foncteur de Rees est exact de la catégorie des espaces vectoriels trifiltrés munie des morphismes strictement compatibles vers la catégorie des faisceaux cohérents équivariants sur le plan projectif privé des trois points origines des cartes affines standard. Le conoyau ainsi défini dans la catégorie des fibrés sur le plan projectif complexe co\"{\i}ncide donc en dehors d'un ensemble de codimension $2$ avec le fibré de Rees associé au conoyau du morphisme de modules de Rees, ainsi les fibrés associés sont isomorphes car les faisceaux réflexifs sont entièrement déterminés par leur restriction
au complémentaire d'un ensemble de codimension $2$.
 
En vue de l'application à la théorie de Hodge, il faut traduire la condition filtrations opposées en terme de fibrés de Rees. Par un calcul de classe de Chern, on sait que les fibrés associés aux structures de Hodge doivent \^etre de degré $0$. On remarquera, en guise de vérification, que toutes les formules de calcul des invariants des fibrés associés à des filtrations opposées $(V,F_{0}^{\bullet },F_{1}^{\bullet },F_{2}^{\bullet })$ sont symétriques en $F_{1}^{\bullet }$ et $F_{2}^{\bullet }$. Supposons que  $(V,F_{0}^{\bullet },F_{1}^{\bullet },F_{2}^{\bullet })$ forment des filtrations opposées. Alors sur chaque espace gradué associé à la filtration $F_{0}^{\bullet }$, celle qui joue un rôle dissymétrique, les deux autres filtrations sont opposées. Ceci se traduit par le fait que le fibré ${\bf G}_{m}$-équivariant obtenu par restriction du fibré de Rees au diviseur associé à la filtration $F_{0}^{\bullet }$, ${\bf P}^{1}_{0}$, soit semistable de pente 0. Réciproquement, si tel est le cas pour un fibré de Rees, alors les filtrations $(F_{1}^{\bullet },F_{2}^{\bullet })$ définissent bien des filtration $n$-opposées sur le gradué à l'étape $n$ associé à $F_{0}^{\bullet }$.
 
Le théorème principal est donc,
\begin{theoremese}{\bf 2} La construction du fibré de Rees sur ${\bf P}^2$ associé à un espace vectoriel trifiltré établit des équivalences de catégories entre :\\
\hspace*{0.1cm}$\bullet$ La catégorie des espaces vectoriels de dimension finie munis de trois filtrations décroissantes et exhaustives et des morphismes strictement compatibles aux filtrations $\calc_{3filtr}$ est équivalente à la catégorie des fibrés vectoriels sur le plan projectif ${\bf P}^2$ munis de l'action de ${\bf T}$ et des morphismes équivariants de fibrés dont les singularités du conoyau sont en codimension $2$, supportées aux points $(1:0:0)$, $(0:1:0)$ et $(0:0:1)$, ${\calf ib}({\bf P}^{2}/{\bf T})$ :
$$
\xymatrix{
\{ \calc_{3filtr} \} \ar@<2pt>[r]^-{\Phi_{R}} &  \{  {\calf ib}({\bf P}^{2}/{\bf T}) \} \ar@<2pt>[l]^-{\Phi_{I}}
}
.$$
\hspace*{0.1cm}$\bullet$ La catégorie des espaces vectoriels de dimension finie munis de trois filtrations décroissantes, exhaustives et opposées et des morphismes strictement compatibles aux filtrations $\calc_{3filtr,opp}$ est équivalente à la catégorie des fibrés vectoriels ${\bf P}^{1}_{0}$-semistables de pente $0$ sur le plan projectif ${\bf P}^2$ munis de l'action de ${\bf T}$ et des morphismes équivariants de fibrés dont les singularités du conoyau sont en codimension $2$, supportées au point $(1:0:0)$, ${\calf ib}_{{\bf P}^{1}_{0}-semistable,\mu=0}({\bf P}^{2}/{\bf T})$ :
$$
\xymatrix{
\{ \calc_{3filtr,opp} \} \ar@<2pt>[r]^-{\Phi_{R}} &  \{  {\calf ib}_{{\bf P}^{1}_{0}-semistable,\mu=0}({\bf P}^{2}/{\bf T}) \} \ar@<2pt>[l]^-{\Phi_{I}}
}
.$$
\end{theoremese}  
Pour l'application à la théorie de Hodge, il faut ajouter de la structure réelle ; on veut trouver la condition sur les fibrés de Rees qui traduise le fait que les filtrations $(F_{1}^{\bullet },F_{2}^{\bullet })$ soient conjuguées par rapport à une structure réelle sous-jacente. Le corps de base est ici celui des complexes, ${\bf C}$. Cette condition est une propriété d'équivariance des fibrés de Rees pour l'action d'une involution anti-holomorphe $\tau$ qui échange les droites associées aux filtrations $F_{1}^{\bullet }$ et $F_{2}^{\bullet }$. Un fibré est ${\bf T}^\tau$-équivariant s'il est à la fois $\tau $ et ${\bf T}$-équivariant. On montre alors que,
\begin{theoremese}{\bf 3}
 La catégorie des espaces vectoriels complexes de dimensions finies munis d'une structure réelle sous-jacente, de trois filtrations décroissantes, exhaustives et opposées telles que deux d'entre elles soient conjuguées, l'autre définie sur la structure réelle, et des morphismes strictement compatibles aux filtrations $\calc_{3filtr,opp, {\bf R}}$ est équivalente à la catégorie des ${\bf T}^{\tau}$-fibrés vectoriels ${\bf P}^{1}_{0}$-semistables de pente $0$ sur le plan projectif ${\bf P}^2$ et des morphismes de ${\bf T}^{\tau}$-fibrés dont les singularités du conoyau sont en codimension au plus $2$, supportées au point $(1:0:0)$, ${\calf ib}_{{\bf P}^{1}_{0}-semistable,\mu=0}({\bf P}^{2}/{\bf T}^{\tau})$ :
$$
\xymatrix{
\{ \calc_{3filtr,opp,{\bf R}} \} \ar@<2pt>[r]^-{\Phi_{R}} &  \{  {\calf ib}_{{\bf P}^{1}_{0}-semistable,\mu=0}({\bf P}^{2}/{\bf T}^{\tau}) \} \ar@<2pt>[l]^-{\Phi_{I}}
}
.$$
\end{theoremese}

Ces descriptions et l'étude faite dans le théorème 1 permettent de montrer, de façon géométrique, que les catégories $\calc_{3filtr,opp} $ et $\calc_{3filtr,opp,{\bf R}} $ sont abéliennes. On introduit alors la catégorie $\calc_{3filtr,opp}^{-}$ dont les objets sont ceux de $\calc_{3filtr,opp} $ mais dans laquelle les morphismes sont supposés compatibles, pas forcément strictement compatibles. On retrouve, comme conséquence des théorèmes précédents, le résultat de Deligne (\cite{del2}, th 1.2.10.) qui dit que tout morphisme de $\calc_{3filtr,opp}^{-} $ est automatiquement strictement compatible aux filtrations et donc que cette catégorie se confond avec $\calc_{3filtr,opp} $.

Il est très important de noter que les catégories abéliennes ${\calf ib}_{{\bf P}^{1}_{0}-semistable,\mu=0}({\bf P}^{2}/{\bf T})$ et ${\calf ib}_{{\bf P}^{1}_{0}-semistable,\mu=0}({\bf P}^{2}/{\bf T}^{\tau})$ qui sont des sous-catégories de la catégorie des faisceaux cohérents sur le plan projectif ne sont pas des sous-catégories abéliennes de cette catégorie. Le foncteur $\Phi_{R}$ n'est pas exact des modules de Rees gradués vers les faisceau cohérents. Les conoyaux dans ${\calf ib}_{{\bf P}^{1}_{0}-semistable,\mu=0}({\bf P}^{2}/{\bf T})$ et ${\calf ib}_{{\bf P}^{1}_{0}-semistable,\mu=0}({\bf P}^{2}/{\bf T}^{\tau})$ sont définis comme les faisceaux réflexifs (donc localement libres puisque l'on est sur une surface) canoniquement associés aux conoyaux faisceautiques. Comme les singularités de ces derniers sont au plus en codimension $2$, prendre le faisceau réflexif associé ne change pas le degré et permet
 d'énoncer des propriétés de ${\bf P}^{1}_{0}$-semistabilité et $\mu$-semistabilité.

La catégorie des fibrés ${\calf ib}_{\mu-semistable,\mu=0}({\bf P}^{2}/{\bf T})$ définie de la même façon mais avec une condition de semistabilité moins forte, la $\mu$-semistabilité, est abélienne (théorème $1$, \cite{pen}). Cette remarque nous permet de définir la catégorie équivalente des espaces vectoriels munis de trois filtrations "$\mu$-opposées", condition symétrique en les trois filtrations.
${}$\\
${}$\\

$\bullet$ Dans la deuxième partie, on veut poursuivre le dictionnaire établi dans la première partie pour les espaces vectoriels filtrés aux familles d'espaces vectoriels filtrés. Par famille d'espaces vectoriels filtrés paramétrisée par une base $S$, on entend fibré vectoriel sur $S$ filtré par des sous-espaces vectoriels stricts. Lorsqu'on fait la construction de Rees associée à un fibré vectoriel filtré par des sous-fibrés stricts et tel que les filtrations soient exhaustives et décroissantes sur une variété algébrique lisse $S$, $(\calv,\calf^{\bullet }_{0},\calf^{\bullet }_{1},\calf^{\bullet }_{2}) \rightarrow S$, on obtient un faisceau de $\calo_{S \times {\bf P}^2}$-modules cohérent $\xi(\calv,\calf^{\bullet }_{0},\calf^{\bullet }_{1},\calf^{\bullet }_{2})$ qui n'est pas localement libre en général, mais toujours réflexif. Le faisceau de Rees relatif $\xi(\calv,\calf^{\bullet }_{0},\calf^{\bullet }_{1},\calf^{\bullet }_{2})$ est ${\bf T}$-équivariant pour l'action de ${\bf T}$ sur le produit $S \times {\bf P}^2$ héritée de l'action par translation sur le deuxième facteur. On ne sait pas en toute généralité retrouver la famille d'espaces vectoriels filtrés à partir du faisceau de Rees équivariant. Par contre si la famille est plus rigide i.e. si les rangs d'intersection deux à deux des différents sous-fibrés ne sautent pas, alors le faisceaux de Rees est un fibré vectoriel équivariant dont la donnée est équivalente à celle de la famille. Cette hypothèse est notée $({\bf H})$. Elle est toujours vérifiée par strates, ainsi le dictionnaire relatif est un dictionnaire par strates ; la condition filtrations opposées est toujours traduite en terme de ${\bf P}^{1}_{0}$-semistabilité des restrictions du fibré de Rees aux $\{s \} \times {\bf P}^{2}$ pour $s \in S$, restrictions permises par l'hypothèse $({\bf H})$. Ainsi, il vient,
\begin{theoremese}${\bf 4}$ Soit $S$ une variété algébrique lisse de type fini sur un corps algébriquement clos de caractéristique nulle $k$. La construction des faisceaux de Rees relatifs sur $S \times {\bf P}^2$ établit les équivalences de catégories entre :\\
$\bullet$ La catégorie $ \calc_{3filtr}(S)$ des fibrés vectoriels sur $S$ trifiltrés par des sous-fibrés stricts tels que les filtrations soient exhaustives et décroissantes munie des morphismes de fibrés strictement compatibles aux filtrations et qui vérifient l'hypothèse $({\bf H})$ et la catégorie $\calf ib(S \times {\bf P}^{2}/{\bf T})$ des fibrés vectoriels ${\bf T}$-équivariants sur $S \times {\bf P}^2$ dont les restrictions ont les propriétés voulues munie des morphismes équivariants de fibrés dont le conoyau est sans torsion :
$$\xymatrix{
 \calc_{3filtr}(S)+({\bf H}) \ar@<2pt>[r]^-{\Phi_{\calr}} &   \calf ib(S \times {\bf P}^{2}/{\bf T}) \ar@<2pt>[l]^-{\Phi_{\cali}}
}.$$ 
$\bullet$ La catégorie $ \calc_{3filtr,opp}(S)$ des fibrés vectoriels sur $S$ trifiltrés par des sous-fibrés stricts tels que les filtrations soient opposées, exhaustives et décroissantes munie des morphismes de fibrés strictement compatibles aux filtrations et qui vérifient l'hypothèse $({\bf H})$ et la catégorie\\
$\calf ib_{{\bf P}^{1}_{0}-semistable/S,\mu=0/S}(S \times {\bf P}^{2}/{\bf T})$ des fibrés vectoriels ${\bf T}$-équivariants sur $S \times {\bf P}^2$ dont les restrictions ont les propriétés voulues et dont les restrictions à $\{s\} \times {\bf P}^2$ pour tout $s \in S$ sont ${\bf P}^{1}_{0}$-semistables de pente $0$ munie des morphismes équivariants de fibrés dont le conoyau est sans torsion : 
$$\xymatrix{
 \calc_{3filtr,opp}(S)+({\bf H}) \ar@<2pt>[r]^-{\Phi_{\calr}} &   \calf ib_{{\bf P}^{1}_{0}-semistable/S,\mu=0/S}(S \times {\bf P}^{2}/{\bf T}) \ar@<2pt>[l]^-{\Phi_{\cali}}
}.$$ 
\end{theoremese}

${}$\\
${}$\\

${\bullet }$ La troisième partie est la traduction de la première en termes de structures de Hodge mixtes. En oubliant la structure du réseau entier et la structure rationnelle, une structure de Hodge mixte est la donnée d'une filtration exhaustive et croissante $W_{\bullet}$ d'un espace vectoriel réel $H_{\bf R}$, la filtration par le poids, et d'une filtration exhaustive et décroissante $F^{\bullet }$ du complexifié $ H_{\bf C}=H_{\bf R}\otimes_{\bf R}{\bf C}$, la filtration de Hodge telle que, ${\overline F}^{\bullet }$ étant la filtration conjuguée à la filtration de Hodge par rapport à la structure réelle sous-jacente, le triplet de filtration $(W_{\bullet},F^{\bullet },{\overline F}^{\bullet })$ forme un système de filtrations opposées sur $H_{\bf C}$. On fait la construction de Rees associée aux données $(H_{\bf C},W^{\bullet},F^{\bullet },{\overline F}^{\bullet })$, où $W^{\bullet}$ est la filtration décroissante canoniquement associée à la filtration $W_{\bullet}$. Alors la construction de Rees et la construction inverse,
   
$$
\xymatrix{
\{  (H_{\bf C},W^{\bullet},F^{\bullet },{\overline F}^{\bullet }) \} \ar@<2pt>[r]^-{\Phi_{R}} &  \{ {\xi}_{{\bf P}^2}(H_{\bf C},W^{\bullet},F^{\bullet },{\overline F}^{\bullet }) \} \ar@<2pt>[l]^-{\Phi_{I}}
}
,$$
établissent l'équivalence de catégorie donnée par le théorème
\begin{theoremese}{\bf 6}
La construction du fibré de Rees sur ${\bf P}^2$ associé à un espace vectoriel trifiltré établit les équivalences de catégories entre :\\
${\bullet}$ La catégorie des structures de Hodge mixtes réelles $\crmhs$ et la catégorie des fibrés vectoriels ${\bf T}^{\tau}$-équivariants ${\bf P}^{1}_{0}$-semistables de pente $\mu=0$ munie des morphismes ${\bf T}^\tau$-équivariants dont les singularités des conoyaux sont supporté part la sous-variété de codimension $2$, $(1:0:0)$ : 
$$
\xymatrix{
\{ \crmhs \} \ar@<2pt>[r]^-{\Phi_{R}} &  \{  {\calf ib}_{{\bf P}^{1}_{0}-semistables,\mu=0}({\bf P}^{2}/{\bf T}^{\tau}) \} \ar@<2pt>[l]^-{\Phi_{I}}
}
.$$

\end{theoremese}

Ce théorème donne une démonstration géométrique des faits suivants :

\begin{corollaire*}$11$(Théorème $(1.3.16)$, \cite{del2})\\ 
$\bullet$ La catégorie des structures de Hodge mixtes réelles $\crmhs$ est abélienne.\\
$\bullet$ Les foncteurs ``oubli des filtrations'', $Gr_{W}$, $Gr_{F}$, $Gr_{\overline F}$ et $Gr_{W}Gr_{F}\cong Gr_{F}Gr_{W} \cong Gr_{W}Gr_{\overline F} \cong Gr_{\overline F}Gr_{W}$ de $\crmhs$ dans la catégorie des espaces vectoriels complexes sont exacts. 
\end{corollaire*}

Explicitons dans le cadre des structures de Hodge mixtes la remarque faite \`a la fin de la présentation de la partie $1$ sur les conoyaux dans la catégorie ${\calf ib}_{{\bf P}^{1}_{0}-semistables,\mu=0}({\bf P}^{2}/{\bf T})$ . Considérons la suite exacte de structures de Hodge mixtes
$$ 0 \rightarrow A \rightarrow H \rightarrow B \rightarrow 0.$$
Cette suite exacte se traduit par le théorème précédent en la suite exacte dans la catégorie\\ $ {\calf ib}_{{\bf P}^{1}_{0}-semistables,\mu=0}({\bf P}^{2}/{\bf T}^{\tau})$

$$  0 \rightarrow \xi_{A} \rightarrow \xi_{H} \rightarrow \xi_{B} \rightarrow 0.$$
Cette suite exacte se transforme en couple de suite exactes dans la catégorie des faisceaux cohérents équivariants
$$\left\lbrace \begin{array}{l}
0 \rightarrow \xi_{A} \rightarrow \xi_{H} \rightarrow \text{Coker} \rightarrow 0,\\    
{}\\
0 \rightarrow \text{Coker} \rightarrow \xi_{B} \rightarrow \calo_{T} \rightarrow 0,
\end{array}
\right.$$
où $T$ est un sous-schéma de codimension $2$ supporté au point $(0:0:1)$. Le foncteur de Rees $\Phi_{R}$ est "exact sur ${\bf P}^{2} \backslash (0:0:1)$". La longueur de $T$ va nous permettre de définir un invariant intéressant des structures de Hodge mixte qui detecte si elles sont ${\bf R}$-scindées ou non.

Avant de décrire cet invariant donnons quelques conséquences du théorème. Il permet de retrouver géométriquement, à l'aide des résultats connus sur les extensions de faisceaux cohérents sur les surfaces, le fait qu'il n'y a pas d'extension non triviales au delà des premiers groupes d'extensions pour les structures de Hodge mixtes. On donne en effet une autre démonstration du fait suivant (\cite{carhai}) : soient $A $ et $B$ deux structures de Hodge mixtes, alors, pour tout $p>1$, $\text{Ext}^{p}(B,A)=0.$

Par le biais des fibrés de Rees, on définit aussi un notion de structure de Hodge mixte indécomposable. Toute structure de Hodge mixte se décompose en sous-structures de Hodge mixtes indécomposables. Ceci permet de comprendre la complexité de la position relative des filtrations par blocs.

Le reste de la partie 3 est voué à la définition et à l'étude du niveau de ${\bf R}$-scindement d'une structure de Hodge mixte. Cet invariant de structures de Hodge mixtes, noté $\alpha$, se définit à partir des invariants discrets des fibrés de Rees, fibrés de degré $0$. Il est donné par la deuxième classe de Chern de ces fibrés. Cette classe de Chern est explicitée par la formule
$$\alpha (H)=\emph{c}_{2}({\xi}_{{\bf P}^2}(H))=\frac{1}{2}\sum_{p,q}(p+q)^{2}(h^{p,q}_{H}-s^{p,q}_{H}),$$
où les entiers $h^{p,q}_{H}$ sont définis par 
$$h^{p,q}_{H}=\text{dim}_{\bf C}Gr_{{\overline F}^{\bullet }}^{q}Gr_{F_{}^{\bullet }}^{p}Gr^{W^{\bullet }}_{-p-q}H_{\bf C}.$$
Ce sont les nombres de Hodge classiques. Les entiers $s^{p,q}_{H}$ sont donnés par $$s^{p,q}_{H}=\text{dim}_{\bf C}Gr_{{\overline F}_{}^{\bullet }}^{q}Gr_{{ F}_{}^{\bullet }}^{p}H_{\bf C}.$$ Ces derniers entiers sont importants pour décrire les sauts de dimension d'intersection entre les sous-espaces vectoriels que définissent la filtration de Hodge et sa conguguée et donc pour mesurer la complexité d'une structure de Hodge mixte en terme de position relative des filtrations. Le niveau de ${\bf R}$-scindement généralise la notion de structure de Hodge mixte ${\bf R}$-scindée dans le sens où il prend ses valeurs dans les entiers naturels et on a l'équivalence suivante $$\alpha(H)=0 \Leftrightarrow H \text{ est }{\bf R}\text{-scindée}.$$ 
Pour une structure de Hodge mixte ${\bf R}$-scindée, ce qui est le cas de toutes les structures de Hodge dont la longueur de la filtration par le poids est strictement inférieure à $2$, on a pour tous les entiers $p,q$, 
$h^{p,q}_{H}=s^{p,q}_{H}$. L'invariant se comporte bien par rapport aux opérations classiques sur les structures de Hodge (cf. Théorème 7.) et est sur-additif par extensions (cf. Théorème 8.).

Dans quelques exemples, pour des courbes complexes de genre $0$ et $1$ qui ne sont pas complètes et ont pour singularités des singularités nodales, on calcule le niveau de ${\bf R}$-scindement des structures de Hodge mixtes sur le premier groupe de cohomologie. L'exemple pertinent le plus simple est celui d'une courbe de genre 0 avec deux points enlevés et deux points recollés. L'espace de modules de telles courbes est isomorphe à ${\bf C}$. Alors l'invariant $\alpha$ est génériquement égal à $1$, la structures de Hodge mixte est génériquement non scindée, et vaut $0$ lorsque la classe d'isomorphisme de la courbe est donnée par un point sur une droite réelle de ${\bf C}$. Les sauts de $\alpha$ se font en codimension réelle.      
 
${}$\\
${}$\\

${\bullet}$ La quatrième partie est d'abord consacrée à des rappels sur les variations de structures de Hodge mixtes et sur la construction des espaces de modules associés. On stratifie ensuite la base d'une variation de structures de Hodge mixtes par l'invariant $\alpha$ qui est semi-continue inferieurement, le niveau de ${\bf R}$-scindement est génériquement le plus haut. On esquisse une idée d'application de l'existence de telles stratifications sans toutefois parvenir à trouver d'exemple. 

L'ambition d'écrire, au moins sous l'hypothèse $({\bf H})$, un dictionnaire
entre variations de structures de Hodge et familles de fibrés sur le plan projectif complexe est vaine. Cette partie ne peut pas être l'application de la partie 2 à la théorie de Hodge comme la 3 l'est de la 1. L'obstacle majeur est que la filtration conjuguée à la filtration de Hodge varie évidemment anti-holorphiquement pour une variation de structures de Hodge mixtes. On propose quelques idées pour pouvoir "déplier" les filtrations de Hodge et conjuguées en deux filtrations holomorphes et ensuite appliquer le travail de la partie 2.

${}$\\
${}$\\

${\bullet}$ Dans l'annexe on rappelle la notion de structure twistorielle mixte d'après \cite{sim2} ainsi que la construction de fibrés sur la sphère holomorphe associés à des paires de filtrations. Les entiers $s^{p,q}$ donnent une expression explicite des structures twistorielles mixtes associées à des structures de Hodge mixtes.

\begin{center}
${\bf ***}$
\end{center}

${}$\\
{\bf Remerciements:} Je tiens \`a remercier C.Sabbah pour m'avoir sugg\' er\'e cette d\'emarche, \`a remercier C.Sorger pour les commentaires sur la condition de semistabilit\'e et je tiens \`a exprimer toute ma gratitude et ma reconnaissance \`a C.Simpson pour avoir encadr\'e ce travail.

%-----------------------------------------
%-----------------------------------------
%---------------------------------------
%\addcontentsline{toc}

\newpage
\section{Construction d'un fibré sur ${{\bf P}^2}$ associé à un espace vectoriel trifiltré}

\subsection{Filtrations}
\subsubsection{Espaces vectoriels filtrés}

Nous rappelons ici, suivant \cite{del2}, des définitions et généralités sur les filtrations. Cette partie est destinée à fixer les notations et à énoncer des propriétés bien connues. On fixe un corps $k$ algébriquement clos et de caractéristique nulle.

Soit $V$ un espace vectoriel de dimension finie sur $k$.\\
\begin{definition}
Une filtration décroissante $F^{\bullet }(V)$ (resp. croissante $F_{\bullet }(V)$) de $V$ est une famille de sous-espaces vectoriels $(F^{p}(V))_{p \in {\bf Z}}$ (resp. $(F_{p}(V))_{p \in {\bf Z}}$) de $V$ tels que :
\begin{center}
$\forall\, m,n,  \,\,\,  n\leq m \Rightarrow F^{m}(V) \subset F^{n}(V)$

\end{center}
(resp. $\forall\, m,n \,\,\,   n\leq m \Rightarrow F_{n}(V) \subset F_{m}(V)$).
Un espace vectoriel $V$ muni d'une filtration $F^{\bullet }$ sera appelé espace vectoriel filtré et noté $(V,F^{\bullet }(V)) $.\\

\end{definition}
On peut associer une filtration décroissante $F^{\bullet }(V)$ à une filtration croissante $F_{\bullet }(V)$ en posant : $\forall\, m \, \,\,  F_{m}(V)= F^{-m}(V)$. Toutes les définitions et propriétés que nous donnerons par la suite pour les filtrations décroissantes se transportent ainsi aux filtrations croissantes.\\ 

\begin{definition}
Une filtration de $V$ est dite exhaustive s'il existe deux entiers $p $ et $q$ tels que : $$F^{p}(V)=\{ 0
 \} \text{ et } F^{q}(V)=V.$$ 
\end{definition}

\begin{definition}
Soient $(V,F^{\bullet }(V))$ et $(V',{F}^{\bullet }(V'))$ deux espaces vectoriels filtrés. Un morphisme d'objets filtrés $f$ est un morphisme de $V$ dans $V'$ tel que :

$$\forall m \in {\bf Z},\,\, f(F^{m}(V)) \subset F^{m}(V').$$
\end{definition}

Les espaces vectoriels de dimension finie sur $k$ forment une catégorie additive, la définition suivante a donc un sens :\\
\begin{definition}
Un morphisme d'objets filtrés $f:(V,F^{\bullet }(V)) \rightarrow (V',F^{\bullet}(V'))$ est strictement compatible aux filtrations si la flèche canonique de $\text{Co\"{\i}m}(f)$ dans $\text{Im}(f)$ est un isomorphisme d'objets filtrés c'est à dire, $V$ et $V'$ étant des $k$-modules comme espaces vectoriels sur $k$ que :

$$\forall m \in {\bf Z},\,\,\, f(F^{m}(V))=Im(f)\cap F^{m}(V').$$
\end{definition}

A tout espace vectoriel filtré $(V,F^{\bullet }(V))$ est associé un espace gradué qui est une famille de $k$-espaces vectoriels indexée par ${\bf Z}$, $(Gr^{n}_{F^{\bullet }}(V))_{n \in {\bf Z}}=(F^{n}(V)/F^{n+1}(V))_{n \in {\bf Z}}$.

Soit  $j: V' \hookrightarrow V$ un morphisme injectif d'espaces vectoriels sur $k$. Supposons que $V$ soit muni d'une filtration $F^{\bullet }(V)$. On peut alors définir sur $V'$ une filtration induite par celle de $V$ : 
\begin{definition}
La filtration $F^{\bullet}(V')$ sur $V'$ induite par  la filtration $F^{\bullet }(V)$ sur $V$ est l'unique filtration qur $V'$ telle que $j$ soit un morphisme strictement compatible aux filtrations c'est à dire que l'on ait :

 $$F^{n}(V')=j^{-1}(F^{n}(V))=V'\cap F^{n}(V).$$

La filtration quotient sur $V/V'$ est l'unique filtration telle que la projection canonique $p:V \rightarrow V/V'$ soit strictement compatible aux filtrations c'est à dire : 

$$F^{n}(V/V')=p(F^{n}(V))\cong(V'+F^{n}(V))/V'\cong F^{n}(V)/(V' \cap F^{n}(V)).$$

\end{definition}

On peut ainsi définir par la stricte compatibilité et la propriété universelle des sommes directes la filtration d'une somme directe finie d'espaces vectoriels filtrés :\\

Soit $(V_{i},F^{\bullet }(V_{i}))_{i \in [1,n]}$ une famille finie d'espaces vectoriels filtrés. On définit une filtration sur $\oplus_{i \in [1,n]}V_{i}$ par :
\begin{center}
$F^{k}(\oplus_{i \in [1,n]}V_{i})=\oplus_{i \in [1,n]} F^{k}(V_{i})$.\\
\end{center}

De façon analogue, pour le produit tensoriel d'un nombre fini d'espaces vectoriels filtrés :\\

Soit $(V_{i},F^{\bullet }(V_{i}))_{i \in [1,n]}$ une famille finie d'espaces vectoriels filtrés. On définit une filtration sur $\otimes_{i \in [1,n]}V_{i}$ par :
\begin{center}
$F^{k}(\otimes_{i \in [1,n]}V_{i})=\sum_{\sum_{i}k_{i}=k} \text{Im}(\otimes_{i \in [1,n]}F^{k_{i}}(V_{i})\rightarrow \otimes_{i \in [1,n]}V_{i})=\sum_{\sum_{i}k_{i}=k} \otimes_{i \in [1,n]}F^{k_{i}}(V_{i})$.
\end{center}
On a aussi, pour deux espaces vectoriels filtrés $(V,F^{\bullet}(V))$ et $(V',F^{\bullet}(V'))$ :
\begin{center}
$F^{k}\text{Hom}(V,V')=\{f: \, V \rightarrow V' \vert \forall n ,f(F^{n}(V) \subset F^{n+k}(V') \}$
\end{center}
d'où :

$$\text{Hom}((V,F^{\bullet}(V)),(V',F^{\bullet}(V')))=F^{0}(\text{Hom}(V,V')).$$

\begin{definition} La filtration décroissante triviale de $V$ notée $Triv^{\bullet }$ est la filtration exhaustive décroissante donnée par $Triv^{i}V=V$ pour $i\leq 0$ et $Triv^{1}V=\{ 0 \}$.\\

Pour $p \in {\bf Z}$, l'opération de décalage $Dec^{p}$ d'une filtration $F^{\bullet }(V)$ associe à cette filtration exhaustive et décroisante la filtration exhaustive et décroissante $Dec^{p}F^{\bullet }(V)$ donnée par $Dec^{p}F^{n}(V)=F^{n+p}(V)$ pour tout $n _in {\bf Z}$.\\
\end{definition}
Cette dernière filtration est notée $F[p]^{\bullet}$ dans \cite{del2}. La filtration décroissante triviale décalée de $p$, notée $Dec^{p}Triv^{\bullet}$ est celle donnée par $Dec^{p}Triv^{i}V=V$ pour $i \leq -p$ et $Dec^{p}Triv^{-p+1}V=\{ 0 \}$.\\ 
\\
{\bf Remarque :}
Pour ne pas alourdir les notations, la filtration triviale sur un espace vectoriel sera encore notée $Triv^{\bullet }$ sur les sous-espaces vectoriels et les espaces vectoriels quotients.
\\

\begin{definition}On notera $\calc_{nfiltr}$ la catégorie des espaces vectoriels de dimension finie sur $k$ munis de $n$ filtrations décroissantes et exhaustives et dont les morphismes sont les morphismes strictement compatibles aux filtrations. 
\end{definition}

\subsubsection{Filtrations opposées}
On introduit ici, toujours suivant \cite{del2} la notion de filtrations opposées. On dira qu'un espace vectoriel $V$ sur $k$ est bifiltré s'il est muni de deux filtrations $F^{\bullet }(V)$ et $G^{\bullet }(V)$, il sera noté $(V,F^{\bullet }(V),G^{\bullet }(V))$. $Gr_{F^{\bullet }}^{n}(V)$ est un quotient d'un sous-objet de $V$ et d'après la section précédente est donc muni d'une filtration induite par la filtration $G^{\bullet }(V)$ sur $V$. On obtient ainsi un objet bigradué $(Gr_{G^{\bullet }}^{n}Gr_{F^{\bullet }}^{m}(V))_{n,m \in {\bf Z}}$. Les objets $(Gr_{G^{\bullet }}^{n}Gr_{F^{\bullet }}^{m}(V))_{n,m \in {\bf Z}}$ et $(Gr_{F^{\bullet }}^{m}Gr_{G^{\bullet }}^{n}(V))_{n,m \in {\bf Z}}$ sont canoniquement isomorphes par le lemme de Zassenhaus.\\

\begin{definition}
Deux filtrations finies $F^{\bullet }(V)$ et $G^{\bullet }(V)$ sur $V$ sont dites $n$-opposées si
 $$ \emph{ pour } p+q\neq n, \, \,Gr_{G^{\bullet }}^{n}Gr_{F^{\bullet }}^{m}(V)=0.$$
\end{definition}
 
Si $(V^{p,q})_{(p,q) \in {\bf Z}\times {\bf Z}}$ est un objet bigradué tel que :\\
$V^{p,q}=0$ sauf pour un nombre fini de couples $(p,q)$ et,\\
$V^{p,q}=0$ pour $p+q \neq n$, alors on définit deux filtrations $n$-opposées de $V=\sum_{p,q}V_{p,q}$ en posant :
\begin{center}
$F^{p}(V)=\sum_{p' \geq p,q'}V^{p',q'}$\\
$G^{q}(V)=\sum_{q' \geq q,p'}V^{p',q'}$
\end{center} 
on obtient :
\begin{center}
$Gr_{G^{\bullet }}^{n}Gr_{F^{\bullet }}^{m}(V)=V^{m,n}$.\\
\end{center}
Par ces constructions, Deligne établit une équivalence de catégories quasi-inverse l'une de l'autre entre la catégorie des espaces vectoriels munis de deux filtrations $n$-opposées et la catégorie des objets bigradués définis ci-dessus. Cette équivalence est donnée par la proposition suivante :\\
\begin{proposition}\cite{del2}
\begin{itemize}\parindent=1cm
\item{(i)} Soient $F^{\bullet }(V)$ et $G^{\bullet }(V)$ deux filtrations finies sur $V$. Ces deux filtrations sont $n$-opposées si et seulement si\\
\begin{center}
$\forall \, p, q ,\,\,\,   p+q=n+1  \Rightarrow  F^{p}(V) \oplus G^{q}(V)\cong V$.
\end{center}
\item{(ii)} Si $F^{\bullet }(V)$ et $G^{\bullet }(V)$ sont $n$-opposées, si l'on pose $V^{p,q}=0$ pour $p+q \neq n$ et $V^{p,q}=F^{p}(V) \cap G^{q}(V)$ pour $p+q=n$, alors $V$ est somme directe des $V^{p,q}$ et $F^{\bullet }$ et $G^{\bullet }$ se déduisent de la bigraduation $(V^{p,q})_{(p,q)\in {\bf Z}\times {\bf Z}}$ de $V$ par le procédé décrit plus haut.\\
\par\end{itemize}
\end{proposition}
Nous aurons surtout à manipuler des triplets de filtrations opposées c'est à dire :\\

\begin{definition}
Trois filtrations finies $W^{\bullet }(V)$, $F^{\bullet }(V)$ et $G^{\bullet }(V)$ d'un espace vectoriel complexe $V$ sont dites opposées si $$ \forall p,q,n, \,\,\, \, , p+q+n\neq0 \Rightarrow \, Gr^{p}_{F^{\bullet }}Gr^{q}_{G^{\bullet }}Gr^{n}_{W^{\bullet }}(V) =0$$
\end{definition}

On peut vérifier aisément que si $F^{\bullet }(V)$ et $G^{\bullet }(V)$ sont des filtrations finies $n$-opposées sur l'espace vectoriel $V$ alors les filtrations finies $Dec^{-n}Triv^{\bullet }$, $F^{\bullet }(V)$ et $G^{\bullet }(V)$ sont opposées sur $V$.\\

Plus généralement, si $W^{\bullet }(V)$, $F^{\bullet }(V)$ et $G^{\bullet }(V)$ sont opposées sur $V$, $F^{\bullet }(V)$ et $G^{\bullet }(V)$ induisent sur $Gr^{n}_{W^{\bullet }}(V) $ des filtrations $-n$-opposées. En posant $V^{p,q}=Gr^{p}_{F^{\bullet }}Gr^{q}_{G^{\bullet }}Gr^{-p-q}_{W^{\bullet }}(V)$ on a d'après la proposition précédente une décomposition en somme directe de $Gr^{n}_{W^{\bullet }}(V) $, $Gr^{n}_{W^{\bullet }}(V)=\oplus_{p+q=-n}V^{p,q}$.\\

Il est très important de remarquer que dans la définition de filtrations opposées, les trois filtrations ne jouent pas un rôle symétrique. Par le lemme de Zassenhaus,  

$$ \forall p,q,n, \,\,\, \,\, Gr^{p}_{F^{\bullet }}Gr^{q}_{G^{\bullet }}Gr^{n}_{W^{\bullet }}(V) \cong Gr^{q}_{F^{\bullet }}Gr^{p}_{G^{\bullet }}Gr^{n}_{W^{\bullet }}(V), $$
ce qui signifie que les filtrations $F^{\bullet }$ et $G^{\bullet }$ jouent des rôles symétriques, mais la filtration $W^{\bullet }$ ne joue pas le même rôle. Par exemple, soit  $V=k^2=<e,f>_{k}$ équipé des filtrations décroissantes et exhaustives $F^{\bullet}_{1}$, $F^{\bullet}_{2}$ et $F^{\bullet}_{3}$ suivantes :\\ 
 \hspace*{.5cm}
$\bullet \ W^{-2}_{}=V, \ W^{-1}_{}=W^{0}=< e>, \  W^{1}_{}=\{0\} $,\\
                \hspace*{.5cm}
$\bullet \ F^{0}_{}=V,\  F^{1}_{}=<f+\lambda e>,\  F^{2}_{}=\{0\}$ où $\lambda \in k$,\\
\hspace*{.5cm}
$\bullet \  G^{0}_{}=V,\  F^{1}_{}=<f+\mu e>,\  F^{2}_{}=\{0\}$ où $\mu \in k$.\\
Alors $Gr^{1}_{F^{\bullet }}Gr^{1}_{G^{\bullet }}Gr^{-2}_{W^{\bullet }}(V)$ est de dimension $1$ et $Gr^{-2}_{W^{\bullet }}Gr^{1}_{F^{\bullet }}Gr^{1}_{G^{\bullet }}(V)$ de dimension $0$.\\

Pour bien noter que l'ordre des filtrations importe dans la définition on remarquera l'équivalence suivante 
$$ (W^{\bullet }(V), F^{\bullet }(V), G^{\bullet }(V)) \text{ opposées } \Leftrightarrow \forall \,\,\, n,\,\, ( F^{\bullet }(V), G^{\bullet }(V)) \text{ sont } n-\text{opposées sur } Gr_{W^{\bullet}}^{n}.$$

Citons un résultat sur la catégorie des espaces vectoriels munis de trois filtrations opposées et des morphismes stricts que nous montrerons par une autre voie, géométrique, que celle suivie dans \cite{del2}.\\

\begin{theoremese}(Th.(1.2.10), \cite{del2}) Soit $\calc_{3filtr,opp}$ la catégorie dont les objets sont les espaces vectoriels sur $k$ munis de trois filtrations finies opposées $W^{\bullet }(V)$, $F^{\bullet }(V)$ et $G^{\bullet }(V)$ et les morphismes sont les morphismes compatibles aux trois filtrations, alors :
\begin{itemize}\parindent=1cm
\item{(i)} $\calc_{3filtr,opp}$ est une catégorie abélienne.
\item{(ii)} Le noyau (resp. conoyau) d'une flèche $f: A \rightarrow B$ dans $\calc_{3filtr,opp}$ est le noyau (resp. conoyau) du morphisme d'espaces vectoriels muni des filtrations induites par celles de $A$ (resp. quotient de celles de $B$).
\item{(iii)} Tout morphisme $f: A \rightarrow B $ dans $\calc_{3filtr,opp}$ est strictement compatible aux filtrations $W^{\bullet }(V)$, $F^{\bullet }(V)$ et $G^{\bullet }(V)$. Le morphisme $Gr_{W}(f)$ est compatible aux bigraduations de $Gr_{W}(A)$ et $Gr_{W}(B)$. Les morphismes $Gr_{F}(f)$ et $Gr_{G}(f)$ sont strictement compatibles à la filtration induite par $W$.
\item{(iv)} Les foncteurs oubli des filtrations de $\calc_{3filtr,opp}$ dans la catégorie des espaces vectoriels, $Gr_{W}$, $Gr_{F}$, $Gr_{G}$ et $Gr_{W}Gr_{F} \cong Gr_{F}Gr_{W} \cong Gr_{F}Gr_{G}Gr_{W} \cong Gr_{G}Gr_{W} \cong Gr_{W}Gr_{G}$ sont exacts.
\par\end{itemize}
\end{theoremese}

\subsubsection{Filtrations et scindements }

Tous les espaces vectoriels sur $k$ munis de filtrations seront supposés de dimension finie. Toutes les filtrations seront supposées exhaustives et sauf avis contraire décroissantes. On note $\text{End}(V)$ l'anneau des endomorphisme de $V$.\\
 
\begin{definition}
Soit $(V,F^{\bullet })$ un espace vectoriel filtré. Un endomorphisme semi-simple $Y \in \emph{End}(V)$ scinde $F^{\bullet }$ si $\forall k \in {\bf Z} \,\,\, F^{k}=F^{k+1}\oplus E_{k}(Y)$ où $E_{k}(Y)$ est le sous-espace propre associé à la valeur propre $k$. $Y$ est une graduation de $(V,F^{\bullet })$.\\
\end{definition}
Toute filtration $F^{\bullet }$ sur un espace vectoriel $V$ peut \^etre scindée. On peut construire un scindement de la façon suivante : $F^{\bullet }$ est décroissante et exhaustive donc il existe un plus grand indice $p$ tel que $F^{p+1}=\{0\}$ et $F^{p} \neq \{0\}$. Soit $\{f^{p}_{i}\}_{i \in I_p}$ une famille d'éléments de $V$ formant une base de $F^p$. On complète cette famille dans $F^{p-1}$ par des vecteurs $\{f^{p-1}_{i}\}_{i \in I_{p-1}}$ pour obtenir une base de $F^{p-1}$. On itére l'opération pour obtenir une base de $V$ formée par la famille finie $\cup_{k \in [q,p]}  \{f^{k}_{i}\}_{i \in I_k}$ où $q$ est un entier tel que $F^{q}=V$. Cette famille est dite compatible avec la filtration. Soit $Y$ l'endomorphisme de $V$ définit par $Y(f^{k}_{i})=k.f^{k}_{i}$ pour $i \in I_{k}$. $Y$ est bien un scindement de $(V,F^{\bullet })$. 

Une tel scindement définit une graduation de $V$ i.e. une décomposition de $V$ en une somme directe compatible avec la filtration : $V=\oplus_{k=q}^{k=p}E_{k}(Y)$ et $F^{n}=\oplus_{k=n}^{k=p}E_{k}(Y)$. Nous noterons ${\mathcal Y}(F^{\bullet })$ l'ensemble des endomorphismes de $V$ qui scindent $F^{\bullet }$. \\

Soit $(V,F^{\bullet }_{1},F^{\bullet }_{2})$ un espace vectoriel muni de deux filtrations. Une bigraduation de $V$, $V=\oplus_{p,q}V^{p,q}$, est dite compatible aux filtrations $F^{\bullet }_{1}$ et $F^{\bullet }_{2}$ si pour tout $p \in {\bf Z}$ :

\begin{center}
 $F^{p}_{1}=\oplus_{a \geq p,q}V^{a,q}$ et $F^{p}_{2}=\oplus_{ p,q \geq b}V^{p,b}$. 
\end{center}
Il existe toujours une graduation de $V$ compatible à $F^{\bullet }_{i}$ associée à un endomorphisme $Y_{F^{\bullet }_{i}}$ qui scinde $F^{\bullet }_{i}$ pour $i \in \{1,2\}$. Le lemme suivant signifie que l'on peut trouver deux tels endomorphismes $Y_{F^{\bullet }_{1}} \in {\mathcal Y}(F^{\bullet }_{1})$ et $Y_{F^{\bullet }_{2}} \in {\mathcal Y}(F^{\bullet }_{2})$ qui commutent.\\

\begin{lemme}\label{bigraduation}
Pour tout espace vectoriel muni de deux filtrations $(V,F^{\bullet }_{1},F^{\bullet }_{2})$ il existe toujours une bigraduation $V=\oplus_{p,q}V^{p,q}$ compatible aux deux filtrations.\\  
\end{lemme}

On dit alors que l'on peut scinder les deux filtrations simultanément. La preuve consiste à exhiber une base de $V$ compatible avec les deux filtrations par complétions de bases successives. 

%\begin{preuve}
%Les deux filtrations étant exhaustives, on peut trouver deux entiers $p$ et $q$ tels que $F^{q }_%{1}=F^{q }_{2}=V$ et $F^{p}_{1}=F^{p }_{2}=\{0\}$. Soit $\{f^{k,l}_{i}\}_{i \in I_{p-1,p-1}}$ une famille de vecteurs de $V$ qui forme une base de $F^{p-1}_{1}\cap F^{p-1}_{2}$. On la complète en une base de $F^{p-1 }_{1}\cap F^{p-2 }_{2}$ par une famille $\{f^{k,l}_{i}\}_{i \in I_{p-1,p-2}}$. On peut procéder ainsi jusqu'à $F^{p-1}_{1}\cap F^{q}_{2}=F^{p-1}$. Complétons la base de $F^{p-1}_{1}$ à une base de $F^{p-1}_{1}\cup (F^{p-2}_{1} \cap F^{p-2})$ par une famille $\{f^{k,l}_{i}\}_{i \in I_{p-2,p-1}}$. Les indices $(k,l) \in [q,p-1]\times[q,p-1]$ peuvent ainsi \^etre "descendus" de $(p-1,p-1)$ à $(q,q)$ par l'algorithme suivant : tant que $l>q$ on complète la base de $F^{k+1}_{1}\cup (F^{k}_{1}\cap F^{l}_{2})$ en une base de $F^{k+1}_{1}\cup (F^{k}_{1}\cap F^{l-1}_{2})$ par une famille de vecteurs $\{f^{k,l}_{i}\}_{i \in I_{k,l-1}}$. Si $l=q$, on complète la base de $F^{k}_{1}$ en une base de $F^{k}_{1}\cup (F^{k-1}_{1}\cap F^{p-1}_{2})$ par une famille de vecteurs $\{f^{k,l}_{i}\}_{i \in I_{k-1,p-1}}$. Soient les endomorphismes $Y_{F^{\bullet }_{1}} \in {\mathcal Y}(F^{\bullet }_{1})$ et $Y_{F^{\bullet }_{2}} \in {\mathcal Y}(F^{\bullet }_{2})$ définis pour tout $i \in I_{k,l}$ par  $Y_{F^{\bullet }_{1}}(f^{k,l}_{i})=k.f^{k,l}_{i}$ et $Y_{F^{\bullet }_{1}}(f^{k,l}_{i})=l.f^{k,l}_{i}$. Ces deux endomorphismes commutent et scindent bien les deux filtrations. 
%\end{preuve}
 Une telle propriété n'est pas vraie en général lorsque $V$ est équipé de plus de deux filtrations. Soit $(V,F^{\bullet}_{1},...,F^{\bullet }_{n})$ un espace vectoriel muni de $n$ filtrations. Il n'est pas toujours possible de trouver $Y_{F^{\bullet }_{i}} \in {\mathcal Y}(F^{\bullet }_{i})$ pour tous $i \in [1,n]$ tels que pour tous $(i,j) \in [1,n]^{2}$, $Y_{F^{\bullet }_{i}}.Y_{F^{\bullet }_{j}}=Y_{F^{\bullet }_{j}}.Y_{F^{\bullet }_{i}}$. 

On peut traduire ceci en terme d'espace multigradué. Un espace vectoriel $V$ de dimension finie est dit multigradué d'ordre $n$ s'il est somme directe de sous-espaces vectoriels indexés par des $n$-uplets : $V=\oplus_{(i_{1},...,i_{n})}V^{i_{1},...,i_{n}}$ où la somme est finie.  Soit $V$ un espace vectoriel multigradué d'ordre $n$ muni de $n$ filtrations $(F^{\bullet}_{1},...,F^{\bullet }_{n})$. On dit que sa multigraduation est compatible avec les filtrations si pour tout $k \in [1,n]$ et $p \in {\bf Z}$ :
$$  F_{k}^{p}V= \oplus_{(i_{1},...,i_{n}),\,i_{k}\geq p}V^{i_{1},...,i_{n}}.$$
Dire que l'on ne peut pas scinder les $n$ filtrations simultanéement équivaut à dire que l'on ne peut pas trouver de multigraduation d'ordre $n$ compatible aux $n$-filtrations.

Ceci peut se produire dès que l'on a trois filtrations, comme le montre l'exemple suivant :\\ 
\\
\textit{\bf Exemple :}
Soit $V=k^2=<e,f>_{k}$ équipé des filtrations décroissantes et exhaustives $F^{\bullet}_{1}$, $F^{\bullet}_{2}$ et $F^{\bullet}_{3}$ suivantes :\\ 
 \hspace*{.5cm}
$\bullet \ F^{0}_{1}=V, \ F^{1}_{1}=<f+\kappa e>, \  F^{2}_{1}=\{0\} $ où $\kappa \in k$,\\
                \hspace*{.5cm}
$\bullet \ F^{0}_{2}=V,\  F^{1}_{2}=<f+\lambda e>,\  F^{2}_{2}=\{0\}$ où $\lambda \in k$,\\
\hspace*{.5cm}
$\bullet \  F^{0}_{3}=V,\  F^{1}_{3}=<f+\mu e>,\  F^{2}_{3}=\{0\}$ où $\mu \in k$.\\
Supposons que $\kappa, \lambda , \mu$ soient distincts deux à deux, sans quoi on est ramené à la situation de scinder au plus deux filtrations. Scinder $F^{\bullet }_{1}$ et $F^{\bullet }_{2}$ simultanément revient à prendre pour bigraduation de $V=V^{1,0} \oplus V^{0,1}$ où $V^{1,0}=<f+\kappa e>$ et $V^{0,1}=<f+\lambda e>$. $V^{1,0}$ est le sous-espace propre associé à la valeur $1$ de tout élément dans ${\mathcal Y}(F^{\bullet }_{1})$, ici, pour que les endomorphismes commutent, on doit choisir pour $Y_{F^{\bullet }_{1}}$ l'élément qui a pour sous-espace propre associé à la valeur $0$, $V^{0,1}$. De m\^eme, $V^{0,1}$ est le sous-espace propre associé à la valeur $1$ de tout élément dans ${\mathcal Y}(F^{\bullet }_{2})$, ici on doit choisir pour $Y_{F^{\bullet }_{2}}$ l'élément de $\text{End}(V)$ qui a pour sous-espace propre associé à la valeur $0$, $V^{1,0}$. Ainsi $Y_{F^{\bullet }_{1}}.Y_{F^{\bullet }_{2}}-Y_{F^{\bullet }_{2}}.Y_{F^{\bullet }_{1}}=0$. Comme $\mu \neq \kappa$ et $\mu \neq \lambda$, on ne peut diagonaliser aucun élément $Y_{F^{\bullet }_{3}} \in {\mathcal Y}(F^{\bullet }_{3})$ dans cette base.

On peut regarder l'espace des configuration des trois filtrations comme $k^3$ où chacunes des coordonnées $\kappa$, $\lambda $ et $\mu$ sont quelconques dans $k$. Alors le lieu où l'on peut simultanéement scinder les trois filtrations est la diagonale généralisée de $k^3$ formée des éléments\\
 $\{ (\kappa, \kappa, \mu ),(\kappa,\lambda ,\lambda ), (\kappa, \lambda , \kappa ) \vert (\kappa,\lambda, \mu) \in {k}^{3} \}$.

On peut de même généraliser au cas où les filtrations, ¨déterminées¨ par le terme d'indice $1$ ici, sont données par des droites quelconques de $k^2$. L'espace des configurations est alors la grassmanienne des droites dans $k^2$, $grass(1,2)={\bf P}^{1}_k$. Le lieu où l'on peut diagonaliser simultanéement les trois filtrations est alors la diagonale généralisée de ${\bf P}^{1}_{k} \times  {\bf P}^{1}_{k} \times {\bf P}^{1}_{k}$. 

Considérons la catégorie dont les objets sont les espaces vectoriels munis de trois filtrations décroissantes, exhaustives et scindées et les morphismes sont les morphismes strictement compatibles aux trois filtrations ${\mathcal C}_{3filtr,sci}$. On peut montrer que :\\

\begin{proposition}
La catégorie ${\mathcal C}_{3filtr,sci}$ est abélienne.\\

\end{proposition}

Nous allons par la suite étudier les propriétés de $\calc_{3filtr,opp}$, la catégorie des espaces vectoriels munis de trois filtrations opposées, puis de ${\mathcal C}_{3filtr,sci}$ et de leur intersection vue dans la catégorie des espaces vectoriels trifiltrés munie des morphismes compatibles aux filtrations.

\subsection{Préliminaires algébriques, modules de Rees}

Pour étudier les espaces vectoriels munis de plusieurs filtrations nous allons leurs associer des objets algébriques puis géométriques. L'objet de cette partie est de définir le module de Rees d'un espace vectoriel multifiltré et d'en dégager les principales propriétés. 

\subsubsection{Modules de Rees}
 Dans cette section, $k$ est toujours un corps algébriquement clos de caractéristique nulle. On rappelle que l'on note $\calc_{nfiltr}$ la catégorie des espaces vectoriels sur $k$ munis de $n$ filtrations exhaustives et décroissantes dont les morphismes sont les morphismes strictement compatibles aux filtrations. On définit le module de Rees de la façon suivante :\\

\begin{definition}
Soit $(V,F^{\bullet}_{1},F^{\bullet}_{2},...,F^{\bullet}_{n})$ un objet de $ \calc_{nfiltr}$. Le module de Rees d'ordre $n$ associé à $(V,F^{\bullet}_{1},F^{\bullet}_{2},...,F^{\bullet}_{n})$, noté $R^{n}(V,F^{\bullet}_{i})$, est le sous-${k}[u_{1},u_{2},...,u_{n}]$-module de\\
 ${k}[u_{1},u_{2},...,u_{n},u_{1}^{-1},u_{2}^{-1},...,u_{n}^{-1}]\otimes_{k} V$ suivant :
$$R^{n}(V,F^{\bullet}_{i})=\sum_{(p_{1},p_{2},...,p_{n}) \in {\bf Z}^{n}} \,u_{1}^{-p_{1}}u_{2}^{-p_{2}}...\,u_{n}^{-p_{n}}\,(F^{p_{1}}_{1} \cap F^{p_{2}}_{2} \cap ... \cap F^{p_{n}}_{n})$$
\end{definition}

La catégorie des espaces vectoriels sur $k$ étant additive, $\calc_{nfiltr}$ est munie de sommes directes.\\

\begin{lemme}\label{isosummodrees}
Soient $(V,F^{\bullet}_{1},F^{\bullet}_{2},...,F^{\bullet}_{n}) \in \calc_{nfiltr}$, $(V',{F^{\bullet}_{1}}',{F^{\bullet}_{2}}',...,{F^{\bullet}_{n}}') \in \calc_{nfiltr}$ et $(V \oplus V',F^{\bullet}_{1}\oplus {F^{\bullet}_{1}}',F^{\bullet}_{2}\oplus{ F^{\bullet}_{2}}',...,F^{\bullet}_{n}\oplus {F^{\bullet}_{n}}')$ l'objet de $ \calc_{nfiltr}$ qui s'en déduit par somme directe, on a l'isomorphisme de ${k}[u_{1},u_{2},...,u_{n}]$-modules suivant :
$$ R^{n}(V,F^{\bullet}_{i}) \oplus R^{n}(V',{F^{\bullet}_{i}}') \cong R^{n}(V\oplus V',F^{\bullet}_{i}\oplus {F^{\bullet}_{i}}')$$
\end{lemme}

\begin{preuve}
Il suffit de constater que par définition de la filtration sur la somme directe :\\
 $(F^{p_1}_{1}\oplus {F^{p_1}_{1}}' \cap F^{p_2}_{2}\oplus {F^{p_2}_{2}}' \cap ... \cap F^{p_n}_{n}\oplus {F^{p_n}_{n}}')=(F^{p_{1}}_{1} \cap F^{p_{2}}_{2} \cap ... \cap F^{p_{n}}_{n}) \oplus ({F^{p_{1}}_{1}}' \cap {F^{p_{2}}_{2}}' \cap ... \cap {F^{p_{n}}_{n}}')$.
\end{preuve}

La catégorie des espaces vectoriels sur $k$ étant munie d'un produit tensoriel, $\calc_{nfiltr}$ est donc une catégorie tensorielle :\\

\begin{lemme}\label{isoprodmodrees}
Soient $(V,F^{\bullet}_{1},F^{\bullet}_{2},...,F^{\bullet}_{n}) \in \calc_{nfiltr}$, $(V',{F^{\bullet}_{1}}',{F^{\bullet}_{2}}',...,{F^{\bullet}_{n}}') \in \calc_{nfiltr}$ et $(V \otimes V',F^{\bullet}_{1}\otimes {F^{\bullet}_{1}}',F^{\bullet}_{2}\otimes {F^{\bullet}_{2}}',...,F^{\bullet}_{n}\otimes {F^{\bullet}_{n}}')$ l'objet de $ \calc_{nfiltr}$ qui s'en déduit par produit tensoriel, on a l'isomorphisme de ${k}[u_{1},u_{2},...,u_{n}]$-modules suivant :
$$ R^{n}(V,F^{\bullet}_{i}) \otimes_{{k}[u_{1},u_{2},...,u_{n}]} R^{n}(V',{F^{\bullet}_{i}}') \cong R^{n}(V\otimes_{k} V',F^{\bullet}_{i}\otimes_{k} {F^{\bullet}_{i}}')$$
\end{lemme}

\begin{preuve}
On remarque que : $(F^{p_1}_{1}\otimes_{k} {F^{p_1}_{1}}' \cap F^{p_2}_{2}\otimes_{k} {F^{p_2}_{2}}' \cap ... \cap F^{p_n}_{n}\otimes_{k} {F^{p_n}_{n}}')=(F^{p_{1}}_{1} \cap F^{p_{2}}_{2} \cap ... \cap F^{p_{n}}_{n}) \otimes_{k} ({F^{p_{1}}_{1}}' \cap {F^{p_{2}}_{2}}' \cap ... \cap {F^{p_{n}}_{n}}')$
ce qui permet de conclure. 
\end{preuve}
{\bf Remarque :}
Le deux lemmes précédents restent vrais pour une famille finie d'éléments de $\calc_{nfiltr}$.\\

\begin{lemme}
Soit $(p_{i}) \in {\bf Z}^{n}$ et $(V,F^{\bullet}_{1},F^{\bullet}_{2},...,F^{\bullet}_{n}) \in \calc_{nfiltr}$, alors l'isomorphisme :
$$(\prod_{i \in [1,n]}u_{i}^{p_i}) .R^{n}(V,F^{\bullet}_{i}) \cong R^{n}(V,Dec^{p_{i}}F^{\bullet}_{i})$$ est un isomorphisme de ${k}[u_{1},u_{2},...,u_{n}]$-modules.\\
\end{lemme}

\begin{preuve}
Il suffit de montrer pour tout $k \in [1,n]$ que $u_{k}. R^{n}(V,F^{\bullet}_{i}) \rightarrow R^{n}(V,Dec^{1_{k}}F^{\bullet}_{i})$, où ¨$1_k$¨ signifie que l'on ne décale que la $k$-ième filtration de $1$, est un isomorphisme de ${k}[u_{1},u_{2},...,u_{n}]$-modules. C'est bien le cas car on envoie l'elément $u_{k}.(u_{1}^{-p_{1}}...u_{k}^{-p_{k}}...\,u_{n}^{-p_{n}})\,v$ où $v \in (F^{p_{1}}_{1} \cap ... \cap F^{p_{k}}_{k} \cap ... \cap F^{p_{n}}_{n})$ sur l'élément $ u_{1}^{-p_{1}}...u_{k}^{-p_{k}+1}...\,u_{n}^{-p_{n}}\,v$ où $v \in (F^{p_{1}}_{1} \cap ...\cap Dec^{1}F^{p_{k}-1}_{k} \cap ... \cap F^{p_{n}}_{n})$ car $Dec^{1}F^{p_{k}-1}_{k}=F^{p_{k}}_{k}$.
\end{preuve}
\begin{center}
{\bf Module de Rees associé à une espace vectoriel filtré $(V,F^{\bullet })$ }
\end{center}

Le module de Rees associé à $(V,F^{\bullet })$, $R(V,F^{\bullet})$, est le sous-${k}[u]$-module de $V$ $\otimes _{k} {k}[u,u^{-1}]$ engendr\'e par les \'el\'ements de la forme $<u^{-p} a_p,a_p \in F^p>\,$ : $R(V,F^{\bullet})=\sum_{p}\,u^{-p}\,F^{p}$. D'après l'étude des scindements associés à un espace vectoriel filtré, on peut choisir une base $\{ v_{i} \}_{i \in I}$ de $V$ adaptée à la filtration c'est à dire telle que pour tout $i \in I$ $v_{i} \in F^{p(i)}-F^{p(i+1)}$. Les éléments $u^{-p(i)}v_{i}$ forment une base de $R(V,F^{\bullet})$.\\ 

\begin{lemme}On a les isomorphismes :\\
${}$\\
$(i)$  $R(V,F^{\bullet })/(u-1)\cong V$.\\
${}$\\
$(ii)$ $R(V,F^{\bullet })[u^{-1}]\cong V\otimes k[u,u^{-1}]$.\\
${}$\\
$(iii)$ $R(V,F^{\bullet })/(u)\cong Gr^{F}V$.
\end{lemme}

\begin{preuve}
L'assertion $(i)$ est évidente.\\
Montrons le $(ii)$. Soit $\{ v_{i} \}_{i \in I}$ une base de $V$ compatible avec la filtration comme exhibé au dessus. L'isomorphisme consiste à envoyer pour tout $i \in I$ l'élément de la base du module de Rees $u^{-p(i)}v_{i} \in R(V,F^{\bullet})$ vers l' élément $v_{i} \in V$.
Pour le (iii), on définit le morphisme qui suit : $\Psi : R(V,F^{\bullet }) \rightarrow Gr^{F}V$ par $ \sum_{i \in I}  u^{p(i)}v_{i} \mapsto (\overline{v_{i}})_{i \in I}$ où $\overline{v_{i}}$ est la projection de $v_{i} \in F^{p(i)} $ sur $Gr_{p(i)}^{F}V$. Soit $\Phi$ la surjection canonique $\Phi :R(V,F^{\bullet }) \rightarrow R(V,F^{\bullet })/(u)$. On a $Ker \Psi =Ker \Phi=u.R(V,F^{\bullet })$. 
\end{preuve}
\begin{center}
{\bf Module de Rees associé à un espace vectoriel muni de deux filtrations $(V,F^{\bullet },G^{\bullet })$ }
\end{center}

Le module double de Rees, noté $RR(V,F^{\bullet },{G}^{\bullet  })$, associé à $(V,F^{\bullet},G^{\bullet})$ est le sous $A$-module $ RR(V,F^{\bullet },{G}^{\bullet }) \subset {k} [{u},{{u}^{-1}},{v},{{v}^{-1}} ] \otimes_{k} V$ engendré par les éléments de la forme ${{u}^{-p}}\,{{v}^{-q}}\,a_{p,q}$ pour $a_{p,q} \in F^{p} \cap G^{q}$ i.e. : $RR(V,F^{\bullet },G^{\bullet })=\sum_{p,q}\,u^{-p}\,v^{-q} (F^{p}\cap G^{q})$. D'après l'étude des scindements associés à un espace vectoriel muni de deux filtrations faites plus haut, on sait qu'il existe une base $\{ v_{i,j} \}_{(i,j) _in I \times J}$ de $V$ qui est adaptée aux deux filtrations c'est à dire telle que l'on ait $F^{p}=< v_{i,j},i \in I , p(i) \leq p, j \in J>$ et de m\^eme $G^{q}=< v_{i,j},j \in J , q(j) \leq q,i\in I>$.

Rappelons que comme pour tout $p$ $Gr^{F}_{p}V$ est le quotient d'un sous-objet $F^pV$ de $V$ par un sous-objet $F^{p+1}V$, il se trouve muni d'une filtration induite par $G^{\bullet}$ (cf section sur les filtrations); ainsi $Gr^{F}V$ est aussi muni d'une filtration induite par $G^\bullet$, que l'on notera $G_{ind}^{\bullet }$. Symétriquement on notera $F_{ind}^{\bullet }$ la filtration induite par $F^{\bullet }$ sur $Gr^{G}V$.\\

\begin{lemme}\label{lemmedebase}
On a les isomorphismes de ${k}[u]$-modules :\\
${}$\\
$(i)$ $RR(V,F^{\bullet },G^{\bullet})/(v-1)\cong R(V,F^{\bullet})$\\
${}$\\
$(ii)$ $RR(V,F^{\bullet },G^{\bullet })/(v)\cong R(Gr^{G}V,F^{\bullet}_{ind})$\\
${}$\\
Ainsi que les isomorphismes de ${k}[v]$-modules :\\
${}$\\
$(iii)$ $RR(V,F^{\bullet },G^{\bullet})/(u-1)\cong R(V,G^{\bullet})$\\
$(iv)$ $RR(V,F^{\bullet },G^{\bullet })/(u)\cong R(Gr^{F}V,G^{\bullet}_{ind})$\\
${}$\\
Et les isomorphismes :\\
${}$\\
$(v)$ $RR(V,F^{\bullet },G^{\bullet })/(u,v)\cong Gr^{G}Gr^{F}V \cong Gr^{F}Gr^{G}V$.\\
${}$\\
$(vi)$ $RR(V,F^{\bullet },G^{\bullet})/(u-1,v-1)\cong V$.\\
${}$\\
$(vii)$ $RR(V,F^{\bullet },G^{\bullet })[u^{-1}]\cong R(V,G^{\bullet})\otimes k[u,u^{-1}]$.\\
${}$\\
$(viii)$ $RR(V,F^{\bullet },G^{\bullet })[v^{-1}]\cong R(V,F^{\bullet})\otimes k[v,v^{-1}]$.\\
${}$\\
$(ix)$ $RR(V,F^{\bullet },G^{\bullet })[u^{-1},v^{-1}]\cong V\otimes k[u,u^{-1},v,v^{-1}]$.\\
\end{lemme}

\begin{preuve}
L'assertion $(i)$ est évidente. Démontrons la deuxième assertion. La preuve est similaire à celle du lemme précédent, on définit le morphisme $\Psi : RR(V,F^{\bullet },G^{\bullet }) \rightarrow R(Gr^{F}V,G^{\bullet}_{ind})$ par  $ \sum_{p,q}  u^{p}v^{q}a_{p,q} \mapsto (\overline{a_{p,q}})_{p,q \leq 0}$ où $\overline{a_{p,q}}$ est la projection de $a_{p,q} \in F^{p} \cap G^{q} $ sur $\frac{ F^{p} \cap G^{q} }{ F^{p-1} \cap G^{q}}$. Le résultat découle du fait que $Ker \Psi= u. RR(V,F^{\bullet },G^{\bullet })$.

$(ii)$ et $(iv)$ sont des énoncés respectivement équivalents à $(i)$ et $(ii)$ (en permutant les deux filtrations en jeu $F^{\bullet }$ et $G^{\bullet })$.

Pour prouver $(v)$, on applique le $(ii)$ du lemme précédent au  ${k}[u]$-module $RR(V,F^{\bullet },G^{\bullet })/(v)\cong R(Gr^{G}V,F^{\bullet}_{ind})$, en découle le premier isomorphisme. Le deuxième isomorphisme $$Gr^{G}Gr^{F}V \cong Gr^{F}Gr^{G}V$$ 
vient du lemme de Zassenhaus. $(vi)$ est direct, tout comme $(vii)$. $(viii)$ est le pendant de $(vii)$. Pour prouver $(ix)$, on applique le lemme précité $(ii)$ à $(vii)$ ou $(viii)$. 
\end{preuve}

\subsubsection{Faisceau cohérent sur $\text{Spec}\,A$ associé à un $A$-module : construction ${}^\sim$ }

Ce paragraphe suit exactement la construction faite dans \cite{har}. Il a pour but de fixer les notations et les propriétés de base de la construction d'un faisceau cohérent sur $\text{Spec}\,A$ associé à un $A$-module auxquelles on se réf\`erera par la suite.\\

\begin{definition}
Soit $A$ un anneau et $M$ un $A$-module. On d\'efinit le faisceau associé à $M$ sur $\emph{Spec}\, A$, noté $\tilde M$ de la façon suivante. Pour tout ideal premier $p \subset A$, soit $M_{p}$ la localisation de $M$ en $p$. Pour tout ouvert $U \subset \emph{Spec}\,A$, on définit le groupe $\tilde M (U)$ comme étant le groupe des fonctions $s: U \rightarrow \cup_{p \in U} M_{p}$ tel que pour tout $p \in U$, $s(p) \in M_{p}$ et tel que localement $s$ est une fraction $m/f$ avec $m \in M$ et $f \in A$.\\

\end{definition}
  
\begin{proposition}\label{tilde1}\cite{har}
Soit $\tilde M$ le faisceau associé au $A$-module $M$ sur $X=\emph{Spec}\,A$. On a :\\
${}$\\
$(i)$ $\tilde M$ est un ${\mathcal O}_{X}$-module.\\
${}$\\
$(ii)$ Pour tout $p \in X$, le germe $({\tilde M})_{p}$ du faisceau $\tilde M$ en $p$ est isomorphe au module localisé $M_{p}$.\\
${}$\\
$(iii)$ Pour tout $f \in A$, le $A_{f}-$module ${\tilde M}(D(f))$ est isomorphe au module localisé $M_f$.\\
${}$\\
$(iv)$ En particulier, $\Gamma(X,{\tilde M})=M$.\\
\end{proposition}

Nous aurons aussi besoin de la proposition :\\

\begin{proposition}\label{tilde2}\cite{har}
Soit $A$ un anneau et $X=\emph{Spec}\,A$. Soit $A \rightarrow B$ un homomorphisme d'anneaux, et $f: \emph{Spec}\,B \rightarrow \emph{Spec}\,A$ le morphisme de spectres correspondant. Alors :\\
${}$\\
$(i)$ L'application $M \rightarrow {\tilde M}$ donne un foncteur pleinement fidèle de la catégorie des $A$-modules vers la catégorie des ${\mathcal O}_X$-modules.\\
${}$\\
$(ii)$ Si $M $ et $N$ sont des $A$-modules, alors $\widetilde {M \otimes_{A} N} \cong {\tilde M} \otimes_{{\mathcal O}_X} {\tilde N}$.\\
${}$\\
$(iii)$ Si $\{M_i\}$ est une famille de $A$-modules, alors $\widetilde{\oplus M_i  } \cong \oplus \tilde{M_i}$.\\
${}$\\
$(iv)$ Pour tout $B$-module $N$, on a $f_{*}({\tilde N}) \cong \tilde{( N_{A})}$ où $N_A$ est $N$ considéré comme un $A$-module.\\
$(v)$ Pour tout $A$-module $M$, on a $f^{*}({\tilde M}) \cong \widetilde{M \otimes_{A} B }$.\\
\end{proposition}

\subsection{Faisceaux cohérents, faisceaux cohérents réflexifs}\label{sectionreflexif}

\subsubsection{Définitions et propriétés}

Nous allons introduire ici des notions utiles à l'étude des faisceaux cohérents. Nous suivrons essentiellement le plan de \cite{oss}, nous ne rappelerons les démonstrations que lorsqu'elles introduisent des notions nécessaires à la suite de notre étude des espaces vectoriels filtrés. Commençons par quelques notions
d'algèbre homologique.\\

Dans cette partie, on entendra par variété algébrique tout schéma projectif lisse et séparé de type fini sur un corps algébriquement clos de caractéristique nulle $k$. 
\begin{center}
{\bf Algèbre homologique, lieu singulier d'un faisceau cohérent}
\end{center}
 
Soit $\calf$ un faisceau cohérent sur une variété algébrique de dimension $n$. Le germe de sections $\calf_{x}$ est un module finiment engendré sur l'anneau local régulier n{\oe}thérien $\calo_{X,x}$.\\
\begin{definition}
La dimension projective de $\calf_{x}$ sur $\calo_{X,x}$, $\emph{pd}(\calf_{x})$, est la longueur minimale d'une résolution projective de $\calf_{x}$.\\
\end{definition} 
D'après \cite{ser}, IV-27, $\text{pd}(\calf_{x})$ est aussi le plus petit entier $k$ tel que pour tout $\calo_{X,x}$-module finiment engendré $M$ et tout $i > k$ on ait $\text{Ext}^{i}_{\calo_{X,x}}(\calf_{x},M)=0$. Elle est parfois appelée dimension homologique de $\calf_{x}$ sur $\calo_{X,x}$ et alors notée $\text{dh}(\calf_{x})$.

Rappelons la définition d'une suite régulière pour un $A$-module $M$. Une suite d'éléments de $A$, $x_{1},x_{2},...,x_{r}$ est une suite régulière pour $M$ si $x_1$ n'est pas un diviseur de zéro dans $M$ et si pour tout $i=2,...,r$, $x_i$ n'est pas un diviseur de zéro dans $M/(x_{1},...,x_{i-1})M$. Si $A$ est un anneau local d'idéal maximal $\mathfrak m$, alors la profondeur de $M$, notée $\text{depth}(M)$, est la longueur maximale d'une suite régulière $x_{1},x_{2},...,x_{r}$ pour $M$ avec pour tout $i$, $x_{i} \in {\mathfrak m }$.\\

\begin{definition}
La profondeur du $\calo_{X,x}$-module $\calf_{x}$, $\emph{depth}(\calf_{x})$, est la longueur maximale d'une suite régulière pour $\calf_x$ dans $\calo_{X,x}$.
\end{definition}
Comme $\calo_{X,x}$ est un anneau local régulier, les deux notions précédentes sont reliées par la formule d'Auslander-Buchsbaum (\cite{ser},\cite{har}) :
$$ \text{pd}(\calf_{x})+\text{depth}(\calf_{x})=\text{dim}(\calo_{X,x}).$$
\begin{lemme}\cite{oss} Soit $\calf$ un faisceau cohérent sur une variété algébrique de dimension $n$. Alors :\\
${}$\\
$(i)$ $\emph{pd}(\calf_{x}) \leq k $ si et seulement si pour tout $i > k $ $({\cale xt}^{i}_{\calo_{X}}(\calf,\calo_{X}))_{x}=0$.\\
${}$\\
$(ii)$ $\emph{dim}(\emph{supp} \,{\cale xt}^{i}_{\calo_{X}}(\calf,\calo_{X})) \leq n-i$.\\
\end{lemme}
$${}$$
\begin{definition}
Soit $\calf $ un faisceau cohérent sur une variété algébrique $X$. Le $m$-ième ensemble de singularité pour la dimension projective de $\calf$ est :
$$ S_{m}(\calf)=\{ x  \in X \vert \emph{depth}(\calf_{x}) \leq m \}.$$
\end{definition}
On a ainsi une suite d'inclusions $X=S_{n}(\calf) \supset S_{n-1}(\calf) \supset ... \supset S_{0}(\calf)$. D'après la formule d'Auslander-Buchsbaum, on peut aussi définir ces ensembles par $ S_{m}(\calf)=\{ x  \in X \vert \text{pd}(\calf_{x}) \geq n-m \}$. Du lemme précédent découle donc l'expression :
$$ S_{m}(\calf)= \cup_{i=n-m}^{n} \text{supp} \,{\cale xt}^{i}_{\calo_{X}}(\calf,\calo_{X}).$$
On peut ainsi montrer :\\
\begin{lemme}(\cite{oss}, II, lemme 1.1.4) 
Les ensembles $ S_{m}(\calf)$ sont des sous-ensembles algébriques fermés et $\emph{dim}( S_{m}(\calf)) \geq n-m$.\\
\end{lemme}

Le germe de $\calf$ en $x \in X$, $\calf_{x}$ est un $\calo_{X,x}$-module libre si et seulement si sa dimension projective est nulle, ainsi l'ensemble singulier de $\calf$ donné par :
$$S(\calf)=\{ x \in X \vert \calf_{x} \text{n'est pas un module libre sur } \calo_{X,x} \}$$
coïncide avec $S_{n-1}(\calf)$.   
Comme corollaire du lemme précédent vient donc :\\
\begin{corollaire}
L'ensemble singulier $S(\calf )$ d'un faisceau cohérent $\calf $ sur une variété algébrique $X$ est de codimension au moins $1$.\\
\end{corollaire}
Ainsi, sur $X \backslash S(\calf)$, $\calf$ est localement libre. Si $X$ est connexe, on peut donc définir le rang du faisceau cohérent $\calf$ par :
$$ \text{rg}(\calf)=\text{rg}(\calf \vert_{X \backslash S(\calf)}).$$ 

\begin{definition}
Un faisceau cohérent $\calf $ sur une variété algébrique $X$ est dit sans torsion si tout germe $\calf_{x}$ est un $\calo_{X,x}$-module sans torsion, i.e. si $f \in \calf_x$ et $a \in \calo_{X,x}$ sont tels que $f.a=0$ alors, ou $f=0$ ou $a=0$.\\
\end{definition}
Les faisceaux localement libres sont sans torsion, les sous-faisceaux de faisceaux sans-torsion sont sans torsion.\\

Le lemme ci-dessus nous donne aussi le résultat

\begin{corollaire}\label{singfaiscohcodim2}
L'ensemble singulier d'un faisceau cohérent sans torsion est au moins de codimension $2$.\\
\end{corollaire}
On en déduit que tout faisceau cohérent sans torsion sur une courbe algébrique est localement libre.

\begin{center}
{\bf Faisceaux réflexifs}
\end{center}

Soit $\calf $ un faisceau cohérent sur une variété algébrique $X$. Le dual de $\calf $ est le faisceau $\calf^{ * } =\calh om (\calf,\calo_{X})$. Il y a un morphisme naturel $\mu : \calf \rightarrow \calf^{**}$.  Le noyau de ce morphisme est le sous-faisceau de torsion $T(\calf )$ de $\calf$ i.e. pour tout $x \in X$ $\text{ker}\, \mu_{x}=\{ f \in \calf_{x} \vert \exists  a \in {\calo_{X,x} \backslash  \{ 0 \}} , fa=0 \}$ (cf \cite{grarem}, p69 pour une preuve dans le cadre analytique).\\

\begin{definition}
Le faisceau cohérent $\calf $ est dit reflexif si le morphisme naturel $\mu $ de $\calf $ vers son bidual $\calf \rightarrow \calf^{**}$ est un isomorphisme.\\
\end{definition}
les faisceaux localement libres sont des faisceaux réflexifs. Comme on l'a vu, pour un faisceau cohérent $\calf $, $\text{ker}\, \mu=T(\calf)$, ainsi faisceaux réflexifs sont sans torsion. Les faisceaux réflexifs sont plus généraux que les faisceaux localement libres mais moins généraux que les faisceaux sans torsion. La présentation des faisceaux réflexifs de la proposition qui suit est utile.\\
\begin{proposition}\cite{har1}\label{har1}
Un faisceau cohérent $\calf$ sur un schéma lisse et séparé $X$ est réflexif si et seulement si il peut-être inclu localement dans une suite exacte $$ 0 \rightarrow \calf \rightarrow \cale \rightarrow \calg \rightarrow 0 ,$$
où $\cale $ est localement libre et $\calg$ est sans torsion.\\
\end{proposition}
\begin{preuve}La question est locale, on omet d'écrire les restrictions à des ouverts. Supposons que $\calf$ est réflexif. $\calf^*$ étant cohérent, on peut trouver une résolution par des faisceaux localements libres $ \call_{1} \rightarrow \call_{0} \rightarrow \calf^{*} \rightarrow 0$. Le foncteur dualisant $\calh om(-,\calo_{X})$ étant exact à gauche, on obtient la suite exacte $0 \rightarrow \calf^{**} \rightarrow  \call_{0}^{*} \rightarrow \call_{1}^{*} $. Notons par $Im\,\phi$ l'image du morphisme de $ \call_{0}^{*}$ vers $ \call_{1}^{*}$. $\calf $ étant réflexif, on a un isomorphisme $\calf \cong \calf^{** }$. En composant par cet isomorphisme on a la suite exacte voulue $ 0 \rightarrow \calf \rightarrow \call^{*}_{0} \rightarrow Im\, \phi \rightarrow 0 $.

Réciproquement, si on a une suite exacte $ 0 \rightarrow \calf \rightarrow \cale \rightarrow \calg \rightarrow 0 $ avec les propriétés voulues, alors $\calf $ est isomorphe à un sous-faisceau d'un faisceau localement libre donc est sans torsion, ainsi l'application naturelle $\mu : \calf \rightarrow \calf^{**}$ est injective. Comme $\cale $ est localement libre, donc réflexif, on peut voir $\calf^{**}$ comme un sous-faisceau de $\cale$. Ainsi le quotient $\calf^{**}/\calf$, qui est un faisceau de torsion, est isomorphe à un sous-faisceau de $\calg$ et est donc sans torsion car $\calg$ est sans torsion, donc nul, d'où l'isomorphisme.  
\end{preuve}
Tout faisceau cohérent $\calf$ a une résolution par des faisceaux localement libre de la forme :
$$ \call_{1} \rightarrow \call_{0} \rightarrow \calf \rightarrow 0$$
Dualisons cette suite exacte, on obtient : 
$$0 \rightarrow  \calf^{*} \rightarrow \call_{0}^{*} \rightarrow \call_{1}^{*}$$
Notons par $\calg$ l'image du morphisme $\call_{0}^{*} \rightarrow \call_{1}^{*}$, $\calg$ est un sous-faisceau de $\call_{1}^{*}$, donc sans torsion. Il vient donc : 
$$ 0 \rightarrow  \calf^{*} \rightarrow \call_{0}^{*} \rightarrow \calg \rightarrow 0.$$
D'où, comme corollaire de la proposition précédente :\\
\begin{corollaire}\label{dualcohref}
Le dual de tout faisceau cohérent est réflexif.\\
\end{corollaire}
La proposition qui suit est la principale caract\'erisation des faisceaux réflexifs dont nous aurons besoin. 
\begin{proposition}\cite{har1}\label{caracrefl}
Un faisceau cohérent $\calf$ sur une variété $X$ est réflexif si et seulement si\\
\hspace*{2cm} $(i)$ $\calf$ est sans torsion et,\\
\hspace*{2cm} $(ii)$ pour tout $x \in X$, $\emph{depth}\,\calf_{x} \geq 2$. \\
\end{proposition}
Nous allons voir que les faisceaux réflexifs sont déterminés par leurs restrictions aux ensembles complémentaires de sous-ensembles de codimension $\geq 2$. Ce sont en fait les faisceaux cohérents qui sont à la fois sans torsion et normaux.\\

\begin{definition}Un faisceau cohérent $\calf$ sur $X$ est normal si pour tout ouvert $\calu \subset X$ et tout sous-ensemble fermé $Y$ de codimension $\geq 2$, le morphisme de restriction $\calf(\calu) \rightarrow \calf(\calu-Y)$ est bijectif.\\  

\end{definition} 

\begin{proposition}\cite{har1}\label{reflnormal} Soit $\calf$ un faisceau cohérent sur $X$, alors les conditions suivantes sont équivalentes :\\ 
\hspace*{2cm} $(i)$ $\calf$ est réflexif,\\
\hspace*{2cm} $(ii)$ $\calf$ est sans torsion et normal. \\
\end{proposition}

Cette proposition sera utilisée sous la forme suivante : si $f : \cale \rightarrow \calf$ est un morphisme de faisceaux réflexifs sur $X$ qui est un isomorphisme en dehors d'un ensemble $Y$ de codimension $\geq 2$, alors $f$ est un isomorphisme sur $X$, en effet, pour tout ouvert $\calu \subset X$, la suite d'isomorphismes suivante permet de conclure : $$  \cale(\calu) \cong \cale(\calu-Y) \cong \calf(\calu-Y) \cong \calf(\calu).$$

\subsubsection{ $\mu$-semistabilité}

Soit $\calf $ un faisceau cohérent sans torsion de rang $r$ sur une variété algébrique $X$. Rappelons que le fibré en droite déterminant associé à $\calf$ est définit par $\text{det}(\calf)=(\Lambda^{r} \calf)^{**}$ (d'après \ref{har1}, un faisceau réflexif de rang $1$ est un fibré en droites). La première classe de Chern de $\calf$ est définie par
$$\text{c}_{1}(\calf)=\text{deg}_{}(\text{det} \calf),$$
et est à valeurs dans ${\bf Z}$ (voir \cite{har} pour la définition et les propriétés du degré).\\

 Nous n'avons donc définit la première classe de Chern que pour les faisceaux cohérents sans torsion comme étant le degré du fibré déterminant qui est un fibré en droite. Supposons qu'un faisceau cohérent (avec ou sans torsion) $\calf $ admette une résolution finie par des faisceaux localement libres $ 0 \rightarrow \calf_{n} \rightarrow \calf_{n-1} \rightarrow ... \rightarrow \calf_{0} \rightarrow \calf \rightarrow 0$. On définit alors $\text{det}(\calf )=\otimes_{i}\,\text{det}(\calf_{i})^{(-1)^{i}}$. On peut montrer que cette définition ne dépend pas de la résolution choisie et permet ainsi de définir $\text{c}_{1}(\calf )$. Cette définition fait sens ici car d'après \cite{har}, III Ex. 6.8. et 6.9., tout faisceau cohérent sur une variété algébrique lisse admet une résolution finie par des faisceaux localement libres.\\

D'après la propriété d'additivité du degré (cf. \cite{har},Ex. 6.12. p149), la première classe de Chern est additive
c'est à dire si l'on a une suite exacte de faisceaux cohérents
$$0  \rightarrow \calf ' \rightarrow \calf \rightarrow \calf'' \rightarrow 0,$$
alors,
$$\text{c}_{1}(\calf)=\text{c}_{1}(\calf')+\text{c}_{1}(\calf'').$$
Définissons la pente de $\calf $ par 
$$\mu(\calf)=\frac{\text{c}_{1}(\calf)}{\text{rg}(\calf)},$$
si $\text{rg}(\calf)=0$, et
$$\mu(\calf)=0,$$
sinon. 

On peut maintenant introduire la notion de semistabilité.\\

\begin{definition}
Un faisceau cohérent sans torsion sur $X$ est dit $\mu$-semistable si pour tout sous-faisceau cohérent tel que $0 \varsubsetneq \cale \subset \calf$ on a
$$\mu(\cale) \leq \mu(\calf).$$ 
\end{definition}

Notons  $\calr efl(\mu)$ la catégorie des faisceaux réflexifs $\mu$-semistables de pente $\mu$ sur une variété algébrique $X$.\\

Rappelons que par additivité de la première classe de Chern, si l'on a une suite exacte courte de faisceaux cohérents sans torsion sur une variété $X$ de la forme
$$ 0 \rightarrow \cale \rightarrow \calf \rightarrow \calg \rightarrow 0,$$
alors,
$$ \mu(\calf)=\frac{\text{c}_{1}(\cale)+\text{c}_{1}(\calg)}{\text{rg}(\cale)+\text{rg}(\calg)},$$
la pente de $\calf$ est donc barycentre des pentes respectives de $\cale$ et $\calg$. De façon explicite :$$\mu(\calf)=\frac{1}{{\text{rg}(\calf)}}(\text{rg}(\cale)\mu(\cale) +\text{rg}(\calg)\mu(\calg)).$$

Le lemme suivant généralise la proposition $5.8$ p$79$ dans \cite{pot}, la démonstration est similaire.\\
\begin{lemme}\label{exact}
Si l'on a une suite exacte de faisceaux sans torsion $0 \rightarrow \cale \rightarrow \calf \rightarrow \calg \rightarrow 0$ telle que $\cale $ et $\calg$ sont $\mu$-semistables de pente $\mu$, alors $\calf $ est $\mu$-semistable de pente $\mu$.\\
\end{lemme}

\begin{preuve}
Comme explicité ci-dessus, $\mu(\calf)$ est barycentre des pentes de $\cale$ et $\calg$, donc $\mu(\calf)=\mu$.

Montrons que $\calf$ est $\mu$-semistable. Soit $\calf'$ un sous-faisceau cohérent non nul de $\calf$. Il est sans-torsion comme sous-faisceau d'un faisceau sans torsion. Soit $\calg'$ l'image de $\calf'$ dans $\calg$ et $\cale'$ l'intersection $\cale \cap \calf'$. $\cale'$ est cohérent comme noyau du morphisme de faisceaux cohérents $ \cale \rightarrow \calf \rightarrow \calf /\calf'$ et sans torsion comme sous-faisceau de $\cale$. $\calg'$ est cohérent comme conoyau du morphisme de faisceaux cohérents $\cale' \rightarrow \calf'$ et sans torsion comme sous-faisceau de $\calg$. On a ainsi la suite exacte courte de faisceaux cohérents sans torsion
$$0 \rightarrow \cale' \rightarrow \calf' \rightarrow \calg' \rightarrow 0$$
donc la pente de $\calf'$ est barycentre à coefficients positifs de $\mu(\cale')$ et $\mu(\calg')$. Si ces faisceaux sont non nuls, on a $\mu(\cale') \leq \mu$ et $\mu(\calg') \leq \mu$ donc $\mu(\calf') \leq \mu$. Si $\cale'=0$, alors $\calf'$ s'dentifie à un sous-faisceau cohérent de $\calg$, si $\calg'=0$, $\calf'$ est un sous-faisceau cohérent de $\cale'$, dans ces deux cas, on a évidement $\mu(\calf') \leq \mu$.   
\end{preuve}

Soit $\calf $ un faisceau cohérent sur une variété algébrique $X$. On a un morphisme canonique $\nu : \calf \rightarrow \calf^{**}$ dont le noyau est la torsion de $\calf$. D'après le corollaire \ref{dualcohref} p.\pageref{dualcohref}, le faisceau dual d'un faisceau cohérent est réflexif donc $\calf^{**}$ est un faisceau reflexif.\\
\begin{definition}
Soit $\calf$ un faisceau cohérent,  le faisceau $\calf^{**}$ est appelé faisceau réflexif associé à $\calf $.\\ 
\end{definition}

\begin{lemme}\label{factorise}
Soit $f: \calf \rightarrow \calg $ un morphisme de faisceaux cohérents, où $\calg$ est réflexif, alors $f$ se factorise de façon unique par $\nu : \calf \rightarrow \calf^{**}$.\\
\end{lemme}

\begin{preuve}
$\calg $ étant réflexif, le morphisme canonique $\nu':\calg \rightarrow \calg^{**}$ est un isomorphisme. Notons par $f^{**}$ l'image de $f$ par le morphisme canonique de $\calh om (\calf, \calg ) \rightarrow \calh om(\calf^{**}, \calg^{**})$. On a le diagramme commutatif suivant :
$$\xymatrix{
\calf \ar[r]^{f} \ar[d]_{\nu} &  \calg \ar[d]_{\nu '}^{\simeq}\\
\calf^{**} \ar[r]_{f^{**}} &  \calg^{**}
}
$$
qui donne la factorisation $f=\nu'{}^{-1}\circ f^{**} \circ \nu$.

\end{preuve}

\begin{lemme}\label{c1pos}
Soit $e : \cale \rightarrow \calf $ un morphisme injectif de faisceaux cohérents sans torsion de même rang sur $X$. Alors les premières classes de Chern de $\cale $ et $\calf $ vérifient l'inégalité $\emph{c}_{1}(\cale) \leq \emph{c}_{1}(\calf)$. De plus le conoyau $Coker(e)$ du morphisme $e$ est de torsion et $\emph{c}_{1}(Coker(e)) \geq 0$.
%\footnote{ Nous n'avons définit la première classe de Chern que pour les faisceaux cohérents sans torsion comme étant le degré du fibré déterminant qui est un fibré en droite. Supposons qu'un faisceau cohérent $\calf $ admette une résolution finie par des faisceaux localement libres $ 0 \rightarrow \calf_{n} \rightarrow \calf_{n-1} \rightarrow ... \rightarrow \calf_{0} \rightarrow \calf \rightarrow 0$. On définit alors $\text{det}(\calf )=\otimes_{i}\,\text{det}(\calf_{i})^{(-1)^{i}}$. On peut montrer que cette définition ne dépend pas de la résolution choisie et permet ainsi de définir $\text{c}_{1}(\calf )$. Cette définition fait sens ici car d'après \cite{har}, III Ex. 6.8. et 6.9., tout faisceau cohérent sur une variété algébrique lisse admet une résolution finie par des faisceaux localement libres.}. \\
\end{lemme}

\begin{preuve}
Prouvons d'abord (cf \cite{oss}, II, lemme 1.1.17) que : un morphisme injectif de faisceau cohérents sans torsion de même rang $e :\cale \rightarrow \calf$ induit un morphisme injectif des fibrés en droite déterminants $\text{det}(e) : \text{det }(\cale ) \rightarrow \text{det }(\calf )$. D'après \cite{oss}, même lemme, $e$ induit un isomorphisme en dehors de l'ensemble algébrique $Y=S(\calf) \cup S(\calf/\cale)$ et donc $\text{det}(e)$ est un isomorphisme en dehors de $Y$. Ainsi $\text{ker}(\text{det}(e))$ est un faisceau de torsion sous-faisceau du faisceau localement libre $\text{det}(\cale)$ donc est nul.

De l'injection de fibrés en droite $\text{det}(e) : \text{det }(\cale) \rightarrow \text{det}(\calf)$ on déduit :
$$ \text{c}_{1}(\cale)=  \text{c}_{1}(\text{det}(\cale))\leq \text{c}_{1}(\text{det}(\calf))= \text{c}_{1}(\calf).$$
Par la formule de Whitney, on déduit de la suite exacte $0 \rightarrow \cale \rightarrow \calf \rightarrow Coker(e) \rightarrow 0$ que $ \text{c}_{1}(Coker(e))=\text{c}_{1}(\calf)-\text{c}_{1}(\cale)$ et donc que $\text{c}_{1}(Coker(e)) \geq 0$.

\end{preuve}

\begin{lemme}\label{codim2} Soit $\calf $ un faisceau cohérent sur une variété algébrique $X$.\\
\hspace*{0.5cm}$(i)$ Si le morphisme canonique $\nu : \calf \rightarrow \calf^{**}$ est injectif, alors c'est un isomorphisme en dehors d'un ensemble algébrique $Y$ de codimension $2$.\\
\hspace*{0.5cm}$(ii)$ Si de plus $X$ est complète, alors, si le morphisme canonique $\nu : \calf \rightarrow \calf^{**}$ est surjectif et que de plus on a $\emph{c}_{1}(\calf)=\emph{c}_{1}(\calf^{**})$, alors c'est un isomorphisme en dehors d'un ensemble algébrique $Y$ de codimension $2$.\\
\end{lemme}

\begin{preuve}$(i)$ Le noyau de $\nu$ est la torsion de $\calf$. Si ce noyau est nul c'est que $\calf $ est sans torsion et donc localement libre en dehors d'un sous-ensemble de $X$ de codimension $2$.\\
$(ii)$ Si $\nu$ est une surjection, on peut alors écrire la suite exacte
$$\xymatrix{
0 \ar[r] & T \ar[r] & \calf  \ar[r]^{\nu } &  \calf^{**} \ar[r] &  0}
$$
où $T$ est le faisceau de torsion de $\calf$. Montrons que la codimension du support de $T$ est au moins $2$ ce qui permettra de conclure. Supposons que $T$ ait une composante $V$ de codimension $1$. Alors, comme $X$ est complète, d'après \cite{ful}, il existe une courbe $i:C \hookrightarrow X$, transverse à $V$. $i^{*}T$ est de torsion sur la courbe $C$ donc supporté par des points qui contribuent chacun positivement à $\text{c}_{1}(i^{*}T)$, donc $\text{c}_{1}(i^{*}T)>0$. Comme $C$ est transverse à $V$, $\text{c}_{1}(i^{*}T)=i^{*}\text{c}_{1}(T)$, donc $i^{*}\text{c}_{1}(T)>0$ ce qui est une contradiction, $\calf$ et $\calf^{**}$ ayant m\^eme première classe de Chern. Donc $\text{Supp}(T)$ est de codimension au moins $2$ ce qui permet de conclure.

\end{preuve}

\begin{proposition}\label{caracref}
Soit $\xymatrix{ 0 \ar[r] &\cale \ar[r]^{e} &\calf \ar[r]^{f}& \calg \ar[r]& 0}$ une suite exacte de faisceaux cohérents sur une variété algébrique $X$ telle que $\calf$ soit réflexif et $\calg $ soit sans torsion, alors $\cale $ est réflexif.\\
\end{proposition}

\begin{preuve} 
Le faisceau cohérent $\cale$ est un sous-faisceau du faisceau réflexif $\calf$ et est donc sans torsion, ainsi le morphisme canonique $\nu :\cale \rightarrow \cale^{**}$ est injectif. On veut montrer que $\nu $ est un isomorphisme. D'après le lemme \ref{factorise} p.\pageref{factorise}, on sait que l'on peut factoriser $e$ par $\nu$. On note $e^{**}$ le morphisme qui s'en déduit. Vient le diagramme commutatif :
 
$$\xymatrix{
{}&  0 \ar[d]\\ 
0 \ar[r] &\cale \ar[r]^{e}\ar[d]^{\nu} &\calf \ar[r]^{f}\ar[d]^{id}& \calg \ar[r] \ar[d]^{\eta}& 0\\
 0 \ar[r] &\cale^{**} \ar[r]^{e^{**}} &\calf \ar[r]^{}& Coker(e^{**})  \ar[r] \ar[d]& 0\\
{} &{} & {} & 0
}
$$
Comme $\cale $ est sans torsion il est localement libre en dehors d'un sous-ensemble algébrique $Y$ de codimension $2$ dans $X$.Ainsi $\nu \vert_{X \backslash Y} : \cale\vert_{X \backslash Y} \rightarrow \cale^{**}\vert_{X \backslash Y}$ soit un isomorphisme. Ainsi, $ \eta\vert_{X \backslash Y} : \calg \rightarrow Coker(e^{**})\vert_{X \backslash Y}$ est un isomorphisme donc le support du noyau $\text{ker}(\eta)$ est de codimension $2$, donc est de torsion, mais $\calg$ est sans torsion d'où $\text{ker}(\eta)=0$. Donc $\eta$ est un isomorphisme, il en va ainsi de même pour $\nu$, ce qui achève la preuve.

\end{preuve}

\begin{theoreme}\label{muab}
Soit $X$ une varièté algébrique lisse et projective sur un corps $k$ algébriquement clos de caractéristique nulle. La catégorie $\emph{Refl}_{\mu}(X)$ des faisceaux réflexifs $\mu$-semistables sur $X$ est abélienne.\\ 
\end{theoreme}

\begin{preuve}
Cette catégorie est une sous-catégorie additive de la catégorie des faisceaux vectoriels, elle a des sommes directes par le lemme \ref{exact} p.\pageref{exact}.

Montrons qu'elle est exacte. Soit $f : \cale \rightarrow \calf$ un morphisme dans cette catégorie :\\
%$\text{Im}(f)$ est un faisceau cohérent sans torsion comme sous-faisceau du faisceau localement libre $\calf$. La suite exacte de faisceau $ 0 \rightarrow \text{Ker}(f) \rightarrow \calf \rightarrow \text{Im}(f) \rightarrow 0$  montre que (cf har1) que $\text{Ker}(f)$ est réflexif donc localement libre. Soit $g : \cald \rightarrow \cale$ un morphisme dans $\text{Sh}_{refl,\mu-semistables,\mu}(X)$ tel que $f \circ g=0$, alors elle se factorise par $i : \text{Ker}(f) \rightarrow \calf$. Donc $\text{Ker}(f)$ est bien un noyau dans cette catégorie. 
%Dans la catégorie des faisceaux $\text{Im}(f)$, $\text{Coker}(f)$ et $\text{Coïm}(f)$ sont des faisceaux cohérents et on a les isomorphismes de faisceaux $\text{Coïm}(f)\cong \cale /\text{Ker}(f) \cong \text{Im}(f)$. Définissons les images, cokernels et coïmages (notées $I$, $CK$ et $CI$) dans $\calc(\nu)$ comme les saturés respectifs $I=(\text{Im}(\calf))^{**}$, $CK=(\text{Coker}(f))^{**}$ et $CI=(\text{Coïm}(f))^{**}$.
%Montrons que $CK$ est un cokernel. Soit $h :\calg \rightarrow \calh $ un morphisme dans $\calc(\nu)$ tel que $h \circ f=0$. Alors comme $\text{Coker}(f)$ est un cokernel dans la catégorie des faisceaux, on peut factoriser $h$ par un morphisme $\text{Coker}(f) \rightarrow \calh$ que l'on note toujours $h$ puis appliquer le lemme \ref{factorise} p.\pageref{factorise}, pour factoriser par $h': CK \rightarrow \calh$.
%Montrons que $I$ est une image...
%Montrons que $CI$ est une coïmage...
%Reste à montrer que $I \cong CI$ :
On dispose dans la catégorie des faisceaux cohérents de noyaux, conoyaux, images et coïmages $\text{Ker}(f)$, $\text{Coker}(f)$, $\text{Im}(f)$ et $\text{Coim}(f)$ qui vérifient les isomorphismes $\text{Im}(f)\cong \cale /\text{Ker}(f) \cong \text{Coim}(f)$ et $\text{Coker}(f) \cong \calg/\text{Im}(f)$. On note $i : \text{Ker}(f) \rightarrow \calf$ et $\pi :
\calg \rightarrow \text{Coker}(f)$.

On définit le noyau de $f$ dans la catégorie $\text{Refl}_{\mu}(X)$, que l'on note $\calk er$, comme étant le faisceau réflexif associé au noyau de $f$ dans la catégorie des faisceaux cohérents $\text{Ker}(f)$. Or $\text{Ker}(f)$ est réflexif. En effet il s'insère dans la suite exacte :
$$ 0 \rightarrow \text{Ker}(f) \rightarrow \cale \rightarrow \text{Coim}(f) \rightarrow 0.$$
Or $\text{Coim}(f)\cong\text{Im}(f)$ et $\text{Im}(f)$ est un sous-faisceau du faisceau réflexif $\calf$, donc sans-torsion, et $\cale $ est réflexif. D'après la proposition \ref{caracref} p.\pageref{caracref}, $\text{Ker}(f)$ est réflexif, ainsi $\calk er=\text{Ker}(f)$. 

$\calk er$ est de pente $\mu$ : $\text{Ker}(f)$ est un sous-faisceau cohérent de $\cale$ donc par semistabilité de $\cale$, $\mu(\text{Ker}(f))\leq \mu$, de même $\text{Im}(f)$ est un sous-faisceau cohérent de $\calf$ donc $\mu(\text{Im}(f))\leq \mu$. Or $\mu $ est barycentre à coefficients positifs de $\mu(\text{Ker}(f))$ et $\mu(\text{Im}(f))$ donc $\mu(\calk er)=\mu(\text{Ker}(f))=\mu(\text{Im}(f))=\mu$. Il est $\mu$-semistable, sinon il aurait un sous-faisceau cohérent de pente $ \geq \mu$ mais ce serait aussi un sous-fibré de $\cale$ ce qui contredirait la semistabilité de $\cale$.  Montrons que $\calk er$ est bien un noyau dans la catégorie $\text{Refl}_{\mu}(X)$. Soit $d : \cald \rightarrow \cale$ un morphisme de faisceaux réflexifs tel que $e \circ d=0$. Alors $d$ se factorise par $\text{Ker}(f)$ dans la catégorie des faisceaux donc par $\calk er$. On a ainsi prouvé que $\calk er$ est un noyau dans $\text{Refl}_{\mu}(X)$.

% En effet soit $e : \cale \rightarrow \calf $ un morphisme de faisceaux réflexifs tel que $f\circ e=0$. On peut alors factoriser $e$ par le noyau dans la catégorie des faisceaux :
%$$
%\xymatrix{
%\cale \ar[r]^{e} \ar[d]_{e'} & \calf \ar[r]^{f} & \calf \\
%\text{Ker}(f) \ar[ru]_{i}
%}
%$$
%est commutatif, ce qui donne , d'après le lemme \ref{factorise}, la factorisation voulue par le diagramme commutatif suivant :
%$$
%\xymatrix{
%\cale \ar[r]^{e} \ar[d]_{e'} & \calf \ar[r]^{f} & \calf \\
%\text{Ker}(f) \ar[r]_-{\nu} &   \calk er=(\text{Ker}(f))^{**} \ar[u]_{i^{**}}
%}
%.$$

On définit le conoyau dans $\text{Refl}_{\mu}(X)$, que l'on note $\calc oker$, comme étant le faisceau réflexif associé au faisceau cohérent $\text{Coker}(f)$, i.e. $\calc oker=(\text{Coker}(f))^{**}$. On note $\pi^{**}: \calf \rightarrow \calc oker$ le morphisme associé à $\pi : \calf \rightarrow \text{Coker}(f)$ ($\pi^{**}=\nu \circ \pi$ où $\nu $ est le morphisme canonique de $\text{Coker}(f)$ vers son bidual). Ainsi définit, $\calc oker$ est bien un conoyau dans la catégorie des faisceaux réflexifs : soit $g :\calf \rightarrow \calg$ un morphisme de faisceaux réflexifs tel que $g \circ f=0$. On peut alors factoriser $g$ ; en effet, dans la catégorie des faisceaux 
$$
\xymatrix{
\cale \ar[r]^{f}  & \calf \ar[r]^{g} \ar[d]_{\pi} & \calg\\
{} & \text{Coker}(f) \ar[ru]_{g'}
}
$$
est commutatif, d'où, par le lemme de factorisation (lemme \ref{factorise}), on obtient le diagramme commutatif :

$$
\xymatrix{
\cale \ar[r]^{f}  & \calf \ar[r]^{g} \ar[dr]^{\pi^{**}} \ar[d]_{\pi} & \calg\\
{} & \text{Coker}(f) \ar[r]_-{\nu} &  \calc oker \ar[u]_{g''}
}
$$
donc $\calc oker$ est bien un conoyau.

Il faut prouver que $\calc oker$ est $\mu$-semistable de pente $\mu$. Considérons le faisceau cohérent $\text{ker}(\pi^{**})$  noyau de $\pi^{**} : \calf  \rightarrow \calc oker$. On a alors un morphisme injectif de $\text{Im}(f)$ vers $\text{ker}(\pi^{**})$. Notons par $T$ le conoyau de ce morphisme et considérons le diagramme commutatif formé de suites exactes courtes : 
$$
\xymatrix{
{} & {} & {} & 0 \ar[d]\\
{} & 0 \ar[d] & {} & T\ar[d] \\
0 \ar[r] &\text{Im}(f) \ar[r]\ar[d] &  \calf \ar[r]^{\pi} \ar[d]^{id} & \text{Coker}(f) \ar[r] \ar[d]^{\nu} & 0 \\
0 \ar[r] & \text{Ker}(\pi^{**})  \ar[r] \ar[d] &  \calf \ar[r]^-{\pi^{**}} & \calc oker =(\text{Coker}(f))^{**} \ar[r] \ar[d] & 0\\
{} & T \ar[d] & {} & 0\\
{} & 0
}
$$ 
Comme $\text{rg}(\text{Coker}(f))=\text{rg}((\text{Coker}(f))^{**})$ on a aussi $\text{rg}(\text{Im}(f))=\text{rg}((\text{Ker}(\pi^{**}))^{**})$. Or,\\
$ \text{c}_{1}(\text{Ker}(\pi^{**}))=\text{c}_{1}(\text{Im}(f))+\text{c}_{1}(T)$, donc $\text{c}_{1}(T) \leq 0 $ sans quoi on aurait $\mu(\text{Ker}(\pi^{**})) > \mu$ ce qui nierait la $\mu$-semistabilité de $\calf$. Mais $\mu(\text{Ker}(\pi^{**})) \leq \mu$ implique $\mu( \calc oker) \geq \mu$ ce qui sachant que $\mu(\text{Coker}(f))=\mu$ implique $\text{c}_{1}(T) \geq 0$ et donc $\text{c}_{1}(T) = 0 $. Ainsi, par le lemme \ref{codim2} p.\pageref{codim2}, $(ii)$, $\nu :  \text{Coker}(f) \rightarrow \calc oker$ est un isomorphisme en dehors d'un sous-ensemble algébrique $Y$ de $X$ de codimension $2$.  On en déduit que $\mu( \calc oker)=\mu( \text{ker}(\pi^{**}))=\mu$ et que le morphisme restreint $\text{Im}(f)\vert_{X \backslash Y} \rightarrow \text{ker}(\pi^{**}) \vert_{X \backslash Y}$ est un isomorphisme.

Montrons que $\calc oker$ est $\mu$-semistable. Soit $E$ un sous-faisceau coh\'erent de $\calc oker$ et soit $F$ le faisceau coh\'erent quotient de $\calc oker$ par $E$, on a la suite exacte $0 \rightarrow E \rightarrow \calc oker \rightarrow F \rightarrow 0$. Comme $\pi^{**}: \calf \rightarrow \calc oker$ est surjective, on a un morphisme surjectif de $\calf$ vers $F$. Notons $K$ son noyau, c'est un sous-faisceau coh\'erent de $\calf$. On a la suite exacte $0 \rightarrow K \rightarrow \calf \rightarrow F \rightarrow 0$. Supposons que $\mu(E) > \mu$, alors comme $\mu( \calc oker)=\mu(\calf)=\mu$, il vient $\mu(F) < \mu$ et donc $\mu(K) > \mu$ ce qui contredit la $\mu$-semistabilit\'e de $\calg$. Donc, si $E$ est un sous-faisceau coh\'erent de $\calc oker$, $\mu(E) \leq \mu$. $\calc oker$ est bien $\mu$-semistable.

On d\'efinit ensuite l'image de $f$ dans $Refl_{\mu}(X)$, $\cali m$, comme \'etant le noyau pour le conoyau, c'est \`a dire le faisceau r\'eflexif associ\'e au faisceau coh\'erent $\text{ker}(\pi^{**})$, noyau faisceautique du morphisme de faisceaux r\'eflexifs $\pi^{**} : \calf  \rightarrow \calc oker$. Comme on l'a montr\'e pour le noyau, $\text{ker}(\pi^{**})$ est d\'ej\`a un faisceau r\'eflexif car $\calg $ est r\'eflexif et $\calc oker$ est sans torsion car r\'eflexif. On a vu de plus que $\text{ker}(\pi^{**})$ est $\mu$-semistable de pente $\mu$. Donc l'image $\cali m\cong\text{ker}(\pi^{**})$ est $\mu$-semistable de pente $\mu$.

La co{\"\i}mage $\calc oim$ de $f$ est le faisceau r\'eflexif associ\'e au conoyau du noyau $i: \calk er \rightarrow \calf$. On a $\calc oim=\text{Coker} (i)^{**}$. D'apr\`es l'\'etude faite sur $\calc oker$, $\calc oim$ est $\mu$-semistable de pente $\mu$. De plus d'apr\`es la proposition \ref{caracref} p.\pageref{caracref}, $\calc oim$ est isomorphe \`a $\text{Coim}(f)$. On a :

$$
\xymatrix{
O \ar[r] & \calk er=\text{Ker}(f)  \ar[r]^{i} &  \cale \ar[r] \ar[rd] & \text{Coker}(i^{}) \ar[r] \ar[d]^{\nu}  & 0\\
{} &  {} & {} &  \calc oim=(\text{Coker}(i))^{**}
}
$$

Reste \`a montrer que l'image et la co{\"\i}mage sont isomorphes. Montrons d'abord que le support de $T$ est de codimension au moins $2$. Supposons qu'il contienne une composante $V$ de codimension $1$. Alors, comme $X$ est compl\`ete, d'apr\`es \cite{ful}, il existe une courbe $i:C \hookrightarrow X$ transverse \`a $V$. $i^{*}T$ est de torsion sur $C$ donc est support\'e par des points de $C$ qui contribuent chacun strictement positivement \`a $\text{c}_{1}(i^{*}T)$ donc $\text{c}_{1}(i^{*}T)>0$. Comme $C$ est transverse \`a $V$, $\text{c}_{1}(i^{*}T)=i^{*}\text{c}_{1}(T)$, donc $i^{*}\text{c}_{1}(T)>0$ ce qui est une contradiction. Donc $\text{supp}(T)$ est de codimension au moins $2$. Ainsi il existe un sous-ensemble alg\'ebrique $Y$ de $X$ de codimension au moins $2$ tel que $\calc oker $ et $\text{Coker}(f)$ soient isomorphes sur $X \backslash Y$. Donc, par le diagramme commutatif pr\'ec\'edent, le morphisme de $\text{Im}(f)$ vers $\cali m$ est un isomorphisme sur $X\backslash Y$. Ainsi la compos\'ee des restrictions \`a $X\backslash Y $ de l'isomorphisme sur $X$ entre $\calc oim$ et $\text{Coim}(f)$ puis de l'isomorphisme sur $X$ entre $\text{Coim}$ et $\text{Im}(f)$ et de l'isomorphisme sur $X\backslash Y$ entre $\text{Im}(f)$ et $\cali m$ fournit un isomorphisme entre $\cali m$ et $\calc oim$ sur $X$ priv\'e d'un ensemble de codimension au moins $2$, $Y$. Ces faisceaux \'etant r\'eflexifs, d'apr\`es \cite{har}, ils sont normaux et donc enti\`erement d\'efinis par leurs restrictions aux compl\'ementaires d'ensembles de codimension au moins $2$, d'o\`u l'isomorphisme $\cali m \cong \calc oim$.

\end{preuve}

Rappelons que les faisceaux réflexifs sur les variétés algébriques de dimension inférieure à $2$ sont les faisceaux localement libres. On retrouve le fait bien connu que (pour les courbes voir \cite{pot} et voir \cite{huyleh} pour les surfaces) :\\
\begin{corollaire}
La catégorie des fibrés vectoriels $\mu$-semistables de pente $\mu$ sur une variété algébrique lisse projective de dimension inférieure à $2$ est abélienne.\\
\end{corollaire}
${}$\\
{\bf Remarque :} Nous appliquerons le théorème \label{muab} et sa démonstration avec $X={\bf P}^2$ pour montrer que la catégorie des fibrés sur ${\bf P}^2$ avec une condition de semistabilité plus forte est abélienne (cf corollaire \ref{geotrifab} p.\pageref{geotrifab}).\\

\subsection{Faisceau cohérent de Rees sur ${\bf A}^n$ associé à un espace vectoriel muni de $n$ filtrations}

\subsubsection{Définition des faisceaux cohérents de Rees sur les espaces affines}
Pour étudier les espaces vectoriels filtrés de $\calc_{nfiltr}$, on leurs associe des faisceaux cohérents sur ${\bf A}^n$ par la construction du faisceau cohérent $\tilde M$ d'un $A$-module de Rees $M$ sur $\text{Spec}(A)$, où $A$ est un anneau. Le faisceau est cohérent car, les filtrations étant exhaustives, le module de Rees est finiment engendré.\\

\begin{definition}
Le faisceau de Rees associé à $(V,F^{\bullet}_{1},F^{\bullet}_{2},...,F^{\bullet}_{n}) \in \calc_{nfiltr}$ est le faisecau cohérent sur ${\bf A}^{n}=\emph{Spec}({k}[u_{1},u_{2},...,u_{n}])$, $\tilde{ R^{n}(V,F^{\bullet}_{i})}$ obtenu à partir du ${k}[u_{1},u_{2},...,u_{n}]$-module de Rees $ R^{n}(V,F^{\bullet}_{i})$. Il sera noté $\xi_{{\bf A}^n}( R^{n}(V,F^{\bullet}_{i}))$. Ainsi :
 $$\xi_{{\bf A}^n}( R^{n}(V,F^{\bullet}_{i}))= R^{n}(V,F^{\bullet}_{i})^{\sim }.$$
\end{definition}
%Le faisceau cohérent sur ${\bf A}^n$ est localement libre car le module de Rees $ R^{n}(V,F^{\bullet}_{i})$ est libre et fini. Pour adopter un language géométrique qui nous conviendra mieux par la suite nous parlerons de fibré vectoriel plut\^ot que de faisceau localement libre.\\
%Cette définition donne une première conséquence directe :\\

Comme $R^{n}(V,F^{\bullet}_{i}$ est un sous ${k}[u_{1},u_{2},...,u_{n}]$-module de ${k}[u_{1},u_{1}^{-1},u_{2},u_{2}^{-1}...,u_{n},u_{n}^{-1}] \otimes V $, $\xi_{{\bf A}^n}( R^{n}(V,F^{\bullet}_{i}))$ est le sous-faisceau de $i^{*}(\calo_{{\bf G}_{m}^{n}}\otimes V)$ engendré par les sections de la forme $ u_{1}^{-p_{1}}u_{2}^{-p_{2}}...u_{n}^{-p_{n}}w$ où $w \in F_{1}^{p_1} \cap F_{1}^{p_2} \cap ... \cap F_{1}^{p_n}$.

\begin{lemme}
La construction du fibré de Rees  à partir des modules de Rees est fonctorielle et le foncteur $ {\xi}_{{\bf A}^n}(.,.,.)$ qui va de $\calc_{nfiltr}$ vers la catégorie des faisceaux cohérents sur ${\bf A}^n$ est exact. \\
\end{lemme}

\begin{preuve}
C'est la proposition \ref{tilde2} p.\pageref{tilde2} (i) appliquée aux modules de Rees d'ordre $n$.
\end{preuve}

Soit $J$ un ensemble de cardinal fini et pour tout $j \in J$, $(V_{j},(F_{ij}^{\bullet })_{i \in [1,n]})$ un élément de $\calc_{nfiltr}$. $(V_{j},(F_{ij}^{\bullet })_{i \in [1,n]})_{j \in J}$ est donc une famille finie d'éléments de $\calc_{nfiltr}$ que l'on notera $(V_{j},F_{ij}^{\bullet })_{j \in J}$. Pour une telle famille :\\

\begin{lemme}\label{isosumfibrees}
$$\xi_{{\bf A}^n}(\otimes_{j \in J}R^{n}(V_{j},F^{\bullet}_{ij})) \cong \otimes_{j \in J}\xi_{{\bf A}^n}(R^{n}(V_{j},F^{\bullet}_{ij}))$$
$$\xi_{{\bf A}^n}(\oplus_{j \in J}R^{n}(V_{j},F^{\bullet}_{ij})) \cong \oplus_{j \in J}\xi_{{\bf A}^n}(R^{n}(V_{j},F^{\bullet}_{ij})).$$
\end{lemme}

\begin{preuve}
Ce sont des applications directes de la proposition \ref{tilde2} p.\pageref{tilde2}, (ii) et (iii). 
\end{preuve}

\subsubsection{Les faisceaux de Rees sont r\'eflexifs}

Pour montrer que les faisceaux de Rees sur les espaces affines sont réflexifs nous allons utiliser la caractérisation donnée dans la proposition \ref{caracrefl} p.\pageref{caracrefl}. Il faut donc prouver que le faisceau de Rees $\xi_{{\bf A}^n}( R^{n}(V,F^{\bullet}_{i}))$ est sans torsion (Prop.\ref{caracrefl} p.\pageref{caracrefl},$(i)$) et qu'en tout point
$x \in {\bf A}^n$ l'inégalité 
 $$\text{depth}(  \xi_{{\bf A}^n}( R^{n}(V,F^{\bullet}_{i}))_{x}) \geq 2,$$
est vérifiée (Prop.\ref{caracrefl} p.\pageref{caracrefl},$(ii)$). La profondeur ici est celle du $\calo_{{{\bf A}^n},x}$-module $\xi_{{\bf A}^n}( R^{n}(V,F^{\bullet}_{i}))_{x}$, où $\calo_{{{\bf A}^n},x}$ est un anneau local d'idéal maximal ${\mathfrak m}_x$.

 La première condition $(i)$ est clairement satisfaite car les modules de Rees $R^{n}(V,F^{\bullet}_{i})$ sont des ${k}[u_{1},u_{2},...,u_{n}]$-modules sans torsion. Pour prouver que la deuxième condition est satisfaite nous allons utiliser un lemme algébrique qui donne une condition suffisante sur un $A$-module $M$ pour que sa profondeur soit supérieure ou égale à deux. Ce lemme nécessite l'introduction de quelques notions de cohomologie locale.

Soit $A$ un anneau local régulier de dimension $n$ et $M$ un $A$-module finiment engendré. Soit $\mathfrak a$ un idéal de $A$. On définit alors le sous-module $\Gamma_{\mathfrak a}(M)$ par
$$ \Gamma_{\mathfrak a}(M)=\{ m \in M \vert {\mathfrak a}^{n.m}=0 \text{ pour } n >> 0 \}.$$
On peut alors montrer (cf \cite{gro}) que le foncteur $ \Gamma_{\mathfrak a}(.)$ qui va de la catégorie des $A$-modules dans elle-même est exact à gauche. On note ses foncteurs dérivés à droite dans la catégorie des $A$-modules par $H^{.}_{\mathfrak a}(.)$. On définit ainsi le $i$-ème groupe de cohomologie locale de $M$ en $\mathfrak a$ par
$$  H^{.}_{\mathfrak a}(M).$$

Un $A$-module sans torsion $M$ est un sous-module d'un $A$-module libre $L$.\\

\begin{definition}
Le saturé $M^{sat}$ de $M$ dans $L$ pour l'idéal $\mathfrak a$ de $A$ est le sous $A$-module de $L$ défini par
$$ M^{sat}=\{ l \in L \vert {\mathfrak a}^{k}.l \in M \text{ pour } k >> 0 \}.$$  
\end{definition}

\begin{lemme}\label{sat}
Soit $A$ un anneau local régulier de dimension $n$, d'idéal maximal $\mathfrak m$, et $M$ un $A$-module finiment engendré. Soit $M^{sat}$ le saturé de $M$ pour $\mathfrak m$ par rapport à un $A$-module libre $L$, alors
$$ M^{sat}=M \Longrightarrow \emph{depth}_{\mathfrak m} M \geq 2.$$

\end{lemme}

\begin{preuve}
Considérons la suite exacte de $A$-modules
$$
\xymatrix{
0 \ar[r] & M \ar[r] & L \ar[r] & L/M \ar[r] & 0
}.$$
Elle mène, d'après \cite{gro}, à la suite exacte longue de cohomologie locale en $\mathfrak m$ :
$$
\xymatrix{
0 \ar[r] & H^{0}_{\mathfrak m}(M) \ar[r] & H^{0}_{\mathfrak m}(L) \ar[r] & H^{0}_{\mathfrak m}(L/M) \ar[r] & {}\\
{} \ar[r] & H^{1}_{\mathfrak m}(M) \ar[r] & H^{1}_{\mathfrak m}(L) \ar[r] & H^{1}_{\mathfrak m}(L/M) \ar[r] & ...
}$$   
Comme $L$ est libre, pour tout $i \geq 0$, $ H^{i}_{\mathfrak m}(L)=0$, donc $ H^{0}_{\mathfrak m}(M)=0$ et pour tout $i \geq 1$ :
$$ H^{i}_{\mathfrak m}(M) \cong H^{i-1}_{\mathfrak m}(L/M).$$
Or dans $L/M$, $ M^{sat}/M \cong H^{0}_{\mathfrak m}(L/M)$. Ainsi :
$$  M^{sat}/M=0 \Longleftrightarrow H^{1}_{\mathfrak m}(M) = H^{0}_{\mathfrak m}(L/M)=0.$$

Pour conclure, rappelons le lien entre la profondeur d'un $A$-module et la cohomologie locale. D'après \cite{har}(III, Ex.3.4.), pour un $A$-module $M$ finiment engendré, pour tous $k \geq 0$ les conditions suivantes sont équivalentes : 
\begin{itemize}\parindent=2cm
\item {(i)} $\text{depth}_{\mathfrak m} M \geq k$.
\item {(ii)} $H^{i}_{\mathfrak m}(M)=0$ pour tout $i<k$.
\par\end{itemize}
Ce qui donne le résultat voulu car $ H^{0}_{\mathfrak m}(M)=0$ et $H^{1}_{\mathfrak m}(M) =0$.

\end{preuve}

Nous sommes maintenant en mesure de prouver que les faisceaux de Rees sont réflexifs.\\ 

\begin{proposition}\label{reesreflexif}
Soit $(V,F^{\bullet}_{1},F^{\bullet}_{2},...,F^{\bullet}_{n}) \in \calc_{nfiltr}$ un espace vectoriel muni de $n$ filtrations exhaustives et décroissantes. Le faisceau de Rees sur ${\bf A}^n$, $\xi_{{\bf A}^n}( R^{n}(V,F^{\bullet}_{i}))$ est un faisceau réflexif.\\

\end{proposition}

\begin{preuve} La propriété est vraie pour les faisceaux de Rees sur la droite affine associés à un espace muni d'une seule filtration car ces faisceaux sont sans torsion donc localement libres sur la droite affine. Supposons que ce soit vrai pour $n-1$ filtrations. Alors pour le prouver pour $n$ filtration il suffit de le prouver que le module de Rees est saturé pour l'idéal maximal associé à L'origine. Dans les autres cas, une coordonnée étant non nulle, on se ramène toujours, par localisation, à $n-1$ filtrations. Rappelons que ${ R^{n}(V,F^{\bullet}_{i})}$, le ${k}[u_{1},u_{2},...,u_{n}]$-module de Rees  associé à $(V,F^{\bullet}_{1},F^{\bullet}_{2},...,F^{\bullet}_{n}) $ est un sous-module du module libre $V \otimes {k}[u_{1},u_{1}^{-1},u_{2},u_{2}^{-1},...,u_{n},u_{n}^{-1}]$.
D'après le lemme précédent, il suffit de montrer qu'il est saturé dans $$V \otimes {k}[u_{1},u_{1}^{-1},u_{2},u_{2}^{-1},...,u_{n},u_{n}^{-1}]$$ pour l'ideal maximal ${\mathfrak m}=(u_{1},u_{2},...,u_{n})$.  

${ R^{n}(V,F^{\bullet}_{i})}^{sat}={ R^{n}(V,F^{\bullet}_{i})}$ signifie que pour tout $ l \in  { R^{n}(V,F^{\bullet}_{i})}^{sat}$ tel que ${\mathfrak m}.l \in  { R^{n}(V,F^{\bullet}_{i})}$ alors $l \in  { R^{n}(V,F^{\bullet}_{i})}$ (${\mathfrak m}.l \in  { R^{n}(V,F^{\bullet}_{i})}$ signifie que pour tout $x \in {\mathfrak m},\,\, x.l \in { R^{n}(V,F^{\bullet}_{i})}$).

Soit $x \in V \otimes {k}[u_{1},u_{1}^{-1},u_{2},u_{2}^{-1},...,u_{n},u_{n}^{-1}] $, alors $x= \sum_{(i_{1},i_{2},...,i_{n})}\,u_{1}^{-i_{1}}.,u_{2}^{-i_{2}}...u_{n}^{-i_{n}} \otimes v_{(i_{1},i_{2},...,i_{n})}$. ${\mathfrak m}.x \in { R^{n}(V,F^{\bullet}_{i})}$ signifie que pour tout $k \in [1,n]$ :
$$ u_{k}.x= \sum_{(i_{1},i_{2},...,i_{n})}\,u_{1}^{-i_{1}}.,u_{2}^{-i_{2}}...u_{k}^{-i_{k}+1}...u_{n}^{-i_{n}} \otimes v_{(i_{1},i_{2},...,i_{n})},$$
c'est à dire que pour tout $k \in [1,n]$ $v_{(i_{1},i_{2},...,i_{n})} \in F_{1}^{i_{1}} \cap F_{2}^{i_{2}} \cap ... \cap F_{k}^{i_{k}-1} \cap ...\cap F_{n}^{i_{n}}$, donc $v_{(i_{1},i_{2},...,i_{n})}$ est dans l'intersection de tous ces ensembles qui est  $F_{1}^{i_{1}} \cap F_{2}^{i_{2}} \cap ...\cap F_{n}^{i_{n}}$ et ainsi $x \in R^{n}(V,F^{\bullet}_{i})$.

Pour conclure, suivant le plan anoncé en début de section, on utilise le lemme précédent et la caractérisations des faisceaux réflexifs : ${ R^{n}(V,F^{\bullet}_{i})}^{sat}={ R^{n}(V,F^{\bullet}_{i})}$ implique par le lemme que la profondeur est plus grande que deux ce qui signifie par la proposition que $\xi_{{\bf A}^n}( R^{n}(V,F^{\bullet}_{i}))$ est un faisceau réflexif. 
\end{preuve}

%{\bf Remarque 1 :} Le lemme \ref{sat} p.\pageref{sat} nous donne la longueur maximale d'une résolution libre de $\xi_{{\bf A}^n}( R^{n}(V,F^{\bullet}_{i}))$ sur ${\bf A}^n$. En effet, d'après l'égalité d'Ausslander-Buchsbaum $$\text{pd}( R^{n}(V,F^{\bullet}_{i})) + \text{depth}( R^{n}(V,F^{\bullet}_{i}))=n,$$ donc $\text{pd}( R^{n}(V,F^{\bullet}_{i})) \leq n-2$. On a donc toujours une résolution par des fibrés vectoriels sur ${\bf A}^n$ (donc triviaux) de la forme :
%$$ 0 \rightarrow \calo^{r_{n-2}} \rightarrow  ... \rightarrow  \calo^{r_{2}} \rightarrow  \calo^{r_{1}} \rightarrow \xi_{{\bf A}^n}( R^{n}(V,F^{\bullet}_{i}))\rightarrow 0.$$

{\bf Remarque :} Nous avons vu plus haut qu'un faisceau réflexif sur une variété de dimension $1$ ou $2$ est localement libre. Un corollaire de la proposition précédente est donc que les faisceaux de Rees sur ${\bf A}^1$ et ${\bf A}^2$ sont des fibrés vectoriels. Ceci peut être vu directement. Pour une filtration c'est immédiat car le module de Rees est libre de dimension finie. Comme on l'a vu plus haut, on peut toujours scinder deux filtrations simultanément. Prenons une base compatible avec les deux filtrations $\{ w_{i,j} \}_{(i,j) \in I \times J}$. Cette base donne une description directe des générateurs du module de Rees $R^{2}(F^{\bullet },G^{\bullet })=< u^{-i}.v^{-j} \otimes w_{i,j}>_{(i,j) \in I \times J}$. Ces générateurs sont sans relation et donnent une trivialisation % en effet $ P(u,v).u^{-i}.v^{-j} \otimes v_{i,j}=Q(u,v) u^{-k}.v^{-l} \otimes v_{i,j}$ où $P$ et $Q$ sont des polynômes implique que $v_{i,j} \in F^{i}\cap    
du faisceau localement libre $\xi_{{\bf A}^2}( R^{2}(V,F^{\bullet},G^{\bullet}))$.\\
\begin{center} 
{\bf Etude du fibré de Rees sur ${\bf A}^1$ d'un espace vectoriel filtré }
\end{center}

Par la construction du faisceau de Rees associé à un espace vectoriel filtré on obtient donc un fibré vectoriel sur ${\bf A}^1$. Le lemme suivant nous donne le comportement du fibré restreint à l'ouvert ${\bf G}_{m}$.\\

\begin{lemme}On a l'isomorphisme :
${\xi}_{{\bf A}^1}(V,F^{\bullet}) \vert_{ {\bf G}_{m}} \cong V \otimes {\mathcal O}_{{\bf G}_{m}}$.\\
\end{lemme}

\begin{preuve}
C'est une application directe de la proposition \ref{tilde1} (iii) p.\pageref{tilde1}.
\end{preuve}

\begin{lemme}On a les isomorphismes :\\
${}$\\
$(i)$ ${\xi}_{{\bf A}^1}(V,F^{\bullet})_{1} \cong V$.\\
${}$\\
$(ii)$ ${\xi}_{{\bf A}^1}(V,F^{\bullet})_{0} \cong Gr^{F}V$.\\
\end{lemme}

\begin{preuve}
Ces isomorphismes proviennent de la proposition \ref{tilde1} (iii) p.\pageref{tilde1} et de l'étude des quotients des modules de Rees faite au lemme \ref{lemmedebase}.
\end{preuve}

\begin{proposition}\label{trivialisationderees1}
Le fibré de Rees est isomorphe au fibré trivial :
$${\xi}_{{\bf A}^1}(V,F^{\bullet}) \cong \oplus_{p}(Gr^{F^{\bullet}}_{p}V \otimes_k {\mathcal O}_{{\bf A}^1}(-p.0)) \cong {\mathcal O}_{{\bf A}^1}^{n}$$ où $n=dim_k\,V$.\\ 
\end{proposition}

\begin{preuve}
Ceci est d\^u au fait que l'on peut toujours trouver un scindement de $V$ compatible à la filtration. L'isomorphisme vient du choix de ce scindement $V \cong \oplus_{p}Gr^{F^{\bullet}}_{p}V$. 
\end{preuve}
\begin{center}
{\bf Etude du fibré de Rees sur ${\bf A}^2$ d'un espace vectoriel muni de deux filtrations }
\end{center}

Notons par $i: {\bf G}_{m} \times {\bf G}_{m} \hookrightarrow {\bf A}^2$, $i_{1}: {\bf G}_{m} \times {\bf A}^{1} \hookrightarrow {\bf A}^2$ et $i_{2}: {\bf A}^{1}\times {\bf G}_{m} \hookrightarrow {\bf A}^2$ les différents morphismes d'inclusion. Le faisceau associé à deux filtrations par la construction de Rees est un fibré vectoriel. Le lemme suivant nous donne le comportement du fibré restreint au groupe multiplicatif ${{\bf G}_{m} \times {\bf G}_{m}}$ et montre comment la construction du fibré de Rees associé à deux filtrations généralise la construction du fibré de Rees sur ${\bf A}^1$ associé à une filtration.\\

\begin{lemme}\label{surlesouverts}On a les isomorphismes :\\
${}$
$(i)$ ${\xi}_{{\bf A}^2}(V,F^{\bullet},G^{\bullet}) \vert_{{\bf A}^{1}\times {\bf G}_{m}} \cong i_{2}^{*}{\xi}_{{\bf A}^2}(V,F^{\bullet},Triv^{\bullet}) \cong {\xi}_{{\bf A}^1}(V,F^{\bullet}) \otimes {\mathcal O}_{{\bf G}_{m}}$.\\
${}$\\
$(ii)$ ${\xi}_{{\bf A}^2}(V,F^{\bullet},G^{\bullet})\vert_{{\bf G}_{m} \times {\bf A}^{1}} \cong i_{1}^{*}{\xi}_{{\bf A}^2}(V,Triv^{\bullet},G^{\bullet}) \cong  {\mathcal O}_{{\bf G}_{m}} \otimes {\xi}_{{\bf A}^1}(V,G^{\bullet})$.\\
${}$\\
$(iii)$ ${\xi}_{{\bf A}^2}(V,F^{\bullet},G^{\bullet})\vert_{{\bf G}_{m} \times {\bf G}_{m}} \cong i_{}^{*}{\xi}_{{\bf A}^2}(V,Triv^{\bullet},Triv^{\bullet})=i^{*}( {\mathcal O}_{{\bf A}^2} \otimes V)= {\mathcal O}_{{\bf G}_{m}}^{2} \otimes V$.\\
\end{lemme}

\begin{preuve}
(i) et (ii) découlent de la proposition \ref{tilde1} (iii) p.\pageref{tilde1}, avec $f=u$ et $f=v$ respectivement, et du lemme \ref{lemmedebase}. Pour (iii) on peut utiliser \ref{tilde1} (iv) et le m\^eme lemme.
\end{preuve}

\begin{lemme}On a les isomorphismes :\\
${}$
$(i)$ ${\xi}_{{\bf A}^2}(V,F^{\bullet},G^{\bullet})_{(1,0)} \cong Gr^{G}V$.\\
${}$\\
$(ii)$ ${\xi}_{{\bf A}^2}(V,F^{\bullet},G^{\bullet})_{(0,1)} \cong Gr^{F}V$.\\
${}$\\
$(iii)$ ${\xi}_{{\bf A}^2}(V,F^{\bullet},G^{\bullet})_{(0,0)} \cong Gr^{F}Gr^{G}V$.\\
${}$\\
$(iv)$ ${\xi}_{{\bf A}^2}(V,F^{\bullet},G^{\bullet})_{(1,1)} \cong V$.\\
\end{lemme}

\begin{preuve}
Ces isomorphismes proviennent de la proposition \ref{tilde1} (iii) p.\pageref{tilde1}, et de l'étude des quotients des modules de Rees faite au lemme \ref{lemmedebase} p.\pageref{lemmedebase}.
\end{preuve}
 Notons $D_{1}$ et $D_{2}$ les diviseurs (droites affines) d'équations respectives $u=0$ et $v=0$.\\
\\
\begin{proposition}\label{trivialisationderees}
Le fibré de Rees est isomorphe au fibré trivial :
$${\xi}_{{\bf A}^2}(V,F^{\bullet},G^{\bullet}) \cong \oplus_{p,q}(Gr^{F^{\bullet}}_{p}Gr^{G^{\bullet}}_{q}V \otimes_k {\mathcal O}_{{\bf A}^2}(-pD_{2}-qD_{1})) \cong (\oplus_{p,q}(Gr^{F^{\bullet}}_{p}Gr^{G^{\bullet}}_{q}V) \otimes_k {\mathcal O}_{{\bf A}^2})={\mathcal O}_{{\bf A}^2}^{n}$$ où $n=dim_k\,V$.\\ 
\end{proposition}

\begin{preuve}
Ceci est d\^u au fait que l'on peut toujours trouver un scindement de $V$ compatible à deux filtrations comme on l'a montré précedement. Le premier isomorphisme vient du choix de $V \cong \oplus_{p,q}Gr^{F^{\bullet}}_{p}Gr^{G^{\bullet}}_{q}V$. Exhibons le deuxième sur chacune des composantes de la somme directe : $$\alpha : s \in (Gr^{F^{\bullet}}_{p}Gr^{G^{\bullet}}_{q}V) \otimes_k {\mathcal O}_{{\bf A}^2}(-pD_{2}-qD_{1}) \rightarrow u^{p}.v^{q}.s \in (Gr^{F^{\bullet}}_{p}Gr^{G^{\bullet}}_{q}V) \otimes_k {\mathcal O}_{{\bf A}^2} .$$ 
\end{preuve}

\subsection{Actions de groupes algébriques linéaires et faisceaux de Rees}\label{groupesalg}

\subsubsection{Goupes algébriques linéaires, $G$-variétés et $G$-faisceaux}

Le but de cette partie est de fixer les notations, de définir les objets que l'on utilisera par la suite ainsi que leurs principales propriétés.

On rappelle qu'un groupe algébrique linéaire $G$ sur un corps algébriquement clos $k$ est une variété algébrique équipée d'une structure de groupe et telle que les lois du groupe soient des morphismes de variétés : le produit $\pi : G \times G \rightarrow G , \, (x,y) \mapsto x.y$ et l'application inverse $i : G \rightarrow G, \, x \mapsto x^{-1}$ sont des morphismes algébriques et il y a un élément neutre $e$. Comme le foncteur contravariant $X \mapsto k[X]$, qui à une variété affine $X$ associe sa $k$-algèbre des fonctions régulières, établit une équivalence de catégories entre la catégorie des variétés affines sur le corps $k$ et la catégorie des $k$-algèbres de type fini, toute variété algébrique affine $X$ peut être vue comme l'ensemble des homomorphismes de $k$-algèbres $\text{Hom}_{k}(k[X],k)$. La loi de groupe sur $G$ vient donc des morphismes de $k$-algèbre suivants :
$$\pi^{*}:k[G] \rightarrow k[G] \otimes k[G] , \, i^{*}:k[G] \rightarrow k[G]$$
l'élément neutre $e$ provient d'un morphisme dans $\text{Hom}_{k}(k[G],k)$, $e: k[G] \rightarrow k$. Notons par $m : k [G] \otimes k[G] \rightarrow k[G]$ la multiplication et par $id$ l'homomorphisme identité dans $k[G]$. Les axiomes de la loi de groupe de $G$ sont exprimés par la commutativité des diagrammes suivants :
 
élément neutre :
$$
\xymatrix{
k[G] &  k[G] \otimes k[G] \ar[l]_{e \otimes id}\\
k[G] \otimes k[G] \ar[u]^{id \otimes e} & k[G] \ar[l]^-{\pi^*} \ar[ul]_{id} \ar[u]_{\pi^*}
}
$$

associativité :

$$
\xymatrix{
k[G] \otimes k[G] \otimes k[G] &  k[G] \otimes k[G] \ar[l]_-{\pi^{*} \otimes id}\\
k[G] \otimes k[G] \ar[u]^{id \otimes  \pi^*} & k[G] \ar[l]^{\pi^*}  \ar[u]_-{\pi^*}
}
$$

inverse à gauche et à droite :

$$
\xymatrix{
k[G] \otimes k[G] \ar[d]_{m} &  k[G] \otimes k[G] \ar[l]_-{i^{*} \otimes id}\\
k[G]  & k[G] \ar[l]^{e}  \ar[u]_{\pi^*}
}
\hspace{1cm}\xymatrix{
k[G] \otimes k[G] \ar[d]_{m} &  k[G] \otimes k[G] \ar[l]_-{id \otimes i^{*} }\\
k[G]  & k[G]. \ar[l]^{e}  \ar[u]_{\pi^*}
}
$$
{\bf Exemple :} Les groupes affines qui nous intéresseront par la suite sont les tores $({\bf G}_{m})^{n}$, où $n \in {\bf N}$, produits du groupe affine ${\bf G}_m=\text{Spec}k[u,u^{-1}]$. Les lois de groupe de ${\bf G}_m$ viennent de la commutativité des diagrammes précédents, les homomorphismes $\pi^*$, $i^*$ et $e$ sont donnés par $\pi^{*}u=u \otimes u$, $i^{*}u=u^{-1}$ et $e(u)=1$.\\    

Introduisons la notion de $G$-variété.\\

\begin{definition}Une variété algébrique $X$ est une $G$-variété s'il existe une action à gauche de $G$ sur l'ensemble $X$ $\sigma : G \times X \rightarrow X,\, \sigma :(g,x) \rightarrow g.x$ qui soit un morphisme de variétés et qui vérifie pour tous $g,h \in G, \, x \in X$ : $(gh).x=g.(h.x)$ et $e.x=x$. \\
\end{definition}

Soient $X$ et $Y$ deux $G$-variétés. Un morphisme équivariant entre les $G$-variétés $X$ et $Y$ est un morphisme $f : X \rightarrow Y$ tel que le diagramme suivant commute (on indexe par $X$ resp. $Y$ le morphisme qui donne l'action de $G$ sur $X$ resp. Y) :

$$
\xymatrix{
G \times X  \ar[r]^{id \times f} \ar[d]_{\sigma_X} &  G \times Y   \ar[d]^{\sigma_Y}\\
X  \ar[r]_{f} &  Y
}
$$

Supposons que $X$ soit une $G$-variété affine. Alors d'après l'équivalence de catégories citée plus haut entre la catégorie des variétés affines et la catégorie des $k$-algèbres finiment engendrées, l'action est donnée par un morphisme $\sigma^{*}: k[X] \rightarrow k[G] \otimes k[X]$. La loi du groupe $G$ est toujours notée par $\pi$. Dire que l'action est une action de groupe signifie que les homomorphismes $\pi^{*} \otimes id$ et $id \otimes \pi^{*}$ coïncident de $k[G]\otimes k[X]$ vers $k[G]\otimes k[G] \otimes k[X]$ et que l'homomorphisme $(e \otimes id) \sigma^{*}$ de $k[X]$ vers lui-même est l'identité.\\
\\
{\bf Exemple :} Soit $G={\bf G}_m$. La donnée d'une ${\bf G}_{m}$-action sur une variété affine $X$ est équivalente à la donnée d'une graduation de l'algèbre $k[X]$, i.e. d'une décomposition en somme directe d'espaces vectoriels $k[X]=\oplus_{m \in {\bf Z}} k[X]_{m}$ tels que pour tous $m,n \in {\bf Z}$, $k[X]_{m}.k[X]_{n}\subset k[X]_{m+n}$. L'action de $\sigma$ (ou de façon équivalente l'homomorphisme $\sigma^*$) se déduisent l'un de l'autre par $\sigma^{*}P=u^{m}\otimes P$ si $P \in k[X]_{m}$.\\
${}$\\

 Définissons la notion d'action d'un groupe algébrique $G$ sur un faisceau cohérent sur une $G$-variété compatible avec l'action sur la variété. Soit $\calf$ un faisceau cohérent sur une $G$-variété $X$.\\
\begin{definition}
Le faisceau cohérent $\calf$ est un $G$-faisceau cohérent s'il existe un isomorphisme de faisceaux de $\calo_{G \times X}$-modules :
$$ \Psi : \sigma^{*}\calf \cong p_{2}^{*}\calf,$$ 
où $p_{2}: G \times X \rightarrow X $ est le morphisme projection sur le deuxième facteur, tel que la condition de cocycle suivante soit vérifiée :
$$ (\pi \otimes id_{X})^{*} \Psi = p_{12}^{*} \circ (id_{G} \otimes \sigma )^{*} \Psi,$$
où $p_{12}$: $G \times G \times X \rightarrow G \times X $ est le morphisme projection sur les deux derniers facteurs.

Un $G$-fibré sur une $G$-variété $X$ est un fibré sur $X$ qui est en plus un $G$-faisceau cohérent.\\
\end{definition}
L'isomorphisme $\Psi$ est appelé $G$-linéarisation du faisceau de $\calo_{X}$-modules cohérents $\calf$. 
Suivant \cite{huyleh}, on donne l'interprétation intuitive qui suit. Pour $(g,x) \in G \times X$, on note toujours $g.x$ à la place de
$\sigma(g,x)$. L'isomorphisme de faisceaux de $\calo_{G \times X}$-modules $\Psi $ donne un isomorphismes entre les fibres de $\calf$ :
$$ \Psi_{g,x}: \calf(g.x) \rightarrow \calf(x).$$
Avec cette notation, la condition de cocyclicité se traduit en :
$$\Psi_{g,x}\circ \Psi_{h,g.x}=\Psi_{h.g,x} : \calf(h.g.x) \rightarrow \calf(x).$$
\begin{definition}Un homomorphisme $\Phi : \calf \rightarrow \calf'$ entre faisceaux de $\calo_{X}$-modules $G$-linéarisés est un homomorphisme de $\calo_{X}$-modules qui commute avec les linéarisations respectives de $\calf$ et $\calf'$, i.e. $\Psi$ et $\Psi'$ ; ceci signifie que
$$\Psi'\circ \sigma^{*}\Phi=p_{2}^{*}\Phi \circ \Psi.$$

\end{definition}

\begin{lemme}\label{gcatab} \cite{thoma}
La catégorie des $G$-faisceaux cohérents est abélienne.
\end{lemme}

\begin{preuve}
Cette catégorie est clairement additive, montrons qu'elle est munie de noyaux et de conoyaux. Soit $f : \calf \rightarrow \calg$ un morphisme de $G$-faisceaux cohérents. Notons $\calk er$ et $\calc oker$ le noyau et conoyau de $f$. La projection $p_{2}$ est un morphisme plat. Le morphisme $\sigma $ est la composée de l'isomorphisme de $G \times X$ donné par la flèche $(g,x) \mapsto (g,g.x)$ et de la projection plate $p_{2}$ et est donc plat. Ainsi les foncteurs $p_{2}^{*}$ et $\sigma^*$ qui vont de la catégorie des faisceaux cohérents sur $X$ vers la catégorie des faisceaux cohérents sur $G \times X$ sont exacts. On en déduit donc le diagramme commutatif ci-dessous, ce qui permet, par le lemme des cinq, de construire les isomorphismes de faisceaux de $\calo_{G \times X}$-modules (les flèches en pointillés) $ \Psi_{\calk er} : \sigma^{*}\calk er \cong p_{2}^{*}\calk er$ et $ \Psi_{\calc oker} : \sigma^{*}\calc oker \cong p_{2}^{*}\calc oker$. La condition de cocycle se vérifie alors immédiatement.

$$\xymatrix{ 0 \ar[r]&  \sigma^{*}(\calk er)  \ar[r]^{\sigma^{*}i} \ar@{.>}[d]^{\Psi_{\calk er}}_{\cong}&  \sigma^{*}(\calf) \ar[r]^{\sigma^{*}f} \ar[d]_{\cong}^{\Psi_{\calf}} & \sigma^{*}(\calg)  \ar[r]^{\sigma^{*}\pi} \ar[d]^{\Psi_{\calg} }_{\cong}&  \sigma^{*}(\calc oker) \ar[r] \ar@{.>}[d]^{\Psi_{\calc oker}}_{\cong}& 0\\
 0 \ar[r]&  p_{2}^{*}(\calk er ) \ar[r]^{p_{2}^{*}i}  &  p_{2}^{*}(\cale) \ar[r]^{p_{2}^{*}f}   & p_{2}^{*}(\calg)  \ar[r]^{p_{2}^{*}\pi}  &  p_{2}^{*}(\calc oker) \ar[r] & 0
}
.$$
L'image $\cali m$, qui est un noyau pour le conoyau, est donc un $G$-faisceau cohérent, de même pour la coïmage $\calc oim$ qui est un conoyau pour le noyau donc un $G$-faisceaux cohérent. La catégorie des faisceaux cohérents sur $X$ est abélienne, donc l'image et la coïmage sont isomorphes, notons par $\Phi$ l'isomorphisme désigné. Comme pour la commutativité du diagramme précédent on a la commutativité de 
$$\xymatrix{ \sigma^{*}(\calc oim) \ar[r]^{\sigma^{*}\Phi}_{\cong} \ar[d]^{\Psi}_{\cong} &  \sigma^{*}(\cali m) \ar[d]^{\Psi'}_{\cong}\\
p_{2}^{*}(\calc oim) \ar[r]^{p_{2}^{*}\Phi}_{\cong} & p_{2}^{*}(\cali m),
}$$
ce qui prouve que $\Phi$ est un morphisme de $G$-faisceaux cohérents.  
\end{preuve}
 
Pour généraliser le théorème \ref{muab} p.\pageref{muab} aux $G$-faisceaux réflexifs $\mu$-semistables de pente $\mu$, nous aurons besoin de savoir si le faisceaux réflexif $\calf^{**}$ associé à un $G$-faisceau $\calf$ est un $G$-faisceau. Le lemme suivant répond par l'affirmative à cette question.\\

Introduisons d'abord la notion d'action plate :

\begin{definition}
Une action $\sigma : G \times X \rightarrow X$ est dite plate si le morphisme de schémas correspondant est plat.
\end{definition}

\begin{lemme}\label{grefl}Soit $G$ un groupe algébrique linéaire sur $k$ et $X$ une $G$-variété lisse et séparée de type fini sur $k$ tels que l'action soit plate.
%\begin{itemize}\parindent=1cm
%\item{(i)}
Soit $\calf$ un $G$-faisceau sur $X$, alors le faisceau réflexif associé $\calf^{**}$ est un $G$-faisceau et le morphisme canonique $\nu$ du faisceau $\calf$ vers son bidual est un morphisme de $G$-faisceaux.
%\item{(ii)} Si $f :\cale \rightarrow \calf$ est un morphisme de $G$-faisceaux réflexifs qui est un isomorphisme en dehors d'un ensemble de codimension $2$, alors c'est un isomorphisme de $G$-faisceaux.\\
%\par\end{itemize}
\end{lemme} 

\begin{preuve} Nous allons montrer que le dual $\calf^{*}$ d'un $G$-faisceau de $\calo_{X}$-modules  coh\'erent $\calf$ est un $G$-faisceaux coh\'erent.

$p_{2}:G\times_{k}X \rightarrow X$ est plate. $\sigma$ est plate comme composée de $p_{2}$ avec l'isomorphisme $G \times_{k}X \cong G\times_{k}X$ donn\'e par $(g,x) \mapsto (g,gx)$. Les foncteurs $p_{2}^{*}$ et $\sigma^{*}$ sont donc exacts. Donc, d'apr\`es \cite{har}, proposition 1.8, on a $p_{2}^{*}(\calf^{*})\cong(p_{2}^{*}\calf)^{*}$ et $\sigma^{*}(\calf^{*})\cong(\sigma^{*}\calf)^{*}$. Ainsi l'isomorphisme image de $\Psi_{\calf}$ dans $\text{Hom}((p_{2}^{*}\calf)^{*},(\sigma^{*}\calf)^{*})$ se transporte en isomorphisme dans $\text{Hom}(p_{2}^{*}(\calf^{*}),\sigma^{*}(\calf^{*}))$ ce qui permet de conclure.

On vient de voir que $\calf^{**}$ est un $G$-faisceau. Soit $\Psi_{\calf}^{**}$ le morphisme entre faisceaux r\'eflexif associ\'e à $\Psi_{\calf}$. Alors $\Psi_{\calf}^{**}\circ \sigma^{*}\nu=p_{2}^{*}\nu\circ \Psi_{\calf}$ et $(\sigma^{*}\calf)^{**}\cong \sigma^{*}(\calf^{**})$ et $(p_{2}^{*}\calf)^{**}\cong p_{2}^{*}(\calf^{**})$ permettent de conclure. 

\end{preuve}

Les lemmes précédents, composés avec le théorème $1$, nous permettent de prouver que la catégorie des faisceaux $G$-équivariants $\mu$-semistables réflexifs sur une $G$ variété $X$ est abélienne. Rappelons que l'on s'est placé sur un corps de caractéristique nulle algébriquement clos. \\

\begin{corollaire}
Soit $G$ un groupe algébrique linéaire sur $k$ et $X$ une $G$-variété lisse et séparée de type fini sur $k$. Alors la catégorie $\emph{Refl}_{\mu-semistable,\mu}(X /G)$ des faisceaux réflexifs sur $X$ $G$-équivariants $\mu$-semistables de pente $\mu$ et des morphismes $G$-équivariants entre ces objets est abélienne.  
\end{corollaire}

\begin{preuve} Soit $f:\cale \rightarrow \calf$ un morphisme dans $\text{Refl}_{\mu-semistable,\mu}(X /G)$. Le théorème \ref{muab} p.\pageref{muab} nous donne l'existence d'un noyau $\calk er $, d'un conoyau $\calc oker$, puis d'une image $\cali m$ et d'une coïmage $\calc oim$ de $f$ qui sont $\mu $-semistables de pente voulue, et tels que l'image et la coïmage soient isomorphes, i.e., tels qu'il exite un isomorphisme $\xymatrix{ g :\calc oim \ar[r]^{\sim}  & \cali m }$. Le lemme \ref{gcatab} p.\pageref{gcatab} affirme que chacun des faisceaux $\calk er $, $\calc oker$, $\cali m$ et $\calc oim$ est un $G$-faisceaux. Il reste donc à vérifier, que l'isomorphisme $g$ dans la catégorie des faisceaux réflexifs $\mu$-semistables de pente $\mu$ sur $X$, $\text{Refl}_{\mu-semistable,\mu}(X)$,
 entre objets de la sous-catégorie des objets équivariants munie des morphismes équivariants $\text{Refl}_{\mu-semistable,\mu}(X /G)$ est $G$-équivariant. Pour cel\`a revenons à la preuve du théorème \ref{muab} p.\pageref{muab}. Pour montrer qu'on avait un isomorphisme $g$ entre la coïmage et l'image, on montrait d'une part que l'on avait un isomorphisme en dehors d'un ensemble de codimension $2$ entre le conoyau faisceautique du noyau noté $\text{Coker}(i)$ et son bidual $\calc oim=(\text{Coker}(i))^{**}$, d'autre part que l'on avait un isomorphisme en dehors d'un ensemble de codimension $2$ entre l'image faisceautique $\text{Im}(f)$ et l'image dans la catégorie $\text{Refl}_{\mu-semistable,\mu}(X)$, $\cali m=\text{Ker}(\pi^{**})$. On en tirait, par composition, les images et coïmages faisceautiques étant isomorphes, un isomorphisme en dehors d'un ensemble de codimension $2$ entre faisceaux réflexifs, donc un isomorphisme entre faisceaux réflexifs. Par le lemme \ref{gcatab} p.\pageref{gcatab}, l'isomorphisme entre $\text{Coim}(f)$ et $\text{Im}(f)$ est un morphisme de $G$-faisceaux ; par le lemme \ref{grefl} p.\pageref{grefl}, le morphisme entre $\text{Coim}$ et $\calc oim$ (son faisceau réflexif associé) est un morphisme de $G$-faisceaux ; le diagramme de la preuve du théorème \ref{muab} suivant
$$
\xymatrix{
0 \ar[r] &\text{Im}(f) \ar[r]\ar[d] &  \calf \ar[r]^{\pi} \ar[d]^{id} & \text{Coker}(f) \ar[r] \ar[d]^{\nu} & 0 \\
0 \ar[r] & \cali m=\text{Ker}(\pi^{**})  \ar[r]  &  \calf \ar[r]^-{\pi^{**}} & \calc oker =(\text{Coker}(f))^{**} \ar[r]  & 0,
}
$$   
montre que le morphisme entre $\text{Im}(f)$ et $\cali m$ est un morphisme de $G$-faisceaux ($\nu$ est un morphisme de $G$-faisceaux par le lemme \ref{grefl} et les morphismes horizontaux sont des morphismes de $G$-faisceaux par le lemme \ref{gcatab}). On en déduit que $g$ est un morphisme de $G$-faisceaux réflexifs qui est un isomorphisme, et donc que $\calc oim$ et $\cali m$ sont des $G$-faisceaux isomorphes, ce qui permet de conclure.

\end{preuve}

Comme on l'a vu plus haut, les faisceaux réflexifs sur une variété de dimension $\leq 2$ sont localement libres, ainsi,\\

\begin{corollaire}\label{gbunab} 
Soit $X$ une $G$-variété du type consid\'er\'e dans le corollaire précédent, avec de plus $\emph{dim}(X) \leq 2$. Alors la catégorie $\text{Bun}_{\mu-semistables,\mu}(X/G)$ des $G$-fibrés vectoriels $\mu$-semistables de pente $\mu$ et des morphismes $G$-équivariants entre ces objets est abélienne.\\ 
\end{corollaire}

 %\begin{proposition}equiv de cat wag+\\

%\end{proposition}

%\begin{preuve}

%\end{preuve}

Terminons cette section en définissant une notion qui nous sera utile pour décrire une certaine classe de fibrés vectoriels sur ${\bf A}^3 \backslash \{0,0,0\}$, la classe de ceux qui peuvent être obtenus comme images inverses de fibrés sur le plan projectif ${\bf P}^2$.\\  

\begin{definition}Soit $G$ un groupe algébrique et $f : X \rightarrow Y$ un morphisme $G$-équivariant de $G$-variétés. On dit que $f$ est $G$-fibré principal s'il existe un morphisme étale surjectif $Y'\rightarrow Y$ et un isomorphisme $G$-équivariant de $G$-variétés $Y'\times G \rightarrow Y'\times_{Y} X$ ce qui signifie que $X$ est localement pour la topologie étale isomorphe comme $G$-variété au produit $Y \times G$.\\ 
\end{definition}
Si $X \rightarrow Y$ est un quotient géométrique et que $f$ est plate et si le morphisme $(\sigma_{X},p_{1}): X \times G \rightarrow X \times_{Y} X$ est un isomorphisme, alors $f$ est un $G$-fibré principal.\\

Nous utiliserons cette définition dans un cadre restreint :

 Soit l'action de ${\bf G}_m$ sur ${\bf A}^{n} \backslash \{0,0,...,0\}$ donnée par $$\sigma : (t,(u_{1},u_{2},...,u_{n})) \mapsto (t.u_{1},t.u_{2},...,t.u_{n}).$$

Cette action a un quotient géométrique qui est l'espace projectif ${\bf P}^{n-1}$. Le morphisme $f$ est bien plat et $(\sigma_{},p_{1}): ({\bf A}^{n} \backslash \{0,0,...,O\}) \times {\bf G}_m \rightarrow ({\bf A}^{n} \backslash \{0,0,...,0\})\times_{{\bf P}^{n-1}}({\bf A}^{n} \backslash \{0,0,...,0\})$ est bien un isomorphisme, donc $f : {\bf A}^{n} \rightarrow {\bf P}^{n-1}$ est bien un ${\bf G}_m$-fibré principal. \\  

\begin{definition}
Soit $f : X \rightarrow Y$ un morphisme de $G$-variétés qui soit un quotient géométrique et soit $\calf $ un $G$-faisceau cohérent. On dit alors que $\calf$ descend à $Y$ s'il existe un faisceau cohérent $\cale$ sur $Y$ tel que l'on ait sur $X$ l'isomorphisme de $G$-faisceaux $\calf \cong f^{*}\cale$.\\
\end{definition} 

La définition des $G$-fibrés principaux sera justifiée ici par le fait suivant :\\
$\,$\\
{\bf Fait :} (\cite{huyleh}, Th.4.2.14 p.87) Soit $G$ un groupe réductif. Si $f : X \rightarrow Y$ est un $G$-fibré principal et $\calf $ un $G$-faisceau, alors $\calf $ descend.\\

Dans les applications par la suite, $G$ sera toujours réductif car on aura toujours $G=({{\bf G}_m})^{n}$, où $n \in {\bf N}$. En particulier, le fait précédent sera utilisé sous la forme suivante : si l'on a un ${\bf G}_m$-faisceau $\calf$ sur ${\bf A}^{n} \backslash \{0,0,...,0\}$, alors ce faisceau descend sur ${\bf P}^{n-1}$ i.e. il existe un faisceau $\cale$ sur ${\bf P}^{n-1}$ tel que $f^{*}\cale \cong_{} \calf$ soit un isomorphisme de ${\bf G}_{m}$-faisceaux. 

%\\
%??????????caractèers
%\\

\subsubsection{Action du groupe multiplicatif $({{\bf G}_m})^n$ sur les faisceaux de Rees associés à $n$ filtrations}
Les faisceaux cohérents de Rees sur ${\bf A}^n$ sont munis d'une structure plus riche par l'action du groupe multiplicatif ${{\bf G}_m}^n$ qui prolonge l'action par translation sur ${\bf A}^n$. 

Notons ${\bf T}^{n}={{\bf G}_m}^{n}$. On a 
\begin{center}
${\bf T}^{n}=\text{Spec}\,k[t_{1},t_{1}^{-1},t_{2},t_{2}^{-1},...,t_{n},t_{n}^{-1}]$ et ${\bf A}^{n}=\text{Spec}\,k[u_{1},u_{2},...,u_{n}]$.
\end{center}
 L'action par translation est donnée par :

\begin{center}
$\sigma : {\bf T}^{n} \times {\bf A}^{n} \rightarrow {\bf A }^{n}$\\
$(t_{1},t_{2},...,t_{n})\times(u_{1},u_{2},...,u_{n}) \mapsto (t_{1}.u_{1},t_{2}.u_{2},...,t_{n}.u_{n})$
\end{center}
Cette action est duale du morphisme d'anneaux (de la coaction)
\begin{center}
$\sigma^{\sharp} :k[u_{1},u_{2},...,u_{n}] \rightarrow  k[u_{1},u_{2},...,u_{n}] \otimes k[t_{1},t_{1}^{-1},t_{2},t_{2}^{-1},...,t_{n},t_{n}^{-1}]$\\
$u_{i} \mapsto u_{i}.t_{i} $
\end{center}
On a aussi le morphisme de projection $p_{2} : {\bf T}^{n} \times {\bf A}^{n} \rightarrow {\bf A }^{n}$ dual du morphisme d'anneaux $p_{2}^{\sharp}:k[u_{1},u_{2},...,u_{n}] \rightarrow  k[u_{1},u_{2},...,u_{n}] \otimes k[t_{1},t_{1}^{-1},t_{2},t_{2}^{-1},...,t_{n},t_{n}^{-1}]$ donné par $p_{2}^{\sharp}: u_{i} \mapsto u_{i} \otimes 1$.\\

\begin{proposition}\label{reesgfib}
Soit $V$ un espace vectoriel sur $k$ muni de $n$ filtrations exhaustives et décroissantes $(F_{i}^{\bullet })_{i \in [1,n]}$. Le fibré de Rees sur ${\bf A}^n$ associé $\xi_{{\bf A}^n}(V,F_{i})$ est un ${\bf T}^n$ fibré pour l'action de ${\bf T}^n$ par translation sur ${\bf A}^n$. 
\end{proposition}

\begin{preuve}
Pour prouver cette assertion, il faut exhiber un isomorphisme entre les faisceaux de $\calo_{{\bf T}^{n} \times {\bf A}^{n}}$-modules $\sigma^{*}\xi_{{\bf A}^n}(V,F_{i})$ et $p_{2}^{*}\xi_{{\bf A}^n}(V,F_{i})$. Pour simplifier l'écriture, on notera $A=k[u_{1},u_{2},...,u_{n}]$ et $T=k[t_{1},t_{1}^{-1},t_{2},t_{2}^{-1},...,t_{n},t_{n}^{-1}]$. On a défini deux morphismes d'anneaux de $A$ vers $A \otimes T$, $p_{2}^{\sharp}$ et $\sigma^{\sharp}$. 

Comme $\xi_{{\bf A}^n}(V,F_{i})={R^{n}(V,F_{i}^{\bullet})}^{\sim}$ ($R^{n}(V,F_{i}^{\bullet})$ est un $A$-module), on a d'après la proposition \ref{tilde2} p.\pageref{tilde2}, 
$$p_{2}^{*}\xi_{{\bf A}^n}(V,F_{i})=({R^{n}(V,F_{i}^{\bullet}) \otimes_{(A,\sigma^{\sharp})}(A\otimes T)})^{\sim} $$
 et, $$\sigma^{*}\xi_{{\bf A}^n}(V,F_{i})=({R^{n}(V,F_{i}^{\bullet}) \otimes_{(A,p_{2}^{\sharp})}(A\otimes T)})^{\sim}.$$ On a ainsi un isomorphisme de $A \otimes G$-modules
$$\Psi^{\sharp}: R^{n}(V,F_{i}^{\bullet}) \otimes_{(A,\sigma_{2}^{\sharp})}(A\otimes T) \rightarrow R^{n}(V,F_{i}^{\bullet}) \otimes_{(A,p_{2}^{\sharp})}(A\otimes T).$$
Rappelons que le foncteur $\sim$ qui à un $A$-module $M$ associe un faisceau $M^{\sim}$ sur $\text{Spec}\,A$ établit une équivalence de catégories, on a donc\\
 
$\text{Hom}_{A\otimes G}( R^{n}(V,F_{i}^{\bullet}) \otimes_{(A,\sigma_{2}^{\sharp})}(A\otimes T) , R^{n}(V,F_{i}^{\bullet}) \otimes_{(A,p_{2}^{\sharp})}(A\otimes T)=$\\
\hspace*{6cm}$\text{Hom}_{\calo_{\text{Spec}\,A \times \text{Spec}\,G}}(\sigma^{*}\xi_{{\bf A}^n}(V,F_{i}),p_{2}^{*}\xi_{{\bf A}^n}(V,F_{i})).$\\

D'où l'isomorphisme de faisceaux de $\calo_{{\bf T}^{n} \times {\bf A}^{n}}$-modules $\Psi$ à droite associé à $\Psi^{\sharp}$ dans le membre de gauche.

\end{preuve}

\begin{lemme}
Les isomorphismes de fibrés de Rees sur ${\bf A}^n$ donnés dans le lemme \ref{isosumfibrees} p.\pageref{isosumfibrees}, sont compatibles avec l'action de ${\bf T}^n$ i.e. sont des isomorphismes de ${\bf T}^{n}$-faisceaux de $\calo_{{\bf A}^n}$-modules.\\
\end{lemme}

\begin{preuve}
Il suffit de remarquer en utilisant la proposition \ref{tilde2} p.\pageref{tilde2} que pour un morphisme d'anneaux $f^{\sharp}:A \rightarrow B$ et une famille de $A$-modules $\{M_{j} \}$ indexée par $j \in J$, on a, si $f: \text{Spec}\,B \rightarrow \text{Spec}\,A$ est le morphisme associé à $f^\sharp$, $f^{*}((\otimes_{j}M_{j})^{\sim})=((\otimes_{j}M_{j})\otimes_{A,f}B)^{\sim}=(\otimes_{j}(M_{j}\otimes_{A,f}B))^{\sim}=\otimes_{j}(M_{j}\otimes_{A,f}B)^{\sim}=\otimes_{j}f^{*}(M_{j}^{\sim})$. On montre de même que $f^{*}((\oplus_{j}M_{j})^{\sim})=\oplus_{j}f^{*}(M_{j}^{\sim})$. On conclut en appliquant ceci à $f^{\sharp}=\sigma^{\sharp}$ et $f^{\sharp}=p_{2}^{\sharp}$.
\end{preuve}
${}$\\
{\bf Remarque :} Une autre fa{\c c}on de pr\'esenter les actions sur les faisceaux de ${\calo}_{X}$-modules est d'utiliser la notion de comodule, cf.e\cite{saa}.

\subsubsection{Construction inverse}
Dans cette section nous donnons un premier outil en vue d'\'etablir d'un dictionnaire entre objets filtrés et fibrés vectoriels. Pour ceci, on asseoit une correspondance entre les espaces vectoriels filtrés et bifiltrés, qui sont des objets algébriques, et les fibrés de Rees sur la droite ou le plan affine, i.e. des objets géométriques. Cette correspondance fonctionne bien puisqu'on exhibe une équivalence de catégories entre une catégorie dont les objets sont des espaces filtrés, les morphismes sont les morphismes strictement compatibles avec les filtrations et une catégorie dont les objets sont des fibrés vectoriels équivariants pour l'action d'un certain groupe sur les espaces affines avec des morphismes que nous spécifierons plus bas. 

Rappelons que deux foncteurs $F: \calc \rightarrow \calc'$ et $G: \calc' \rightarrow \calc$ définissent une équivalence de catégories entre les catégories $\calc$ et $\calc'$ si $G\circ F$ est naturellement équivalent à l'identité de $\calc $ et $F\circ G$ est naturellement équivalent à l'identité de $\calc' $. Rappelons de plus que pour qu'un foncteur entre catégories $F: \calc \rightarrow \calc'$ définisse une équivalence de catégories entre les catégories $\calc$ et $\calc'$, il faut et il suffit que $F$ soit essentiellement surjectif et pleinement fidèle. Essentiellement surjectif signifie que pour tout objet $C'\in \calc'$ il existe un objet $C \in  \calc$ tel que $F(C) \cong C'$. Pleinement fidèle signifie que pour toute paire d'objets $(A,B)$ de $\calc$ l'application définie par $F$ induit une bijection de $\text{Hom}_{\calc}(A,B)$ vers $\text{Hom}_{\calc'}(F(A),F(B))$.  

Notons par exemple, et puisque cet exemple sera utilisé pour la preuve de nos équivalences, que, un anneau $A$ quelconque (resp. noetherien) étant fixé, le foncteur $\sim$ d\'ecrit plus haut qui à un $A$-module $M$ associe le faisceau $\tilde M$ sur $\text{Spec}(A)$ établit une équivalence de catégories entre la catégorie des $A$-modules (resp. $A$-modules finiment engendrés) et la catégorie des $\calo_{X}$-modules quasi-cohérents (resp. cohérents) (cf \cite{har}, Cor.5.5, p. 113).   

On désignera par $\Phi_{R}$ le foncteur qui à un espace vectoriel muni de $n$ filtrations\\ $(V,F^{\bullet}_{1},F^{\bullet}_{2},...,F^{\bullet}_{n}) \in \calc_{nfiltr}$ associe le faisceau cohérent sur ${\bf A}^{n}$, $$\xi_{{\bf A}^n}( R^{n}(V,F^{\bullet}_{i}))={\tilde R^{n}(V,F^{\bullet}_{i})}=\Phi_{R}((V,F^{\bullet}_{1},F^{\bullet}_{2},...,F^{\bullet}_{n})).$$ 

Nous allons montrer que l'on peut définir un foncteur inverse , que l'on notera $\Phi_{I}$, qui définit avec $\Phi_{R}$ une équivalence de catégorie dans les cas $n=1,2$.\\ 

\begin{center}
{\bf Construction inverse de la construction du fibré de Rees associée à un espace vectoriel filtré }
\end{center}
Soit $(V,F^{\bullet })$ un espace vectoriel filtré. Comme on l'a vu plus haut, $\xi_{{\bf A}^1}(V,F^{\bullet })$ est un ${\bf G}_m$-fibré vectoriel sur ${\bf A}^1$ pour l'action qui prolonge l'action de ${\bf G}_m$ sur ${\bf A}^1$ par translation. L'image de $\Phi_R$ est donc bien dans la catégorie voulue. 

Réciproquement, soit $\calf$ un faisceau localement libre (i.e. un fibré vectoriel) sur ${\bf A}^1$ muni d'une action de ${\bf G}_m$ prolongeant l'action par translation. Soit $V=\calf_1$, la fibre de ce faisceau au point $1 \in {\bf A}^1$. ${\bf G}_m$ est abélien, ainsi on peut décomposer son action suivant ses caractères, et donc, en regardant l'ordre des zéros des sections invariantes, avoir une graduation et donc une filtration $F^{\bullet }$ de la fibre au dessus de $1$, $V$. Exprimons ceci en terme de $k[u]$-modules. D'après l'équivalence donnée dans l'introduction de cette section, $\calf$ correspond à un $k[u]$-module libre $B=\calf({\bf A}^{1})$. La restriction de l'action par translation à ${\bf G}_{m} \subset {\bf A}^1$ étant transitive, elle fournit une trivialisation de $\calf \vert_{{\bf G}_m}$ qui en terme de modules donne l'isomorphisme :
$$B \otimes_{k[u]}k[u,u^{-1}] \cong W \otimes_{k} k[u,u^{-1}] . $$
Prenons le quotient à gauche par l'idéal de $k[u,u^{-1}]$, $(u-1)$, on obtient $B/(u-1)B$. A droite on obtient $W$. Cela donne une identification canonique entre l'espace vectoriel $W$ et la fibre en $1$ de $\calf$ (par $B/(u-1)B$, proposition \ref{tilde1}). D'où l'inclusion :
$$ B \subset V \otimes_{k}k[u,u^{-1}].$$
On obtient ainsi une filtration de $V$ décroissante et exhaustive en posant $$F^{p}V=\{ v \in V , \, u^{-p}v \in B \}.$$ 

On pose alors :
$$\Phi_{I}(\calf)=(V,F^{\bullet }).$$ 

\begin{proposition}\cite{sim1}\label{inverse1}
Les foncteurs $\Phi_R$ et $\Phi_I$ établissent une équivalence de catégories entre la catégorie des espaces vectoriels de dimension finie munis d'une filtration exhaustive et décroissante et dont les morphismes sont les morphismes strictement compatibles aux filtrations et la catégorie des fibrés vectoriels sur ${{\bf A}^{1}}$ munis de l'action de ${{\bf G}_{m}} $ prolongeant l'action standard dont les morphismes sont les morphismes de fibrés vectoriels dont les conoyaux sont sans torsion donc des faisceaux localement libres sur la droite affine et qui commutent avec l'action :
$$
\xymatrix{
\{ \calc_{1filtr} \} \ar@<2pt>[r]^-{\Phi_{R}} &  \{  {\text Bun}({\bf A}^{1}/{\bf G}_{m}) \} \ar@<2pt>[l]^-{\Phi_{I}}
}
.$$

\end{proposition}

\begin{preuve} $\Phi_I$ établit une construction inverse de la construction du fibré de Rees $\Phi_R$. En effet, si l'on part d'un espace vectoriel filtré $(V,F^{\bullet })$ et que l'on applique $\Phi_I$ à $\xi_{{\bf A}^1}(V,F^{\bullet })=\Phi_{R}((V,F^{\bullet }))$ on obtient la même filtration sur le même espace vectoriel $V$ (la fibre en $1$ de $\xi_{{\bf A}^1}(V,F^{\bullet })$). 

Réciproquement partons d'un fibré vectoriel $\calf$ muni de l'action du groupe multiplicatif. Par la construction $\Phi_I$ faite au-dessus, on obtient un espace vectoriel filtré $(V,F^{\bullet })$ associé au module de Rees $B$ définit par $B=\calf ({\bf A}^{1})$. Soit $B'$ le module de Rees associé à cet espace vectoriel filtré : $B'=\Phi_{R}(V,F^{\bullet })$. $F^{p}$ a été défini comme l'ensemble des vecteurs $v \in V$ tels que $u^{-p}v \in B$, donc $B'\subset B$. Réciproquement, soit $u^{-p}.v \in B$, on a $v \in F^{p}V$ et donc d'après la définition de $B'$, $u^{-p}.v \in B'$ donc $B=B'$. 

Reste à montrer que la correspondance est pleinement fidèle. Considérons l'application 
\begin{center} $\text{Hom}_{ \calc_{1filtr} }((V,F^{\bullet}),(W,G^{\bullet})) \rightarrow \text{Hom}_{  {\text Bun}({\bf A}^{1}/{\bf G}_{m}) }(\Phi_{R}(V,F^{\bullet}),\Phi_{R}(W,G^{\bullet}))$ \end{center}
définie comme il suit : à tout morphisme $f$ dans $\calc_{1filtr}$ on associe un morphisme de modules de Rees qui à $u^{-p}.v$ associe $u^{-p}.f(v)$. La compabilité du morphisme d'objets filtrés $f$ assure que $u^{-p}.f(v)$ est bien un élément du module de Rees associé à $(W,G^{\bullet})$ et donc l'existence du morphisme de module. L'image de $f$ dans $\text{Hom}_{  {\text Bun}({\bf A}^{1}/{\bf G}_{m}) }(\Phi_{R}(V,F^{\bullet}),\Phi_{R}(W,G^{\bullet}))$ est alors le morphisme image du morphisme de $k[u]$-modules par la construction fonctorielle des faisceaux cohérents associés aux modules.\\   

Pour montrer que la correspondance est pleinement fidèle montrons que cette application est bijective. Le terme de source de l'application concerne l'ensemble des morphismes compatibles d'espace vectoriels filtrés et le terme de droite l'ensemble des morphismes ${\bf G}_m$-équivariants de ${\bf G }_m$-fibrés. Commençons par montrer la surjectivité de cette application. Soit \\
$f \in \text{Hom}_{  {\text Bun}({\bf A}^{1}/{\bf G}_{m}) }(\Phi_{R}(V,F^{\bullet}),\Phi_{R}(W,G^{\bullet}))$. D'après la proposition \ref{tilde2}, pour un anneau $A$, l'application $M \rightarrow M^{\sim}$ qui va des $A$-modules vers les faisceaux cohérents sur $\text{Spec} \,A$ est un foncteur pleinement fidèle. On déduit de $f$ un morphisme de $k [u]$-modules $f' \in \text{Hom}_{k[u]}(B_{F},B_{G})$ où $B_{F}$ et $B_{G}$ sont les modules donnés par les sections globales des deux fibrés $\Phi_{R}(V,F^{\bullet})({\bf A}^{1})$ et $\Phi_{R}(W,G^{\bullet})({\bf A}^{1})$. Ces modules sont munis de l'action de ${\bf G}_m$ héritée de l'action sur les fibrés et $f'$ définit un morphisme d'espaces vectoriels qui respecte les filtration $f''$. En effet, soit $v \in F^{p}(V)$ alors par définition $u^{-p}.v \in B_{F}$. Alors, comme $f'$ respecte l'action $f'(u^{-p}.v)=u^{-p}.P(u).w$ où $w \in W$ et $P \in k[u]$. Donc $f''(v)=w$ et $f''(F^{p}(V)) \subset G^{p}(W) $ car $u^{-p}.P(u).w \in B_{G} \Longrightarrow w \in G^{p}$. Pour voir que $f''$ est strictement compatible aux filtrations, on exprime le fait que $f$ est un morphisme de ${\bf G}_m$-fibré i.e. que le noyau, l'image de $f$ sont des fibrés. On en déduit un isomorphisme de ${\bf G}_m$-fibrés entre la coïmage et l'image de $f$ qui descend aux espaces vectoriels filtrés. L'injectivité est claire à chacune des étapes : de  $\text{Hom}_{\calc_{1filtr} }((V,F^{\bullet}),(W,G^{\bullet}))$ vers $\text{Hom}_{ k[u] }(R(V,F^{\bullet}),R(W,G^{\bullet}))$ puis vers $\text{Hom}_{  {\text Bun}({\bf A}^{1}/{\bf G}_{m}) }(\Phi_{R}(V,F^{\bullet}),\Phi_{R}(W,G^{\bullet}))$. On arrive bien dans les morphismes de fibrés vectoriels car la stricte compatibilité des morphismes d'espaces filtrés fait que les images et noyaux des morphismes induits sur les $k[u]$-modules sont sans torsion, les images et noyaux des morphismes dans ${\text Bun}({\bf A}^{1}/{\bf G}_{m})$ sont donc des sous-faisceaux cohérents sans torsion, donc singuliers sur des ensembles de codimension au moins $2$, donc localement libres sur la droite affine.

\end{preuve}

\begin{center}
{\bf Construction inverse de la construction du fibré de Rees associée à un espace vectoriel muni de deux filtrations }
\end{center}
Soit $(V,F^{\bullet },G^{\bullet })$ un espace vectoriel bifiltré. Le fibré de Rees $\Phi( V,F^{\bullet },G^{\bullet })=\xi_{{\bf A}^2}(V,F^{\bullet},G^{\bullet})$ est bien d'après la proposition \ref{reesgfib} p.\pageref{reesgfib} un ${{\bf G}_m}^2$-fibré vectoriel sur le plan affine ${\bf A}^2$. 

Réciproquement, soit $\calf $ un ${{\bf G}_m}^2$-fibré vectoriel sur ${\bf A}^2$. Soit $V=\calf_1$ la fibre en $(1,1) \in {\bf A}^2$. On a vu que $\calf \vert_{{{\bf A}^1} \times \{1\}}$ est un ${\bf G}_m$-fibré sur la droite affine de fibre naturellement isomorphe à $V$ en $1$. On en déduit donc, comme dans la section précédente, par le morphisme $\Phi_I$ une filtration $F^{\bullet }$. De même on obtient une filtration $G^\bullet$ en restreignant le fibré à $\{1\} \times {\bf A}^1$ par $(V,G^{\bullet })=\Phi_{I}(\calf \vert_{{ \{ 1 \} \times {\bf A}^1} })$. On en déduit donc un morphisme que l'on notera toujours $\Phi_I$ de la catégorie des fibrés vectoriels sur ${{\bf A}^{2}}$ munis de l'action de ${{\bf G}_{m}} \times {{\bf G}_{m}}$ prolongeant l'action standard dont les morphismes sont les morphismes de fibrés vectoriels qui commutent avec l'action, notée ${\text Bun}({\bf A}^{2}/({\bf G}^{m})^{2})$, vers la catégorie $\calc_{2filtr}$.

\begin{proposition}\label{inverse2}
On a une équivalence de catégories entre la catégorie des espaces vectoriels de dimension finie munis de deux filtrations exhaustives et décroissantes et dont les morphismes sont les morphismes strictement compatibles aux filtrations et la catégorie des fibrés vectoriels sur ${{\bf A}^{2}}$ munis de l'action de ${{\bf G}_{m}} \times {{\bf G}_{m}}$ prolongeant l'action standard dont les morphismes sont les morphismes de fibrés vectoriels dont les conoyaux sont sans torsion et qui commutent avec l'action :
$$
\xymatrix{
\{ \calc_{2filtr} \} \ar@<2pt>[r]^-{\Phi_{R}} &  \{  {\text Bun}({\bf A}^{2}/({\bf G}_{m})^{2}) \} \ar@<2pt>[l]^-{\Phi_{I}}
}
.$$
\end{proposition}

\begin{preuve} Partons d'un espace vectoriel bifiltré $(V,F^{\bullet },G^{\bullet })$. La construction du fibré de Rees équivariant restreinte à ${\bf A}^{1} \times \{1\}$ et la proposition \ref{inverse1} p.\pageref{inverse1}, nous permettent de voir que par $\Phi_{I} \circ \Phi_{R}$ on récupère la même filtration $F^{\bullet }$ sur $V$. On fait de même pour la deuxième filtration. Ainsi $\Phi_{I}\circ \Phi_{R}((V,F^{\bullet },G^{\bullet }))=(V,F^{\bullet },G^{\bullet })$.

Réciproquement, soit $\calf $ un objet de ${\text Bun}({\bf A}^{2}/({\bf G}_{m})^{2})$ et soit $B=\calf ({\bf A}^{2})$ le $k[u,v]$-module associé. Comme dans la construction inverse dans le cas d'un espace vectoriel filtré, on a une inclusion $B \subset V \otimes_{k} k[u,u^{-1},v,v^{-1}]$ et les filtrations précédement définies par $\Phi_{I}$ sont entièrement caractérisées par $F^{p}\cap G^{q}=\{ w \in V \vert u^{-p}v^{-q}w \in B \}$. Soit $B'$ le module de Rees associé à $(V,F^{\bullet },G^{\bullet })$ et $\calf '=\Phi_{R}(V,F^{\bullet },G^{\bullet })=\Phi_{R}\circ \Phi_{I}(\calf )$ le fibré vectoriel associé. D'après la construction des filtrations $B'$ est un sous-module de $B$. On en déduit un morphisme de fibrés vectoriels $f : \calf '\rightarrow \calf$. D'après l'équivalence montrée dans la section précédente ce morphisme est un isomorphisme en dehors de $(0,0) \in {\bf A}^2$. En effet, pour tout $x \neq 0 \in {\bf A}^1$ : $f \vert_{\{x \} \times {\bf A}^{1}}: \calf ' \vert_{\{x \} \times {\bf A}^{1}} \rightarrow \calf \vert_{\{x \} \times {\bf A}^{1}}$ est l'isomorphisme de ${\bf G}_m$-fibrés sur ${\bf A}^1$ de $\Phi_{R}\circ \Phi_{I} (\calg) $ vers $\calg$ donné dans la section précédente où $\calg=(B/(u-x)B)^{\sim }$ et si l'on note par $A$ le $k[v]$-module $B/(u-x)B$, $\Phi_{R}\circ \Phi_{I} (\calg)=(A')^{\sim }$ c'est à dire le ${\bf G}_m$-fibré vectoriel sur ${\bf A}^1$ associé à l'espace vectoriel filtré $\Phi_{I}(\calg)$. Donc pour tous $(x,y) \neq (0,0)$, $f_{(x,y)}:{\calf '}_{(x,y)} \rightarrow \calf_{(x,y)}$ est un isomorphisme. C'est donc un isomorphisme de fibré vectoriels sur ${\bf A}^{2}- \{(0,0)\}$ or $\calf '$ et $\calf $ sont des fibrés vectoriels sur ${\bf A}^2$, donc $f$ est un isomorphisme : $\calf'\cong \calf$.

Montrons que cette correspondance est pleinement fidèle. On considère l'application
\begin{center} $ \text{Hom}_{\calc_{2filtr}}( (V,F^{\bullet },G^{\bullet }),(W,H^{\bullet },I^{\bullet })) \rightarrow \text{Hom}_{k[u,v]}(R^{2}(V,F^{\bullet },G^{\bullet }),R^{2}(W,H^{\bullet },I^{\bullet }))$, \end{center}
qui est bien définie car les morphismes à gauche sont compatibles aux filtrations, composée à la filtration

\begin{center} $  \text{Hom}_{k[u,v]}(R^{2}(V,F^{\bullet },G^{\bullet }),R^{2}(W,H^{\bullet },I^{\bullet })) \rightarrow \text{Hom}_{\text{Bun}({{\bf A}^2}/({{\bf G}_m})^{2})}(\Phi_{R}(V,F^{\bullet },G^{\bullet }),\Phi_{R}(W,H^{\bullet },I^{\bullet }))$, 
\end{center}
bien définie par la construction fonctorielle $\sim$.\\

Montrons que la composée de ces deux flèches est bijective. L'injectivité des deux flèches est claire. De plus un morphisme strictement compatible entre espaces vectoriels bifiltrés donne lieu à un morphisme de modules de Rees dont noyau et image sont sans torsion ce qui prouve que le morphisme induit entre $\xi_{{\bf A}^2}(V,F^{\bullet },G^{\bullet })$ et $\xi_{{\bf A}^2}(W,H^{\bullet },I^{\bullet })$ a son noyau et son image sans torsion. La surjectivité en termes de morphismes d'espaces vectoriels filtrés est directe, reste à vérifier que les morphismes sont strictement compatibles aux filtrations. On vérifie cela pour $F^{\bullet }$ et $H^{\bullet }$ (resp. $G^{\bullet }$ et $I^{\bullet }$) en se plaçant sur une droite affine de ${\bf A}^2$ de la forme ${\bf A}^{1} \times \{v_{0} \}$ (resp. $\{ u_{0} \} \times {\bf A}^{1} $) qui évite le lieu singulier des noyaux et images du morphisme de faisceaux ce qui est toujours possible car ces faisceaux étant sans torsion, le lieu singulier est de codimension $2$ par le corollaire \ref{singfaiscohcodim2} p.\pageref{singfaiscohcodim2}.

\end{preuve}

\subsection{Construction du fibré de Rees sur ${\bf P}^2$ à partir de trois filtrations}

\subsubsection{Construction}
Recouvrons ${\bf P}^{2}=\text{Proj} \ k[u_{0},u_{1},u_{2}]$ par les trois cartes affines standard, $U_{k}=\{(u_{0},u_{1},u_{2}) \in {\bf P}^{2}, u_{k}\neq 0\}={\bf A}^{2}_{ij}=\text{Spec} \ k[\frac{u_{i}}{u_k},\frac{u_{j}}{u_k}]$ où $i,j,k$ sont deux à deux distincts, $\{i,j,k\}=\{0,1,2\}$ et $i<j$. Sur les trois cartes, on effectue la construction du fibré de Rees associé à deux filtrations. Pour les trois triplets $(i,j,k)$ considérés au-dessus on construit le fibré de Rees $\xi_{{\bf A}^{2}_{ij}}(V,F^{\bullet}_{i},F^{\bullet }_{j})$ sur $U_{k}={\bf A}^{2}_{ij}$. Nous allons montrer que l'on peut recoller ces faisceaux obtenus sur chacune des cartes affine pour obtenir un faisceau sur ${\bf P}^{2}$. En fait le recollement est simple puisque deux faisceaux sur deux des cartes affines coïncident sur l'intersection de ces cartes. Il suffit alors pour décrire la structure de fibré d'exhiber les trivialisations locales. D'après la description des fibrés faite dans la proposition \ref{trivialisationderees} p.\pageref{trivialisationderees}, sur chacune des cartes affines les fibrés sont isomorphes au fibré trivial ${\mathcal O}_{U_k}^{\text{dim}(V)}$ sur $U_k$ par les isomorphismes $\phi_{k}:{\mathcal O}_{U_k}^{\text{dim}(V)} \rightarrow \xi_{{\bf A}^{2}_{ij}}(V,F^{\bullet}_{i},F^{\bullet }_{j})$. \\

\begin{propdef}\label{recolledef}
Soit $(V,F_{0}^{\bullet },F_{1}^{\bullet },F_{2}^{\bullet })$ un objet de $\calc_{3filtr}$ On peut recoller les fibrés de Rees $\xi_{{\bf A}^{2}_{ij}}(V,F^{\bullet}_{i},F^{\bullet }_{j})$ obtenus à partir de deux filtrations sur chacun des ouverts affines ${U_k}$ où $\{i,j,k\}=\{1,2,3\}$ et $i<j$. On en déduit un fibré vectoriel sur ${\bf P}^2$. Le fibré ainsi obtenu est appelé fibré de Rees associé à l'espace vectoriel trifiltré  $(V,F_{0}^{\bullet },F_{1}^{\bullet },F_{2}^{\bullet })$ et est noté ${\xi}_{{\bf P}^2}(V,F_{0}^{\bullet },F_{1}^{\bullet },F_{2}^{\bullet })$.\\

\end{propdef}

\begin{preuve} On utilise les restrictions des isomorphismes trivialisants des fibrés de Rees sur les plans affines décrits plus haut pour former des cocycles définissant un fibré. Introduisons quelques notations. Les intersections deux à deux des ouverts affines standards $U_{k}$ où $k \in \{0,1,2 \}$ seront notées $U_{kl}=U_{k}\cap U_{l}$ où $k,l \in \{0,1,2 \}$ et $k < l$. $j_{kl}$ sera l'inclusion $j_{kl}: U_{kl} \hookrightarrow U_{k}$ et $j_{lk}$ sera l'inclusion $j_{lk}: U_{kl} \hookrightarrow U_{l}$. De m\^eme, on a les trois inclusions $j_{ijk} : U_{ijk} \hookrightarrow U_{ij}$. 
Pour $l \neq k$, on peut supposer par exemple que $l=j$ où $(i,j,k)$ est un triplet de la forme décrite plus haut. D'après le lemme \ref{surlesouverts} p.\pageref{surlesouverts}, $ \xi_{{\bf A}^{2}_{ij}}(V,F^{\bullet}_{i},F^{\bullet }_{j})\vert_{U_{jk}}= \xi_{{\bf A}^{2}_{ik}}(V,F^{\bullet}_{i},F^{\bullet }_{k})\vert_{U_{jk}}$ comme sous faisceaux de $j^{*}_{jk}({\mathcal O}_{U_{ijk}}\otimes V)$. En effet ces deux faisceaux sont égaux au faisceau $\xi_{{\bf A}^1}(V,F_{i}^{\bullet })\otimes {\mathcal O}_{{\bf G}_{m}}$ sur $U_{jk}=\text{Spec} \ k[u,v,v^{-1}]$, où ${\bf G}_{m}=\text{Spec} \ k[v,v^{-1}]$. Cela suffit donc à recoller les faisceaux locallement libres sur les cartes affines donc triviaux, les isomorphismes de transition sont les identités sur les intersections des cartes affines, les conditions de cocycles sont évidement vérifiées. 
\end{preuve}

Décrivons la structure de fibré vectoriel en donnant des trivialisations. Sur $U_k$ on a l'isomorphisme $\phi_{k}:{\mathcal O}_{U_k}^{\text{dim}(V)} \rightarrow \xi_{{\bf A}^{2}_{ij}}(V,F^{\bullet}_{i},F^{\bullet }_{j})$. On peut ainsi définir $\phi_{kl}=\phi_{l}^{-1}\circ \phi_{k}\vert_{U_{kl}}$ comme fonction de transition entre les fibrés triviaux sur les cartes $U_{k} $ et $U_{l}$. La condition de cocycle est automatiquement vérifiée car $\phi_{k}\vert_{U_{ijk}}(\xi_{{\bf A}^{2}_{ij}}(V,F^{\bullet}_{i},F^{\bullet }_{j}))={\mathcal O}_{U_{ijk}}\otimes V$. Ainsi sur $U_{ijk}$ : $\phi_{ij} \circ \phi_{jk} \circ \phi_{kl}=id$.
% Alors, comme $U_{i}=Spec\, k[u,v]$, $U_{j}=Spec\, k[u,v]$ et que l'on recolle ces deux ouverts
%par l'intermédiaire de l'automorphisme de $U_{ij}=Spec\, k[u,u^{-1},v]$, $(u,v) \mapsto (u^{-1},v/u)$, on a les données suivantes : les isomorphismes $\alpha_{i}: {\xi}_{U_{i}}(V,F_{j}^{\bullet },F_{k}^{\bullet }) \cong {\mathcal O}_{U_{i}}^{n}$ pour tout $i$ sur $U_{i}$ (avec les conditions précédentes pour $j$ et $k$), ${\xi}_{U_{ij}}(V,F_{j}^{\bullet },F_{k}^{\bullet })={\xi}_{U_{ij}}(V,F_{i}^{\bullet },F_{k}^{\bullet })$ (c'est la construction $^{\sim}$ restreinte à $Spec\,k[u,u^{-1},v]$). On peut appliquer le lemme précédent en prenant ${\mathcal O}_{U_{i}}^n$ pour ${\mathcal E}_{i}$, et $(\alpha_{j}\vert_{U_{ij}}) \circ (\alpha_{i}\vert_{U_{ij}})^{-1}$ pour $\phi_{ij}: {\mathcal E}_{i}\vert_{U_{ij}} \mapsto {\mathcal E}_{j}\vert_{U_{ij}}$. La condition de cocycle est ainsi automatiquement vérifiée. Le fibré est unique à isomorphisme prés.
Reste à montrer que la construction ne dépend pas des bigraduations $V=\oplus_{p,q}V^{p,q}$ choisies pour chaque construction $\xi_{{\bf A}^{2}_{ij}}(V,F^{\bullet}_{i},F^{\bullet }_{j})$ sur les cartes affines $U_{k}$. Supposons que l'on ait une autre graduation $V=\oplus_{p,q}U^{p,q}$. Alors prendre le bigradué associé aux filtrations $F^{\bullet }_{i}$ et $F^{\bullet}_{j}$, donne des isomorphismes $V^{p,q} \cong U^{p,q}$. Ces isomorphismes induisent des isomorphismes au niveau des modules de Rees puis des fibrés de Rees sur les cartes affines $U_{k}$ pour tous les triplets $i,j,k$. Ainsi, changer les bigraduations associées aux paires de filtrations revient à changer de trivialisations locales du fibré localement libre sur ${\bf P}^2$ et conduit donc à un fibré isomorphe. Le fibré de Rees ${\xi}_{{\bf P}^2}(V,F_{0}^{\bullet },F_{1}^{\bullet },F_{2}^{\bullet })$ est bien déterminé à isomorphisme prés à partir de $(V,F_{0}^{\bullet },F_{1}^{\bullet },F_{2}^{\bullet })$.\\     
%ecrire avec Schaf

%\begin{figure}[ht] 

%\begin{center}

%\scalebox{1}{\input{petitp2.pstex_t}}

%\end{center}

%\caption {Construction de ${\xi}_{{\bf P}^2}(V,F_{0}^{\bullet },F_{1}^{\bullet },F_{2}^{\bullet })$}

%\end{figure}
${}$\\
{\bf Exemple :}
  Donnons la description explicite du fibré de Rees associé à un espace vectoriel de dimension $1$ muni de trois filtrations. Cet exemple est important car on essaiera toujours par la suite de se ramener à une somme directe d'espaces vectoriels de dimension 1 munis de trois filtrations, la construction précédente se comportant bien pour les sommes directes. Soit $V=k$ muni de trois filtrations. On a forcément $(V,F_{0}^{\bullet },F_{1}^{\bullet },F_{2}^{\bullet })=(V ,Dec^{r}Triv^{\bullet},Dec^{p}Triv^{\bullet},Dec^{q}Triv^{\bullet})$ où $r,p,q$ sont les rangs respectifs auxquels $F^{\bullet}_{0}$, $F^{\bullet}_{1}$ et $F^{\bullet}_{2}$ sautent.

 ${\xi}_{{\bf P}^2}(V,Dec^{p}Triv^{\bullet},Dec^{q}Triv^{\bullet},Dec^{r}Triv^{\bullet})$ est un fibré en droites sur ${\bf P}^2$ donc de la forme $ O_{{\bf P}^2}(D)$. Explicitons le diviseur $D$. Pour cela notons $D_{k}$ le complémentaire de $U_{k}$ dans ${\bf P}^2$ (la droite a l'infini associée à l'ouvert affine $U_k$), c'est à dire : pour $k \in \{0,1,2\}$ : $U_{k}=\{ (u_{0}:u_{1}:u_{2}) \in {{\bf P}^2} \vert u_{k}\neq 0 \} \cong {\bf A}^{2}_{ij}$ et $D_{k}=\{ (u_{0}:u_{1}:u_{2}) \in {{\bf P}^2} \vert u_{k}= 0 \} \cong {\bf P}^{1}$. La construction dans le cas d'un fibré de rang $1$ est simplifiée par le fait que l'on peut évidemment trouver des trivialisations sur les ouverts compatibles à toutes les filtrations. Une section méromorphe $s$ du fibré est donnée par (on note $v$ un vecteur engendrant la droite $V$) : $ (\frac{u_0}{u_2})^{-r}.(\frac{u_1}{u_2})^{-p}v $ sur $U_{2}$, $ (\frac{u_1}{u_0})^{-p}.(\frac{u_2}{u_0})^{-q}v $ sur $U_{0}$ et $ (\frac{u_0}{u_1})^{-r}.(\frac{u_2}{u_1})^{-q}v $ sur $U_{1}$. D'où : $D=(s)=-rD_{0}-pD_{1}-qD_{2}$.\\
${}$\\
Nous noterons ce fibré\footnote{On écrit $(r,p,q)$ et non $(p,q,r)$ car dans les applications \`a la th\'eorie de Hodge $F_{1}^{\bullet}$ et $F_{2}^{\bullet }$ seront la filtration de Hodge et sa filtration opposée qui sont décrites classiquement par les indices $p$ et $q$} $$\xi_{{\bf P}^2}^{r,p,q}.$$ 
On a décrit l'isomorphisme
$$\xi_{{\bf P}^2}^{r,p,q} \cong \calo_{{\bf P}^2}(-r-p-q).$$

\begin{lemme}\label{plusfois}
Soit $(V_{i},{F_{0}^{\bullet }}_{i},{F_{1}^{\bullet }}_{i},{F_{2}^{\bullet }}_{i})$ une famille finie d'espaces vectoriels munis de trois filtrations, alors :\begin{center}
 ${\xi}_{{\bf P}^2}(\oplus_{i}(V_{i},{F_{0}^{\bullet }}_{i},{F_{1}^{\bullet }}_{i},{F_{2}^{\bullet }}_{i})) \cong \oplus_{i}{\xi}_{{\bf P}^2}(V_{i},{F_{0}^{\bullet }}_{i},{F_{1}^{\bullet }}_{i},{F_{2}^{\bullet }}_{i}) $. 
  ${\xi}_{{\bf P}^2}(\otimes_{i}(V_{i},{F_{0}^{\bullet }}_{i},{F_{1}^{\bullet }}_{i},{F_{2}^{\bullet }}_{i})) \cong \otimes_{i}{\xi}_{{\bf P}^2}(V_{i},{F_{0}^{\bullet }}_{i},{F_{1}^{\bullet }}_{i},{F_{2}^{\bullet }}_{i}) $.\\
\end{center}
\end{lemme}
 
\begin{preuve}
Rappelons (\cite{har}) que si $M$ et $N$ sont deux $A$-modules, et que l'on note par une tilde la construction du faisceau cohérent sur $\text{Spec} \ A$ associée à un $A$-module, alors : $ (M \oplus N)^{\sim} \cong M^{\sim } \oplus N^{\sim } $ et $(M \otimes_{A} N)^{\sim}  \cong M^{\sim} \otimes_{\text{Spec} \ A} N^{\sim}$. Les constructions des fibrés de Rees sur les plans affines sur les ouverts standards associées à la somme directe et au produit tensoriel se déduisent donc des constructions sur chacun des facteurs. Les trivialisations sur ces ouverts sont données par chacunes des trivialisations correspondantes aux espaces vectoriels trifiltrés. Les fonctions de transitions du fibré associé à la somme directe et du fibré associé au produit tensoriel sont donc données par des blocs correspondants aux somme directes et par produit tensoriel des fonctions de transition. Ainsi le fibré de Rees de la somme directe est bien la somme de Whitney des fibrés de Rees associés aux espaces vectoriels en somme directe et le fibré de Rees associé au produit tensoriel est le produit tensoriel des fibrés de Rees de chacun des facteurs.
\end{preuve}
                                        
{\bf Construction équivalente du fibré de Rees associé à trois filtrations}\\

Rappelons que la construction du faisceau de Rees associé à trois filtrations $(V_{},{F_{0}^{\bullet }}_{},{F_{1}^{\bullet }}_{},{F_{2}^{\bullet }}_{})$ sur l'espace affine ${\bf A}^3$ donne un faisceau cohérent réflexif $$\xi_{{\bf A}^3}(V_{},{F_{0}^{\bullet }}_{},{F_{1}^{\bullet }}_{},{F_{2}^{\bullet }}_{})=(R^{3}(V,F^{0},F_{1}^{\bullet },F_{2}^{\bullet })_{u_0})^{\sim}.$$
 
Ce faisceau est en fait un $({\bf G}_{m})^3$-faisceau pour l'action standard de $({\bf G}_{m})^3$ sur ${\bf A}^3$ par translation. Considérons la diagonale de $({\bf G}_{m})^3$. La restriction de l'action à la diagonale est l'action décrite dans la section \ref{groupesalg} p.\pageref{groupesalg}. $\xi_{{\bf A}^3}(V_{},{F_{0}^{\bullet }}_{},{F_{1}^{\bullet }}_{},{F_{2}^{\bullet }}_{})\vert_{{\bf A}^{3}\backslash \{ (0,0,0) \}}$ est bien sur un ${\bf G}_m$-fibré pour cette action. Comme on l'a vu dans l'étude de cette action, le morphisme entre la ${\bf G}_m$-variété ${\bf A}^{3}\backslash \{ (0,0,0) \}$ et ${\bf P}^2$ est un ${\bf G}_m$-fibré principal. Donc le\\
${\bf G}_m$-fibré $\xi_{{\bf A}^3}(V_{},{F_{0}^{\bullet }}_{},{F_{1}^{\bullet }}_{},{F_{2}^{\bullet }}_{})\vert_{{\bf A}^{3}\backslash \{ (0,0,0) \}}$ descend sur ${\bf P}^2$, i.e. il existe un faisceau cohérent $\cale $ sur ${\bf P}^2$ tel que 
$$ \xi_{{\bf A}^3}(V_{},{F_{0}^{\bullet }}_{},{F_{1}^{\bullet }}_{},{F_{2}^{\bullet }}_{})\vert_{{\bf A}^{3}\backslash \{ (0,0,0) \}} \cong f^{*}\cale$$
soit un isomorphisme de ${\bf G}_m$-fibrés sur ${\bf A}^{3}\backslash \{ (0,0,0) \}$ ($f$ est le morphisme de l'espace affine privé de l'origine vers l'espace projectif).

Comme l'action de ${\bf G}_m$ est sans point fixe, $\cale $ est unique à isomorphisme prés (cf \cite{huyleh}). Donc, à isomorphisme près il existe un unique faisceau cohérent $\cale $ sur le plan projectif tel que $ \xi_{{\bf A}^3}(V_{},{F_{0}^{\bullet }}_{},{F_{1}^{\bullet }}_{},{F_{2}^{\bullet }}_{})\vert_{{\bf A}^{3}\backslash \{ (0,0,0) \}} \cong f^{*}\cale$ soit un isomorphisme de ${\bf G}_m$-fibrés. \\ 
 
\begin{proposition}
$$ \xi_{{\bf A}^3}(V_{},{F_{0}^{\bullet }}_{},{F_{1}^{\bullet }}_{},{F_{2}^{\bullet }}_{})\vert_{{\bf A}^{3}\backslash \{ (0,0,0) \}} \cong f^{*} \xi_{{\bf P}^2}(V_{},{F_{0}^{\bullet }}_{},{F_{1}^{\bullet }}_{},{F_{2}^{\bullet }}_{})$$
comme isomorphisme de ${\bf G}_m$-fibrés sur ${\bf A}^{3}\backslash \{ (0,0,0) \}$.\\
\end{proposition}

\begin{preuve}
Plaçons nous d'abord sur un ouvert affine $U_{k}=\text{Spec}k[\frac{u_{i}}{u_k},\frac{u_{j}}{u_k}]$ de ${\bf P}^{2}=\text{Proj} \ k[u_{1},u_{2},u_{3}]$ décrit en début de section, $U_{0}$ par exemple. $f^{-1}(U_{0})=\{(u_{0},u_{1},u_{2}) \in {\bf A}^{3}\vert\, u_{0}\neq 0 \}$. En terme d'anneaux cette restriction de $f$ est associée au morphisme $g : A=k[\frac{u_{1}}{u_0},\frac{u_{2}}{u_0}] \rightarrow B=k[u_{0},u_{0}^{-1},u_{1},u_{2}]=k[u_{0},u_{1},u_{2}]_{u_0}$ ($\text{Spec}\,A={\bf A}^{2}$ et $\text{Spec}\,B=f^{-1}(U_{0})$). Alors,

\begin{eqnarray*}   f^{*} \xi_{{\bf P}^2}(V_{},{F_{0}^{\bullet }}_{},{F_{1}^{\bullet }}_{},{F_{2}^{\bullet }}_{})\vert_{\text{Spec}\,B}& \cong &(R^{2}(V,F_{1}^{\bullet },F_{2}^{\bullet }) \otimes_{A}B)^{\sim}\\
{}&\cong &(( R^{3}(V,F^{0},F_{1}^{\bullet },F_{2}^{\bullet })_{u_0})^{\sim}\\
{} &\cong &\xi_{{\bf A}^3}(V_{},{F_{0}^{\bullet }}_{},{F_{1}^{\bullet }}_{},{F_{2}^{\bullet }}_{})\vert_{\text{Spec}\,B}.
\end{eqnarray*}
 On obtient ainsi trois isomorphismes dans les images réciproques des trois cartes affines du plan projectif. Ces isomorphismes se recollent bien car sur le plan projectif se placer sur l'intersection de deux cartes affines revient à ne considérer qu'une seule filtration (celle qui est en commun aux constructions sur chacunes des deux cartes affines) et il en va de même ¨au dessus¨ : l'image réciproque de l'intersection est isomorphe à ${\bf A}^{1}\times {\bf G}_{m}\times {\bf G}_{m}$, et le fibré restreint à cet ouvert fait apparaître la construction de Rees avec la filtration voulue dans la direction de la droite affine et est trivial dans les autres directions. On obtient ainsi l'isomorphisme voulu sur ${\bf A}^{3}\backslash \{ (0,0,0) \}$. 

La compatibilité de cet isomorphisme avec l'action de ${\bf G}_m$ induite par $f^*$ est claire.   
\end{preuve}

{\bf Les fibrés de Rees sur ${\bf P}^2$ sont des $G$-fibrés}\\

Le fibré de Rees ${\xi}_{{\bf P}^2}(V,F_{0}^{\bullet },F_{1}^{\bullet },F_{2}^{\bullet })$ est un ${{\bf G}_m}^3$-fibré sur le plan projectif pour l'action de ${{\bf G}_m}^3$ héritée de l'action par translation sur l'espace affine ${\bf A}^3$. L'action de la diagonale de ${{\bf G}_m}^3$, $\Delta ( {\bf G}_{m}^{3})$ est triviale (où $\Delta({\bf G}_{m}^{3})$ est l'image du morphisme ${\bf G}_{m} \hookrightarrow {\bf G}_{m}^3$ donné par $t \mapsto (t,t,t)$). Sur chaque ouvert affine $U_{k}$, on a une action de ${{\bf G}_m}\times {\bf G}_m$ qui est le quotient non canonique ${\bf G}_{m}^{3}/\Delta({\bf G}_{m}^{3})$. Cette action, restreinte aux ouverts standard, correspond \`a l'action de ${{\bf G}_m}\times {\bf G}_m$ que l'on a décrite sur les fibrés de Rees associés à deux filtrations sur les ouverts affines. L'action de ${{\bf G}_m}^3$ sur ${\bf A}^{3}$ est donnée par :\\
\begin{center}
$k[u_{0},u_{1},u_{2}] \rightarrow k[u_{0},u_{1},u_{2}] \otimes k[t_{0},t_{0}^{-1},t_{1},t_{1}^{-1},t_{2},t_{2}^{-1}]$\\
$u_{i} \mapsto u_{i}\otimes t_{i}$, où $i \in \{0,1,2\}$.
\end{center}
Elle induit, sur $U_{0}=\text{Spec}k[\frac{u_{1}}{u_0},\frac{u_{2}}{u_0}]$ par exemple, l'action de ${{\bf G}_m}\times {\bf G}_{m}=\text{Spec}\,k[\frac{t_{1}}{t_0},\frac{t_{2}}{t_0},{\frac{t_{1}}{t_0}}^{-1},{\frac{t_{2}}{t_0}}^{-1}]$ :\\
\begin{center}
$k[\frac{u_{1}}{u_0},\frac{u_{2}}{u_0}] \rightarrow k[\frac{u_{1}}{u_0},\frac{u_{2}}{u_0}] \otimes k[\frac{t_{1}}{t_0},\frac{t_{2}}{t_0},{\frac{t_{1}}{t_0}}^{-1},{\frac{t_{2}}{t_0}}^{-1}]$\\
$\frac{u_{i}}{u_0}  \mapsto \frac{u_{i}}{u_0} \otimes \frac{t_{i}}{t_0}$, où $i \in \{ 1,2 \}$.
\end{center}
Le quotient ${{\bf G}_m}\times {\bf G}_{m}$ de ${\bf G}_{m}^{3}$ vient du morphisme :\\
\begin{center}
$k[t_{0},t_{0}^{-1},t_{1},t_{1}^{-1},t_{2},t_{2}^{-1}] \rightarrow k[t_{0},t_{0}^{-1},t_{1},t_{1}^{-1},t_{2},t_{2}^{-1}] \otimes k[t,t^{-1}]$\\
$t_{i} \mapsto t_{i}\otimes t $ et $t_{i}^{-1} \mapsto t_{i}^{-1}\otimes t^{-1} $.
\end{center} 
Posons ${\bf T}={\bf G}_{m}^{3}/\Delta({\bf G}_{m}^{3})$, nous avons montré que\\

\begin{proposition} Le fibré de Rees sur ${\bf P}^{2}$ associé à un espace vectoriel muni de trois filtrations exhaustives et décroissantes $(V,F_{0}^{\bullet },F_{1}^{\bullet },F_{2}^{\bullet })$, ${\xi}_{{\bf P}^2}(V,F_{0}^{\bullet },F_{1}^{\bullet },F_{2}^{\bullet })$ est un ${\bf T}$-fibré pour l'action de ${\bf T}$ par translation sur ${\bf P}^2$.

\end{proposition}
% ${\bf P}^{2}=\text{Proj} \ k[u_{1},u_{2},u_{3}]$ par les trois cartes affines standard $U_{k}=\{(u_{1},u_{2},u_{3}) \in {\bf P}^{2}, u_{k}\neq 0\}={\bf A}^\{2}_{ij}=\text{Spec} \ k[\frac{u_{i}}{u_k},\frac{u_{j}}{u_k}]$ où $i,j,k$ sont deux à deux distincts, $\{i,j,k\}=\{1,2,3\}$ et $i<j$
${}$\\
{\bf Exemple :} Revenons à l'exemple des fibrés en droite de Rees $\xi_{{\bf P}^2}^{r,p,q}$. Décrivons l'action de ${\bf T}$ sur un des ouverts standards, $U_{0}$ par exemple. Soit $w$ un vecteur qui engendre l'espace filtré et $s$ une section équivariante sur $U_{0}$ telle que $s(1)=w$. Alors pour tout $(u_{1},u_{2}) \in {\bf A}^2$, $s(u_{1},u_{2})=u_{1}^{-p}u_{2}^{-q}w$. $p$ et $q$ sont les caractères de l'action restreinte au plan affine.\\

\subsubsection{Restrictions des fibrés de Rees aux droites standards}
La construction du fibr\'e de Rees \'etant sym\'etrique par rapport aux trois filtrations, il suffit de d\'ecrire la restriction du fibr\'e à l' une des droites standards. Nous \'etudions ici la restriction à la droite ${\bf P}^{1}_{0}$ qui correspond aux pôles apport\'es par la filtration $F^{\bullet }_{0}$.

\begin{proposition}\label{resfibree}
La restriction du fibr\'e de Rees peut être d\'ecrite en terme de fibr\'e de Rees sur ${\bf P}^{1}$ (cf annexe A), en effet,
$${\xi}_{{\bf P}^2}(V,F_{0}^{\bullet },F_{1}^{\bullet },F_{2}^{\bullet })\cong \oplus_{r}{\xi}_{{\bf P}^1}(Gr_{F_{0}^{\bullet }}^{r},F_{1}^{\bullet }{}_{ind},F_{2}^{\bullet }{}_{ind})\otimes \calo_{{\bf P}^1}(-r),$$
où $F_{1}^{\bullet }{}_{ind}$ et $F_{2}^{\bullet }{}_{ind}$ sont les filtrations induites par $F_{1}^{\bullet }$ et $F_{2}^{\bullet }$ sur l'espace gradu\'e associ\'e à la première filtration.
\end{proposition}

\subsubsection{R\'esolution par des fibr\'es scind\'es}

On montre que l'on peut trouver une r\'esolution de tout fibr\'es de Rees par des fibr\'es somme de fibr\'es en droite (des $\xi_{{\bf P}^2}^{r,p,q}$). La r\'esolution n'est pas a priori finie et les termes de la r\'esolution d'un fibr\'e de Rees associ\'e \`a des filtrations oppos\'ees ne proviennent pas forc\'ement de filtrations oppos\'ees.

\begin{proposition}
Soit $\xi$ un ${\bf T}$-fibr\'e \'equivariant sur ${\bf P}^2$. Il existe alors un r\'esolution par des ${\bf T}$-fibr\'es scind\'es $\xi_{i}$, $i \geq 0$ :
$$
\xymatrix{
...\ar[r]&   \xi_{1} \ar[r] &  \xi_{0} \ar[r] &  \xi \ar[r] & 0.}
$$

\end{proposition}

\begin{preuve}Nous allons montrer l'assertion plus g\'en\'erale :

Soit $X$ une vari\'et\'e projective lisse munie de l' action plate d'un tore ${\bf T}$, $\sigma$, et soit $\calf$ un faisceau de $\calo_{X}$-modules ${\bf T}$-\'equivariant, alors il existe une r\'esolution de $\calf$ par des fibr\'es scind\'es ${\bf T}$-\'equivariants.

Notons ${\calo}_{X}(1)$ le fibr\'e ample de $X$. D'apr\`es \cite{har}, il existe $n_0$ tel que pour tout $n \geq n_0$ $\calf(n)=\calf \otimes \calo_{X}(n)$ soit engendr\'e par ses sections. D'o\`u
$$\xymatrix{
H^{0}(X,\calf(n))\otimes \calo_{X} \ar[r] &  \calf(n) \ar[r] & 0. }
$$

C'est un morphisme \'equivariant de faisceaux coh\'erents \'equivariants. En effet, notons par $p$ la projection de $X$ sur le point $*$. Les foncteurs $p^*$ et $p_*$ sont adjoints de la cat\'egorie des ${\bf T}$-faisceaux coh\'erents sur $X$ vers la catégorie des ${\bf T}$-faisceaux coh\'erents sur le point i.e. des espaces vectoriels munis de l'action du tore. Le morphisme exhib\'e est donc l' adjoint du morhisme $id$ par
$$
\text{Hom}_{{\bf T},X}(p^{*}V,W)=\text{Hom}_{{\bf T},*}(V,p_{*}W),
$$
o\`u $V=p_{*}\calf(n)$ et $W=\calf(n)$.

Donc le tore ${\bf T}$ op\`ere sur l' espace vectoriel $H^{0}(X,\calf(n))$. L' action se d\'ecompose donc suivant les caract\`eres $\chi$ :
$$
H^{0}(X,\calf(n))=\oplus_{i \in \chi}L_{i},$$
o\`u le tore agit sur $L_i$ par le caract\`ere $\chi$. D'o\`u le morphisme \'equivariant surjectif qui permet de d\'efinir $\xi_{0}$ somme de fibr\'es en droites \'equivariants,
$$\xymatrix{
\oplus_{i \in\chi} 
((\xi^{p_{i},q_{i},r_{i}}_{{\bf P}^2})^{\text{dim}_{\bf C}L_i})\otimes \calo_{X}(-n)  \ar[r] &  \calf \ar[r] & 0. }
$$
Le noyau de ce morphisme est un faisceau coh\'erent \'equivariant sans torsion, on peut donc par la m\^eme m\'ethode trouver une fl\`eche surjective \'equivariante qui permet de d\'efinir $\xi_{1}$. On continue ainsi pour trouver la r\'esolution voulue.
\end{preuve}

\subsubsection{Etude du fibré de Rees}
Cette partie est préparatoire à la section dans laquelle on calcule le caractère de Chern. L'idée générale est d'essayer de ramener l'étude du fibré à celle de fibrés scindés, somme de fibrés en droite de Rees de type  $\xi_{{\bf P}^2}^{r,p,q}$. Pour cela on ``sépare'' les filtrations de façon à n'avoir à scinder que deux filtrations simultanément, ce qui est toujours possible.

L'étude du fibré de Rees associé à un espace vectoriel muni de trois filtrations (décroissantes et exhaustives) et de ses invariants va nous permettre de déceler à quel point les filtrations sont loins d'\^etre dans la position la plus simple, celle où elles sont simultanéement scindées. Rappelons que les trois filtrations $(F_{0}^{\bullet },F_{1}^{\bullet },F_{2}^{\bullet })$ sont simultanéement scindées s'il existe des sous-espaces vectoriels $V^{p,q,r}$ tels que :
\begin{center}$V =\oplus_{p,q,r}V^{p,q,r}$ (où la somme est finie) et 
$\left\lbrace \begin{array}{l}
         F_{0}^{r}=\oplus_{\{(p',q',r') \vert r'\leq  r \}}\,V^{p',q',r'},\\
         F_{1}^{p}=\oplus_{\{(p',q',r') \vert p'\leq  p \}}\,V^{p',q',r'},\\
         F_{2}^{q}=\oplus_{\{(p',q',r') \vert q'\leq  q \}}\,V^{p',q',r'}.
\end{array}
\right.$.
\end{center}
$${}$$
\begin{lemme}\label{scin}
Soit $V$ un espace vectoriel muni de trois filtrations $(F_{0}^{\bullet },F_{1}^{\bullet },F_{2}^{\bullet })$ simultanéement scindées et de scindement $V=\oplus_{p,q,r}V^{p,q,r}$. Alors :
\begin{center}
${\xi}_{{\bf P}^2}(V,F_{0}^{\bullet },F_{1}^{\bullet },F_{2}^{\bullet })=\oplus_{p,q,r}{\xi}_{{\bf P}^2}(V^{p,q,r},Dec^{r}Triv^{\bullet },Dec^{p}Triv^{\bullet },Dec^{q}Triv^{\bullet })$.
\end{center}
\end{lemme}

\begin{preuve}
Il suffit d'écrire que $(V,F_{0}^{\bullet },F_{1}^{\bullet },F_{2}^{\bullet })=\oplus_{p,q,r}(V^{p,q,r},{F_{0}^{\bullet }}_{ind},{F_{1}^{\bullet }}_{ind},{F_{2}^{\bullet }}_{ind})$ et que sur chaque "bloc" $V^{p,q,r}$, ${F_{i}^{\bullet }}_{ind}$ est une filtration de niveau $1$, c'est à dire qu'il existe $k_{i}$ tel que $F^{k_i}_{i\  ind}=V^{p,q,r}$ et $F^{k_i}_{i\  ind}=\{0\}$. Sur $V^{p,q,r}$ on a $k_{i}=p$. Et $F^{\bullet }_{i\, ind}$ est donc une filtration triviale décalée $Dec^{p}Triv^{\bullet }$. On applique ensuite le lemme \ref{plusfois} p.\pageref{plusfois}. 
\end{preuve}

Nous n'avons pas de moyens, pour étudier le fibré de Rees associé à trois filtrations\\
 ${\xi}_{{\bf P}^2}(V,F_{0}^{\bullet },F_{1}^{\bullet },F_{2}^{\bullet })$, d'écrire des suites exactes en le scindant par l'une des filtrations pour faire appara\^{\i}tre des fibrés de rang plus bas construits à partir de filtrations qui sont en positions plus simples (on pense ici à la filtration par le poids pour une structure de Hodge mixte et aux structures pures sur chacun des gradués par le poids). Ceci vient du fait que, comme on l'a vu plus haut dans la section sur les filtrations, si l'on peut toujours scinder deux filtrations simultanément, ce n'est pas en général possible pour trois filtrations. Pour pouvoir étudier ${\xi}_{{\bf P}^2}(V,F_{0}^{\bullet },F_{1}^{\bullet },F_{2}^{\bullet })$ on va lui associer un fibré que l'on pourra démonter suivant les sous-espaces associés aux différentes filtrations.

Eclatons ${\bf P}^2$ en $(0:0:1)$, on notera $e: \widetilde{{\bf P}^2} \rightarrow {\bf P}^2$ le morphisme correspondant à l'éclatement et $E=e^{-1}(0)$ le diviseur exceptionnel. L'éclatement est décrit dans l'ouvert $U_{0}={\bf A}^{2}_{12}=\text{Spec}\,k[u,v]$ par les morphismes $e_{i}: \text{Spec}\,B_{i} \rightarrow \text{Spec}\,A$ où $i \in \{1,2\}$ correspondants aux morphismes d'anneaux suivants (on note les morphismes d'anneaux de la m\^eme façon) : $e_{0}: A=k[u,v] \rightarrow B_{0}=k[x,y]$ tel que $e_{0}(u,v)=(uv,v)$ et  $e_{1}: A=k[u,v] \rightarrow B_{1}=k[z,t]$ tel que $e_{0}(u,v)=(u,uv)$. $\text{Spec}\,B_{i}$ est noté $U_{0}^{i}$. On fait la construction du fibré de Rees sur $\text{Spec}\,B_{i}$ associée au $B_{i}$-module $R^{2}(V,F_{i}^{\bullet },Triv^{\bullet })$.\\
$\,$\\
{\bf Remarque :}
Sur l'ordre des filtrations. Le but du travail sur le fibré ${\xi}_{{\bf P}^2}(V,F_{0}^{\bullet },F_{1}^{\bullet },F_{2}^{\bullet })$ est de conna{\^\i}tre les positions relatives des filtrations dans le cas où l'espace vectoriel trifiltré est un structure de Hodge mixte. Dans ce cas, on prendra $F_{0}^{\bullet }=(W_{\bullet })^{\bullet }$ la filtration croissante associée à la filtration par le poids, et $F_{1}^{\bullet }=F_{}^{\bullet }$, $F_{2}^{\bullet }={\overline F}_{}^{\bullet }$ la filtration de Hodge et sa conjuguée par rapport à la structure réelle sous-jacente. On va s'intéresser aux quotients par la filtration par le poids qui sont des structures de Hodge pures. C'est pourquoi dans la construction du fibré associé à ${\xi}_{{\bf P}^2}(V,F_{0}^{\bullet },F_{1}^{\bullet },F_{2}^{\bullet })$, on distingue $F_{0}^{\bullet }$ pour pouvoir quotienter par des sous-espaces vectoriels associés à cette filtration alors que les trois filtrations jouent des r\^oles symétriques dans la construction du fibré de Rees.\\ 

\begin{propdef}
On peut utiliser les trivialisations locales des fibrés de Rees associés aux paires de filtrations sur les quatre plans affines qui recouvrent $\widetilde{{\bf P}^2}$ pour obtenir un fibré vectoriel unique à isomorphisme prés. Il sera appelé le fibré vectoriel sur $\widetilde{{\bf P}^2}$ associé à  ${\xi}_{{\bf P}^2}(V,F_{0}^{\bullet },F_{1}^{\bullet },F_{2}^{\bullet })$, il sera noté ${\xi}_{\widetilde{{\bf P}^2}}(V,F_{0}^{\bullet },F_{1}^{\bullet },F_{2}^{\bullet },Triv^{\bullet })$.\\
\end{propdef}

\begin{preuve}
La preuve est la m\^eme que pour la construction de ${\xi}_{{\bf P}^2}(V,F_{0}^{\bullet },F_{1}^{\bullet },F_{2}^{\bullet })$, on trivialise localement sur les quatre ouverts par l'intermédiaire des isomorphismes $\alpha_{i}$ exhibés plus haut ce qui permet d'avoir automatiquement les conditions de cocycles et d'utiliser de recoller.
\end{preuve}

%\begin{figure}[ht] 

%\begin{center}

%\scalebox{1}{\input{petitecl.pstex_t}}

%\end{center}

%\caption {Construction de ${\xi}_{\widetilde{{\bf P}^2}}(V,F_{0}^{\bullet },F_{1}^{\bullet },F_{2}^{\bullet },Triv^{\bullet })$.}

%\end{figure}

\begin{lemme}\label{dirsumxitilde}Soit $(V_{i},{F_{0}^{\bullet }}_{i},{F_{1}^{\bullet }}_{i},{F_{2}^{\bullet }}_{i})$ une famille finie d'espaces vectoriels munis de trois filtrations, alors :\begin{center}
${\xi}_{\widetilde{{\bf P}^2}}(\oplus_{i}(V_{i},{F_{0}^{\bullet }}_{i},{F_{1}^{\bullet }}_{i},{F_{2}^{\bullet }}_{i},Triv^{\bullet })) \cong \oplus_{i}{\xi}_{\widetilde{{\bf P}^2}}(V_{i},{F_{0}^{\bullet }}_{i},{F_{1}^{\bullet }}_{i},{F_{2}^{\bullet }}_{i},Triv^{\bullet})$. \\
${\xi}_{\widetilde{{\bf P}^2}}(\otimes_{i}(V_{i},{F_{0}^{\bullet }}_{i},{F_{1}^{\bullet }}_{i},{F_{2}^{\bullet }}_{i},Triv^{\bullet })) \cong \otimes_{i}{\xi}_{\widetilde{{\bf P}^2}}(V_{i},{F_{0}^{\bullet }}_{i},{F_{1}^{\bullet }}_{i},{F_{2}^{\bullet }}_{i},Triv^{\bullet})$.  
\end{center}
\end{lemme}
\begin{preuve}
La preuve est similaire à celle du lemme \ref{isosumfibrees}.   

\end{preuve}

D\'ecrivons le fibr\'e $e^{*}{\xi}_{{\bf P}^2}(V,F_{0}^{\bullet },F_{1}^{\bullet },F_{2}^{\bullet })$. Il suffit de se placer dans la carte $U_0$, des coordonn\'ees $(u,v)$ d\'ecrites pr\'ec\'edemmemt. Dans la direction de la coordonn\'ee $uv$, le fibr\'e image inverse du fibr\'e est le fibr\'e de Rees sur un affine associ\'e aux filtrations $F^{\bullet}_1$ et $G^{\bullet }$ o\`u $G^{\bullet}$ est la filtration associ\'ee \`a $F^{\bullet }_{1}$ et $F^{\bullet}_{2}$ et d\'efinie par
$$
G^{r}=\sum_{p+q\geq r}F^{p}_{1}\cap F^{q}_{2}.
$$
Cette filtration est appel\'ee la convol\'ee de $F^{\bullet}_{1}$ et $F^{\bullet}_{2}$ suivant \cite{stezuc}. 

 \begin{definition}
Une filtration décroissante et exhaustive $F^{\bullet }$ est dite positive si le plus petit entier  $p$ tel que $F^{k}=V$ pour tout $k \leq p$ est positif ou nul. Il est équivalent de dire qu'il existe un morphisme d'objets filtrés de $(V,Triv^{\bullet })$ vers $(V,F^{\bullet })$.
\end{definition}
La filtration $G^{\bullet}$ est donc positive. Les fibrés ainsi définis sont liés par le lemme suivant :\\

\begin{lemme}Soit $(V,F^{\bullet }_{0},F^{\bullet }_{1},F^{\bullet }_{2})$ un espace vectoriel muni de trois filtrations. Si $F_{1}^{\bullet}$ et $F_{2}^{\bullet}$ sont positives on a un morphisme injectif de ${\xi}_{\widetilde{{\bf P}^2}}(V,F_{0}^{\bullet },F_{1}^{\bullet },F_{2}^{\bullet },Triv^{\bullet })$ vers $e^{*}{\xi}_{{\bf P}^2}(V,F_{0}^{\bullet },F_{1}^{\bullet },F_{2}^{\bullet })$.\\
\end{lemme}

\begin{preuve}
La filtration $G^{\bullet}$ est positive donc il existe un morphisme d'objets filtr\'es de $(V,Triv^{\bullet})$ vers $(V,G^{\bullet})$. On en d\'eduit l'existence d'un morphisme de ${\xi}_{\widetilde{{\bf P}^2}}(V,F_{0}^{\bullet },F_{1}^{\bullet },F_{2}^{\bullet },Triv^{\bullet })$ vers\\
$e^{*}{\xi}_{{\bf P}^2}(V,F_{0}^{\bullet },F_{1}^{\bullet },F_{2}^{\bullet })$. Ce morphisme est un isomorphisme sur un ouvert donc est injectif.
\end{preuve}

En terme de coordonn\'ees, sur chacun des ouverts de $\widetilde{{\bf P}^2}$, on a la description qui suit : notons  $U_{1}'\cong U_{1}$ et $U_{2}'\cong U_{2}$ les images inverses de $U_{1}$ et $U_{2}$ par $e$. Sur les ouverts $U_{1}'$ et $U_{2}'$ il n'y a rien \`a montrer. Tout se passe au dessus de la carte $U_{0}$ dans ${\bf P}^2$ (rappelons qu'un recouvrement affine de $e^{-1}(U_{0})$ est donné par $U_{0}^{0}$ et $U_{0}^{1}$). Comme on l'a vu plus haut, pour un $A$-module $M$ et un morphisme d'anneaux $A \rightarrow B$, on a $f^{*}({\tilde M}) \cong \widetilde{M \otimes_{A} B }$. Par exactitude de $^\sim$ il suffit donc de montrer qu'on a une flèche de $B_{i}$-modules $RR(V,F_{i}^{\bullet },Triv^{\bullet }) \rightarrow  RR(V,F_{i}^{\bullet },F_{2}^{\bullet }) \otimes_{A} B_{i}$ pour $i \in \{1,2\}$ injective. Il suffit de le montrer pour $i=1$.

$RR(V,F_{1}^{\bullet },F_{2}^{\bullet }) \otimes_{A} B_{1}$ est engendré par les éléments de la forme $u^{-p}.v^{-q}.a_{p,q} \otimes_{A}P(x,y)$ où $P \in B_{1}$ et $a_{p,q} \in F_{1}^{p } \cap F_{2}^{q }$, or $u^{-p}.v^{-q}.a_{p,q} \otimes_{A}P(x,y)=a_{p,q} \otimes_{A} e_{1}(u^{-p}).e_{1}(v^{-q}).P(x,y)=a_{p,q} \otimes_{A} (x.y)^{-p}.y^{-q}.P(x,y)=a_{p,q} \otimes_{A} x^{-p}.y^{-p-q}.P(x,y)$. Prenons comme flèche le morphisme injectif : $R^{2}(V,F_{1}^{\bullet },Triv^{\bullet }) \rightarrow  R^{2}(V,F_{1}^{\bullet },F_{2}^{\bullet }) \otimes_{A} B_{1}$ donnée par $x^{-p}.a_{p} \mapsto a_{p} \otimes_{A}x^{-p}$ où $a_{p} \in  F_{1}^{p }$.\\
${}$\\
{\bf Remarque :}
Lorsque le fibré sera associé à une struture de Hodge $F_{1}^{\bullet }=F_{}^{\bullet }$ et $F_{2}^{\bullet }={\overline F}_{}^{\bullet }$ et les filtrations au rang $p$ sont données par des $p$-formes et sont donc positives.   
Notons $\mathcal F $ le faisceau quotient de $e^{*}{\xi}_{{\bf P}^2}(V,F_{0}^{\bullet },F_{1}^{\bullet },F_{2}^{\bullet })$ par ${\xi}_{\widetilde{{\bf P}^2}}(V,F_{0}^{\bullet },F_{1}^{\bullet },F_{2}^{\bullet },Triv^{\bullet })$, on a la suite exacte : 
\begin{center}
$0 \rightarrow {\xi}_{\widetilde{{\bf P}^2}}(V,F_{0}^{\bullet },F_{1}^{\bullet },F_{2}^{\bullet },Triv^{\bullet }) \rightarrow e^{*}{\xi}_{{\bf P}^2}(V,F_{0}^{\bullet },F_{1}^{\bullet },F_{2}^{\bullet }) \rightarrow {\mathcal F} \rightarrow 0$.
\end{center}
La construction de ${\xi}_{\widetilde{{\bf P}^2}}(V,F_{0}^{\bullet },F_{1}^{\bullet },F_{2}^{\bullet },Triv^{\bullet })$ est légitimée par le lemme suivant. Il permet de construire des suites exactes de fibrés sur $\widetilde{{\bf P}^2}$ associés aux sous-espaces et espaces quotients de la filtration $F^{\bullet }_{0}$ et donc de réduire l'étude de ${\xi}_{\widetilde{{\bf P}^2}}(V,F_{0}^{\bullet },F_{1}^{\bullet },F_{2}^{\bullet },Triv^{\bullet })$ puis de ${\xi}_{{\bf P}^2}(V,F_{0}^{\bullet },F_{1}^{\bullet },F_{2}^{\bullet })$ à celle des fibrés de Rees qui sont des sommes directes de fibrés en droite et donc de calculer les invariants topologiques du fibré de Rees associé à trois filtrations. \\

\begin{lemme}\label{estilde}%exactsequ.tilde
Soit $(V,F_{0}^{\bullet },F_{1}^{\bullet },F_{2}^{\bullet },Triv^{\bullet })$ un espace vectoriel muni de trois filtrations exhaustives et décroissantes, si $V'$ est un sous-espace vectoriel de $V$ donné par un des termes de la filtration $F_{0}^{\bullet }$(i.e. il existe $p$ tel que $V'=F_{0}^{p}$), alors, en notant par des $'$ les sous-objets et les filtrations
induites sur $V'$ et par des $''$ les objets quotients et les filtrations quotient sur $V''=V/V'$, on a la suite exacte :
\begin{center}
\scalebox{0.9}[1]{$0 \rightarrow {\xi}_{\widetilde{{\bf P}^2}}(V',{F'}_{0}^{\bullet },{F'}_{1}^{\bullet },{F'}_{2}^{\bullet },Triv^{\bullet }) \rightarrow {\xi}_{\widetilde{{\bf P}^2}}(V,F_{0}^{\bullet },F_{1}^{\bullet },F_{2}^{\bullet },Triv^{\bullet }) \rightarrow {\xi}_{\widetilde{{\bf P}^2}}(V'',{F''}_{0}^{\bullet },{F''}_{1}^{\bullet },{F''}_{2}^{\bullet },Triv^{\bullet }) \rightarrow 0$.} 
\end{center}
\end{lemme}

\begin{preuve}
La suite est évidement exacte sur chacun des quatre ouverts affines decrits plus haut. Le recollement est rendu possible par le fait que sur chaque intersection triple deux filtrations différentes sont en jeu et donc toutes les trivialisations sont compatibles.
\end{preuve}

{\bf Remarque :}
L'idée donnée dans la démonstration du lemme nous donne m\^eme un résultat plus fort : tous les quotients de $(V,F_{0}^{\bullet },F_{1}^{\bullet },F_{2}^{\bullet },Triv^{\bullet })$ gradués par $F^{\bullet }_{0}$ forment des sous-fibrés de ${\xi}_{\widetilde{{\bf P}^2}}(V,F_{0}^{\bullet },F_{1}^{\bullet },F_{2}^{\bullet },Triv^{\bullet })$ et leurs quotients sont de m\^eme des sous-fibrés. \\

Le lemme suivant est l'analogue du lemme \ref{scin} p.\pageref{scin} :\\

\begin{lemme}
Soit $V$ un espace vectoriel muni de trois filtrations $(F_{0}^{\bullet },F_{1}^{\bullet },F_{2}^{\bullet })$ scindées. Alors :
\begin{center}
${\xi}_{\widetilde{{\bf P}^2}}(V,F_{0}^{\bullet },F_{1}^{\bullet },F_{2}^{\bullet },Triv^{\bullet})=\oplus_{p,q}{\xi}_{\widetilde{{\bf P}^2}}(V^{p,q},Dec^{r}Triv^{\bullet },Dec^{p}Triv^{\bullet },Dec^{q}Triv^{\bullet },Triv^{\bullet})$.
\end{center}
\end{lemme}

\begin{preuve}
La preuve se calque sur celle du lemme \ref{scin}.
\end{preuve}

\subsection{Calcul explicite du caractère de Chern des fibrés de Rees}

Le but de cette section est de calculer le caractère de Chern du fibré ${\xi}_{{\bf P}^2}(V,F_{0}^{\bullet },F_{1}^{\bullet },F_{2}^{\bullet })$ associé à un
espace vectoriel trifiltré $(V,F_{0}^{\bullet },F_{1}^{\bullet },F_{2}^{\bullet }) \in \calc_{3filtr}$ afin de voir en quoi ce fibré diffère du fibré trivial obtenu pour la construction associée avec trois filtrations simultanéement scindées provenant des gradués donnés par les trois filtrations. Pour cel\`a, on utilise la décomposition de la section précédente pour calculer les caractères de Chern $\text{ch}({\mathcal F})$ et $\text{ch}({\xi}_{\widetilde{{\bf P}^2}}(V,F_{0}^{\bullet },F_{1}^{\bullet },F_{2}^{\bullet },Triv^{\bullet }))$.\\

Notons ${\widetilde{{D}_{i}}}$ la transformée stricte dans ${\widetilde{{\bf P}^2}}$ de $D_{i}$ donnée par l'application de l'éclatement $e$, et $E$ le diviseur exceptionnel.

\begin{proposition}
 Le calcul du caractère de Chern $\emph{ch}({\xi}_{\widetilde{{\bf P}^2}}(V,F_{0}^{\bullet },F_{1}^{\bullet },F_{2}^{\bullet },Triv^{\bullet }))$ peut toujours se ramener au calcul des caractères de Chern $\emph{ch}({\xi}_{\widetilde{{\bf P}^2}}(k,Dec^{r}Triv^{\bullet },Dec^{p}Triv^{\bullet },Dec^{q}Triv^{\bullet },Triv^{\bullet }))$ pour $(r,p,q) \in {\bf Z}^3$. De plus, en notant  $\eta_{D} \in H^{2}({\widetilde{{\bf P}^2}},{\bf Z})={\bf Z}$ le dual de Poincaré du diviseur $D$ et ${\tilde w}^{4}$ une forme qui engendre $H^{4}({\widetilde{{\bf P}^2}},{\bf Z})={\bf Z}$ :\\
$\emph{ch}({\xi}_{\widetilde{{\bf P}^2}}(k,Dec^{r}Triv^{\bullet },Dec^{p}Triv^{\bullet },Dec^{q}Triv^{\bullet },Triv^{\bullet }))=1+r {\eta}_{{\tilde D}_{0}}+ p{\eta}_{{\tilde D}_{1}}+q{\eta}_{{\tilde D}_{2}}+{\frac{1}{2}}(r^2+2rp+2rq){\tilde w}^{4}$
On a ainsi la relation :\\
$\emph{ch}({\xi}_{\widetilde{{\bf P}^2}}(V,F_{0}^{\bullet },F_{1}^{\bullet },F_{2}^{\bullet },Triv^{\bullet }))={\emph{dim}}_kV+ \sum_{p,q,r}{\delta }_{p,q,r}.(r {\eta}_{{\tilde D}_{0}}+ p{\eta}_{{\tilde D}_{1}}+q{\eta}_{{\tilde D}_{2}}+{\frac{1}{2}}(r^2+2rp+2rq){\tilde w}^{4})$
où ${\delta }_{p,q,r}={\emph{dim}}_kGr_{F_{2}^{\bullet }}^{q}Gr_{F_{1}^{\bullet }}^{p}Gr_{F_{0}^{\bullet }}^{r}V$.\\
\end{proposition}

Avant de démontrer cette proposition, établissons le lemme suivant :\\

\begin{lemme}$ $\\
${\xi}_{{\bf P}^2}(k,Dec^{r}Triv^{\bullet },Dec^{p}Triv^{\bullet },Dec^{q}Triv^{\bullet })\cong O_{{\bf P}^2}(rD_{0}+pD_{1}+qD_{2})$.\\ 
$c_{1}({\xi}_{{\bf P}^2}(k,Dec^{r}Triv^{\bullet },Dec^{p}Triv^{\bullet },Dec^{q}Triv^{\bullet }))=r \eta_{D_{0}}+p\eta_{D_{1}}+q\eta_{D_{2}}$.\\${\xi}_{\widetilde{{\bf P}^2}}(k,Dec^{r}Triv^{\bullet },Dec^{p}Triv^{\bullet },Dec^{q}Triv^{\bullet },Triv^{\bullet })\cong O_{\widetilde {{\bf P}^2}}(r{\tilde D}_{0}+p{\tilde D}_{1}+q{\tilde D}_{2})$.\\
$c_{1}({\xi}_{\widetilde{{\bf P}^2}}(k,Dec^{r}Triv^{\bullet },Dec^{p}Triv^{\bullet },Dec^{q}Triv^{\bullet },Triv^{\bullet }))=r {\eta}_{{\tilde D}_{0}}+ p{\eta}_{{\tilde D}_{1}}+q{\eta}_{{\tilde D}_{2}}.$\\

\end{lemme}

\begin{preuve} La première assertion a été démontrée dans l'exemple qui suit la proposition-definition \ref{recolledef} p.\pageref{recolledef}. Pour distinguer plus bas leurs transformées strictes dans l'éclaté $\widetilde{{\bf P}^2}$ de ${\bf P}^2$, on garde les notations des $D_{i}$ bien que tous soient homologues dans ${\bf P}^2$ et  définissent des fibrés isomorphes. D'après \cite{gh} par exemple, sur une variété compacte complexe $M$, la classe de Chern d'un fibré en droite de la forme $\calo_{M}(D)$ pour $D \in \text{Div}(M)$ est donnée par $c_{1}(\calo_{M}(D))= \eta_{D}$ où $\eta_{D} \in H^{2}_{DR}(M) $ est le dual de Poincaré de $D$. Ainsi :\\
$$c_{1}({\xi}_{{\bf P}^2}(V,Dec^{r}Triv^{\bullet },Dec^{p}Triv^{\bullet },Dec^{q}Triv^{\bullet }))= r.\eta_{D_{0}}+p.\eta_{D_{1}}+q.\eta_{D_{2}}.$$

Le fibré en droites ${\xi}_{\widetilde{{\bf P}^2}}(k,Dec^{r}Triv^{\bullet },Dec^{p}Triv^{\bullet },Dec^{q}Triv^{\bullet },Triv^{\bullet }))$ sur ${\widetilde{{\bf P}^2}}$ est de la forme $\calo_{\widetilde {{\bf P}^2}}(D)$ pour $[D] \in \text{Pic}({\widetilde {{\bf P}^2}})={\bf Z}\oplus {\bf Z}$. Un calcul analogue à celui fait dans l'exemple précité montre l'égalité $D=r{\tilde D}_{0}+p{\tilde D}_{1}+q{\tilde D}_{2}$ et donc : $$c_{1}({\xi}_{\widetilde{{\bf P}^2}}(k,Dec^{r}Triv^{\bullet },Dec^{p}Triv^{\bullet },Dec^{q}Triv^{\bullet },Triv^{\bullet }))=r {\eta}_{{\tilde D}_{0}}+ p{\eta}_{{\tilde D}_{1}}+q{\eta}_{{\tilde D}_{2}}.$$

\end{preuve}
Démontrons la proposition :
\begin{preuve} Démontrons la première assertion. On notera toujours par $Triv^{\bullet }$ la filtration triviale induite par la filtration triviale sur des sous-espaces ou des espaces quotients de $V$. Soit $V$ un espace vectoriel muni de trois filtrations opposées $(F_{0}^{\bullet },F_{1}^{\bullet },F_{2}^{\bullet })$. Rappelons que pour une suite exacte de faisceaux cohérents de la forme
 $0 \rightarrow {\mathcal G}' \rightarrow {\mathcal G} \rightarrow {\mathcal G}'' \rightarrow 0$ la relation suivante est vérifiée : $ \text{ch}( {\mathcal G})= \text{ch}({\mathcal G}')+ \text{ch}({\mathcal G}'')$. Ainsi en filtrant $V$ par $F_{0}^{\bullet }$ et d'après la suite exacte exhibée dans le lemme \ref{estilde} :

$$\text{ch}({\xi}_{\widetilde{{\bf P}^2}}(V,F_{0}^{\bullet },F_{1}^{\bullet },F_{2}^{\bullet },Triv^{\bullet }))=\sum_{r}\text{ch}({\xi}_{\widetilde{{\bf P}^2}}(Gr_{F_{0}^{\bullet }}^{r}V,F_{0,ind_{n}}^{\bullet },F_{1,ind_{r}}^{\bullet },{F_{2,ind_{}r}^{\bullet }},Triv^{\bullet })).$$

Comme $F_{0,ind_{r}}^{\bullet }$ est de longueur $1$ sur $Gr_{F_{0}^{\bullet }}^{r}V$ c'est la filtration décalée de $r$ par rapport à la filtration triviale, donc :

\hspace{1cm}$ \text{ch}({\xi}_{\widetilde{{\bf P}^2}}(Gr_{F_{0}^{\bullet }}^{r}V,F_{0,ind_{r}}^{\bullet },F_{1,ind_{r}}^{\bullet },{F_{2,ind_{r}}^{\bullet }},Triv^{\bullet }))=$

\hspace{2cm}  $ \text{ch}({\xi}_{\widetilde{{\bf P}^2}}(Gr_{F_{0}^{\bullet }}^{r}V,Dec^{r}Triv^{\bullet},F_{1,ind_{r}}^{\bullet },{F_{2,ind_{}r}^{\bullet }},Triv^{\bullet })).$

Pour tout $r$, filtrons $Gr_{F_{0}^{\bullet }}^{r}$ par la filtration induite par $F_{1,ind_{r}}^{\bullet }$ sur cet espace. Il vient alors :
 
\hspace{1cm}$ \text{ch}({\xi}_{\widetilde{{\bf P}^2}}(Gr_{F_{0}^{\bullet }}^{r}V,Dec^{r}Triv^{\bullet},F_{1,ind_{r}}^{\bullet },{F_{2,ind_{}r}^{\bullet }},Triv^{\bullet }))=$ 

\hspace{2cm}$\sum_{p}  \text{ch}({\xi}_{\widetilde{{\bf P}^2}}(Gr_{F_{1,ind_{r}}^{\bullet }}^{p}Gr_{F_{0}^{\bullet }}^{r}V,Dec^{r}Triv^{\bullet},F_{1,ind_{r},ind_{p}}^{\bullet },{F_{2,ind_{r},ind_{p}}^{\bullet }},Triv^{\bullet })).$

Or $F_{1,ind_{r},ind_{p}}^{\bullet }$ est de longueur $1$ sur $Gr_{F_{1,ind_{r}}^{\bullet }}^{p}Gr_{F_{0}^{\bullet }}^{r}V=Gr_{F_{1}^{\bullet }}^{p}Gr_{F_{0}^{\bullet }}^{r}V$, ainsi :

\hspace{1cm}$ \text{ch}({\xi}_{\widetilde{{\bf P}^2}}(Gr_{F_{1,ind_{r}}^{\bullet }}^{p}Gr_{F_{0}^{\bullet }}^{r}V,Dec^{r}Triv^{\bullet},F_{1,ind_{r},ind_{p}}^{\bullet },{F_{2,ind_{r},ind_{p}}^{\bullet }},Triv^{\bullet }))=$

\hspace{2cm}$ \text{ch}({\xi}_{\widetilde{{\bf P}^2}}(Gr_{F_{1}^{\bullet }}^{p}Gr_{F_{0}^{\bullet }}^{r}V,Dec^{r}Triv^{\bullet},Dec^{p}Triv^{\bullet },{F_{2,ind_{r},ind_{p}}^{\bullet }},Triv^{\bullet })).$

Le m\^eme argument appliqué à la dernière filtration permet de montrer que :

$ \text{ch}({\xi}_{\widetilde{{\bf P}^2}}(V,F_{0}^{\bullet },F_{1}^{\bullet },F_{2}^{\bullet },Triv^{\bullet }))=$

\hspace{1cm} $\sum_{p,q,r}  \text{ch}({\xi}_{\widetilde{{\bf P}^2}}(Gr_{F_{2}^{\bullet }}^{q}Gr_{F_{1}^{\bullet }}^{p}Gr_{F_{0}^{\bullet }}^{r}V,Dec^{r}Triv^{\bullet },Dec^{p}Triv^{\bullet },Dec^{q}Triv^{\bullet },Triv^{\bullet })).$

Puis finalement que :\\

$ \text{ch}({\xi}_{\widetilde{{\bf P}^2}}(V,F_{0}^{\bullet },F_{1}^{\bullet },F_{2}^{\bullet },Triv^{\bullet }))=$

\hspace{0.5cm}$ \sum_{p,q,r} \,  \text{dim}_k( Gr_{F_{2}^{\bullet }}^{q}Gr_{F_{1}^{\bullet }}^{p}Gr_{F_{0}^{\bullet }}^{r}V).\, \text{ch}({\xi}_{\widetilde{{\bf P}^2}}(k,Dec^{r}Triv^{\bullet },Dec^{p}Triv^{\bullet },Dec^{q}Triv^{\bullet },Triv^{\bullet })).$\\

Ce qui démontre la première assertion.\\
La deuxième découle du lemme précédent. Pour un fibré en droite $\mathcal L$, on a la relation $ \text{ch}({\mathcal L})=\exp(c_{1}({\mathcal L}))$ c'est à dire sur une surface : $ \text{ch}({\mathcal L})=1+c_{1}({\mathcal L})+{\frac{1}{2}}c_{1}({\mathcal L})^{2}$. Les produits d'intersection dans ${\bf {P^2}}$ sont donnés par ${{\tilde D}_{0}}^{2}=1,\,{{\tilde D}_{1}}^{2}=0,\,{{\tilde D}_{2}}^{2}=0,\,E^{2}=-1$ et ${{\tilde D}_{0}}.{{\tilde D}_{1}}=1,\,{{\tilde D}_{0}}.{{\tilde D}_{2}}=1,\,E.{{\tilde D}_{1}}=1,\,E.{{\tilde D}_{2}}=1,\,{{\tilde D}_{1}}.{{\tilde D}_{2}}=0$. D'où :\\

$ \text{ch}({\xi}_{\widetilde{{\bf P}^2}}(k,Dec^{r}Triv^{\bullet },Dec^{p}Triv^{\bullet },Dec^{q}Triv^{\bullet },Triv^{\bullet }))=$\\
\hspace{0.5cm} $1+r {\eta}_{{\tilde D}_{0}}+ p{\eta}_{{\tilde D}_{1}}+q{\eta}_{{\tilde D}_{2}}+{\frac{1}{2}}(r^2+2rp+2rq){\tilde w}^{4}.$
%où $w$ engendre $H^{4}({\widetilde{{\bf P}^2}},{\bf Z})={\bf Z}$.

%$F_{1,ind_{n}}^{\bullet }$ et ${F_{2,ind_{n}}^{\bullet }}$ sont $(-n)$-opposées (à rappeler...cf \cite{del2}) sur $Gr_{F_{0}^{\bullet }}^{n}V$ et définissent donc une décomposition de $Gr_{F_{0}^{\bullet }}^{n}V$ en somme directe: $Gr_{F_{0}^{\bullet }}^{n}V=\oplus_{p+q=-n}V_{n}^{p,q}$ et les sous espaces $V_{n}^{p,q}$ sont tels que $F_{1,ind_{n}}^{p }=\oplus_{p'\leq p}V_{n}^{p',-n+p'}$ et ${F_{2,ind_{n}}^{q }}=\oplus_{q'\leq q}V_{n}^{-n+q',q'}$. Sur chacun des sous-espaces en somme directe $F_{1,ind_{n},ind}^{\bullet }$ et ${F_{2,ind_{n},ind}^{\bullet }}$ sont de longueur $1$, donc:
%\begin{center}
%$\sum_{n}ch({\xi}_{\widetilde{{\bf P}^2}}(Gr_{F_{0}^{\bullet }}^{n}V,Dec^{n}Triv^{\bullet},F_{1,ind_{n}}^{\bullet },{F_{2,ind_{}n}^{\bullet }},Triv^{\bullet }))=\sum_{n}\sum_{p+q=-n}ch({\xi}_{\widetilde{{\bf P}^2}}(V_{n}^{p,q},Dec^{n}Triv^{\bullet},Dec^{p}Triv^{\bullet },Dec^{q}Triv^{\bullet },Triv^{\bullet }))=\sum_{n}\sum_{p+q=-n}\,dim_kV_{n}^{p,q}\,ch({\xi}_{\widetilde{{\bf P}^2}}(k,Dec^{n}Triv^{\bullet},Dec^{p}Triv^{\bullet },Dec^{q}Triv^{\bullet },Triv^{\bullet }))$
%\end{center}
%Prouvons la deuxième assertion:\\
%Nous allons prouver l'assertion plus générale:
%$$ch({\xi}_{\widetilde{{\bf P}^2}}(k,Dec^{p}Triv^{\bullet },Dec^{q}Triv^{\bullet },Dec^{r}Triv^{\bullet }))=...$$
%pour $p,q,r \in {\bf Z}$. 

\end{preuve}
Reste à déterminer $ \text{ch}(\mathcal F)$. L'idée pour le faire est de voir que $\mathcal F $ ne dépend pas de la filtration $F_{0}^{\bullet }$. En effet le support de $\mathcal F $ est dans l'éclaté de la carte affine $U_{0}$ et $\calf $ peut \^etre déterminé par les restrictions des faisceaux  ${\xi}_{{\bf P}^2}(V,F_{0}^{\bullet },F_{1}^{\bullet },F_{2}^{\bullet })$ et ${\xi}_{{\bf P}^2}(V,F_{0}^{\bullet },F_{1}^{\bullet },F_{2}^{\bullet })$ à l'éclaté de cet ouvert affine. Or ces restrictions sont construites à partir uniquement des filtrations $F^{\bullet }_{1}$, $F^{\bullet }_{2}$ et $Triv^{\bullet }$ et ne dépendent pas de la filtration $F_{0}^{\bullet }$. Donc on peut écrire une suite exacte courte :
 $$0 \rightarrow {\xi}_{\widetilde{{\bf P}^2}}(V,F^{\bullet },F_{1}^{\bullet },F_{2}^{\bullet },Triv^{\bullet }) \rightarrow e^{*}{\xi}_{{\bf P}^2}(V,F^{\bullet },F_{1}^{\bullet },F_{2}^{\bullet }) \rightarrow {\mathcal F} \rightarrow 0$$ 
pour n'importe quelle filtration exhaustive décroissante $F^{\bullet }$ de $V$ avec le m\^eme faisceau $\calf$. Ainsi, on va choisir une filtration $F^{\bullet }$ telle que $ \text{ch}({\xi}_{{\bf P}^2}(V,F^{\bullet },F_{1}^{\bullet },F_{2}^{\bullet }))$ soit facilement calculable. Pour cela il faut qu'il y ait un scindement commun de $V$ par les trois filtrations $F^{\bullet }$, $F_{1}^{\bullet }$ et $F_{2}^{\bullet }$ ce qui est toujours le cas lorsqu'il n'y a que deux filtrations, comme dans le cas où $F^{\bullet }=F_{1}^{\bullet }$ ou $F^{\bullet }=F_{2}^{\bullet }$ par exemple.\\

\begin{proposition}
Le caractère de Chern de $\calf$ est $ \emph{ch}(\calf)=\frac{1}{2}(p+q)^{2}{\tilde w}^{4}$.\\
\end{proposition}

\begin{preuve}
  On a le choix sur $F^{\bullet }$ pour effectuer le calcul. Celui-ci est facilité si l'on prend par exemple $F^{\bullet }=F_{1}^{\bullet }$. 

La suite exacte courte de faisceaux explicitée plus haut nous donne :

$$ \text{ch}(\calf)= \text{ch}(e^{*}{\xi}_{{\bf P}^2}(V,F_{1}^{\bullet },F_{1}^{\bullet },F_{2}^{\bullet })) -\text{ch}({\xi}_{\widetilde{{\bf P}^2}}(V,F_{1}^{\bullet },F_{1}^{\bullet },F_{2}^{\bullet },Triv^{\bullet })).$$
D'après la proposition précédente, en posant ${\delta }_{p,q}= \text{dim}_kGr_{F_{2}^{\bullet }}^{q}Gr_{F_{1}^{\bullet }}^{p}V$ :

 $  \text{ch}({\xi}_{\widetilde{{\bf P}^2}}(V,F_{1}^{\bullet },F_{1}^{\bullet },F_{2}^{\bullet },Triv^{\bullet }))= \text{dim}_kV+ \sum_{p,q}{\delta }_{p,q}.(p {\eta}_{{\tilde D}_{0}}+ p{\eta}_{{\tilde D}_{1}}+q{\eta}_{{\tilde D}_{2}}+{\frac{1}{2}}(p^2+2p^{2}+2pq)w)$.\\
D'un isomorphisme $V \cong Gr_{F_{2}^{\bullet }}^{q}Gr_{F_{1}^{\bullet }}^{p}V$ (rappelons qu'il en existe toujours un) et du lemme \ref{isosumfibrees} on tire :
%              \begin {eqnarray*}
%                |v(z)|^2 & \le &
%                        C_1\left( \ep^2 \|\db v\|^2_{L^{\infty}B(z,\ep)} +
%                        \ep^{-2n} \|v\|^2_{L^2(B(z,\ep))}\right)\\
%                 & \le & C_1\left( N
%                ^2\ep^2 \|\db h\
%                h^{N-1}\|^2_{L^\infty(B(z,\ep))} + \ep^{-2n} \|v e^{-
%                \frac{\psi}{2}} 
%              e^{ \frac{\psi}{2} } \|^2_{L^2(B(z,\ep))} \right).\\
%                                \end{eqnarray*}
\begin{eqnarray*}
{\xi}_{{\bf P}^2}(V,F_{1}^{\bullet },F_{1}^{\bullet },F_{2}^{\bullet }) & \cong & {\xi}_{{\bf P}^2}(\oplus_{p,q}Gr_{F_{2}^{\bullet }}^{q}Gr_{F_{1}^{\bullet }}^{p}V,\oplus_{p,q}Dec^{p}Triv^{\bullet }_{p,q},\oplus_{p,q}{Dec^{p}Triv^{\bullet }}_{p,q},\oplus_{p,q}Dec^{q}Triv^{\bullet }_{p,q})\\
& \cong &  \oplus_{p,q}{\xi}_{{\bf P}^2}(k,Dec^{p}Triv^{\bullet},Dec^{p}Triv^{\bullet},Dec^{q}Triv^{\bullet})^{dim_kGr_{F_{2}^{\bullet }}^{q}Gr_{F_{1}^{\bullet }}^{p}V} .
\end{eqnarray*}
Donc, si l'on note $w^{4}$ le générateur de $H^{4}({{\bf P}^2},{\bf Z})$ tel que $e^{*}(w^{4})={\tilde w}^{4}$ :\\
${}$\\
$ \text{ch}(e^{*}{\xi}_{{\bf P}^2}(V,F_{1}^{\bullet },F_{1}^{\bullet },F_{2}^{\bullet }))$
\begin{eqnarray*}
&=&e^{*} \text{ch}({\xi}_{{\bf P}^2}(V,F_{1}^{\bullet },F_{1}^{\bullet },F_{2}^{\bullet }))\\
&=&e^{*} \text{ch}(\oplus_{p,q}{\xi}_{{\bf P}^2}(k,Dec^{p}Triv^{\bullet},Dec^{p}Triv^{\bullet},Dec^{q}Triv^{\bullet})^{ \text{dim}_kGr_{F_{2}^{\bullet }}^{q}Gr_{F_{1}^{\bullet }}^{p}V})\\
&=&e^{*}( \text{dim}_kV+\sum_{p,q}{\delta}_{p,q}(p {\eta}_{{ D}_{0}}+ p{\eta}_{{ D}_{1}}+q{\eta}_{{ D}_{2}}+{\frac{1}{2}}(p^2+p^{2}+q^{2}+2(p^{2}+pq+pq))w^{4})\\
&=& \text{dim}_kV+\sum_{p,q}{\delta}_{p,q}(p {\eta}_{{\tilde D}_{0}}+ p{\eta}_{{\tilde D}_{1}}+q{\eta}_{{\tilde D}_{2}}+{\frac{1}{2}}(p^{2}+p^{2}+q^{2}+2(p^{2}+pq+pq)){\tilde w}^{4}).
\end{eqnarray*}
D'après la formule : $ \text{ch}({\xi}_{{\bf P}^2}(k,Dec^{r}Triv^{\bullet},Dec^{p}Triv^{\bullet},Dec^{q}Triv^{\bullet})=1+r \eta_{D_{0}}+p\eta_{D_{1}}+q\eta_{D_{2}}$

\hspace{6.5cm} $+{\frac{1}{2}} (r^{2}+p^{2}+q^{2}+2(rp+rq+pq))w^{4}$).\\
\\
Ainsi :
$$ \text{ch}(\calf)=\frac{1}{2}\sum_{p,q}{\delta}_{p,q}(p+q)^{2}{\tilde w}^{4}.$$

\end{preuve}
Finalement, le caractère de Chern du fibré ${\xi}_{{\bf P}^2}(V,F_{0}^{\bullet },F_{1}^{\bullet },F_{2}^{\bullet })$ est donnée par :\\
\begin{proposition}\label{cherncaracter}%\hspace{1cm}
$ \emph{ch}({\xi}_{{\bf P}^2}(V,F_{0}^{\bullet },F_{1}^{\bullet },F_{2}^{\bullet })=\emph{dim}_kV$
\begin{eqnarray*}
\hspace{2cm}&+&\sum_{p,q,r}\Bigl[(\emph{dim}_kGr_{F_{2}^{\bullet }}^{q}Gr_{F_{1}^{\bullet }}^{p}Gr_{F_{0}^{\bullet }}^{r}V).\Bigl((r+p+q){w}^{2}+{\frac{1}{2}}(r^2+2rp+2rq){ w}^{4}\Bigr)\Bigr]\\
\hspace{2cm}&+&\sum_{p,q}\Bigr[(\emph{dim}_kGr_{F_{2}^{\bullet }}^{q}Gr_{F_{1}^{\bullet }}^{p}V).({\frac{1}{2}}(p+q)^{2}{ w}^{4})\Bigr]
\end{eqnarray*}
\end{proposition}

\begin{preuve}
C'est une conséquence directe du calcul de $ \text{ch}(\calf)$ et de $ \text{ch}\,{\xi}_{\widetilde{{\bf P}^2}}(V,F^{\bullet }_{1},F_{1}^{\bullet },F_{2}^{\bullet },Triv^{\bullet })$.
\end{preuve}
%\textit{Remarque :}
%stabilite des fibre obtenus en fonction des classes de Chern et du rang d'après Drezet et Le Potier cite?
Nous allons effectuer ce calcul dans le cas qui nous intéressera par la suite c'est à dire celui où les trois filtrations sont opposées. Soit $V$ un espace vectoriel muni de trois filtrations $(F_{0}^{\bullet },F_{1}^{\bullet },F_{2}^{\bullet })$ exhaustives,finies et décroissantes. On rappelle que ces filtrations sont dites opposées si $Gr_{F_{1}^{\bullet }}^{p}Gr_{F_{2}^{\bullet }}^{q}Gr_{F_{0}^{\bullet }}^{n}=0$ pour $p+q+n \neq 0$. On note $\calc_{3filtr,opp}$ la catégorie des espaces vectoriels munis de trois filtrations décroissantes, exhaustives et opposées munie des morphismes compatibles aux filtrations.\\ 

%\textit{Remarque :}
%Si trois filtrations sur $V$ sont scindées alors elles sont opposées. En effet ...???\\                    
\begin{corollaire}
\label{filtrscin}
Soit $(V,F_{0}^{\bullet },F_{1}^{\bullet },F_{2}^{\bullet }) \in \calc_{3filtr,opp}$ alors :\\
$\,$\\
\hspace{1cm} $ \emph{ch}({\xi}_{{\bf P}^2}(V,F_{0}^{\bullet },F_{1}^{\bullet },F_{2}^{\bullet })= \emph{dim}_kV$

\hspace{3cm} $+{\frac{1}{2}}\sum_{p,q}\Bigl[ \emph{dim}_kGr_{F_{2}^{\bullet }}^{q}Gr_{F_{1}^{\bullet }}^{p}V- \emph{dim}_kGr_{F_{2}^{\bullet }}^{q}Gr_{F_{1}^{\bullet }}^{p}Gr_{F_{0}^{\bullet }}^{-p-q}V \Bigr].(p+q)^{2}.w^{4}$.\\
\end{corollaire}

\begin{preuve}
On applique la formule de la proposition précédente. Tous les coefficients sont nuls si $p+q+n \neq 0$.  
\end{preuve}
${}$\\
{\bf Remarque :} Lorsque $(F_{0}^{\bullet },F_{1}^{\bullet },F_{2}^{\bullet })$ sont trois filtrations opposées sur $V$, on a :

$$\text{c}_{1}({\xi}_{{\bf P}^2}(V,F_{0}^{\bullet },F_{1}^{\bullet },F_{2}^{\bullet }) )=-\text{ch}_{1}({\xi}_{{\bf P}^2}(V,F_{0}^{\bullet },F_{1}^{\bullet },F_{2}^{\bullet }) )=0.$$

\subsection{Suites exactes de faisceaux associ\'ees aux suites exactes d'espaces trifiltr\'es}

Considérons une suite exacte courte dans la catégorie abélienne des espaces vectoriels munis de trois filtrations opposées $\calc_{3filtr,opp}$ :

$$  0 \rightarrow (V',F_{0}^{\bullet }{}',F_{1}^{\bullet }{}',F_{2}^{\bullet }{}')  \rightarrow (V,F_{0}^{\bullet }{},F_{1}^{\bullet }{},F_{2}^{\bullet }{})  \rightarrow (V'',F_{0}^{\bullet }{}'',F_{1}^{\bullet }{}'',F_{2}^{\bullet }{}'')  \rightarrow 0$$
et son image par la construction des modules de Rees associés :
$$  0 \rightarrow R^{3}(V',F_{0}^{\bullet }{}',F_{1}^{\bullet }{}',F_{2}^{\bullet }{}')  \rightarrow R^{3}(V,F_{0}^{\bullet }{},F_{1}^{\bullet }{},F_{2}^{\bullet }{})  \rightarrow R^{3}(V'',F_{0}^{\bullet }{}'',F_{1}^{\bullet }{}'',F_{2}^{\bullet }{}'')  \rightarrow 0 $$
ou de façon équivalente son image par le foncteur qui à tout $k[u,v,w]$-module gradué fini associe un faisceau cohérent sur ${\bf P}^2$ : 
$$  0 \rightarrow \xi_{{\bf P}^2}(V',F_{0}^{\bullet }{}',F_{1}^{\bullet }{}',F_{2}^{\bullet }{}')  \rightarrow \xi_{{\bf P}^2}(V,F_{0}^{\bullet }{},F_{1}^{\bullet }{},F_{2}^{\bullet }{})  \rightarrow \xi_{{\bf P}^2}(V'',F_{0}^{\bullet }{}'',F_{1}^{\bullet }{}'',F_{2}^{\bullet }{}'')  \rightarrow 0 $$
que nous noterons pour simplifier :
$$ 0 \rightarrow \xi '\rightarrow \xi \rightarrow \xi '' \rightarrow 0 .$$
Cette suite n'est pas exacte en général, mais on a seulement la suite exacte :
$$ 0 \rightarrow \xi'\rightarrow \xi \rightarrow \xi''.$$
Notons $(\xi/\xi')^{faisc}$ le quotient faisceautique des fibrés $\xi $ et $\xi'$ qui s'inscrit dans la suite exacte courte $ 0 \rightarrow \xi '\rightarrow \xi \rightarrow \xi/\xi' \rightarrow 0 $. D'après l'étude de la section \ref{sectionreflexif} p.\pageref{sectionreflexif} sur les suite exactes de faisceaux réflexifs (i.e. de fibrés car on est sur une surface), on a le diagramme commutatif dans lequel toutes les suites sont exactes 
$$
\xymatrix{
{}&{}&{} & 0 \ar[d]\\
0 \ar[r] & \xi' \ar[r] &  \xi  \ar[r] \ar[rd] &  (\xi/\xi')^{faisc} \ar[d] \ar[r]  & 0\\
{} & {  }   & {}  &   \xi''\ar[d]\\
{} & {} &   {}    &  \cals \ar[d] \\
{}  &  {}  &  {}  &  0
}
$$

Les caractéristiques de Chern des fibrés de Rees associés aux espaces vectoriels trifiltrés qui forment la suite exacte verifient ainsi l'égalité

$$\text{ch}( \xi_{{\bf P}^2}(V',F_{0}^{\bullet }{}',F_{1}^{\bullet }{}',F_{2}^{\bullet }{}')) + \text{ch} (\xi_{{\bf P}^2}(V'',F_{0}^{\bullet }{}'',F_{1}^{\bullet }{}'',F_{2}^{\bullet }{}''))= \text{ch} (\xi_{{\bf P}^2}(V,F_{0}^{\bullet }{},F_{1}^{\bullet }{},F_{2}^{\bullet }{})) + \text{ch}(\cals) .$$

\begin{proposition} Si $(V,F_{0}^{\bullet }{},F_{1}^{\bullet }{},F_{2}^{\bullet }{}) $, $(V',F_{0}^{\bullet }{}',F_{1}^{\bullet }{}',F_{2}^{\bullet }{}')$ et $(V'',F_{0}^{\bullet }{}'',F_{1}^{\bullet }{}'',F_{2}^{\bullet }{}''))$ sont des éléments de $\calc_{3filtr,opp}$, alors    $$\emph{ch}_{2}( \xi_{{\bf P}^2}(V',F_{0}^{\bullet }{}',F_{1}^{\bullet }{}',F_{2}^{\bullet }{}'))+\emph{ch}_{2}(\xi_{{\bf P}^2}(V'',F_{0}^{\bullet }{}'',F_{1}^{\bullet }{}'',F_{2}^{\bullet }{}'')) \leq \emph{ch}_{2}(\xi_{{\bf P}^2}(V,F_{0}^{\bullet }{},F_{1}^{\bullet }{},F_{2}^{\bullet }{})).$$
\end{proposition}
La preuve de cette proposition se déduit de l'étude du signe de $\text{ch}_{2}(\cals)$ faite dans le lemme ci-dessous.\\

\begin{lemme}
$$\emph{ch}_{2}(\cals) \geq 0.$$
 \end{lemme}

\begin{preuve}(du lemme) Il faut prouver que si $\calf$ est un faisceau cohérent sans torsion et si $\cals$ est le conoyau du morphisme canonique $\nu : \calf \rightarrow \calf^{**}$ (qui est injectif car $\calf$ est sans torsion), alors $\text{ch}_{2}(\cals) \leq 0$, ou de façon équivalente, puisque $\text{ch}_{1}(\cals)=0$, $\text{c}_{2}(\cals)\leq 0$. 

\end{preuve}

La suite 
$$  0 \rightarrow \xi_{{\bf P}^2}(V',F_{0}^{\bullet }{}',F_{1}^{\bullet }{}',F_{2}^{\bullet }{}')  \rightarrow \xi_{{\bf P}^2}(V',F_{0}^{\bullet }{},F_{1}^{\bullet }{},F_{2}^{\bullet }{})  \rightarrow \xi_{{\bf P}^2}(V'',F_{0}^{\bullet }{}'',F_{1}^{\bullet }{}'',F_{2}^{\bullet }{}'')  \rightarrow 0 $$
est en fait exacte dans la catégorie des fibrés vectoriels dans le sens où $\xi_{{\bf P}^2}(V'',F_{0}^{\bullet }{}'',F_{1}^{\bullet }{}'',F_{2}^{\bullet }{}'')$ est le conoyau de $ \xi_{{\bf P}^2}(V',F_{0}^{\bullet }{}',F_{1}^{\bullet }{}',F_{2}^{\bullet }{}')  \rightarrow \xi_{{\bf P}^2}(V,F_{0}^{\bullet }{},F_{1}^{\bullet }{},F_{2}^{\bullet }{})$ dans la catégorie des fibrés vectoriels obtenu en prenant le bidual du conoyau faisceautique $((\xi/\xi')^{faisc})^{**}$, comme on l'a vu dans la démonstration du théorème \ref{muab} p.\pageref{muab}.\\ 
${}$\\
{\bf Remarque :} D'après la remarque qui suit le corollaire \ref{filtrscin} p.\pageref{filtrscin}, comme le $\text{ch}_1$ de chacun des fibrés est nul parce que les filtrations sont opposées, la proposition peut s'écrire
 $$\emph{c}_{2}( \xi_{{\bf P}^2}(V',F_{0}^{\bullet }{}',F_{1}^{\bullet }{}',F_{2}^{\bullet }{}'))+\emph{c}_{2}(\xi_{{\bf P}^2}(V'',F_{0}^{\bullet }{}'',F_{1}^{\bullet }{}'',F_{2}^{\bullet }{}'')) \leq \emph{c}_{2}(\xi_{{\bf P}^2}(V,F_{0}^{\bullet }{},F_{1}^{\bullet }{},F_{2}^{\bullet }{})).$$

\begin{corollaire}
Soit $(V,F_{0}^{\bullet }{},F_{1}^{\bullet }{},F_{2}^{\bullet }{}) \in \calc_{3filtr,opp}$. Alors
$$\text{c}_{2}(\xi_{{\bf P}^2}(V,F_{0}^{\bullet }{},F_{1}^{\bullet }{},F_{2}^{\bullet }{})) \geq 0.$$
\end{corollaire}

\begin{preuve} ``Coupons'' $V$ par l'une de ses trois filtrations, $F^{\bullet}_{0}$ par exemple de façon à avoir une suite exacte dans $\calc_{3filtr,opp}$ de la forme 
$$ 0 \rightarrow (F^{p}_{0}V,F_{0}^{\bullet }{}_{inj},F_{1}^{\bullet }{}_{inj},F_{2}^{\bullet }{}_{inj}) \rightarrow (V,F_{0}^{\bullet }{},F_{1}^{\bullet }{},F_{2}^{\bullet }{}) \rightarrow (V/F^{p}_{0}V,F_{0}^{\bullet }{}_{ind},F_{1}^{\bullet }{}_{ind},F_{2}^{\bullet }{}_{ind}) \rightarrow 0.$$
D'après la proposition précédente, on a l'inégalité 
$$\text{c}_{2}( \xi_{{\bf P}^2}(F^{p}_{0}V,F_{0}^{\bullet }{}_{inj},F_{1}^{\bullet }{}_{inj},F_{2}^{\bullet }{}_{inj}))+\text{c}_{2}(\xi_{{\bf P}^2}(V/F^{p}_{0}V,F_{0}^{\bullet }{}_{ind},F_{1}^{\bullet }{}_{ind},F_{2}^{\bullet }{}_{ind})) \leq \text{c}_{2}(\xi_{{\bf P}^2}(V,F_{0}^{\bullet }{},F_{1}^{\bullet }{},F_{2}^{\bullet }{}).$$
On peut à nouveau minorer les deux termes de gauche en coupant les deux filtrations obtenus à partir de la première. Comme cette filtration est exhaustive en coupant un nombre fini de fois, on peut arriver à ce que les filtrations soit de longueur $1$ i.e. triviales. La filtration triviale se scinde évidemment simultanéement avec les deux autres filtrations, les fibrés obtenus sont donc tous triviaux, donc de deuxième classe de Chern nulle. D'où l'inégalité.
\end{preuve}
${}$\\
{\bf Exemple :} Donnons un exemple sur un espace vectoriel $V=<e,f>$ sur $k$ de dimension $2$. On définit les filtrations (décroissantes et exhaustives) par\\
\hspace{0.5cm}${\bullet}$ $F^{-2}_{0}=V$, $F^{-1}_{0}=F^{0}_{0}=<e>$ et $F^{1}_{0}=\{ 0\}$.\\
\hspace{0.5cm}${\bullet}$ $F^{0}_{1}=V$, $F^{1}_{1}=<f+\lambda. e>$ et $F^{2}_{1}=\{ 0\}$ où $\lambda \in k$.\\
\hspace{0.5cm}${\bullet}$ $F^{0}_{2}=V$, $F^{1}_{2}=<f+\kappa. e>$ et $F^{2}_{2}=\{ 0\}$ où $\kappa \in k$.\\
On vérifie alors que ces trois filtrations sont opposées sur $V$ et que\\ 
\hspace{0.5cm} $\bullet$ $\text{dim}_{k}Gr_{F_{2}^{\bullet }}^{q}Gr_{F_{1}^{\bullet }}^{p}Gr_{F_{0}^{\bullet }}^{-p-q}V=\left \lbrace \begin{array}{l}
1 \text{ si } p=q=1 \\
1 \text{ si } p=q=0 \\
0 \text{ sinon ,}
\end{array} \right.$\\
si $\lambda=\kappa$,\\
\hspace{0.5cm}$\bullet$ $\text{dim}_{k}Gr_{F_{2}^{\bullet }}^{q}Gr_{F_{1}^{\bullet }}^{p}V=\left \lbrace \begin{array}{l}
1 \text{ si } p=q=1 \\
1 \text{ si } p=q=0 \\
0 \text{ sinon ,}
\end{array} \right.$\\
et si $\lambda \neq \kappa$,\\
\hspace{0.5cm}$\bullet$ $\text{dim}_{k}Gr_{F_{2}^{\bullet }}^{q}Gr_{F_{1}^{\bullet }}^{p}V=\left \lbrace \begin{array}{l}
1 \text{ si } p=0, q=1 \\
1 \text{ si } p=1, q=0 \\
0 \text{ sinon .}
\end{array} \right.$\\
D'où le calcul de la deuxième classe de Chern\\
$$\text{c}_{2}(\xi_{{\bf P}^2}(V,F_{0}^{\bullet }{},F_{1}^{\bullet }{},F_{2}^{\bullet }{}))=\left|\begin{array}{l}
0 \text{ si } \lambda=\kappa\\
1 \text{ si } \lambda \neq \kappa.
\end{array}\right.
$$ 
On observe que sur l'espace des configurations des filtrations de cet exemple, $k^2$ (le choix de $(\lambda,\kappa) \in k^2$), la classe de Chern du fibré est génériquement $1$ et qu'elle est zéro lorsque qu'on peut scinder simultanéement les filtrations, i.e. sur la strate fermée d'équation $\lambda=\kappa$.    
%{\frac{1}{2}}\sum_{p,q}\Bigl[ \emph{dim}_kGr_{F_{2}^{\bullet }}^{q}Gr_{F_{1}^{\bullet }}^{p}V- \emph{dim}_kGr_{F_{2}^{\bullet }}^{q}Gr_{F_{1}^{\bullet }}^{p}Gr_{F_{0}^{\bullet }}^{-p-q}V \Bigr].(p+q)^{2}.w^{4}$.\\

\subsection{${\bf P}^1$-semistabilité et $\mu$-semistabilit\'e}

La donn\'ee ``filtrations oppos\'ees'' pour les espaces vectoriels trifiltr\'es se traduit pour les fibr\'es de Rees par une condition de semistabilit\'e que nous allons expliciter dans cette section. Rappelons que la propriété que les trois filtrations de l'espace vectoriel $V$, $F_{0}^{\bullet }{},F_{1}^{\bullet }{}$ et $F_{2}^{\bullet }{}$ soient opposées est équivalente au fait que pour tout $n \in {\bf Z}$, $F_{1}^{\bullet }{}$ et $F_{2}^{\bullet }{}$ induisent des filtrations $n$-opposées sur $Gr_{F^{\bullet }_{0}}^{-n}$. La condition de semistabilité des fibrés de Rees associée aux filtrations opposées est une semistabilité de la restriction du fibré au diviseur qui correspond à la filtration $F^{\bullet }_{0}$, droite que nous noterons ${\bf P}^{1}_{0}$. Notons $j$ le morphisme
$$j : {\bf P}^{1}_{0} \hookrightarrow {\bf P}^{2}.$$

\begin{definition}
Un fibré $\calf$ sur ${\bf P}^2$ est dit ${\bf P}^{1}_{0}$-semistable si pour tout sous-faisceau cohérent $\cale$ sur ${\bf P}^{1}_{0}$ tel que $\cale \subset j^{*}\calf$ on a 
$$\mu(\cale) \leq \mu(j^{*}\calf).$$
\end{definition}    
Comme $j^{*}$ induit un isomorphisme de $H^{2}({\bf P}^{2},{\bf Z})$ vers $H^{2}({\bf P}^{1}_{0},{\bf Z})$ on a pour tout faisceau cohérent $\calf$ sur ${\bf P}^2$, $\text{c}_{1}(j^{*}\calf)=j^{*}\text{c}_{1}(\calf)$ et donc, 
$$ {\bf P}^{1}_{0}\text{-semistable } \Rightarrow \mu\text{-semistable }.$$

Rappelons (cf \cite{gro1}) que tout fibré $\calf$ sur ${\bf P}^1$ se décompose de manière unique à l'ordre des termes près en somme de fibrés en droite $\calf \cong \oplus_{i}\,\calo_{{\bf P}^1}(a_{i})^{p_{i}}$. Un tel fibré est donc semistable seulement si tous les poids $a_{i}$ de la d\'ecomposition sont les mêmes.

\begin{proposition}\label{p1sstable}
 Les fibrés de Rees associés à des filtrations opposées sont ${\bf P}^{1}_{0}$-semistables.
\end{proposition}     
 
\begin{preuve}
Soit $\calf=\xi_{{\bf P}^2}(V,F_{0}^{\bullet }{},F_{1}^{\bullet }{},F_{2}^{\bullet }{}))$ un tel fibré. Il suffit de montrer que $\calf \cong \calo_{{\bf P}^1}(a)^{\text{rang}\calf}$. Ce qui est immédiat par la proposition \ref{resfibree} car la restriction est alors isomorphe à un fibr\'e de Rees sur ${\bf P}^1$ qui est trivial donc semistable puisque $F_{1}^{\bullet }{}$ et $F_{2}^{\bullet }{}$ induisent des filtration $r$-oppos\'ees sur $Gr_{F^{\bullet }_{0}}^{r}$.
\end{preuve}

On en déduit, d'après la relation d'ordre entre semistabilités donnée au dessus, les fibrés de Rees associés à des filtrations opposées sont $\mu$-semistable.

On peut aussi démontrer directement la $\mu$-semistabilité des fibrés de Rees associés à des filtrations opposées directement sans passer par la proposition précédente qui utilise la proposition \ref{resfibree}.

 Rappelons que sur toute vari\'et\'e lisse projective $X$, si $\calf$ est un faisceau coh\'erent de pure dimension $\text{dim}_{k}(X)$ alors il existe un sous-faisceau coh\'erent de $\calf$, $\cale_{max}$ tel que pour tout sous-faisceau coh\'erent $\cale$ de $\calf$, $\mu(\cale)\leq \mu(\cale_{max})$ et en cas d'\'egalit\'e $\cale \subset \cale_{max}$. Le faisceau $\cale_{max}$ est appel\'e le sous-faisceau d\'estabilisant maximal de $\calf$.

\begin{proposition}\label{semistab}
Les fibr\'es de Rees associ\'es aux espaces vectoriels trifiltr\'es dont les filtrations sont oppos\'ees sont $\mu$-semistables de degr\'e $0$.
\end{proposition}

\begin{preuve}
Soit $\calf=\xi_{{\bf P}^2}(V,F_{0}^{\bullet }{},F_{1}^{\bullet }{},F_{2}^{\bullet }{}))$ le fibr\'e de Rees associ\'e \`a l'espace vectoriel muni de trois filtrations oppos\'ees $(V,F_{0}^{\bullet }{},F_{1}^{\bullet }{},F_{2}^{\bullet }{}))$. Par le calcul du caract\`ere de Chern fait au dessus, nous savons que $\calf$ est de degr\'e $0$.

Pour montrer que $\calf$ est $\mu$-semistable nous allons proc\'eder en deux \'etapes. D'abord nous allons d\'emontrer qu'il est $\mu$-semistable dans la cat\'egorie des ${\bf T}$-faisceaux \'equivariants i.e. que tout sous-faisceau ${\bf T}$-\'equivariant a un degr\'e inferieur \`a z\'ero. Nous montrerons ensuite que le faisceau d\'estabilisant peut \^etre muni d'une action de ${\bf T}$ ce qui permettra de conclure.

Soit $\cale$ un sous-faisceau coh\'erent ${\bf T}$-\'equivariant de $\calf$. Alors la d\'ecomposition suivant les caract\`eres de l'action du tore sur $\cale$ est incluse dans la d\'ecomposition pour $\calf$, ce qui prouve que les filtrations associ\'ees \`a cette action sur $\cale$ sont oppos\'ees, et donc que le degr\'e de $\cale$ est $0$.

Soit $\cale_{max}$ le sous-faisceau d\'estabilisant maximal de $\calf$. Montrons qu'il peut-\^etre muni d'une structure de ${\bf T}$-faisceau telle que le morphisme
$$\xymatrix{0 \ar[r]&\cale_{max}\ar[r]&\calf}$$
soit un ${\bf T}$-morphisme. Consid\'erons pour tout $t \in {\bf T}$ le morphisme
$$\xymatrix{f_{t},g_{t}:\{ t \} \times {{\bf P}^2} \ar[r]^-{i} & {\bf T} \times {{\bf P}^2} \ar[r]^-{p_{2},\sigma} & {{\bf P}^2}},
$$
o\`u pour $f$, la deuxi\`eme fl\`eche est $p_2$ et pour $g$, elle est $\sigma$. Alors pour tout $t$, $f_{t}^{*}\cale_{max}$ est le sous-faisceau d\'estabilisant maximal de $f_{t}^{*}\calf$ ($p_{2}$ correspond \`a l'action triviale). De m\^eme comme d'apr\`es \cite{huyleh}, lemme $1.3.5$, le sous-faisceau d\'estabilisant maximal est unique, donc $g_{t}^{*}\cale_{max}$ est le sous-faisceau d\'estabilisant maximal de $g_{t}^{*}\calf$. Or $\calf $ est un ${\bf T}$-faisceau donc si $\Psi$ est l'isomorphisme sur ${\bf T}\times {\bf P}^2$ entre l'image r\'eciproque par le morphisme de l'action et l'image r\'eciproque par le morphisme de projection, il vient,
$$
\xymatrix{
\Psi_{t}:  g_{t}^{*}\calf \ar[r] &f_{t}^{*}\calf
}.$$
Ce qui permet d'exhiber pour tout $t \in {\bf T}$, par unicit\'e du sous-faisceau d\'estabilisant maximal, l'isomorphisme
$$
\xymatrix{
\Psi_{t}\vert_{g_{t}^{*}\cale_{max}}:  g_{t}^{*}\cale_{max} \ar[r] &f_{t}^{*}\cale_{max}
}.$$
$\cale_{max}$ est donc un ${\bf T}$-faisceau.

Tout sous-faisceau $\cale $ de $\calf$ v\'erifie $\mu(\cale)\leq \mu(\cale_{max})=0$ ce qui permet de conclure.

\end{preuve}
${}$\\
{\bf Remarque :} L'action du tore l'ensemble des classes d'isomorphismes de faisceaux coh\'erents sans torsion fixe par unicit\'e la classe d'isomorphisme du sous-faisceau d\'estabilisant maximal. Pour montrer que le faisceau est muni d'une action on utilise le fait que $H^{2}({\bf T },\text{Aut}, d'après (\cale_{max}))=0$(\cite{demgab}).

\subsection{Dictionnaire : Espaces vectoriels filtrés-Fibrés vectoriels }

Dans cette section nous rassemblons et complétons les résultats des sections précédentes afin d'écrire un dictionnaire entre les espaces vectoriels $n$-filtrés et les fibrés vectoriels. Nous avons déjà établi à travers les propositions \ref{inverse1} p.\pageref{inverse1} et \ref{inverse2} p.\pageref{inverse2} des équivalences de catégories entre les catégories des espaces vectoriels filtrés ou bifiltrés munies des morphismes strictement compatibles aux filtrations et les catégories des fibrés ${\bf G}_m$ ou $ {\bf G}_{m}^2$-équivariants munies des morphismes de fibrés vectoriels dont les conoyaux sont sans torsion. 

Nous voulons établir ce type d'équivalence entre \index{espaces} vectoriels trifiltrés et fibrés vectoriels sur ${\bf P}^2$. Rappelons que l'on note :\\
\hspace*{0.5cm}$\bullet$ $\calc_{nfiltr}$ la catégorie des espaces vectoriels munis de $n$-filtrations exhaustives et décroissantes et dont les morphismes sont les morphismes strictement compatibles aux filtrations.\\
\hspace*{0.5cm}$\bullet$ $\calc_{nfiltr,scind}$ la sous-catégorie pleine de $\calc_{nfiltr}$ dont les objets sont les espaces vectoriels $n$-filtrés dont les filtrations sont simultanéement scindées.\\
\hspace*{0.5cm}$\bullet$ $\calc_{3filtr,opp}$ la sous-catégorie pleine de $\calc_{3filtr}$ dont les objets sont les espaces vectoriels trifiltrés dont les filtrations sont opposées.\\
\hspace*{0.5cm}$\bullet$ $\calc_{3filtr,opp,scind}$ la sous-catégorie pleine de $\calc_{3filtr,opp}$ dont les objets sont les espaces vectoriels trifiltrés dont les filtrations opposées sont simultanéement scindées.\\

On rappelle que ${\bf T}={\bf G}_{m}^{3}/\Delta({\bf G}_{m}^{3})$ est le tore quotienté par la diagonale. ${\bf P}^{2}=\text{Proj}[u_{0},u_{1},u_{2}]$ est fix\'e et on note ${\bf P}^{1}_{0}$ le diviseur qui correspond à $\text{Proj}[u_{1},u_{2}]$. $P$ est le point $(1:0:0)$.

\begin{theoreme}\label{th2}La construction du fibré de Rees sur ${\bf P}^2$ associé à un espace vectoriel trifiltré établit des équivalences de catégories entre :\\
\hspace*{0.1cm}$\bullet$ La catégorie des espaces vectoriels de dimension finie munis de trois filtrations décroissantes et exhaustives et des morphismes strictement compatibles aux filtrations $\calc_{3filtr}$ est équivalente à la catégorie des fibrés vectoriels sur le plan projectif ${\bf P}^2$ munis de l'action héritée de l'action par translation sur ${\bf A}^{3} \backslash \{(0,0,0)\}$ de ${{\bf G}_m}^3$ dont la diagonale agit trivialement et des morphismes équivariants de fibrés dont les singularités du conoyau sont en codimension $2$, supportées aux points $(1:0:0)$, $(0:1:0)$ et $(0:0:1)$, ${\calf ib}({\bf P}^{2}/{\bf T})$ :
$$
\xymatrix{
\{ \calc_{3filtr} \} \ar@<2pt>[r]^-{\Phi_{R}} &  \{  {\calf ib}({\bf P}^{2}/{\bf T}) \} \ar@<2pt>[l]^-{\Phi_{I}}
}
.$$
\hspace*{0.1cm}$\bullet$ La catégorie des espaces vectoriels de dimension finie munis de trois filtrations décroissantes, exhaustives et simultanéement scindées et des morphismes strictement compatibles aux filtrations $\calc_{3filtr,scind}$ est équivalente à la catégorie des fibrés vectoriels sommes de fibrés en droite sur le plan projectif ${\bf P}^2$ munis de l'action héritée de l'action par translation sur ${\bf A}^{3} \backslash \{(0,0,0)\}$ de ${{\bf G}_m}^3$ dont la diagonale agit trivialement et des morphismes équivariants de fibrés, ${\calf ib}_{scind}({\bf P}^{2}/{\bf T})$ :
$$
\xymatrix{
\{ \calc_{3filtr,scind} \} \ar@<2pt>[r]^-{\Phi_{R}} &  \{  {\calf ib}_{scind}({\bf P}^{2}/{\bf T}) \} \ar@<2pt>[l]^-{\Phi_{I}}
}
.$$
\hspace*{0.1cm}$\bullet$ La catégorie des espaces vectoriels de dimension finie munis de trois filtrations décroissantes, exhaustives et opposées et des morphismes strictement compatibles aux filtrations $\calc_{3filtr,opp}$ est équivalente à la catégorie des fibrés vectoriels semistables de pente $0$ sur le plan projectif ${\bf P}^2$ munis de l'action héritée de l'action par translation sur ${\bf A}^{3} \backslash \{(0,0,0)\}$ de ${{\bf G}_m}^3$ dont la diagonale agit trivialement et des morphismes équivariants de fibrés dont les singularités du conoyau sont en codimension $2$, supportées au point $P=(1:0:0)$, ${\calf ib}_{\mu-semistable,\mu=0}({\bf P}^{2}/{\bf T})$ :
$$
\xymatrix{
\{ \calc_{3filtr,opp} \} \ar@<2pt>[r]^-{\Phi_{R}} &  \{  {\calf ib}_{{\bf P}^{1}_{0}-semistable,\mu=0}({\bf P}^{2}/{\bf T}) \} \ar@<2pt>[l]^-{\Phi_{I}}
}
.$$

\end{theoreme}  

\begin{preuve}
\hspace*{0.1cm}$\bullet$ Montrons la première équivalence. Pour vérifier que les foncteurs $\Phi_{R}$ et $\Phi_{I}$ établissent une équivalence de catégorie, commençons par prouver que $\Psi_{R}$ est essentiellement surjective. Soit $\xi$ un élément de ${\calf ib}_{}({\bf P}^{2}/{\bf T})$. Les restrictions de $\xi$ aux ouverts affines standards $\calu_{k}=\{ (u_{0},u_{1},u_{2}) \in {\bf P}^{2} \vert u_{k}\neq 0 \}$ pour tout $k \in \{0,1,2 \}$ sont des fibrés vectoriels $({\bf G}_{m})^2$-équivariants, noté $\xi_k$. D'après la proposition \ref{inverse2}, chacun de ces fibrés $\xi_k$ permet de définir un couple de filtrations $(F^{\bullet }_{k},G^{\bullet }_{k})$ sur le même espace vectoriel donné par la fibre en $(1:1:1)$ de $\xi$, $V=\xi_{(1:1:1)}$. Comme l'action sur $\xi_k$ et $\xi_{k'}$ est la même sur les restrictions de ces fibrés à $\calu_{k} \cap \calu_{k'}$ ($\xi_{k}\vert_{\calu_{k'}}=\xi_{k'}\vert_{\calu_{k}}=\xi\vert_{\calu_{k}\cap \calu_{k'}}$), les filtrations données par l'action sur les fibrés restreints sont les mêmes. Ainsi $G^{\bullet}_{0}=F^{\bullet }_{1}$, $G^{\bullet}_{1}=F^{\bullet }_{2}$ et $G^{\bullet}_{2}=F^{\bullet }_{0}$. $\xi $ est donc isomorphe au fibré obtenu en recollant par l'action les fibrés
$({\bf G}_{m})^2$-équivariants sur les cartes affines obtenus à partir des couples de filtrations, donc par la proposition-définition \ref{recolledef} p.\pageref{recolledef}, $\xi \cong \xi_{{\bf P}^2}(V,F^{\bullet }_{0},F^{\bullet }_{1},F^{\bullet}_{2})$ comme fibrés vectoriels sur le plan projectif puis de façon directe comme ${\bf T}$-fibrés équivariants. D'où l'essentielle surjectivité de $\Phi_{R}$.    

Vérifions que l'image d'un morphisme $f$ entre deux objets $(V,F^{\bullet}_{0},F^{\bullet }_{1},F^{\bullet }_{2})$ et $(W,G^{\bullet}_{0},G^{\bullet }_{1},G^{\bullet }_{2})$ dans $\calc_{3filtr}$ par $\Phi_{R}$ est bien un morphisme dans ${\calf ib}_{}({\bf P}^{2}/{\bf T})$. $\Phi_{R}(f)$ est bien par construction un morphisme ${\bf T} $-équivariant. Reste à vérifier que son conoyau est sans torsion. La question est locale, plaçons-nous sur les ouverts affines standards $\calu_k$ décrits plus haut, sur $\calu_0$ par exemple. Alors le conoyau de $\Phi_{R}(f)\vert_{\calu_0}$ est sans torsion car ce morphisme de fibrés vectoriels ${{\bf G}_m}^2$-équivariants est l'image par le foncteur de Rees sur les espaces bifiltrés, déjà noté $\Phi_{R}$, du morphisme strictement compatible aux filtrations induit par $f$ en oubliant la première filtration  $\tilde{f}:(V,F^{\bullet }_{1},F^{\bullet }_{2}) \rightarrow (W,G^{\bullet }_{1},G^{\bullet }_{2})$ dans $\calc_{2filtr}$. Or la proposition \ref{inverse2} affirme que le conoyau de $\Phi_{R}(\tilde{f})$, isomorphe au conoyau de $\Phi_{R}(f)\vert_{\calu_{0}}$, est sans torsion. 

Reste à voir que le foncteur $\Phi_R$ établit une correspondance pleinement fidèle. Soient $(V,F^{\bullet}_{0},F^{\bullet }_{1},F^{\bullet }_{2})$ et $(W,G^{\bullet}_{0},G^{\bullet }_{1},G^{\bullet }_{2})$ deux objets dans $\calc_{3filtr}$. Soient $f,f'$ deux éléments de \begin{center}$\text{Hom}_{\calc_{3filtr}}( (V,F^{\bullet}_{0},F^{\bullet }_{1},F^{\bullet }_{2}),(W,G^{\bullet}_{0},G^{\bullet }_{1},G^{\bullet }_{2}))$ tels que $\Phi_{R}(f)=\Phi_{R}(f')$\end{center} dans
\begin{center}
 $\text{Hom}_{{\calf ib}_{}({\bf P}^{2}/{\bf T})}(\Phi_{R}((V,F^{\bullet}_{0},F^{\bullet }_{1},F^{\bullet }_{2})),\Phi_{R}((W,G^{\bullet}_{0},G^{\bullet }_{1},G^{\bullet }_{2})))$.
\end{center}
En se restreignant aux ouverts affines standards $\calu_{k}$, $\Phi_{R}(f)\vert_{\calu_k}=\Phi_{R}(f')\vert_{\calu_k}$ signifie d'après la proposition \ref{inverse2} que les restrictions de $f$ et $f'$ aux espaces vectoriels bifiltrés par oubli de la filtration $F^{\bullet }_{k}$ sont égales, ceci sur les trois ouverts, ce qui signifie que $f$ et $f'$ coïncident sur chacunes des filtrations, i.e. pour tout $k \in \{0,1,2 \}$ et tout $p \in {\bf Z}$, $f(F_{k}^{p})=f'(F^{p}_{k})$) et donc $f$ et $f'$ sont égaux, d'où l'injectivité. Montrons la surjectivité. Soit $g \in  \text{Hom}_{{\calf ib}_{}({\bf P}^{2}/{\bf T})}(\xi,\xi'))$. Les restrictions $g_k$ aux ouverts affines $\calu_k$ de $g$ sont des morphismes $({\bf G}_{m})^2$-équivariants de $({\bf G}_{m})^2$-fibrés vectoriels sur les ouverts affines dont le conoyau est sans torsion puisque $g$ est sans torsion. Ainsi pour tout $k$, $g_{k} \in \text{Hom}_{{\calf ib}({{\bf A}^2}/({{\bf G}_m})^{2})}(\xi\vert_{\calu_{k}},\xi'\vert_{{\calu}_k})$ a un antécédent $f_{k} \in \text{Hom}_{\calc_{2filtr}}(\Phi_{I}(\xi\vert_{\calu_{k}}),\Phi_{I}(\xi'\vert_{{\calu}_k}))$ par la proposition \ref{inverse2}. Les morphismes d'espaces vectoriels bifiltrés $f_k$ ont deux à deux une filtration en commun, ce qui permet de définir un morphisme d'espace vectoriel trifiltré $f$ qui coïncide avec les $f_k$ lorsqu'il est restreint aux espaces bifiltrés. On  a bien $\Phi_{R}(f)=g$ ce qui permet de conclure. (On aurait pu aussi utiliser la proposition \ref{inverse1} appliquée aux fibrés ${\bf G}_m$-équivariants sur les intersections deux à deux des ouverts standards $\calu_{k} \cap \calu_{k'}$).

\hspace*{0.1cm}$\bullet$ Deuxième équivalence. L'image par $\Phi_R$ d'un espace vectoriel trifiltré scindé est clairement un fibré équivariant sur le plan projectif somme directe de fibrés en droite équivariants.

Soit $\xi \in   {\calf ib}_{scind}({\bf P}^{2}/{\bf T})$. $\xi$ est somme de fibrés en droite ${\bf T}$-équivariants $\xi=\oplus_{i \in I}\, \xi_i$. Par l'équivalence précédente, il existe un espace vectoriel trifiltré $(V_{i},F^{\bullet }_{i}{}_{0},F^{\bullet }_{i}{}_{1},F^{\bullet }_{i}{}_{2})$ de dimension $1$ tel que $\xi_{i}\cong \Phi_{R}((V_{i},F^{\bullet }_{i}{}_{0},F^{\bullet }_{i}{}_{1},F^{\bullet }_{i}{}_{2}))$. Posons $( V,F^{\bullet }{}_{0},F^{\bullet }{}_{1},F^{\bullet }{}_{2})=\oplus_{i \in I}\,(V_{i},F^{\bullet }_{i}{}_{0},F^{\bullet }_{i}{}_{1},F^{\bullet }_{i}{}_{2})$. On alors, par le lemme \ref{isosumfibrees} p.\pageref{isosumfibrees} l'isomorphisme de fibrés équivariants sur le plan projectif $\xi \cong \oplus_{i \in I} \xi_{i} \cong \oplus_{i \in I}  \Phi_{R}((V_{i},F^{\bullet }_{i}{}_{0},F^{\bullet }_{i}{}_{1},F^{\bullet }_{i}{}_{2})) \cong \Phi_{R}( V,F^{\bullet }{}_{0},F^{\bullet }{}_{1},F^{\bullet }{}_{2})$ qui prouve l'essentielle surjectivité de $\Phi_R$.

La fidélité de $\Psi_{R}$ est immédiate. En effet, si $f$ et $f'$ sont deux morphismes d'espaces vectoriels trifiltrés scindés tels que les morphismes de fibrés sommes de fibrés en droites équivariants $\Phi_{R}(f)$ et $\Phi_{R}(f')$ coïncident, alors $f$ et $f'$ sont égaux dans $\calc_{3filtr}$, donc dans $\calc_{3filtr,scind}$. Pour la surjectivité, on peut remarquer que si $\xi $ et $\xi'$ sont des éléments de $ {\calf ib}_{scind}({\bf P}^{2}/{\bf T})$, alors, comme $\xi =\oplus_{i \in I}\, \xi_{i}$ et $\xi' =\oplus_{j \in J}\, \xi_{j}'$, on a \begin{center}$\text{Hom}_{ {\calf ib}_{scind}}({\bf P}^{2}/{\bf T})( \xi, \xi')= \oplus_{(i,j) \in I \times J}\text{Hom}_{ {\calf ib}_{scind}({\bf P}^{2}/{\bf T})}( \xi_{i}, \xi_{j}').$ \end{center} Soit $g$ un morphisme dont on veut trouver un antécédent. $g$ se décompose par blocs en $g^{ij}$. Par la première équivalence de catégories exhibée au dessus, chacun des termes $g^{ij}$ de cette décomposition admet un antécédent par $\Phi_{R}$, que nous notons $f^{ij}$. Le morphisme d'espaces vectoriels trifiltrés scindés obtenu par la famille des $f^{ij}$ a clairement $g$ pour image, ce qui montre la surjectivité. D'où la pleine fid\'elité.

\hspace*{0.1cm}$\bullet$ Troisième équivalence. Soit $\xi \in  {\calf ib}_{{\bf P}^{1}_{0}-semistable,\mu=0}({\bf P}^{2}/{\bf T})$. %\marginpar{Je crois qu'il faut prendre une action tq $F_{1},F_{2}$ filtr positives par rapport $F_{0}$....?????????}

Considérons le foncteur $\Phi_R $ de la catégorie des filtrations opposées $ \calc_{3filtr,opp}$ vers la catégorie des fibrés équivariants semistables de pente nulle $ {\calf ib}_{{\bf P}^{1}_{0}-semistable,\mu=0}({\bf P}^{2}/{\bf T})$. Ce foncteur est bien d\'efinit par la proposition \ref{p1sstable} p.\pageref{p1sstable}, en effet les fibr\'es de Rees associ\'es aux espaces vectoriels munis de trois filtrations oppos\'ees sont ${\bf P}^{1}_{0}$-semistables par cette proposition. Ils sont de pente nulle d'apr\`es les calculs effectu\'es dans la section ``calcul explicite du caract\`ere de Chern des fibr\'es de Rees''.  On voit directement par les caractérisations précédentes que $\Phi_R$ est injectif. La surjectivité de $\Phi_R$ est également immédiate.    

Reste \`a voir que le foncteur $\Phi_{I}$ a bien pour images des \'el\'ements de $\calc_{3filtr,opp}$. Par la premi\`ere équivalence de cat\'egories, \`a partir d'un fibr\'e ${\bf T}$-\'equivariant $\xi_{{\bf P}^2}$ on obtient par le foncteur $\Phi_{I}$ un espace vectoriel trifiltr\'e $(V,F_{0}^{\bullet },F_{1}^{\bullet },F_{2}^{\bullet })$. Montrons que ces filtrations sont oppos\'ees i.e. que $Gr_{F_{1}}^{p}Gr_{F_{2}}^{q}Gr_{F_{0}}^{r}V\neq0$ implique $p+q+r=0$. C'est \'equivalent \`a montrer que $Gr_{F_{1}}^{p}Gr_{F_{2}}^{q}Gr_{F_{0}}^{r}V=0$ pour tout $p+q+r\neq0$ ce que nous allons d\'ecomposer en deux parties, la partie positive (demi-espace strictement positif), $p+q+r>0$, et la partie n\'egative $p+q+r<0$. Supposons qu'il existe un triplet $(p_{0},q_{0},r_{0})$ dans le demi-espace strictement positif tel que l'espace trigradu\'e associ\'e ne soit pas nul. Alors, on a un morphisme injectif d'espaces vectoriels trifiltr\'es 
$$\xymatrix{0 \ar[r]  & (k,Dec^{r_0}Triv^{\bullet},Dec^{p_0}Triv^{\bullet},Dec^{q_0}Triv^{\bullet}) \ar[r] & (V,F_{0}^{\bullet },F_{1}^{\bullet },F_{2}^{\bullet }).}$$
 Ce morphisme m\`ene \`a un morphisme injectif de ${\bf T}$-fibr\'es sur le plan projectif$$
\xymatrix{
0 \ar[r] & \xi_{{\bf P}^2}^{(r_{0},p_{0},q_{0})} \ar[r] &  \xi_{{\bf P}^2} }.$$
Or $\mu(\xi_{{\bf P}^2}^{(r_{0},p_{0},q_{0})})=r_{0}+p_{0}+q_{0}>0$ ce qui contredit la $\mu$-semistabilit\'e de $\xi_{{\bf P}^2}$ donc la ${\bf P}^{1}_{0}$-semistabilité de ce fibré. Les espaces trigradu\'es correspondants \`a des indices dans le demi-espace strictement positif sont donc nuls.

Revenons \`a la formule du degr\'e du fibr\'e $\xi_{{\bf P}^2}$. On a
$$ \text{c}_{1}({\xi}_{{\bf P}^2})=\sum_{p,q,r}(\emph{dim}_kGr_{F_{2}^{\bullet }}^{q}Gr_{F_{1}^{\bullet }}^{p}Gr_{F_{0}^{\bullet }}^{r}V).(r+p+q){w}^{2}.$$
Comme pour tout triplet $(r,p,q)$, $p+q+r \leq 0$, la somme est nulle si et seulement si tous les termes sont nuls ce qui permet de conclure. Les filtrations sont bien oppos\'ees.   

\end{preuve}
%\section{Constructions de fibrés associés à un nombre fini de filtrations}
%Généralisons la construction précédente. Soit $V$ un espace vectoriel muni de $n$ filtrations décroissantes et exhaustives :$(V,F_{0}^{\bullet },F_{1}^{\bullet },...,F_{n-1}^{\bullet })$ ... \\q
%Justification de la construction des fibrés...\\
%\\
%Fibrés sur ${\bf P}^1$ et cohomologie quantique()

%\xymatrix{
%\{ \calc_{3filtr,scind} \} \ar@<2pt>[r]^-{\Phi_{R}} &  \{  {\text Bun}_{split}({\bf P}^{2}/({{\bf G}_m}^{3}/\Delta({{\bf G}_m})) \} \ar@<2pt>[l]^-{\Phi_{I}}
%}
%.$$
%\hspace*{0.1cm}$\bullet$ La catégorie des espaces vectoriels de dimension finie munis de trois filtrations décroissantes, exhaustives et opposées et des morphismes strictement compatibles aux filtrations $\calc_{3filtr,opp}$ est équivalente à la catégorie des fibrés vectoriels semistables de pente $0$ sur le plan projectif ${\bf P}^2$ munis de l'action héritée de l'action par translation sur ${\bf A}^{3} \backslash \{(0,0,0)\}$ de ${{\bf G}_m}^3$ dont la diagonale agit trivialement et des morphismes équivariants de fibré dont le conoyau est sans torsion $\text{Bun}_{split}({\bf P}^{2}/({{\bf G}_m}^{3}/\Delta({{\bf G}_m}))$ :
%$$
%\xymatrix{
%\{ \calc_{3filtr,opp} \} \ar@<2pt>[r]^-{\Phi_{R}} &  \{  {\text Bun}_{semistable,\mu=0}({\bf P}^{2}/({{\bf G}_m}^{3}/\Delta({{\bf G}_m})) \} \ar@<2pt>[l]^-{\Phi_{I}}
%}
%.$$
Dégageons une première conséquence de ce théorème :

\begin{corollaire}\label{geotrifab}
La catégorie $\calc_{3filtr,opp}$ est abélienne.
\end{corollaire}

\begin{preuve}
Ce corollaire se déduit du fait que ${\calf ib}_{{\bf P}^{1}_{0}-semistable,\mu=0}({\bf P}^{2}/{\bf T})$ est abélienne. Montrons cette affirmation ; on se place dans le contexte du théorème \ref{muab} avec $X={\bf P}^2$, les faisceaux réflexifs sont donc localement libres. Ce dernier théorème et son corollaire \ref{gbunab} montrent que ${\calf ib}_{\mu-semistable,\mu=0}({\bf P}^{2}/{\bf T})$ est abélienne. Pour montrer que la sous-catégorie\\ ${\calf ib}_{{\bf P}^{1}_{0}-semistable,\mu=0}({\bf P}^{2}/{\bf T})$ est abélienne il suffit de voir que le noyau et le conoyau d'un morphisme $f:\cale \rightarrow \calf$ dans cette catégorie sont aussi dans cette catégorie i.e. que le noyau et le conoyau sont ${\bf P}^{1}_{0}$-semistables. 

Considérons la suite exacte dans ${\calf ib}_{\mu-semistable,\mu=0}({\bf P}^{2}/{\bf T})$
$$0 \rightarrow \calk er \rightarrow \cale \rightarrow \calf \rightarrow \calc oker \rightarrow 0.$$

La restriction au diviseur ${\bf P}^{1}_{0}$ est un foncteur exact car le seul point de torsion est le point $(0:0:1)$. D'où la suite exacte dans la catégorie des faisceaux cohérents sur la doite projective (avec $j:{\bf P}^{1}_{0} \rightarrow {\bf P}^2$)
$$0 \rightarrow j^{*}\calk er \rightarrow j^{*}\cale \rightarrow j^{*}\calf \rightarrow j^{*}\calc oker \rightarrow 0.$$

La semistabilité de $j^{*}\calk er$ se déduit directement de la semistabilité de $j^{*}\cale$. Soit $\calh$ un sous-faisceau cohérent de $j^{*}\calc oker $ on peut alors écrire les deux suites exactes $0 \rightarrow \calh \rightarrow j^{*}\calc oker \rightarrow j^{*}\calc oker/\calh \rightarrow 0$ et $0 \rightarrow \calk \rightarrow j^{*}\calf \rightarrow j^{*}\calc oker/\calh \rightarrow 0$ où $\calk $ est le noyau du morphisme canonique $ j^{*}\calf \rightarrow j^{*}\calc oker/\calh $. Comme $j^{*}\calc oker$ et $j^{*}\calf$ sont de degré $0$, $\mu(\calk) \leq 0$ implique $\mu(\calh) \leq 0$ ce qui permet de conclure.    
% Le corollaire \ref{gbunab} montre que ${\calf ib}_{{\bf P}^{1}_{0}-semistable,\mu=0}({\bf P}^{2}/{\bf T})$ est abélienne ce qui permet de conclure par l'équivalence de catégorie précédente. 

\end{preuve}

${\bullet}$ Définissons la catégorie $\calc_{3filtr,opp}^{-}$ dont les objets sont les espaces vectoriels trifiltrés dont les filtrations sont opposées et dont les morphismes sont les morphismes compatibles aux trois filtrations. Cette dernière catégorie est la catégorie étudiée par Deligne dans \cite{del2}.\\

La condition sur les morphismes est donc plus faible que celle exigée pour les morphismes de $\calc_{3filtr,opp}^{}$ ; cette dernière catégorie est une sous-catégorie de $\calc_{3filtr,opp}^{-}$. Nous allons voir que la condition de stricte compatibilité est automatiquement vérifiée pour les morphismes de $\calc_{3filtr,opp}^{-}$, autrement dit que tout morphisme dans $\calc_{3filtr,opp}^{-}$ est strictement compatible aux filtrations, c'est à dire $\calc_{3filtr,opp}^{-}=\calc_{3filtr,opp}^{}$, ce qui est le point essentiel du théorème de Deligne (\cite{del2}, Th. (1.2.10) p.12), que le théorème précédent nous permet de retrouver de façon géométrique. 

\begin{corollaire}\label{thdel} Soit $\calc_{3filtr,opp}^{-}$ la catégorie définie ci-dessus.\\
\hspace*{1cm}$(i)$ $\calc_{3filtr,opp}^{-}$ est une catégorie abélienne.\\
\hspace*{1cm}$(ii)$ Le noyau (resp. conoyau) d'un morphisme $f:V \rightarrow W$ dans $\calc_{3filtr,opp}^{-}$ est le noyau (resp. conoyau) du morphisme d'espaces vectoriels sous-jacent muni des filtrations induites par celles de $V$ (resp. celles de $W$).\\
\hspace*{1cm}$(iii)$ Tout morphisme $f: V \rightarrow W$  dans $\calc_{3filtr,opp}^{-}$ est strictement compatible aux trois filtrations des triplets de filtrations ; le morphisme $Gr_{F_{0}^{\bullet }}(f)$ est compatible aux bigraduations de $Gr_{F_{0}^{\bullet }}(V)$ et $Gr_{F_{0}^{\bullet }}(W)$ ; les morphismes $Gr_{F_{1}^{\bullet }}(f)$ et $Gr_{F_{2}^{\bullet }}(f)$ sont strictement compatibles à la filtration induite par $F^{\bullet }_{0}$.\\
\hspace*{1cm}$(iv)$ Les foncteurs oubli des filtrations, $Gr_{F_{0}^{\bullet }}$, $Gr_{F_{1}^{\bullet }}$, $Gr_{F_{2}^{\bullet }}$, et $$Gr_{F_{0}^{\bullet }}Gr_{F_{1}^{\bullet }} \cong  Gr_{F_{1}^{\bullet }}Gr_{F_{0}^{\bullet }} \cong Gr_{F_{2}^{\bullet }}Gr_{F_{1}^{\bullet }}Gr_{F_{0}^{\bullet }} \cong Gr_{F_{2}^{\bullet }}Gr_{F_{0}^{\bullet }} \cong Gr_{F_{0}^{\bullet }}Gr_{F_{2}^{\bullet }}$$ sont exacts de la catégorie $\calc_{3filtr,opp}^{-}$ vers la catégorie des $k$-espaces vectoriels.
\end{corollaire}

\begin{preuve}$(iii)$ Commençons par démontrer la première assertion de $(iii)$. Considérons le morphisme $f: V \rightarrow W$ dans $\calc_{3filtr,opp}^{-}$ et $\Phi_{R}(f)$ le morphisme entre fibrés vectoriels qui lui est associé. Ce morphisme de fibrés vectoriels est de rang constant sauf éventuellement en un nombre fini de points (i.e. un ensemble de codimension $2$). En effet, si ce n'est pas le cas, le conoyau a de la torsion donc a une première classe de Chern strictement positive et donc son noyau aurait une première classe de Chern strictement positive ce qui contredit la $\mu$-semistabilité et donc la ${\bf P}^{1}_{0}$-semistabilité de la source du morphisme entre fibrés vectoriels. On en déduit en utilisant le proposition $14$ que les morphismes d'objets filtrés sont strictement compatibles avec les filtrations puisque le lieu singulier du conoyau est inclus dans le lieu $P=(1:0:0)$.
%Le morphisme de fibrés vectoriels restriction de $f$ au diviseur ${\bf P}^{1}_{0}$ est donc de rang constant, ce qui permet de conclure pour la deuxième assertion et la troisième assertion.\\
$(i)$ Tous les morphismes dans $\calc_{3filtr,opp}^{-}$ sont strictement compatibles donc les catégories $\calc_{3filtr,opp}^{}$ et $\calc_{3filtr,opp}^{-}$ sont les mêmes ce qui permet de conclure car la catégorie $\calc_{3filtr,opp}^{}$ est abélienne par le corollaire précédent.\\
%$(iv)$  

$(ii)$ Comme $\calc_{3filtr,opp}^{-}=\calc_{3filtr,opp}^{}$, on peut ramener l'assertion à $\calc_{3filtr,opp}^{}$ qui est alors évidente.

$(iv)$ Ces foncteurs sont par le dictionnnaire les foncteurs fibre en respectivement les points $(1:1:1)$, $(0:1:1)$, $(1:0:1)$, $(1:1:0)$, au diviseur (droite projective) qui correspond à $F^{\bullet}_{1}$ puis à l'image de $(0:0:1)$ sur ce diviseur, au diviseur (droite projective) qui correspond à $F^{\bullet}_{0}$ puis à l'image de $(0:0:1)$ sur ce diviseur, au diviseur (droite projective) qui correspond à $F^{\bullet}_{2}$ puis à l'image de $(0:1:0)$ sur ce diviseur, au diviseur (droite projective) qui correspond à $F^{\bullet}_{0}$ puis à l'image de $(0:1:0)$ sur ce diviseur. Toutes ces restrictions se font sur des sous-ensembles de ${\bf P}^{2}\backslash P$ et donc préserve les suites exactes car les morphismes sont de rang constant sur cet ensemble.

\end{preuve}

\subsection{Structures réelles}
Dans cette section on se place sur le corps ${\bf C}$.\\
$\,$\\
Supposons que le ${\bf C}$-espace vectoriel $V$ muni de trois filtrations $(F_{0}^{\bullet }{},F_{1}^{\bullet }{},F_{2}^{\bullet }{})$ ait une structure réelle sous-jacente $V=V_{\bf R} \otimes_{\bf R} {\bf C}$. Il est ainsi muni d'une conjugaison $\tau : V \rightarrow V$. Nous voulons étendre la comparaison entre espaces vectoriels munis de trois filtrations et fibrés équivariants sur ${\bf P}^2$ pour l'action d'un certain groupe $G$ aux cas où les filtrations sont munies de certaines structures réelles. Plus précisement lorsque les filtrations $F_{1}^{\bullet }{}$ et $F_{2}^{\bullet }{}$ sont conjuguées vis-à-vis de la structure réelle de $V$. Lorsque $V$ est muni d'une structure réelle, on peut parler de la filtration conjuguée $\overline {F^{\bullet }_1}$ à une filtration ${F^{\bullet }_1}$. Pour tout $p$, $\overline {F^{p}_1}$ est le sous-espace $\tau({F^{p }_1})$. Si ${F^{\bullet }_1}$ est exhaustive et décroissante il en est de même pour $\overline {F^{\bullet }_1}$. 

Notons $\calc_{3filtr,opp,{\bf R}}$ la catégorie dont les objets sont les espaces vectoriels qui ont une structure réelle sous-jacente et qui sont munis de trois filtrations (ordonnées) exhaustives décroissantes et opposées $(F_{0}^{\bullet }{},F_{1}^{\bullet }{},F_{2}^{\bullet }{})$ telle que les filtrations $F_{1}^{\bullet }{}$ et $F_{2}^{\bullet }{}$ sont conjuguées et $F_{0}^{\bullet }$ est déjà définie au niveau de $V_{\bf R}$ c'est à dire que pour tout $p$, $F^{p}_{0}=\overline{F^{p}_{0}}$ ($F^{\bullet}_{0}$ est une filtration de $V_{\bf R}$ et on note de la même façon, par abus de notation, la filtration induite sur $V=V_{\bf R} \otimes_{\bf R} {\bf C}$) ; les morphismes sont les morphismes de ${\bf C}$-espaces vectoriels strictement compatibles aux filtrations.    

Nous voulons ici compléter le dictionnaire de la section précédente à $\calc_{3filtr,opp,{\bf R}}$. Il faut pour cela traduire la structure réelle en terme de structure sur le fibré, i.e. traduire le fait que deux des filtrations sont conjuguées, en terme de structure des fibrés sur ${\bf P}^2$ associés aux objets de $\calc_{3filtr,opp,{\bf R}}$ qui sont des objets trifiltrés.

Considérons l'involution antiholomorphe de ${\bf P}^2$ donnée par 

$$ \tau_{{\bf P}^2} : {{\bf P}^2} \rightarrow {{\bf P}^2},\,\, [u_{0}:u_{1}:u_{2}] \mapsto [{\overline u_{0}}:{\overline u_{2}}:{\overline u_{1}}].$$
Un fibré $\cale $ sur ${\bf P}^2$ est un $\tau$-fibré s'il est muni d'une involution antilinéaire $\tau$ qui induit $ \tau_{{\bf P}^2} $ sur la base. Ceci s'explicite de façon faisceautique comme on peut le voir dans \cite{sim2} pour les fibrés sur ${\bf P}^1
$. Consid\'erons le fibré $\cale$ comme un faisceau localement libre de $\calo_{{\bf P}^2}$-modules, on peut alors définir le faisceau $\tau^{*}\cale$ par
$$ \tau^{*}(\cale)(\calu)=\cale( \tau_{{\bf P}^2}(\calu))$$
pour tout ouvert $\calu$ de ${\bf P}^2$.

 Remarquons que cette définition a bien un sens car pour tout ouvert de Zariski $\calu$, $\tau_{{\bf P}^2}(\calu)$ est aussi un ouvert de Zariski. En effet, $\calu={\bf P}^{2}-\{[u_{0}:u_{1}:u_{2}] \vert P_{i}(u_{0},u_{1},u_{2})=0\,\, \text{pour tout}\,\,i \in I \} $ où $I$ est un ensemble fini qui indexe des polynômes homogènes $P_{i}(u_{0},u_{1},u_{2})$, et donc $\tau_{{\bf P}^2}(\calu)$ est l'ouvert de Zariski complémentaire du lieu des zéros de la famille des polynômes ${\overline P}_{i}(u_{0},u_{2},u_{1})$ indexée par $I$.

 $\tau^{*}\cale$ est aussi un $\calo_{{\bf P}^2}$-module localement libre. En effet $\tau^{*}(\cale)(\calu)$ a une structure de $\calo_{{\bf P}^2}(\calu)$-module donnée par : pour tout $e \in \cale( \tau_{{\bf P}^2}(\calu))$ et $f \in \calo_{{\bf P}^2}(\calu)$ par $f.e={\overline {\tau_{*}(f)}}e$. 

Le foncteur $\tau^*$ qui va de la catégorie des faisceaux localement libres sur ${\bf P}^2$ est antilinéaire c'est à dire que l'application induite par $\tau^*$ entre les espaces vectoriels complexes $\text{Hom}(\cale,\cale')$ et $\text{Hom}(\tau^{*}(\cale),\tau^{*}(\cale'))$  est antilinéaire. Il en va ainsi de même pour l'application induite en cohomologie $H^{i}(\cale) \rightarrow H^{i}(\tau^{*}(\cale))$. 

Le foncteur $\tau^*$ est une involution. Un $\tau$-fibré, fibré muni d'une involution antilinéaire au dessus de $ \tau_{{\bf P}^2}$, peut donc être défini comme suit :\\
 
\begin{definition}Un $\tau$-fibré vectoriel sur ${\bf P}^2$ est la donnée d'un fibré vectoriel $\cale $ et d'un morphisme $f: \cale \rightarrow \tau^{*}(\cale)$ tel que $\tau^{*}(f) \circ f=\text{id}_{\cale}$.\\
\end{definition}  

Soit ${\bf T}^{\tau}$ le sous-groupe des automorphismes de ${\bf P}^2$ engendré par $G={{\bf G}_m}^{3}/\Delta({{\bf G}_m})$ et $G.\tau_{{\bf P}^2}$. ${\bf T}^{\tau}$ est engendré par $G$ et $\tau_{{\bf P}^2}$. Une fibré vectoriel $\cale $ est un ${\bf T}^{\tau}$-fibré s'il est à la fois un $G$-faisceau cohérent et un $\tau$-fibré.

\begin{lemme}
Le noyau et le conoyau d'un morphisme de ${\bf T}^{\tau}$-fibrés sont des ${\bf T}^{\tau}$-fibrés.
\end{lemme} 

\begin{preuve}
Le fait que le noyau et conoyau d'un morphisme de $G$-fibrés soient des $G$-fibrés a été démontré dans la section \ref{groupesalg} p.\pageref{groupesalg}. Il reste donc à montrer que l' assertion est vraie pour les $\tau$-fibrés. 

Soit $h : \cale \rightarrow \calf$ un morphisme de $\tau$-fibrés. Notons par $\calk er$ et $\calc oker$ le noyau et le conoyau de $h$. Pour tout $\tau$-fibré $\cale$ vu comme faisceau localement libre, et tout point $P \in {\bf P}^{2}$ on a 
$$(\tau^{*}\cale)_{P}=\lim_{\xymatrix{{} \ar[r]_{ \calu \ni P} & {}}}\tau^{*}\cale(\calu)=\lim_{\xymatrix{{} \ar[r]_{ \calu \ni P} & {}}}\cale({\tau_{{\bf P}^2}}(\calu))=(\cale_{\tau_{{\bf P}^2}})_{P}.$$
Ainsi $\tau^* $ est un foncteur exact.%\footnote{\cite{eis}}. 

Notons par $f_{\cale}$, $f_{\calg}$ les morphismes définissant respectivement $\cale $ et $\calg $ comme des $\tau$-fibrés. L'exactitude du foncteur $\tau^{*}$ donne le diagramme commutatif suivant

$$\xymatrix{ 0 \ar[r]&  \calk \ar[r]^{i} \ar[d]^{f_{\cale } \circ i}&  \cale \ar[r]^{h} \ar[d]^{f_{\cale }} & \calf  \ar[r]^{\pi} \ar[d]^{f_{\calf} }&  \calc oker \ar[r] \ar[d]^{\varphi}& 0\\
 0 \ar[r]&  \tau^{*}(\calk) \ar[r]^{\tau^{*}i} \ar@/^/@{.>}[u]^{\tau^{*}( f_{\cale } \circ i)} &  \tau^{*}(\cale) \ar[r]^{\tau^{*}h}  \ar@/^0.5pc/[u]^{\tau^{*} f_{\cale }} & \tau^{*}(\calf)  \ar[r]^{\tau^{*}\pi} \ar@/^0.5pc/[u]^{\tau^{*} f_{\calf }} &  \tau^{*}(\calc oker) \ar[r] \ar[r] \ar@/^/@{.>}[u]^{\tau^{*} \varphi}& 0}$$
où $\varphi \circ \pi= \tau^{*}\pi \circ g$. On vérifie aisément que $\tau^{*}( f_{\cale } \circ i) \circ ( f_{\cale } \circ i) = id_{\calk}$. Pour finir \begin{eqnarray*} \tau^{*}\varphi \circ \varphi \circ \pi&  = &\tau^{*} \varphi \circ \tau^{*} \pi \circ g\\
{} & = & \tau^{*} ( \varphi \circ \pi ) \circ g\\
 {} & = & \tau^{*} ( \tau^{*}(\pi ) \circ g ) \circ g \\
{} & = & \pi \circ \tau^{*}g \circ g \\
{} & = & \pi.
\end{eqnarray*}
et donc, $\pi $ étant surjective, $\tau^{*}\varphi \circ \varphi= id_{\calc oker}$. 

\end{preuve}

Remarquons que $\tau^*$ est aussi une involution antiholomorphe du diviseur ${\bf P}^{1}_{0}$. Notons $\tau_{{\bf P}^1}$ la restriction de $\tau$ à ${\bf P}^{1}_{0}$. La restriction d'un $\tau$-fibré est un $\tau_{{\bf P}^1}$-fibré au sens de l'annexe A.

\begin{theoreme}\label{th3}
 La catégorie des espaces vectoriels complexes de dimensions finies munis d'une structure réelle sous-jacente, de trois filtrations décroissantes, exhaustives et opposées telles que deux d'entre elles soient conjuguées, l'autre définie sur la structure réelle, et des morphismes strictement compatibles aux filtrations $\calc_{3filtr,opp, {\bf R}}$ est équivalente à la catégorie des ${\bf T}^{\tau}$-fibrés vectoriels ${\calf ib}_{{\bf P}^{1}_{0}}$-semistables de pente $0$ sur le plan projectif ${\bf P}^2$ et des morphismes de $H$-fibrés dont les singularités du conoyau sont supportées en codimension $2$, par $(0:0:1)$, ${\calf ib}_{{\bf P}^{1}_{0}-semistable,\mu=0}({\bf P}^{2}/{\bf T}^{\tau})$ :
$$
\xymatrix{
\{ \calc_{3filtr,opp,{\bf R}} \} \ar@<2pt>[r]^-{\Phi_{R}} &  \{  {\calf ib}_{{\bf P}^{1}_{0}-semistable,\mu=0}({\bf P}^{2}/{\bf T}^{\tau}) \} \ar@<2pt>[l]^-{\Phi_{I}}
}
.$$

\end{theoreme}

\begin{preuve}
  Montrons d'abord qu'avec une telle source et un tel but $\Phi_{R}$ et $\Phi_{I}$ sont bien définis. Commençons par $\Phi_{R}$. Soit $(V,F^{\bullet}_{0},F^{\bullet}_{1},F^{\bullet }_{2})$ un objet de $\calc_{3filtr,opp,{\bf R}}$. On définit les fibrés de Rees suivants : $R^{2}(V,{\overline F}_{0}^{\bullet},{\overline F}_{2}^{\bullet })^{\sim}=R^{2}(V,{F}_{0}^{\bullet},{\overline F}_{2}^{\bullet })^{\sim}$ sur ${\bf A}^{2}=U_{2}=\text{Spec}[u_{0}/u_{2},u_{1}/u_{2}]$,  $R^{2}(V,{\overline F}_{0}^{\bullet},{\overline F}_{2}^{\bullet })^{\sim}=R^{2}(V,{F}_{0}^{\bullet},{\overline F}_{1}^{\bullet })^{\sim}$ sur ${\bf A}^{2}=U_{1}=\text{Spec}[u_{0}/u_{1},u_{2}/u_{1}]$ et $R^{2}(V,{\overline F}_{2}^{\bullet},{\overline F}_{1}^{\bullet })^{\sim}$ sur ${\bf A}^{2}=U_{0}=\text{Spec}[u_{1}/u_{0},u_{2}/u_{0}]$. Le fibré de Rees obtenu par la construction habituelle est alors $\tau^{*}(\xi_{{\bf P}^2}(V,F^{\bullet}_{0},F^{\bullet}_{1},F^{\bullet }_{2}))$. Le morphisme $f$ entre $\xi_{{\bf P}^2}(V,F^{\bullet}_{0},F^{\bullet}_{1},F^{\bullet }_{2})$ et $\tau^{*}(\xi_{{\bf P}^2}(V,F^{\bullet}_{0},F^{\bullet}_{1},F^{\bullet }_{2}))$ est donné sur tout affine par le morphisme de modules de Rees $u^{-p}v^{-q}w \mapsto u^{-p}v^{-q}\overline w$. Ainsi $\Phi_{R}$ est bien défini. Dans l'autre direction, soit $\xi$ un objet de ${\calf ib}_{{\bf P}^{1}_{0}-semistable,\mu=0}({\bf P}^{2}/{\bf T}^{\tau})$. Par restriction aux droites affines données par une des coordonnées égale à $1$ sur les plans affines standards qui recouvrent le plan projectif on voit que l'action de $\tau$ donne un ismorphisme entre $\xi_{{\bf A}^1}(V_{(1:1:1)},\calf^{\bullet }_{1})$ et $\xi_{{\bf A}^1}(V_{(1:1:1)},{\overline \calf}^{\bullet }_{1})$ et entre $\xi_{{\bf A}^1}(V_{(1:1:1)},\calf^{\bullet }_{1})$ et $\xi_{{\bf A}^1}(V_{(1:1:1)},{\overline \calf}^{\bullet }_{2})$ ce qui permet de conclure par l'équivalence de catégorie entre celle des fibrés de Rees sur la droite affine et celle des espavecs vectoriels filtrés. Les correspondances sont donc bien définies. 

L'essentielle surjectivité de $\Phi_{R}$ est claire. $\Phi_{R}$ est bien fidèle : si deux flèches de $ \calc_{3filtr,opp,{\bf R}} $ $f$ et $f'$ ont mêmes images, elles sont égales dans $ \calc_{3filtr,opp} $ et donc dans $ \calc_{3filtr,opp,{\bf R}} $. La surjectivité de $\Phi_{R}$ est de même évidente.

\end{preuve}
\begin{corollaire}
$\calc_{3filtr,opp,{\bf R}} $ est une catégorie abélienne.

\end{corollaire}

\begin{preuve}
Le lemme précédent nous dit que la catégorie des ${\bf T}^{\tau}$-fibrés est abélienne ce qui associé au fait que la catégorie ${\calf ib}_{{\bf P}^{1}_{0}-semistable,\mu=0}({\bf P}^{2}/{\bf T})$ est abélienne nous permet de conclure par l'équivalence de catégorie donnée dans le théorème.
\end{preuve}

\subsection{K-théorie équivariante de ${\bf P}^2$ et invariants d'espaces vectoriels filtrés}

Dans cette section on calcule des invariants discrets des fibrés de Rees associés aux espaces vectoriels trifiltrés. Ces invariants sont plus fins que le caractère de Chern qui fait l'objet de la section 1.7, ``Calcul explicite d'un invariant du fibré de Rees''. Pour ce faire, on calcule la classe des fibrés $G$-équivariants sur ${\bf P}^2$ pour l'action standard par translation du tore ${\bf T}$ dans la $K$-théorie équivariante de ${\bf P}^2$, ou du moins dans $K_{0}({\bf P}^{2},{\bf T})$.

Dans une première partie, on rappelle la définition de $K_{0}(X,G)$ pour une $G$-variété $X$, puis on introduit le matériel nécessaire pour donner le ``théorème de reconstruction à partir des strates'' de \cite{vevi} qui nous permettra de le décrire explicitement.  

\subsubsection{Définition de $K_{*}(X,G)$, reconstruction à partir des strates}
On se place toujours sur un corps $k$ de caractéristique nulle et algébriquement clos. Soit $G$ un groupe réductif et $X$ une $G$-variété. D'après \ref{gcatab} p.\pageref{gcatab}, la catégorie des $G$-fibrés vectoriels sur $X$ est abélienne. On note $K_{0}(X,G)$ la $k$-théorie de cette  catégorie abélienne.

$K_{0}(X,G)$ est construit de la façon suivante (cf. \cite{hir}) : soit $C(X,G)$ l'ensemble engendré par les classes d'isomorphismes $[\cale]$ de $G$-fibrés vectoriels sur $X$. La somme de Whitney $\oplus$ muni $C(X,G)$ d'une structure de semi-groupe. Soit $F(X,G)$ le groupe abélien engendré par $C(X,G)$. Notons $R(X,G)$ le sous-groupe engendré par les éléments de la forme $[\cale]-[\cale']-[\cale'']$ lorsqu'on a la suite exacte courte de $G$-fibrés
 $0 \rightarrow \cale' \rightarrow \cale \rightarrow \cale'' \rightarrow 0 $. Alors 
$$K_{0}(X,G)=\dfrac{F(X,G)}{R(X,G)}.$$
A partir d'ici on supposera que $G$ est diagonalisable. Pour tout entier positif $s$, notons par $X_{s} \subset X$ le lieu des points où les stabilisateurs sont de dimension exactement $s$.

Pour tout $s$, on considère le fibré normal à $X_s$ dans $X$, $N_s$. Il est naturellemnt muni d'une action de $G$ pour laquelle les stabilisateurs sont de dimension au plus $s$. Soit $N_{s,s-1}$ le lieu des points où la dimension du stabilisateur est exactement $s-1$. On a alors un morphisme composé $\xymatrix{N_{s,s-1} \ar@{^{(}->}[r]^{i}&  N_{s} \ar[r]^{\pi} & X_{s}}$, dont l'image inverse induit un morphisme en $K$-théorie 
$$ (\pi \circ i)^{*}:K_{*}(X_{s},G) \rightarrow K_{*}(N_{s,s-1},G).$$
Nous définirons plus bas (\cite{vevi}, section 3) un morphisme de spécialisation par une déformation du fibré normal
$$ \text{Sp}^{s-1}_{X,s} :  K_{*}(X_{s-1},G) \rightarrow K_{*}(N_{s,s-1},G).$$
Ces morphismes de spécialisation permettent de creconstruire exactement la $K$-théorie équivariante de $X$ pour l'action de $G$ à partir de la $K$-théorie équivariantes des strates.\\
\begin{theoremese}\cite{vevi} Soit $G$ un groupe diagonalisable de dimension $n$ qui agit sur une variété $X$ séparée de type fini sur $k$. Alors les morphismes de restriction $ K_{*}(X,G) \rightarrow K_{*}(X_{s},G)$ induisent un isomorphisme\\
\begin{center}
$ K_{*}(X,G) \cong K_{*}(X_{n},G) \times_{ K_{*}(N_{n,n-1},G)}   K_{*}(X_{n-1},G) \times_{ K_{*}(N_{n-1,n-2},G)}...$
\end{center}
\hspace*{2cm}\begin{center}$...  \times_{ K_{*}(N_{2,1},G)}   K_{*}(X_{1},G) \times_{ K_{*}(N_{1,0},G)}   K_{*}(X_{0},G) .$\\   
\end{center}
${}$\\
\end{theoremese}
En d'autres termes les homomorphismes de restriction $ K_{*}(X,G) \rightarrow K_{*}(X_{s},G)$ induisent un homomorphisme injectif $ K_{*}(X,G) \rightarrow \prod_{s} K_{*}(X_{s},G)$ et un élément $(a_{n},...,a_{0})$ du produit $\prod_{s} K_{*}(X_{s},G)$ est dans l'image de  $K_{*}(X,G)$ si et seulement si l'image inverse de $a_{s} \in  K_{*}(X_{s},G)$ dans  $K_{*}(N_{s,s-1},G)$ coïncide avec $\text{Sp}^{s-1}_{X,s} (a_{s-1}) \in K_{*}(N_{s,s-1},G)$ pour tout $s \in \{1,...,n\}$.\\
${}$\\
Explicitons les morphismes d'image inverse et de spécialisation. Notons d'abord que si $X$ et $Y$ sont des $G$-variétés algébriques de type fini séparées sur $k$ et que si le $G$-morphisme $i : Y \rightarrow X $ est une immersion fermée de $Y$ dans $X$ et que le $G$-morphisme $j:X\backslash Y \rightarrow X$ est une immersion ouverte, on a alors une suite exacte de localisation (\cite{thoma}, Th 2.7),
$$\xymatrix{ ... \ar[r]&   K_{n}(Y,G) \ar[r]^{i_{*}} &  K_{n}(X,G) \ar[r]^{j^{*}} &  K_{n}(X\backslash Y,G) \ar[r]^-{\delta} & ...}$$
Soit $s$ un entier positif. Considérons l'immersion fermée $\xymatrix{ X_{s} \ar@{^{(}->}[r] & X_{\leq s}}$ et $N_{s}$ le fibré normal qui s'en déduit. Considérons, suivant \cite{ful}, la déformation du fibré normal $\pi : M_{s} \rightarrow {\bf P}^1$. Alors le morphisme $\pi$ est plat et $G$-équivariant et\\
\hspace*{1cm}$\bullet$ $\pi^{-1}({\bf A}^{1})=X_{\leq s} \times {\bf A}^{1}$ et,\\ 
\hspace*{1cm}$\bullet$ $\pi^{-1}(\infty)=N_{s}$ .\\
Considérons la restriction $\pi^{0}$ de $\pi$ à l'ouvert $M_{s}^{0}=(M_{s})_{<s}$. On a alors\\
\hspace*{1cm}$\bullet$ $(\pi^{0})^{-1}({\bf A}^{1})=X_{< s} \times {\bf A}^{1}$ et,\\ 
\hspace*{1cm}$\bullet$ $(\pi^{0})^{-1}(\infty)=N_{s}^{0}=(N_{s})_{<s}$ .\\
On peut alors définir le morphisme de spécialisation par la composition de l'image inverse $ K_{*}(X_{<s},G) \rightarrow K_{*}(X_{<s} \times {\bf A}^{1},G)=K_{*}(M_{s}^{0}\backslash N_{s}^{0}, G)$ et du morphisme de localisation $K_{*}(M_{s}^{0}\backslash N_{s}^{0}, G) \rightarrow K_{*}( N_{s}^{0}, G)$. On obtient alors le morphisme voulu en se restreignant à la strate fermée sur laquelle les stabilisateurs ont pour dimension exactement $s-1$, le morphisme voulu est celui qui fait commuter le diagramme
$$
\xymatrix{
K_{*}(X_{s-1},G) \ar@{-->}[r]^{\text{Sp}^{s-1}_{X,s}} \ar[d]^{\text{image inverse}} &   K_{*}(N_{s,s-1},G) \\
K_{*}(X_{s-1}\times {\bf A}^{1},G) \ar@{=}[r] &   K_{*}((M_{s}^{0})_{s-1} \backslash N_{s,s-1},G). \ar[u]^{\text{localisation}}}$$ 
\subsubsection{Calcul explicite de $K_{0}({{\bf P}^2},{\bf T})$ et classes des fibrés de Rees}

Explicitons à l'aide de la section précédente $K_{0}({{\bf P}^2},{\bf T})$ où ${\bf T}$ est le groupe diagonal de dimension $2$, ${{{\bf G}_m}^3}/\Delta({{\bf G}_m})$, qui agit par translation sur $X={{\bf P}^2}$ donnée dans les sections précédentes. Donons tout d'abord la décomposition en strates du plan projectif pour cette action :\\
Le lieu où les stabilisateurs sont de dimension $2$ est\\
\hspace*{1cm}$\bullet$ $X_{2}=(1:0:0) \coprod (0:1:0) \coprod (0:0:1) =: A_{0} \coprod A_{1} \coprod A_{2}$.\\
Le lieu où les stabilisateurs sont de dimension $2$ est\\
\hspace*{1cm}$\bullet$ $X_{1}=\{(0:u_{1}:u_{2}) \vert u_{1}.u_{2}\neq 0 \} \coprod \{(u_{0}:0:u_{2}) \vert u_{0}.u_{2}\neq 0 \} \coprod \{(u_{0}:u_{1}:0) \vert u_{0}.u_{1}\neq 0 \} =: {{\bf G}_m}_{12} \coprod {{\bf G}_m}_{02} \coprod {{\bf G}_m}_{01}$.\\
Le lieu où les stabilisateurs sont de dimension $0$ est\\
\hspace*{1cm}$\bullet$ $X_{0}=\{(u_{0}:u_{1}:u_{2}) \vert u_{0}.u_{1}.u_{2}\neq 0 \}=: {{\bf G}_m}^{2}$.\\ 

D'après le théorème on a $K_{0}({{\bf P}^2},{\bf T})\cong K_{0}(X_{2},{\bf T})\times_{K_{0}(N_{2,1},{\bf T})} K_{0}(X_{1},{\bf T})\times_{K_{0}(N_{1,0},{\bf T})} K_{0}(X_{0},{\bf T})$. Explicitons chacun des termes de ce produit fibré itéré.\\
${ }$\\
$\bullet$ Calcul de $K_{0}(X_{2},{\bf T})$ :

On a $K_{0}(X_{2},{\bf T}) \cong K_{0}(A_{0},{\bf T})\oplus K_{0}(A_{1},{\bf T}) \oplus K_{0}(A_{2},{\bf T})$. Chacun des termes du membre de droite se calcule de la même façon. Calculons $K_{0}(A_{0},{\bf T})$.

Montrons que l'on a $K_{0}(A_{0},{\bf T})={\bf Z}[u_{0},u_{0}^{-1},v_{0},v_{0}^{-1}]$. Notons par $V_{1,0}$ le fibré vectoriel en droite sur $A_0$ muni de l'action linéaire suivant le premier facteur de ${\bf T}$ et triviale suivant le deuxième facteur puis symétriquement par $V_{1,0}$ le fibré muni de l'action linéaire suivant le deuxième facteur de ${\bf T}$ et linéaire suivant le premier facteur. Pour tout entier $m \geq 0$ on pose $V_{m,0}=V_{1,0}^{\otimes m}$, de même pour tout entier $n \geq 0$, on pose $V_{0,n}=V_{0,n}^{\otimes n}$. On généralise aux entiers négatifs en posant si $m \leq 0$, $V_{m,0}=V_{-m,0}^{*}$, si $n \leq 0$, $V_{0,n}=V_{0,-n}^{*}$. Puis finalement pour tous entiers $m,n$ $V_{m,n}=V_{m,0}\otimes V_{0,n}$. On note toujours par $[\cale]$ la classe d'un fibré $\cale$ dans le $K_0$. Considérons le morphisme

\begin{center}
$\mu :  {\bf Z}[u_{0},u_{0}^{-1},v_{0},v_{0}^{-1}]  \rightarrow  K_{0}(A_{0},{\bf T})$\\
\end{center}
\begin{center} 
$ \sum_{i,j}\,a_{ij}u_{0}^{i}v_{0}^{j}    \mapsto      \sum_{i,j } \,a_{ij} [V_{i,j}]$\\
\end{center}
$\mu$ est clairement injectif. La surjectivité vient du fait que tout fibré vectoriel sur le point se décompose en somme directe sur $(i,j)$ de fibrés de rang $a_{i,j}$ sommes directes de $V_{i,j}$, i.e. $\cale \cong \oplus_{i,j} \oplus_{r(i,j)} \, V_{i,j}$ où $r(i,j)$ est le rang du sous-fibré de $\cale$ sur lequel ${\bf T}$ agit par le caractère $(i,j)$.

Ainsi 
$$K_{0}(X_{2},{\bf T})={\bf Z}[u_{0},u_{0}^{-1},v_{0},v_{0}^{-1}]  \oplus {\bf Z}[u_{1},u_{1}^{-1},v_{1},v_{1}^{-1}] \oplus {\bf Z}[u_{2},u_{2}^{-1},v_{2},v_{2}^{-1}].$$  
$\bullet$ Calcul de $K_{0}(X_{1},{\bf T})$ :

On a $K_{0}(X_{1},{\bf T}) \cong K_{0}({{\bf G}_m}_{12},G)\oplus K_{0}({{\bf G}_m}_{02},G) \oplus K_{0}({{\bf G}_m}_{01},{\bf T})$. Chacun des termes du membre de droite se calcule de la même façon. Calculons $K_{0}({{\bf G}_m}_{12},{\bf T})$.

Rappelons pour ce calcul le fait suivant (Morita): soit $H$ un sous-groupe de ${\bf T}$ qui agit sur $X$ ($X$ est une $H$-variété) et soit $X\times_{H} G /H :=(X \times G/H)/H$. Alors, on a l'isomorphisme\\
\begin{center} $K_{0}(X,H)\cong K_{0}( X\times_{H} G /H,G)$.\\
\end{center}
En appliquant Morita à $X=\{*\}$ et $H={\bf G}_m$ (un des deux facteurs de $G={\bf G}_m$), il vient $K_{0}({{\bf G}_m}_{12},{\bf T})=K_{0}(*,{{\bf G}_m})={\bf Z}[u_{12},u_{12}^{-1}]$. D'où
$$K_{0}(X_{1},{\bf T})={\bf Z}[u_{12},u_{12}^{-1}] \oplus {\bf Z}[u_{02},u_{02}^{-1}] \oplus {\bf Z}[u_{01},u_{01}^{-1}].$$
$\bullet$ Calcul de $K_{0}(X_{0},{\bf T})$ :

Appliquons à nouveau Morita à $X=* $ et $H=e$ l'élément neutre de ${\bf T}$. Il vient
$$K_{0}(X_{0},{\bf T})=K_{0}(*)={\bf Z}.$$ 
 
Nous avons aussi à calculer ces anneaux pour les fibrés normaux.\\
${}$\\
$\bullet$ Calcul de $K_{0}(N_{2,1},{\bf T})$ :

Rappelons que $N_{2}$ est le fibré normal sur $X_{2}$ qui vient de l'immersion fermée $X_{2} \rightarrow X_{\leq 2}={{\bf P}^2}$. Il a trois composantes connexes $N_{2,i}$, correspondantes aux trois points. Il est donc ici isomorphe (comme ${\bf T}$-variété), en chaque point $A_{i}$ constituant $X_2$, au plan affine restriction du fibré tangent au plan projectif en chacun de ces trois points, l'action étant la même. Ainsi, pour la restriction du normal à  $A_0$ par exemple, le lieu où la dimension des stabilisateurs est $1$ est le produit de droites multiplicatives noté ${{\bf G}_m}_{0,01} \times {{\bf G}_m}_{0,02}$. On a donc
$$N_{2,1}={{\bf G}_m}_{0,01} \times {{\bf G}_m}_{0,02} \coprod {{\bf G}_m}_{1,01} \times {{\bf G}_m}_{0,12} \coprod {{\bf G}_m}_{2,02} \times {{\bf G}_m}_{2,12}.$$      
Dont le $K_{0}$ se déduit des calculs précédents.
%$$K_{0}(N_{2,1},{\bf T})=...$$ 
%$\bullet$ Calcul de $K_{0}(N_{1,0},{\bf T})$ :

$N_{1}$ est le fibré normal sur $X_{1}$ qui vient de l'immersion fermée $X_{1} \rightarrow X_{\leq 1}={{\bf P}^2}\backslash \{A_{0},A_{1},A_{2}\}$. Il a trois composantes connexes $N_{1,i}$ où $i \in \{0,1,2\}$, correpondantes aux trois droites multiplicatives. Comme chacune des composantes connexes de $X_{1}$ est dans un ouvert affine du plan projectif, le fibré normal à chacune des composantes
${{\bf G}_{m}}_{01}$, ${{\bf G}_{m}}_{02}$ et ${{\bf G}_{m}}_{12}$ est ${{\bf G}_{m}}_{01} \times {{\bf A}^1} $, resp. ${{\bf G}_{m}}_{02} \times {{\bf A}^1}$ puis ${{\bf G}_{m}}_{12} \times {{\bf A}^1}$. L'action se déduisant de celle sur les ouverts affines, il vient, pour $N_{1,0}$
$$N_{1,0}= {{\bf G}_{m}}_{01} \times {{\bf G}_{m}^1} \coprod {{\bf G}_{m}}_{02} \times {{\bf G}_{m}^1} \coprod {{\bf G}_{m}}_{12} \times {{\bf G}_{m}^1}.$$
Dont le $K_{0}$ se déduit par les calculs précédents.
%$$K_{0}(N_{1,0},{\bf T})=...$$
%Explicitons les morphismes qui permettent d'exprimer le produit fibré % \marginpar{à écrire...}

\begin{proposition}\label{ko} Soit $(V,F^{\bullet}_{0},F^{\bullet }_{1},F^{\bullet }_{2})$ un objet de $\calc_{3filtr}$, alors la classe de $\xi_{{\bf P}^2}(V,F^{\bullet}_{0},F^{\bullet }_{1},F^{\bullet }_{2})$ dans $K_{0}({\bf P}^{2},{\bf T})$ est donnée par
$$[\xi_{{\bf P}^2}(V,F^{\bullet}_{0},F^{\bullet }_{1},F^{\bullet }_{2})]_{K_{0}({\bf P}^{2},{\bf T})}=((P_{A_{0}},P_{A_{1}},P_{A_{2}}),(P_{{\bf G}_{m12}},P_{{\bf G}_{m02}},P_{{\bf G}_{m01}}),P_{{\bf G}_{m}^{2}}),$$
où,\\
$\hspace*{1cm}$ $\bullet$ $P_{A_{0}}=\sum_{p,q}\emph{dim}_{k}Gr^{p}_{F^{\bullet}_{0}}Gr^{q}_{F^{\bullet}_{1}}V\,\,\,u_{0}^{p}v_{0}^{q}$,\\
$\hspace*{1cm}$ $\bullet$ $P_{A_{1}}=\sum_{p,q}\emph{dim}_{k}Gr^{p}_{F^{\bullet}_{0}}Gr^{q}_{F^{\bullet}_{2}}V\,\,\,u_{1}^{p}v_{1}^{q}$,\\
$\hspace*{1cm}$ $\bullet$ $P_{A_{2}}=\sum_{p,q}\emph{dim}_{k}Gr^{p}_{F^{\bullet}_{1}}Gr^{q}_{F^{\bullet}_{2}}V\,\,\,u_{2}^{p}v_{2}^{q}$,\\
$\hspace*{1cm}$ $\bullet$ $P_{{\bf G}_{m12}}=\sum_{p}\emph{dim}_{k}Gr^{q}_{F^{\bullet}_{0}}V\,\,\,u_{12}^{p}$,\\
$\hspace*{1cm}$ $\bullet$ $P_{{\bf G}_{m02}}=\sum_{p}\emph{dim}_{k}Gr^{q}_{F^{\bullet}_{1}}V\,\,\,u_{02}^{p}$,\\
$\hspace*{1cm}$ $\bullet$ $P_{{\bf G}_{m01}}=\sum_{p}\emph{dim}_{k}Gr^{q}_{F^{\bullet}_{2}}V\,\,\,u_{01}^{p}$,\\
$\hspace*{1cm}$ $\bullet$ $P_{{\bf G}_{m}^{2}}=\emph{dim}_{k}V$.\\ 
\end{proposition}
\newpage
$\,$
\vspace{4cm}
$\,$\\
{\bf Dictionnaire : Espaces vectoriels trifiltrés-Fibrés équivariants.}
${}$\\
${}$\\
\begin{center}
%\begin{tabular}{|c||c|}
\begin{tabular}
{|>{\centering}m{2in}
||>{\centering}m{2in}|}
\hline

Catégorie d'espaces multifiltrés & Catégorie de fibrés équivariants sur ${\bf P}^2$\\  %\hline 

\end{tabular}

\begin{tabular}
{|>{\centering}m{2in}
||>{\centering}m{2in}|}
\hline

 $\calc_{3filtr} $ &  $  {\calf ib}({\bf P}^{2}/{\bf T}) $ \\

\end{tabular}

\begin{tabular}
{|>{\centering}m{2in}
||>{\centering}m{2in}|}
\hline

$ \calc_{3filtr,scind} $ &  $  {\calf ib}_{scind}({\bf P}^{2}/{\bf T}) $ \\

\end{tabular}

\begin{tabular}
{|>{\centering}m{2in}
||>{\centering}m{2in}|}
\hline

$ \calc_{3filtr,opp} $ &  $  {\calf ib}_{{\bf P}^{1}_{0}-semistable,\mu=0}({\bf P}^{2}/{\bf T}) $ \\
\end{tabular}

\begin{tabular}
{|>{\centering}m{2in}
||>{\centering}m{2in}|}
\hline

$ \calc_{3filtr,opp,{\bf R}} $ &  $  {\calf ib}_{{\bf P}^{1}_{0}-semistable,\mu=0}({\bf P}^{2}/{\bf T}^{\tau}) $ \\
\end{tabular}

\begin{tabular}
{>{\centering}m{2in}
>{\centering}m{2in}}
\hline
$\,$ & $\,$ \\
\end{tabular}

%\caption{Dictionnaire : Espaces vectoriels trifiltrés-Fibrés équivariants.}
\end{center}

%$ \calc_{3filtr} $ &  $  {\text Bun}({\bf P}^{2}/{\bf T}) $ \\

%\hline

%$\calc_{3filtr,scind} $ & $   {\text Bun}_{split}({\bf P}^{2}/{\bf T}) $ \\

%\hline

%$ \calc_{3filtr,opp} $ & $  {\text Bun}_{semistable,\mu=0}({\bf P}^{2}/{\bf T}) $ \\

%\hline

%$ \calc_{3filtr,opp,{\bf R}} $ &  $  {\text Bun}_{semistable,\mu=0}({\bf P}^{2}/H) $ \\

%\hline

\clearpage

\section{Construction relative}
\subsection{Faisceau de Rees relatif}

Dans cette section, nous voulons étendre les constructions précédentes à des ``familles d'espaces vectoriels filtrés'' paramétrées par une variété algébrique, nous et ainsi décrire une construction relative du fibré de Rees sur le plan projectif de la partie $1$.\\
${}$\\
En guise d'introduction, pour présenter l'idée de la construction, donnons la construction relative du fibré de Rees sur la droite affine ${\bf A}^1$ associée à un fibré vectoriel filtré.\\
${}$\\
Par variété algébrique on entendra schéma lisse de type fini sur un corps algébriquement clos de caractéristique nulle $k$.

\subsubsection{Idée de la construction}

Soit $S$ une variété algébrique et $\calv$ un fibré vectoriel sur $S$. Supposons que $\calv$ soit muni d'une filtration décroissante et exhaustive par des sous-fibrés vectoriels stricts $\calf^{\bullet}$. Exhaustive signifie qu'il existe deux entiers $m$ et $n$ tels que $\calf^{m}=\calv$ et $\calf^{n}=S$. Le fait que les sous-fibrés soient stricts signifie que pour tout $p$, $\calf^p$ est un sous-fibré vectoriel de $\calv $ (le faisceau de $\calo_{S}$-module localement libre associé à $\calf^p$ est plus qu'un sous-faisceau du faisceau localement libre associé à $\calv$). En chaque point $s \in S$, la fibre $\calv_{s}$ de $\calv$ est filtrée par des sous-espaces vectoriels $\calf^{p}_{s}$, 
$$\calv_{s}=\calf^{m}_{s} \supset \calf^{m-1}_{s} \supset ... \supset \calf^{n}_{s}=S.$$ 

On obtient donc en chaque point $s \in S$ un espace vectoriel filtré par une filtration exhaustive et décroissante $(\calv_{s},\calf^{\bullet}_{s})$. On peut faire la construction du fibré de Rees sur la droite affine pour cette donnée $\xi_{{\bf A}^1}( \calv_{s},\calf^{\bullet}_{s})$. On veut construire un fibré vectoriel $\xi(\calv,\calf^{\bullet})$ sur $S \times {{\bf A}^1}$ tel que pour tout $s \in S$ :
$$ \xi(\calv,\calf^{\bullet}) \vert_{\{s \} \times {{\bf A}^1}}\cong \xi_{{\bf A}^1}( \calv_{s},\calf^{\bullet}_{s}).$$    
Rappelons qu'un fibré vectoriel de rang $n$ sur $S$ est la donnée d'un schéma $\calv$ et d'un morphisme de schémas $\pi : \calv \rightarrow S$ ainsi que d'un recouvrement de $S$ par des ouverts $S_{i}, i \in I$ de $S$ et des isomorphismes $\Psi_{i}:{\pi}^{-1}(S_{i}) \rightarrow {\bf A}^{n}_{S_{i}}$, tels que pour tout $i,j \in I$ et pour tout ouvert affine $V=\text{Spec}A \subset S_{i} \cap S_{j}$, l'automorphisme $\Psi=\Psi_{j} \circ \Psi_{i}^{-1}$ de ${\bf A}^{n}_{V}=\text{Spec}A[x_{1},...,x_{n}]$ soit donné par un automorphisme linéaire $\theta$ de $A[x_{1},...,x_{n}]$ i.e. tel que l'on ait $\theta(a)=a$ pour $a \in A$ et $\theta(x_{i})=\sum_{i,j}\,a_{ij}\,x_{j}$ où les $a_{ij}$ sont des éléments de $A$. La matrice donnée par les $a_{ij}$ est inversible. 

Dans ce cadre la donnée d'un fibré vectoriel filtré sur $S$ par une suite décroissante exhaustive de sous-fibrés vectoriels stricts se traduit par le fait que pour tout $i,j$ et pour tout ouvert affine $V$ comme décrit préc\'edement, la matrice représentant $\theta$, automorphisme linéaire de $A[x_{1},...,x_{n}]$, peut être choisie dans une base adaptée triangulaire supérieure par blocs. Notons $f^{p}=\text{dim}_{k}\calf^{p}-\text{dim}_{k}S$, le rang du fibré $\pi \vert_{\calf^{p}}:\calf^{p} \rightarrow S$. Supposons de plus que $\calf^{0}=\calv$ et que $\calf^{m} \varsupsetneq \calf^{m+1}=S$. Alors on peut écrire 
\begin{center}$ {\bf A}^{n}_{V}=\text{Spec}A[x_{1},...,x_{f^{m}-1},x_{f^{m}},x_{f^{m}+1},...,x_{f^{m-1}-1},,x_{f^{m-1}},,x_{f^{m-1}+1},...(...)...,x_{f^{0}}]$.\end{center}
et la matrice qui représente le changement de carte sur l'ouvert affine $V \subset S_{i} \cap S_{j}$ pour le fibré vectoriel $\calv=\calf^{0}$ dans la base des $x_{i}, \, i \in [1,f^{0}]$, $\theta$, est de la forme
$$\left(
\begin{array}{cccc}
\theta_{m} & * & ... & *\\
0 & \theta_{m-1} &... &* \\
...   & ...&...&...\\
0 & 0 & ... & \theta_{0}\\
\end{array}
\right),$$
les éléments diagonaux $\theta_{i}$ sont des matrices inversibles de tailles $(f^{i}-f^{i+1}). (f^{i}-f^{i+1})$, à coefficients dans $A$, la matrice est de taille $f^{0}.f^{0}$.\\
Ainsi la matrice qui représente l'automorphisme de changement de carte sur l'ouvert affine $V \subset S_{i} \cap S_{j}$ pour le fibré vectoriel $\calf^{p}$ dans la base des $x_{i}, \, i \in [1,f^{0}]$, se déduit par troncature de la matrice qui représente $\theta$
$$\left(
\begin{array}{cccc}
\theta_{m} & * & ... & *\\
0 & \theta_{m-1} &... &* \\
...   & ...&...&...\\
0 & 0 & ... & \theta_{p}\\
\end{array}
\right),$$
c'est une matrice $f^{p}.f^{p}$.

Chacun des ouverts $S_i$ dans les donn\'ees définissant le fibré peut être recouvert par des ouverts affines. On peut donc supposer que les ouverts $S_{i}=\text{Spec}A_{i}$ sont affines. Nous allons tout d'abord faire la construction relative du fibré de Rees associé à un espace vectoriel filtré sur les ouverts affines $S_{i}$ pour obtenir des fibrés vectoriels sur les $S_{i}\times {\bf A}^{1}$, $\xi_{i}(\calv,\calf^{\bullet })$. Nous recollerons ensuite ces fibrés pour obtenir le fibré voulu sur $S\times{\bf A}^{1}$, que nous noterons $\xi(\calv,\calf^{\bullet})$.  

L'anneau $A_{i}[x_{i1},...,x_{if^{m}-1},x_{if^{m}},x_{if^{m}+1},...,x_{if^{m-1}-1},,x_{if^{m-1}},,x_{if^{m-1}+1},...(...)...,x_{if^{0}}]$ issu de la description précédente du fibré vectoriel sur la carte affine $S_i$  sera noté $A_{i}[x_{ik}]$ où $k \in [1,f^{0}]$ ; on note aussi ${\bf A}^{1}=\text{Spec}k[u]$. Considérons le $A_{i}[u][x_{ik}]$-module $B_{i}$ tel que
$$B_{i}=<u^{-p}.x_{ik}\vert \,k \leq f^{p}>_{A_{i}[u][x_{ik}]}.$$ 
Définissons le fibré vectoriel $\xi_{i}(\calv,\calf^{\bullet })$ sur $S_{i}\times {\bf A}^{1}$ par
$$ \xi_{i}(\calv,\calf^{\bullet })=\text{Spec}B_{i}.$$
Le morphisme d'inclusion $\pi_{i}^{\sharp}:A_{i}[u] \rightarrow B_{i}$ donne $\pi_{i}:\xi_{i}(\calv,\calf^{\bullet }) \rightarrow S\times{\bf A}^{1}$ (on notera toujours $f^{\sharp}:A \rightarrow B$ le morphisme d'anneaux associé au morphisme de schémas $f: \text{Spec}B \rightarrow \text{Spec}A$). De plus on a un isomorphisme de $A_{i}[u]$-modules $\varphi_{i}^{\sharp}$ entre $B_{i}$ et $A_{i}[u][x_{ik}]$ défini par 
$$\varphi_{i}^{\sharp}: B_{i} \rightarrow  A_{i}[u][x_{ik}], \,\, \varphi_{i}^{\sharp}(x_{ik})=u^{p}.x_{ik},\, \text{ pour } k \in [f^{p+1}+1,f^{p}],$$
 d'où l'isomorphisme explicite sur $S_{i}\times {{\bf A}^1}$ 

$$ \varphi_{i} :\text{Spec}A_{i}[u][x_{ik}] \rightarrow  \xi_{i}(\calv,\calf^{\bullet }) .$$ 
Comme nous avons pris un recouvrement de définition du fibré par des affines, les ouverts $S_{i}\cap S_{j}$ sont affines, $S_{i}\cap S_{j}=\text{Spec}A_{ij}$ : nous dénoterons l'automorphisme $A_{ij}[u]$-linéaire de $A_{ij}[u][x_{k}]$ où $k \in [1,f^{0}]$ par $\theta_{ij}^{\sharp}$. Le morphisme de schéma $\theta_{ij}$ qui s'ne déduit n'est autre que le morphisme de changement de cartes $\Psi=\Psi_{j}\circ\Psi_{i}^{-1}$ pour l'image inverse sur $S \times {\bf A}^1$ du fibré $\calv$ sur $S$ par la projection sur le premier facteur.  

Les morphismes d'anneaux $\varphi_{i}^{\sharp}$, $\varphi_{i}^{\sharp}$ et $\theta_{ij}$ permettent donc de définir des isomorphismes d'anneaux :
$$\varphi_{ij}^{\sharp}:\xymatrix{  (\varphi_{j}^{\sharp})^{-1}(A_{ij}[u][x_{k}]) \ar[r]^-{\varphi_{j}^{\sharp}} &  A_{ij}[u][x_{k}] \ar[r]^-{\theta_{ij}^{\sharp}} &  A_{ij}[u][x_{k}] \ar[r]^-{(\varphi_{i}^{\sharp})^{-1}} &  (\varphi_{i}^{\sharp})^{-1}(A_{ij}[u][x_{k}]) },$$
où l' on omet d'écrire les morphismes de projection de $A_{i}$ ou $A_{j}$ vers $A_{ij}$. On en déduit les isomorphismes 

$$\varphi_{ij}: \xi_{i}\vert_{S_{i}\cap S_{j}} \cong \xi_{i}\vert_{S_{i}\cap S_{j}}, \,\,\text{ pour tout } \,i,j.$$
Rappelons que\\
{\bf Fait :}(\cite{har}, exercice 1.22, p.69)\\
Soit $S$ un espace topologique et $\{S_{i}\}_{i \in I}$ un recouvrement de $S$ par des ouverts. Supposons que pour tout $i \in I$ on ait un faisceau $\xi_{i}$ et pour tout couple $i,j$ un isomorphisme $\varphi_{ij}: \xi_{i}\vert_{S_{i}\cap S_{j}} \cong  \xi_{i}\vert_{S_{i}\cap S_{j}}$ tel que : $(1)$ pour tout $i$, $\varphi_{ii}=id$ , et $(2)$ pour tout $i,j,k$, $\varphi_{ik}=\varphi_{jk}\circ\varphi_{ij}$ sur $S_{i}\cap S_{j}\cap S_{k}$. Alors il existe un unique faisceau $\xi$ sur $S$ et des isomorphismes $\Psi_{i}:\xi\vert_{S_{i}}\cong \xi_{i}$ tels que pour tous $i,j$, $\Psi_{j}=\varphi_{ij}\circ\Psi_{i}$ sur $S_{i} \cap S_{j}$. $\xi$ est le faisceau obtenu en collant les faisceaux $\xi_{i}$ par les isomorphismes $\varphi_{ij}$. 

Ce fait énoncé dans la catégorie des faisceaux reste vrai dans la catégorie des faisceaux cohérents et dans la catégorie des fibrés vectoriels sur une vari\'et\'e alg\'ebrique.\\
${}$

Appliquons ceci aux données $(\{S_{i}\}_{i \in I}, \{\xi_{i}\}_{i \in I},\{\varphi_{ij}\}_{i,j \in I})$ de notre construction. La condition $(1)$ est évidement vérifiée. Pour vérifier la condition $(2)$, il suffit de voir pour tout $i,j,k$ que $\varphi_{ik}^{\sharp}=\varphi_{ij}^{\sharp} \circ \varphi_{jk}^{\sharp}$. Ce qui se ramène à
$$   (\varphi_{i}^{\sharp})^{-1}\circ \theta_{ik}^{\sharp} \circ \varphi_{k}^{\sharp}= ( (\varphi_{i}^{\sharp})^{-1}\circ \theta_{ij}^{\sharp} \circ \varphi_{j}^{\sharp})\circ( (\varphi_{j}^{\sharp})^{-1}\circ \theta_{jk}^{\sharp} \circ \varphi_{k}^{\sharp}),$$
qui est équivalent à
$$ \theta_{ik}^{\sharp}= \theta_{ij}^{\sharp}\circ \theta_{jk}^{\sharp},$$
ce qui est vrai puisque $\calv$ est un fibré. D'où le fibré $\pi:\xi \rightarrow S\times {{\bf A}^1}$.\\
${}$\\
Le fibré vectoriel ainsi obtenu a bien les propriétés voulues.\\
\begin{proposition} Soit $\calv \rightarrow S$ un fibré vectoriel sur un schéma $S$ muni d'une filtration exhaustive $\calf^{\bullet }$ par des sous-fibrés vectoriels stricts, alors le fibré de Rees $\xi(\calv,\calf^{\bullet }) \rightarrow S \times {\bf A}^{1}$ de la construction précédente vérifie\\
\hspace*{1cm} $\bullet$ Pour tout $u \in {\bf A}^{1},\,\,u\neq 0$, $\xi(\calv,\calf^{\bullet})\vert_{S\times\{u\}} \cong \calv$,\\ 
\hspace*{1cm} $\bullet$  $\xi(\calv,\calf^{\bullet})\vert_{S\times\{0\}} \cong Gr_{\calf^{\bullet}}\calv$,\\
\hspace*{1cm} $\bullet$  Pour tout $s \in S$, $\xi(\calv,\calf^{\bullet})\vert_{\{s\}\times{{\bf A}^1}} \cong \xi_{{\bf A}^1}(\calv_{s},\calf^{\bullet}_{s})$.\\
\end{proposition}
$Gr_{\calf^{}}\calv$ est le fibré $\oplus_{p}\calf^{p-1}/\calf^{p}$.\\
${}$\\
Nous démontrerons ces propriétés plus bas, sous certaines hypothèses restreintes, dans le cadre plus général d'un fibré muni de trois filtrations par des sous-fibrés stricts. 

\subsubsection{Construction de Rees associée à un fibré muni de trois filtrations}

Cette partie donne la version relative de la construction du fibré de Rees sur ${\bf P}^2$ associée à un espace vectoriel
muni de trois filtrations. Supposons que l'on ait, dans un sens que l'on précisera, une famille d'espaces vectoriels trifiltrés $(\calv_{s},\calf_{0}^{\bullet }{}_{s},\calf_{1}^{\bullet }{}_{s},\calf_{2}^{\bullet }{}_{s})$ paramétrée par un schéma $S$. En chaque point de $S$ on peut effectuer la construction du fibré de Rees sur le plan projectif $\xi_{{\bf P}^2}(\calv_{s},\calf_{0}^{\bullet }{}_{s},\calf_{1}^{\bullet }{}_{s},\calf_{2}^{\bullet }{}_{s}) \rightarrow {\bf P}^{2}$. Nous allons construire un faisceau cohérent $\xi \rightarrow S \times {\bf P}^{2}$ tel que pour tout point $s \in S$ on ait l'isomorphisme de fibré vectoriels
$$ \xi\vert_{\{s\}\times {\bf P}^{2}} \cong \xi_{{\bf P}^2}(\calv_{s},\calf_{0}^{\bullet }{}_{s},\calf_{1}^{\bullet }{}_{s},\calf_{2}^{\bullet }{}_{s}).$$ 

Considérons un fibré vectoriel $\calv$ sur un schéma $S$ équipé de trois filtrations exhaustives et décroissantes par des sous-fibrés stricts, $\calf_{0}^{\bullet }$, $\calf_{1}^{\bullet }$ et $\calf_{2}^{\bullet }$. Rappelons que dans partie $1$ pour construire un fibré sur $ {\bf P}^2$ associé à un espace vectoriel trifiltré $(V,F^{\bullet}_{0},F^{\bullet}_{1},F^{\bullet}_{2})$ nous construisions tout d'abord les fibré de Rees associés aux couples de filtrations sur les ouverts affines standards du plan projectif pour ensuite recoller ces construtions. Les cartes affines standards étaient et seront notées $U_{0}$, $U_{1}$ et $U_{2}$. Les fibrés sur les cartes affines étaient noté $\xi_{U_{i}}(V,F^{\bullet}_{j},F^{\bullet}_{k})$ pour $\{i,j,k\}=\{0,1,2\}$. Nous allons ici suivre la même démarche et donc construire des faisceaux cohérents $\xi_{i}(\calv,\calf_{j}^{\bullet},\calf_{k}^{\bullet })$ sur chacun des produits $S  \times U_{i}$. Nous recollerons ensuite les trois fibrés pour obtenir le fibré de Rees sur $S \times {\bf P}^{2}$, $\xi_{}(\calv_{},\calf_{0}^{\bullet }{}_{},\calf_{1}^{\bullet }{}_{},\calf_{2}^{\bullet }{}_{})$. Explicitions la construction pour deux filtrations.
\begin{center} 
{\bf Construction de $\xi(\calv,\calf^{\bullet},\calg^{\bullet}) \rightarrow S \times {\bf A}^{2}$} 
\end{center}
Soit $\calv \rightarrow S$ un fibré vectoriel muni de deux filtrations exhaustives et décroissantes $\calf^{\bullet}$ et $\calg^{\bullet }$ par des sous-fibrés stricts.

Etudions tout d'abord les propriétés des intersections des sous-fibrés stricts donnés par les filtrations. Pour tous $(p,q)$, on définit le faisceau cohérent $\calf^{p}\cap \calg^{q}$ comme étant le noyau du morphisme canonique de fibrés vectoriels vus comme des faisceaux localement libres
$$\calf^{p}\cap \calg^{q}:= \text{Ker}(\calf^{p}\rightarrow  \calv/\calg^{q}).$$
$\calf^{p}\cap \calg^{q}$ est ainsi un sous-fibré de $\calf^{p}$ et de $ \calg^{q}$. Par l'inclusion de $\calf^{p}$ dans $\calf^{p-1}$ on voit que $\calf^{p-1}\cap \calg^{q}$ est un sous-faisceau cohérent de $\calf^{p}\cap \calg^{q}$. D'où la suite exacte
$$ 0 \rightarrow \calf^{p-1}\cap \calg^{q} \rightarrow \calf^{p}\cap \calg^{q} \rightarrow Gr_{\calf}^{p}\calv\cap \calg^{q} \rightarrow 0,$$
qui définit le faisceau $Gr_{\calf}^{p}\calv\cap \calg^{q}$. C'est un faisceau cohérent comme conoyau d'un morphisme de faisceaux cohérents. On prenda garde au fait que  $Gr_{\calf}^{p}\calv\cap \calg^{q}$ n'est pas un sous-faisceau de $\calg^{q}$ contrairement à ce que suggère cette notation pratique.

De telles considérations suivant les deux indices mènent au diagramme commutatif de faisceaux cohérents sur $S$ dont les lignes et les colonnes sont exactes
$$\xymatrix{
{}  &  0 \ar[d]  &  0 \ar[d]   &  0 \ar[d]  & {}\\
0  \ar[r]  &   \calf^{p}\cap \calg^{q}  \ar[d] \ar[r] &  \calf^{p-1}\cap \calg^{q} \ar[d] \ar[r] &  Gr_{\calf}^{p}\calv\cap \calg^{q}   \ar[d] \ar[r] &  0 \\
0  \ar[r]  &   \calf^{p}\cap \calg^{q-1}  \ar[d] \ar[r] & \calf^{p-1}\cap \calg^{q-1}  \ar[d] \ar[r] &   Gr_{\calf}^{p}\calv\cap \calg^{q-1}  \ar[d] \ar[r] &  0 \\
0  \ar[r]  &  \calf^{p}\cap Gr_{\calg}^{q}\calv   \ar[d] \ar[r] &   \calf^{p-1}\cap Gr_{\calg}^{q}\calv \ar[d] \ar[r] &  Gr_{\calf}^{p}Gr_{\calg}^{q}\calv   \ar[d] \ar[r] &  0 \\
{} &  0   & 0   & 0  &  {}
}$$
Ce diagramme permet de définir le faisceau cohérent $Gr_{\calf}^{p}Gr_{\calg}^{q}\calv$.\\
${} $\\

 Dans un premier temps nous allons ajouter une hypothèse sur ces filtrations : nous supposerons 
\begin{center} ${\bf (H)}$ pour tous $(p,q)$, les faisceaux cohérents $ \calf^{p}\cap \calg^{q}$ sont des sous-fibrés stricts de $\calv$.
\end{center}
 Cela signifie que l'une des filtrations induit une filtration par des sous-espaces vectoriels stricts sur chacune des composantes du fibré vectoriel gradué associé à l'autre filtration et réciproquement.\\
${}$\\
$\bullet $  Construction de $\xi(\calv,\calf^{\bullet},\calg^{\bullet}) \rightarrow S \times {\bf A}^{2}$ sous l'hypothèse ${\bf (H)}$ :\\  
${}$\\
La construction se calque sur celle qui est faite dans l'introduction de cette partie pour un fibré vectoriel muni d'une seule filtration. Pour être cohérent avec les notations de cette introduction, on en gardera les notations pour tout ce qui concerne la filtration $\calf^{\bullet }$ et on mettra des ' pour tout ce qui concerne la filtration $\calg^{\bullet }$. On supposera toujours que le recouvrement d'ouverts $S=\cup_{i \in I}\,S_{i}$ qui permet de définir le fibré $\calv \rightarrow S$ est formé d'affines i.e. que pour tout $i \in I$, $S_{i}=\text{Spec}A_{i}$.

 Notons $f^{p}=\text{dim}_{k}\calf^{p}-\text{dim}_{k}S$, le rang du fibré $\pi \vert_{\calf^{p}}:\calf^{p} \rightarrow S$. Supposons de plus que $\calf^{0}=\calv$ et que $\calf^{m} \varsupsetneq \calf^{m+1}=S$. On note aussi $g^{p}=\text{dim}_{k}\calg^{p}-\text{dim}_{k}S$, le rang du fibré $\pi \vert_{\calg^{p}}:\calg^{p} \rightarrow S$ et on suppose de plus que $\calg^{0}=\calv$ et que $\calg^{n} \varsupsetneq \calg^{n+1}=S$. On écrit aussi $g^{p,q}=\text{dim}_{k}Gr_{\calf}^{p}\calv\cap \calg^{q}-\text{dim}_{k}S$.

Pour tous $(i,j) \in I$ et tout ouvert affine $V=\text{Spec}A \subset S_{i}\cap S_{j}$ (en particulier pour $U_{ij}=U_{i}\cap U_{j}=\text{Spec}A_{ij}$) les automorphismes $\Psi_{}=\Psi_{j}\circ\Psi_{i}^{-1}$ de ${\bf A}^{n}_{V}$ sont $A$-linéaires et compatibles aux filtrations puisque les sous-fibrés vectoriels qu'elles définissent sont stricts et ont pour changement de cartes des restrictions de ces automorphismes. On peut donc écrire ${\bf A}^{n}_{V}$ comme le spectre d'un anneau de polynômes qui fasse ``apparaître'' l'une et l'autre des filtrations. L'hypothèse ${\bf (H)}$ nous permet de procéder de la même façon pour deux filtrations. Dans la matrice associée à l'automorphisme $\Psi$ dans l'introduction chaque bloc matrice inversible $\theta_{p}$ où $p \in [1,f^{0}]$ devient elle-même une matrice diagonale par blocs dont la diagonale est composée de $n+1$ blocs inversibles de la forme précédente.\\

%\scalebox{0.65}[0.65]{

%$$\left(
%\begin{array}{cccc}
%\theta_{m}=\left(
%\begin{array}{cccc}
%\theta_{m,n}' & * & ... & *\\
%0 & \theta_{m,n-1}' &... &* \\
%...   & ...&...&...\\
%0 & 0 & ... & \theta_{m,0}'\\
%\end{array}
%\right) & * & ... & *\\
%0 & \theta_{m-1}=\left(
%\begin{array}{cccc}
%\theta_{m-1,n}' & * & ... & *\\
%0 & \theta_{m-1,n-1}' &... &* \\
%...   & ...&...&...\\
%0 & 0 & ... & \theta_{m-1,0}'\\
%\end{array}
%\right) &... &* \\
%...   & ...&...&...\\
%0 & 0 & ... & \theta_{0}=\left(
%\begin{array}{cccc}
%\theta_{0,n}' & * & ... & *\\
%0 & \theta_{0,n-1}' &... &* \\
%...   & ...&...&...\\
%0 & 0 & ... & \theta_{0,0}'\\
%\end{array}
%\right)\\
%\end{array}
%\right).$$
%}
${}$\\
On peut donc écrire que $\theta$ (sur $S_{ij}$) est un automorphisme $A_{ij}$-linéaire de $A_{ij}[x_{kl_{k}}]$ où $k \in [1,f^{0}]$ et, comme il existe $p$ tel que $k \in [f^{p+1}+1,f^{p}]$ et $l_{k} \in [g^{p,n},g^{p,0}]$. 
\begin{definition}Le module de Rees associé à la trivialisation $\Psi_{i}$ sur $S_{i}=\emph{Spec}A_{i}$ du fibré $\calv$ muni de deux filtrations vérifiant $(H)$ est le $A_{i}[u,v]$-module
$$B_{i}:=<u^{-p}v^{-q}x_{kl_{k}}\vert k \leq f^{p},\,l_{k}\leq g^{p,q}>_{A_{i}[u,v]}.$$
Le fibré de Rees $\xi_{i}(\calv,\calf^{\bullet },\calg^{\bullet})$ sur $ S_{i}\times {\bf A}^{2}$ est le fibré associé au $A_{i}[u,v]$-module de Rees
$$\xi_{i}(\calv,\calf^{\bullet },\calg^{\bullet}):=\text{Spec}B_{i},$$
où le morphisme $\xi_{i}(\calv,\calf^{\bullet },\calg^{\bullet}) \rightarrow  S_{i}\times {\bf A}^{2}$ est donné par le morphisme injectif d'anneaux $\pi_{i}^{\sharp}:A_{i}[u,v] \rightarrow B_{i}$.
\end{definition}   

On a un isomorphisme de $A_{i}[u,v]$-modules $\varphi_{i}^{\sharp}$ entre $B_{i}$ et $A_{i}[u,v][x_{ik}]$ défini par 
$$\varphi_{i}^{\sharp}: B_{i} \rightarrow  A_{i}[u,v][x_{ik}], \,\, \varphi_{i}^{\sharp}(x_{ik})=u^{p}v^{q}x_{kl_{k}},\, \text{ pour } k \in [f^{p+1}+1,f^{p}],\,\,l_{k}\in [g^{p,q+1}+1,g^{p,q}]$$
 d'où l'isomorphisme explicite sur $S_{i}\times {{\bf A}^2}$ 
$$ \varphi_{i} :\text{Spec}A_{i}[u,v][x_{kl_{k}}] \rightarrow  \xi_{i}(\calv,\calf^{\bullet },\calg^{\bullet}).$$ 
Comme nous avons pris un recouvrement de définition du fibré par des affines, les ouverts $S_{i}\cap S_{j}$ sont affines, $S_{i}\cap S_{j}=\text{Spec}A_{ij}$ : nous dénoterons l'automorphisme $A_{ij}[u]$-linéaire de $A_{ij}[u,v][x_{k}]$ où $k \in [1,f^{0}]$ par $\theta_{ij}^{\sharp}$. Le morphisme de schéma $\theta_{ij}$ qui s'ne déduit n'est autre que le morphisme de changement de cartes $\Psi=\Psi_{j}\circ\Psi_{i}^{-1}$ pour l'image inverse sur $S_{i} \times {\bf A}^2$ du fibré $\calv$ sur $S_{i}$ par la projection sur le premier facteur.  

Les morphismes d'anneaux $\varphi_{i}^{\sharp}$, $\varphi_{i}^{\sharp}$ et $\theta_{ij}$ permettent donc de définir des isomorphismes d'anneaux :
$$\varphi_{ij}^{\sharp}:\xymatrix{  (\varphi_{j}^{\sharp})^{-1}(A_{ij}[u,v][x_{k}]) \ar[r]^-{\varphi_{j}^{\sharp}} &  A_{ij}[u,v][x_{k}] \ar[r]^-{\theta_{ij}^{\sharp}} &  A_{ij}[u,v][x_{k}] \ar[r]^-{(\varphi_{i}^{\sharp})^{-1}} &  (\varphi_{i}^{\sharp})^{-1}(A_{ij}[u,v][x_{k}]) },$$
où l' on omet d'écrire les morphismes de projection de $A_{i}$ ou $A_{j}$ vers $A_{ij}$. On en déduit les isomorphismes 

$$\varphi_{ij}: \xi_{i}(\calv,\calf^{\bullet },\calg^{\bullet})\vert_{S_{i}\cap S_{j}} \cong \xi_{i}(\calv,\calf^{\bullet },\calg^{\bullet})\vert_{S_{i}\cap S_{j}}, \,\,\text{ pour tous} \,i,j.$$
et donc, à l'aide de la discussion faite dans l'introduction, par recollement, d'un fibré $\xi(\calv,\calf^{\bullet },\calg^{\bullet})$ sur $S\times {\bf A}^2$.\\
${}$\\
$\bullet $ Construction de $\xi(\calv,\calf^{\bullet},\calg^{\bullet}) \rightarrow S \times {\bf A}^2$ dans le cas général :\\
${}$\\
Chaque faisceau $\calf^{p}\cap \calg^{q}$ est cohérent sur $S$, en particulier sur chacun des ouverts affines $S_{i}=\text{Spec}A_{i}$, on a  $\calf^{p}\cap \calg^{q}=M_{i,pq}^{\sim}$ où $M_{i,pq}$ est un $A_i$-module finiment engendré. Comme pour tout $(p,q)$ on a des morphismes injectifs de faisceaux cohérents $\calf^{p}\cap \calg^{q} \rightarrow \calf^{p-1}\cap \calg^{q}$ et  $\calf^{p}\cap \calg^{q} \rightarrow \calf^{p}\cap \calg^{q-1}$ (qui commutent sur les carrés d'indices), on déduit du fait que le foncteur covariant $M \mapsto M^{\sim}$ établisse une équivalence de catégories entre la catégorie des $A$-modules finiment engendrés (où $A$ est n{\oe}therien) et la catégorie des $\calo_{\text{Spec}A}$-modules cohérents un diagramme commutatif de morphismes de $A_{i}$-modules dans le quel chaque flèche est l'injection du sous-module source dans le module but$$\xymatrix{
M_{i,mn} \ar[r]\ar[d] & M_{i,m-1n} \ar[r]\ar[d] & ... \ar[r]\ar[d]  &  M_{i,0n} \ar[d]\\
M_{i,mn-1} \ar[r]\ar[d] & M_{i,m-1n-1} \ar[r]\ar[d] & ... \ar[r]\ar[d]  &  ...\ar[d]\\
... \ar[r]\ar[d] & ... \ar[r]\ar[d] & ... \ar[r]\ar[d]  &  M_{i,01} \ar[d]\\
M_{i,m0} \ar[r] & ... \ar[r]  & M_{i,10} \ar[r]  &  M_{i,00},}
$$
où $k \in [1,\text{rang}(\calv)]$. On a $M_{i,m-1n} \cap M_{i,mn-1}=M_{i,mn}$.

Ceci décrit pour tous $(p,q)$ $M_{i,pq}$ comme un sous-$A_{i}$-module du module libre $M_{i,00}$, noté $M_{i}$. Définissons alors le $A_{i}[u,v]$-module de Rees associé à la restriction du fibré bifiltré sur $\text{Spec}A_i$
$$B_{i}=\sum_{p,q}\,\,u^{-p}v^{-q}M_{i,pq}.$$ 
On étend la définition précédente dans ce cadre
\begin{definition}Le $A_{i}[u,v]$-module $B_{i}$ est le module de Rees associé à la restriction de $\calv$ à $S_{i}$.

Le faisceau cohérent $\xi_{i}(\calv,\calf^{\bullet},\calg^{\bullet}):=B_{i}^{\sim}$ sur $S_{i} \times {\bf A}^2$ qui s'en déduit est le faisceau de Rees associé à la restriction de $\calv$ à $S_{i}$.\\

\end{definition}

Les restrictions de  $\xi_{i}(\calv,\calf^{\bullet},\calg^{\bullet})$ et  $\xi_{j}(\calv,\calf^{\bullet},\calg^{\bullet})$ à $(S_{i} \cap S_{j})\times {\bf A}^2$ sont isomorphes. On peut donc recoller les faisceaux cohérents sur les ouverts $(S_{i} \cap S_{j})\times {\bf A}^2$ qui recouvrent $S \times {\bf A}^2$, d'où un faisceau cohérent sur $S \times {\bf A}^2$. Le faisceau de $\calo_{S \times {\bf A}^2}$-modules $\xi(\calv,\calf^{\bullet },\calg^{\bullet}) \rightarrow S\times {\bf A}^2$ ainsi obtenu est le faisceau cohérent de Rees associé au fibré vectoriel muni de deux filtrations exhaustives par des sous-fibrés stricts.\\
${}$\\
{\bf Remarque :} Les deux constructions précédentes coïncident lorsque les données de départ satisfont l'hypothèse ${\bf (H)}$. Nous avons cependant distingué les deux cas pour mettre en valeur le fait que la construction sous cette hypothèse mène à un fibré vectoriel, alors que sans ${\bf (H)}$ on obtient un faisceau cohérent (qui est en fait réflexif comme nous le verrons plus bas).

Cette distinction est aussi préparatoire aux applications à la théorie de Hodge de cette construction dans la partie 3. En effet, supposons que les filtrations proviennent d'une variation de structures de Hodge mixtes. Alors la filtration de Hodge et la filtration par le poids qui filtrent le fibré vectoriel associé à la variation vérifient l'hypothèse ${\bf (H)}$ alors que ce n'est pas le cas pour la filtration par le poids et sa conjuguée (on précisera plus loin comment la filtration conjuguée par rapport à une structure réelle se retrouve à partir d'une filtration par des sous-fibrés vectoriels stricts algébriques). Ainsi, $({\bf H})$ est satisfaite pour la construction sur les cartes $S \times {\bf A}^2$ concernant la filtration par le poids et l'une des deux autres filtrations mais pas vérifiée en général sur la carte concernant la filtration par le poids et sa conjuguée.

\begin{center}
$({\bf H})$ {\bf ``mauvaise hypothèse''}

\end{center} 
En fait, pour poursuivre la remarque ci-dessus, l'hypothèse $({\bf H})$ est illusoire en général. Le travail de la recherche d'invariants associés aux espaces vectoriels trifiltrés a pour but de mesurer à quel point elle n'est pas vérifiée ou d'étudier les incidences des strates sur lesquelles elle est vérifiée. $({\bf H})$ est une bonne hypothèse par strates. 
 
${} $\\
%Nous avons ainsi défini un fon 

 \begin{center}
{\bf Construction de $\xi(\calv,\calf^{\bullet}_{0},\calf^{\bullet}_{1},\calf^{\bullet }_{2}) \rightarrow S \times {\bf P}^{2}$} 
 
\end{center}
 
Considérons un fibré vectoriel $\calv $ sur une variété algébrique $S$ muni des trois filtrations décroissantes et exhaustives par des sous-fibrés stricts $(\calf^{\bullet}_{0},\calf^{\bullet}_{1},\calf^{\bullet}_{2})$. Nous venons de voir comment construire des faisceaux de Rees associés à deux filtrations par des sous-fibrés stricts sur le produit de $S$ par le plan affine. Notons $\xi_{j}(\calv,\calf_{k}^{\bullet },\calf_{l}^{\bullet })$ le fibré sur $S \times U_{j}$ obtenu par cette construction pour $\{j,k,l\}=\{0,1,2\}$ et par $\xi_{i_j}(\calv,\calf_{k}^{\bullet },\calf_{l}^{\bullet })$ sa restriction à l'ouvert affine $S_{i}$, où $i \in I$.

\begin{proposition}On peut recoller les faisceaux cohérents de Rees $\xi_{j}(\calv,\calf_{k}^{\bullet },\calf_{l}^{\bullet })$ sur $S \times U_{j}$ où $j \in \{0,1,2\}$ pour obtenir le faisceau cohérent de Rees associé au fibré vectoriel muni de trois filtrations exhaustives et décroissantes par des sous-fibrés stricts $(\calv,\calf^{\bullet}_{0},\calf^{\bullet}_{1},\calf^{\bullet }_{2})$, que l'on note
$$\xi(\calv,\calf^{\bullet}_{0},\calf^{\bullet}_{1},\calf^{\bullet }_{2}) \rightarrow S \times {\bf P}^{2}.$$  

\end{proposition}

\begin{preuve}%Les ouverts affines $(S_{i} \times U_{j})_{i \in I,j \in \{0,1,2\} }$ recouvrent $S \times {\bf P}^2$.
De la proposition \ref{proplocales} p.\pageref{proplocales}, on trouve pour tout $(j,j') \in \{0,1,2\}$ tels que $j \neq j'$, avec $n$ l'entier tel que $\{j,j',n\}=\{0,1,2\}$, un isomorphisme  
\begin{eqnarray*}\varphi_{jj'}=\Psi_{j'}^{-1}\circ\Psi_{j}:\xi_{j}(\calv,\calf_{k}^{\bullet },\calf_{l}^{\bullet })\vert_{S \times (U_{j}\cap U_{j'})}&\cong_{\Psi_j} &\xi_{}(\calv,\calf_{n}^{\bullet },\calt riv^{\bullet })\vert_{S \times (U_{j}\cap U_{j'})}\\
{}& \cong_{\Psi_{j'}^{-1}}&\xi_{j'}(\calv,\calf_{k}^{\bullet },\calf_{l}^{\bullet })\vert_{S \times (U_{j}\cap U_{j'})}.
\end{eqnarray*}
Pour définir les $\varphi_{jj}$ on prend les isomorphismes identités. Pour $\{j,j',j''\}=\{0,1,2\}$, la relation $\varphi_{jj''}=\varphi_{j'j''}\circ\varphi_{jj'}$ est bien vérifiée sur $S \times (U_{j} \cap U_{j'} \cap U_{j''})=S \times (U_{0} \cap U_{1} \cap U_{2})$. Toutes les hypothèses pour recoller de façon unique les faisceaux de Rees relatifs sur les affines sont donc réunies. On en déduit le faisceau de Rees $\xi(\calv,\calf^{\bullet}_{0},\calf^{\bullet}_{1},\calf^{\bullet }_{2}) \rightarrow S \times {\bf P}^{2}$.
\end{preuve}

\begin{proposition}
Le faisceau cohérent de Rees $\xi(\calv,\calf_{0}^{\bullet },\calf_{1}^{\bullet },\calf_{2}^{\bullet })$ sur $S \times {\bf P}^2$ à partir du fibré vectoriel $\calv$ sur $S$ muni de trois filtrations exhaustives et décroissantes par des sous-fibrés stricts est réflexif.

\end{proposition}
\begin{preuve} D'après la proposition \ref{proplocales} il suffit de regarder en l'origine des $U_{j}$ i.e. aux trois points $S \times \{(1:0:0)\}$, $S \times \{(0:1:0)\}$ et $S \times \{(0:0:1)\}$. En effet, sur les $S \times {\bf G}_{m}^2$ le faisceau est localement libre ainsi que sur $S \times {\bf G}_{m} \times {\bf A}^1$ et $S \times {\bf A}^{1} \times {\bf G}_m$. 
Nous allons utiliser la caractérisation des faisceaux réflexifs donnée dans la proposition \ref{reesreflexif} de la partie 1 : un faisceau cohérent $\xi$ sur $S$ est réflexif s'il est sans torsion et pour tout $s \in S$ $\text{depht}\,\xi_{s}\geq 2$. La preuve copie celle de la proposition \ref{reesreflexif} pour les faisceaux de Rees sur les espaces affines. L'acception est locale, plaçons-nous sur les ouverts affines $S_ {i}\times U_{j}$ de la construction de $\xi(\calv,\calf_{0}^{\bullet },\calf_{1}^{\bullet },\calf_{2}^{\bullet })$ qui recouvrent $S\times {\bf P}^2$ ($i \in I$ et $j \in \{0,1,2\}$).

Etudions la torsion. Pour tout $i$, $S_{i}=\text{Spec}A_{i}$ et le module $M_{i}$ est un $A_{i}$-module libre. Les modules de Rees $B_{i}$ sont des $A_{i}[u,v]$-modules sans torsion et donc $\xi(\calv,\calf_{0}^{\bullet },\calf_{1}^{\bullet },\calf_{2}^{\bullet })$ est sans torsion.

Etudions la profondeur. Pour tout $i$, $B_{i}$ est un sous-$A_{i}[u,v]$-module du $A_{i}[u,v]$-module libre $M_{i}[u,u^{-1},v,v^{-1}]$. Considérons l'idéal maximal ${\mathfrak m}=A_{i}(u,v)$ de $A_{i}[u,v]$. On veut montrer que $B_{i}$ est égal à son saturé $B_{i}^{sat}$. Soit $x  \in B_{i}^{sat} \subset M_{i}[u,u^{-1},v,v^{-1}]$. Alors $x$ s'écrit $x=\sum_{(p,q)}u^{-p}v^{-q}m_{i,pq}$. Le fait que $u.x \in B_{i}$ montre que pour pour tous $(p,q)$, $m_{i,pq} \in M_{i,p-1q}$. De même $v.x \in B_{i}$ montre que pour pour tous $(p,q)$, $m_{i,pq} \in M_{i,pq-1}$. Ainsi pour tous $(p,q)$,  $m_{i,pq} \in M_{i,pq}$ ce qui prouve que $x \in B_{i}$. Donc $B_{i}=B_{i}^{sat}$. Par le lemme \ref{sat} p.\pageref{sat}, on en déduit que $\text{depht}_{\mathfrak m}B_{i} \geq 2$. Ce qui permet de conclure.   
\end{preuve}

\subsubsection{Action du tore ${\bf T}^3$ sur le faisceau de Rees}

Soit, comme dans le dictionnaire de la partie $1$, ${\bf T}={\bf G}_{m}^{3}/\Delta({\bf G}_{m}^{3})$. Considérons l'action par translation de ${\bf T}$ sur ${\bf P}^{2}$ : 
$$\sigma : {\bf T} \times {\bf P}^{2} \rightarrow {\bf P}^{2}.$$

Soit $S$ une variété algébrique, on considère l'action de ${\bf T}$ sur $S \times {\bf P}^{2}$ triviale dans la direction de $S$ et induite par $\sigma$ dans la direction de ${\bf P}^{2}$; on la note aussi $\sigma$. Elle est donnée explicitement par
\begin{center}
$\sigma : {\bf T} \times S  \times {\bf P}^{2} \rightarrow S \times {\bf P}^{2}$\\
$(g,s,x) \mapsto (s,\sigma(g,x))$
\end{center}

Les faisceaux cohérents de Rees sont munis de l'une action de $G$. En effet

\begin{proposition}
Le faisceau de Rees $\xi(\calv,\calf_{0}^{\bullet },\calf_{1}^{\bullet },\calf_{2}^{\bullet })$ sur $S \times {\bf P}^{2}$ associé à un fibré vectoriel muni de trois filtrations par des sous-fibrés stricts $(\calv,\calf_{0}^{\bullet },\calf_{1}^{\bullet },\calf_{2}^{\bullet })$ sur $S$ est équivariant pour l'action de ${\bf T}$ sur $S \times {\bf P}^2$ donnée par $\sigma$.
\end{proposition}

\begin{preuve}Notons $p_{2}$ le morphisme de projection de $ {\bf T} \times S  \times {\bf P}^{2}$ sur $S \times {\bf P}^{2}$. Pour montrer que $\xi(\calv,\calf_{0}^{\bullet },\calf_{1}^{\bullet },\calf_{2}^{\bullet })$ est ${\bf T}$-équivariant pour $\sigma$, il faut montrer que l'on a un isomorphisme de faisceaux cohérents sur $ {\bf T} \times S  \times {\bf P}^{2}$, $\Psi : \sigma^{*}\xi(\calv,\calf_{0}^{\bullet },\calf_{1}^{\bullet },\calf_{2}^{\bullet }) \cong p_{2}^{*}\xi(\calv,\calf_{0}^{\bullet },\calf_{1}^{\bullet },\calf_{2}^{\bullet })$. L'action de ${\bf T}$ sur le plan projectif préserve chacun des ouverts standards le recouvrant, donc l'action sur le produit conserve chaque $S_{i}\times U_{j}$ où $i \in I$ et $j \in \{0,1,2\}$. On peut donc se restreindre à ces ouverts pour montrer l'isomorphisme. Les morphismes $\sigma$ et $p_{2}$ sont donnés sur les ouverts affines ${\bf T} \times S_{i} \times U_{j}$ par les morphismes d'anneaux
\begin{center}
$ \sigma^{\sharp}: A_{i}[u,v] \rightarrow A_{i}[u,v]\otimes k[x,x^{-1},y,y^{-1}] $\\
$\sum_{p,q}a_{i,pq}u^{p}v^{q} \mapsto  \sum_{p,q}a_{i,pq}u^{p}v^{q}\otimes x^{p}y^{q}$\\
\end{center}
et,
\begin{center}
$ p_{2}^{\sharp}: A_{i}[u,v] \rightarrow A_{i}[u,v]\otimes k[x,x^{-1},y,y^{-1}] $\\
$\sum_{p,q}a_{i,pq}u^{p}v^{q} \mapsto  \sum_{p,q}a_{i,pq}u^{p}v^{q}\otimes 1$\\
\end{center}
où ${\bf T}=\text{Spec} k[x,x^{-1},y,y^{-1}]$.

L'isomorphisme entre\\ 
$\sigma^{*}\xi(\calv,\calf_{0}^{\bullet },\calf_{1}^{\bullet },\calf_{2}^{\bullet })\vert_{{\bf T} \times S_{i} \times U_{j}}\cong (B_{i}\otimes_{A_{i}[u,v],\sigma^{\sharp}}A_{i}[u,v]\otimes k[x,x^{-1},y,y^{-1}])^{\sim}$\\ 
et,\\ 
$ p_{2}^{*}\xi(\calv,\calf_{0}^{\bullet },\calf_{1}^{\bullet },\calf_{2}^{\bullet })\vert_{{\bf T} \times S_{i} \times U_{j}}\cong (B_{i}\otimes_{A_{i}[u,v],p_{2}^{\sharp}}A_{i}[u,v]\otimes k[x,x^{-1},y,y^{-1}])^{\sim}$ vient de l'isomorphisme de $A_{i}[u,v][x,x^{-1},y,y^{-1}]$ donné par $\sigma^{\sharp}$ entre $B_{i}\otimes_{A_{i}[u,v],\sigma^{\sharp}}A_{i}[u,v]\otimes k[x,x^{-1},y,y^{-1}]$ et $B_{i}\otimes_{A_{i}[u,v],p_{2}^{\sharp}}A_{i}[u,v]\otimes k[x,x^{-1},y,y^{-1}]$ et de l'équivalence de catégories donnée par le foncteur ${}^{\sim}$ qui à un isomorphisme de modules associe un isomorphisme de faisceaux cohérents. 
\end{preuve}

\subsubsection{Morphismes stricts de fibrés filtrés, morphismes de faisceaux de Rees}
Soient $(\calv,\calf^{\bullet}) \rightarrow S$ et $(\calv',\calf'{}^{\bullet}) \rightarrow S$ deux fibrés vectoriels sur une variété algébrique $S$, filtrés par des sous-fibrés vectoriels stricts, et $f:\calv \rightarrow \calv'$ un morphisme de fibrés vectoriels.% de rang constant. Alors l'image de $\calv$, $Im(f)$ est un fibré vectoriel, sous-fibré de $\calv'$. 

\begin{definition} On dit que $f:(\calv,\calf^{\bullet}) \rightarrow (\calv',\calf'{}^{\bullet})$ est un morphisme de fibrés vectoriels filtrés si pour tout $p$
$$ f(\calf^{p}) \subset \calf'{}^{p},$$
i.e. $ f(\calf^{p})$ est un sous faisceau de $ \calf'{}^{p}$.\\
Lorsque $f$ est un morphisme entre fibrés vectoriels munis de $n$ filtrations ordonnées, on dit que c'est un morphisme d'espaces vectoriels filtrés si c'est un morphisme d'espaces vectoriels filtrés pour les fibrés vectoriels munis d'une seule filtration correspondant à chacun des $n$ indices.   
\end{definition}
Remarquons qu'alors pour tout $s \in S$,  $f_{s}(\calf^{p}_{s}) \subset \calf_{s}'{}^p$ i.e. que $f_{s} $ est un morphisme d'espaces vectoriels filtrés.

Par analogie avec la notion de morphisme d'espaces vectoriels filtrés strictement compatibles, introduisons la notion de morphisme de fibrés vectoriels filtrés strictement compatibles. 
\begin{definition}Le morphisme de fibrés vectoriels filtrés sur $S$,  $f:(\calv,\calf^{\bullet}) \rightarrow (\calv',\calf'{}^{\bullet})$ est dit strictement compatible aux filtrations si, pour tout $p$
$$ f(\calf^{p}) =Im(f) \cap \calf'{}^{p}.$$
La définition s'étend aux morphismes de fibrés vectoriels munis de plusieurs filtrations.
\end{definition}
Si $f$ est un morphisme de fibrés vectoriels filtrés sur $S$ strictement compatible aux filtrations, on vérifie que l'on a bien, pour tout $s \in S$,
$$  f_{s}(\calf^{p}_{s}) =Im(f_{s}) \cap \calf_{s}'{}^{p},$$
et donc qu'un morphisme de fibrés vectoriels strictement compatible aux filtrations par des sous-fibrés vectoriels induit bien sur les fibres un morphisme d'espaces vectoriels strictement compatible aux filtrations.

Donnons-nous un morphisme strictement compatible de fibrés vectoriels bifiltrés sur une variété $S$, $f :(\calv,\calf^{\bullet},\calg^{\bullet }) \rightarrow (\calv',\calf'{}^{\bullet},\calg'{}^{\bullet})$ où $\calv $ est de rang $n$ et $\calv'$ de rang $n'$. Recouvrons $S$ par des ouverts affines trivialisant les deux fibrés $S_{i}=\text{Spec}A_{i}$ où $i \in I$. C'est toujours possible quitte à prendre les intersections deux à deux des ouverts trivialisant respectivement l'un et l'autre des fibrés. La construction du faisceau de Rees sur chacun des ouverts $S_ {i} \times {\bf A}^2$ pour chacun des fibrés donne pour tout $i$, les fibrés $\xi_{i}(\calv,\calf^{\bullet},\calg^{\bullet }) \rightarrow S_{i} \times {\bf A}^2$ et $\xi_{i}(\calv',\calf'{}^{\bullet},\calg'{}^{\bullet}) \rightarrow S_{i} \times {\bf A}^2$. Nous allons associer à $f$ un morphisme entre les faisceaux de Rees. Exhibons ce morphisme localement sur les $S_{i} \times {\bf A}^2$. Nous gardons ici les notations pour les fibrés vectoriels donnée en début de partie, on ajoute des $'$ pour tout ce qui concerne $\calv'$. Le morphisme de fibrés vectoriels $f$ est décrit par des morphismes de schémas $f_{i}: {\bf A}^{n}_{S_i} \rightarrow {\bf A}^{n'}_{S_{i}}$ dont les morphismes d'anneaux associés $f_{i}^{\sharp}:A_{i}[x_{1},...,x_{n'}] \rightarrow A_{i}[x_{1},...,x_{n}]$ sont $A_{i}$-linéaires. Ces morphismes s'expriment en fonction des trivialisations et de $f$ par le diagramme qui commutatif
$$
\xymatrix{
\pi^{-1}(S_{i}) \ar[r]^-{\Psi_{i}}_{\cong} \ar[d]^{f}&  {\bf A}^{n}_{S_i} \ar@{.>}[d]^{f_i}\\
\pi^{-1}{}'(S_{i}) \ar[r]^-{\Psi_{i}'}_{\cong}&  {\bf A}^{n'}_{S_i} \\  
}.$$   
La stricte compatibilité du morphisme de fibrés bifiltrés signifie que pour tout $i \in I$ et tout $(p,q)$,
$$f_{i}(M_{i,pq})=f_{i}(M) \cap M_{i,pq}'.$$
Le morphisme entre modules de Rees relatifs
\begin{center}
$ f_{i}^{r}{}^{\sharp}: B_{i}=\sum_{p,q}u^{-p}v^{-q}M_{i,pq} \rightarrow B_{i}'=\sum_{p,q}u^{-p}v^{-q} M_{i,pq}'$\\
$x=\sum_{p,q}u^{-p}v^{-q}m_{i,pq} \mapsto f_{i}^{r}(x)=\sum_{p,q}u^{-p}v^{-q} f_{i}(m_{i,pq})$
\end{center} 
est donc bien défini. D'où le morphisme associé à $f_{i}^{r}$ par la construction ${}^\sim$ (on garde la même notation),
$$f_{i}^{r}: \xi_{i}(\calv,\calf^{\bullet}{},\calg^{\bullet }{}) \rightarrow \xi_{i}(\calv',\calf^{\bullet}{}',\calg^{\bullet }{}') .$$

\begin{proposition} $(i)$ Les $f_{i}^{r}$ sont des morphismes de faisceaux cohérents ${\bf G}_{m}^{2}$-équivariants pour l'action de ${\bf G}_{m}^{2}$ sur $S_{i}\times {\bf A}^2$ triviale dans la direction de $S_{i}$ et par translation dans la direction du plan affine.\\
$(ii)$ Les $f_{i}^{r}$ se recollent pour former un morphisme ${\bf G}_{m}^{2}$-équivariant entre les ${\bf G}_{m}^{2}$-faisceaux de Rees $f^{r}: \xi_{}(\calv,\calf^{\bullet},\calg^{\bullet }) \rightarrow \xi_{}(\calv',\calf'{}^{\bullet},\calg'{}^{\bullet}) $.\\
$(iii)$ Le conoyau du morphisme $f^r$ est un faisceau cohérent sans torsion.
\end{proposition}

\begin{preuve}$(i)$ Il suffit de voir que l'action (la coaction) est compatible avec le morphisme entre modules de Rees relatifs, i.e. que les diagrammes suivants commutent
$$
\xymatrix{
B_{i} \ar[r]^-{p_{2}^{\sharp}} \ar[d]_-{f_{i}^{\sharp}} &  B_{i} \otimes k[u,u^{-1},v,v^{-1}] \ar[d]^-{f_{i}^{\sharp} \otimes id}\\
B_{i}' \ar[r]^-{p_{2}^{\sharp}}  &  B_{i}' \otimes k[u,u^{-1},v,v^{-1}]\\
}\,\,\,\,\,\,\,\,\,\,\,\,\,\,\xymatrix{
B_{i} \ar[r]^-{\sigma^{\sharp}} \ar[d]_-{f_{i}^{\sharp}} &  B_{i} \otimes k[u,u^{-1},v,v^{-1}] \ar[d]^-{f_{i}^{\sharp} \otimes id}\\
B_{i}' \ar[r]^-{\sigma^{\sharp}}  &  B_{i}' \otimes k[u,u^{-1},v,v^{-1}]\\
},$$
ce qui est bien le cas par la stricte compatibilité des $f_{i}^{r}$.\\ 
$(ii)$ On vérifie que pour tout $i,j \in I$, les restrictions à $(S_{i}\cap S_{j}) \times {\bf A}^2$, les $f_{i}^r$, $f_{j}^r$ commutent aux isomorphismes de changement de cartes des deux fibrés $\varphi_{ij}$ et $\varphi_{ij}'$, ce qui permet de recoller les morphismes pour obtenir un morphisme $f^{r}$ sur le faisceau de Rees obtenu dans la proposition \ref{recolledef} en collant les descriptions locales. L'action étant triviale dans la direction de $S_{i}$ sur tous les $ S_{i}\times {\bf A}^2$, on peut recoller.\\
$(iii)$ La question est locale. Plaçons-nous sur les $S_{i} \times {\bf A}^2$. Il suffit de montrer que le conoyau du morphisme de $A_{i}[u,v]$-modules $B_{i}\rightarrow B_{i}'$ est sans torsion.   
\end{preuve}

On en déduit pour la construction du faisceau de Rees relatif sur le plan projectif :
\begin{corollaire} Soit sur une variété $S$ un morphisme de fibrés vectoriels filtrés par des sous-fibrés stricts strictement compatibles aux filtrations  $f:(\calv,\calf^{\bullet}_{0},\calf^{\bullet }_{1},\calf^{\bullet }_{2}) \rightarrow (\calv',\calf_{0}'{}^{\bullet },\calf_{1}'{}^{\bullet },\calf_{2}'{}^{\bullet })$. Alors le morphisme de faisceaux de Rees ${\bf T}$-équivariants sur $S \times {\bf P}^2$ qui s'en déduit
$$f^{r}:\xi(\calv,\calf^{\bullet}_{0},\calf^{\bullet }_{1},\calf^{\bullet }_{2}) \rightarrow \xi(\calv',\calf_{0}'{}^{\bullet },\calf_{1}'{}^{\bullet },\calf_{2}'{}^{\bullet }),$$
est un morphisme ${\bf T}$-équivariant pour l'action $\sigma$. De plus, le conoyau de $f^{r}$ est un faisceau sans torsion.
\end{corollaire}

%%%%%%%%%%%%%%%%%%%%%%%%%%%%%%%%%%%%%

\subsection{Dictionnaire sous l'hypothèse $({\bf H})$}

{\bf ATTENTION : dans toute cette section les filtrations par des fibrés vectoriels doivent vérifier l'hypoth\`ese }$({\bf H})$.

\subsubsection{ Propriétés de $\xi(\calv,\calf^{\bullet},\calg^{\bullet}) \rightarrow S \times {\bf A}^{2}$}

Etudions les propriétés du faisceau cohérent de Rees associé à un fibré vectoriel filtré par des sous-fibrés vectoriels stricts. 

\begin{lemme}
Soit $A$ un anneau, $\mathfrak a$ un idéal de $A$ et $B$ un $A$-module. Alors l'application $A/{\mathfrak a}\times B \rightarrow B/{\mathfrak a}$ induite par  $\forall a \in A,\,b \in B, \,\,(a,b) \mapsto a.b \,\emph{mod}\,  aB$ est bilinéaire et induit un isomorphisme 
$$ A/{\mathfrak a} \otimes B \cong B/ {\mathfrak a}B.$$
\end{lemme}

\begin{preuve} Remarquons d'abord que l'application $A/{\mathfrak a}\times B \rightarrow B/{\mathfrak a}$ est bilinéaire. D'où le morphisme $A/{\mathfrak a}\otimes B \rightarrow B/{\mathfrak a}$. La surjectivité est évidente. On construit un morphisme inverse avec la flèche $ B \rightarrow  A/{\mathfrak a}\otimes B$ qui à $b $ associe ${\overline 1} \otimes b$ dont le noyau est ${\mathfrak a}B$.

\end{preuve}

\begin{proposition}\label{proplocales} Le faisceau de Rees $\xi(\calv,\calf^{\bullet },\calg^{\bullet})$ sur $ S\times {\bf A}^2$ obtenu à partir du fibré vectoriel $\calv$ sur $S$ filtré par des sous-fibrés vérifie les propriétés\\
\hspace*{1cm}$(i)$ Pour tout $s \in S$, $\xi(\calv,\calf^{\bullet },\calg^{\bullet})\vert_{\{s\}\times {{\bf A}^2}}\cong  \xi_{{\bf A}^2}(\calv_{s},\calf^{\bullet }_{s},\calg^{\bullet}_{s})$ est un isomorphisme de fasiceaux localement libres sur ${\bf A}^2$.\\
\hspace*{1cm}$(ii)$ Pour tout $(u_{0},v_{0}) \in {\bf A}^{2}$ tels que $u_{0}.v_{0}\neq 0$, $\xi(\calv,\calf^{\bullet },\calg^{\bullet})\vert_{S \times (u_{0},v_{0}) }\cong \calv$ est un isomorphisme de fibrés vectoriels sur $S$.\\
\hspace*{1cm}$(iii) $ Pour tout $u_{0}\in {\bf A}^{1}$, $\xi(\calv,\calf^{\bullet },\calg^{\bullet})\vert_{S \times \{u_{0}\} \times {\bf A}^1}\cong \xi(Gr_{\calg}\calv,\calf^{\bullet })$ est un isomorphisme de faisceaux cohérents sur $S \times {\bf A}^1$ (resp. un isomorphisme de fibrés vectoriels sur $S \times {\bf A}^1$ si les filtrations vérifient ${\bf (H)}$).\\ 
\hspace*{1cm}$(iv) $ Pour tout $v_{0}\in {\bf A}^{1}$, $\xi(\calv,\calf^{\bullet },\calg^{\bullet})\vert_{S \times {\bf A}^{1} \times \{v_{0}\}}\cong \xi(Gr_{\calf}\calv,\calg^{\bullet })$  est un isomorphisme de faisceaux cohérents sur $S \times {\bf A}^1$ (resp. un isomorphisme de fibrés vectoriels sur $S \times {\bf A}^1$ si les filtrations vérifient ${\bf (H)}$).\\ 
\hspace*{1cm}$(v) $ Pour tout $u_{0}\in {\bf A}^{1}$, $\xi(\calv,\calf^{\bullet },\calg^{\bullet})\vert_{S \times (u_{0},0)}\cong \xi(Gr_{\calg}\calv)$ est un isomorphisme de faisceaux cohérents sur $S $ (resp. un isomorphisme de fibrés vectoriels sur $S $ si les filtrations vérifient ${\bf (H)}$).\\ 
\hspace*{1cm}$(vi) $ Pour tout $v_{0}\in {\bf A}^{1}$, $\xi(\calv,\calf^{\bullet },\calg^{\bullet})\vert_{S \times (0,v_{0})}\cong \xi(Gr_{\calf}\calv,\calg^{\bullet })$  est un isomorphisme de faisceaux cohérents sur $S $ (resp. un isomorphisme de fibrés vectoriels sur $S $ si les filtrations vérifient ${\bf (H)}$).\\
\hspace*{1cm}$(vii)$ $\xi(\calv,\calf^{\bullet },\calg^{\bullet})\vert_{S \times (0,0)}\cong Gr_{\calf}Gr_{\calg}\calv$  est un isomorphisme de faisceaux cohérents sur $S$ (resp. un isomorphisme de fibrés vectoriels sur $S$ si les filtrations vérifient ${\bf (H)}$).\\
Ainsi que les propriétés\\
\hspace*{1cm}$(viii) $ $\xi(\calv,\calf^{\bullet },\calg^{\bullet})\vert_{S \times {\bf A}^{1}\times {\bf G}_{m}}\cong \xi(\calv,\calf^{\bullet },\calt riv^{\bullet})\vert_{S \times  {\bf A}^{1}\times {\bf G}_{m}}$ est un isomorphisme de fibrés vectoriels sur $ S \times  {\bf A}^{1}\times {\bf G}_{m}$.\\
\hspace*{1cm}$(ix) $ $\xi(\calv,\calf^{\bullet },\calg^{\bullet})\vert_{S \times {\bf G}_{m}\times {\bf A}^{1}}\cong \xi(\calv,\calg^{\bullet },\calt riv^{\bullet})\vert_{S \times  {\bf G}_{m}\times {\bf A}^{1}}$ est un isomorphisme de fibrés vectoriels sur $ S \times  {\bf G}_{m}\times {\bf A}^{1}$.\\
\hspace*{1cm}$(x) $ $\xi(\calv,\calf^{\bullet },\calg^{\bullet})\vert_{S \times {\bf G}_{m}\times {\bf G}_{m}}\cong \xi(\calv,\calt riv^{\bullet },\calt riv^{\bullet})\vert_{S \times  {\bf G}_{m} \times {\bf G}_{m}}$ est un isomorphisme de fibrés vectoriels sur $ S \times  {\bf G}_{m}\times {\bf G}_{m}$.\\

\end{proposition}

\begin{preuve}
$(i)$ Soit $i \in I$ tel que $s \in S_{i}$. Il faut montrer que $\xi_{i}(\calv,\calf^{\bullet },\calg^{\bullet})\vert_{\{s\}\times {{\bf A}^2}}\cong  \xi_{{\bf A}^2}(\calv_{s},\calf^{\bullet }_{s},\calg^{\bullet}_{s})$. Soit ${\mathfrak a}_s$ l'idéal maximal de $A_{i}$ correspondant au point $s \in S_{i}$ et ${\mathfrak a}_s[u,v]$ l'idéal de $A_{i}[u,v]$ correspondant à $\{s\} \times {\bf A}^2$. Soit $j: \{s\} \times {\bf A}^{2}  \rightarrow S_{i} \times {\bf A}^2$ le morphisme d'inclusion. Alors $\xi(\calv,\calf^{\bullet },\calg^{\bullet})\vert_{\{s\}\times {{\bf A}^2}}=j^{*}\xi(\calv,\calf^{\bullet },\calg^{\bullet})\cong(B_{i} \otimes_{A_{i}[u,v]}  A_{i}[u,v]/{\mathfrak a}_s[u,v])^{\sim} \cong(B_{i}/{\mathfrak a}_s[u,v].B_{i})^{\sim}$; le dernier isomorphisme vient du lemme précédent. Pour tout $(p,q)$, $M_{i,pq}/{\mathfrak a}_{s}$ est un $k$-espace vectoriel de dimension $\text{dim}_{k}((\calf^{p}\cap\calg^{q})_{s})$ que nous noterons $V_{pq}$. Du diagramme commutatif de $A_{i}$-modules dans lequel les flèches sont injectives on peut déduire par localisation le diagramme de $k$-espaces vectoriels de même forme dans lequel toutes les flèches sont injectives. On retrouve ainsi les filtrations $F^{\bullet }=\calf^{\bullet}_{s}$ et $G^{\bullet }=\calg^{\bullet}_{s}$ du $k$-espace vectoriel $V=M_{i}/{\mathfrak a}_{s}$. D'où l'isomorphisme $B_{i}/{\mathfrak a}_s[u,v].B_{i} \cong R^{2}(V,F^{\bullet },G^{\bullet })=R^{2}(\calv_{s},\calf^{\bullet }_{s},\calg^{\bullet }_{s})$ qui donne le résultat.

$(ii)$ Montrons d'abord qu'il y a isomorphisme sur chacun des ouverts affines $S_{i}$ qui recouvrent $S$. Pour tout $i$, le morphisme d'inclusion $k:S_{i}\times (u_{0},v_{0}) \rightarrow S_{i} \times {\bf A}^{2}$ vient du morphisme surjectif d'anneaux $A_{i}[u,v]\rightarrow A_{i}[u,v]/(u-u_{0},v-v_{0}).A_{i}[u,v]\cong A_{i}$. $\xi(\calv,\calf^{\bullet },\calg^{\bullet})\vert_{S \times (u_{0},v_{0})}=k^{*}\xi(\calv,\calf^{\bullet },\calg^{\bullet})\cong (B_{i}\otimes_{A_{i}[u,v]}A_{i}[u,v]/(u-u_{0},v-v_{0}).A_{i}[u,v])^{\sim}\cong(B_{i}/(u-u_{0},v-v_{0}).B_{i})^{\sim}\cong M^{\sim}$ ce qui permet de conclure pour chacun des ouverts affines. Il suffit ensuite de remarquer que les isomorphismes de recollement sur les $(S_{i}\cap S_{j})\times (u_{0},v_{0})$ et les $S_{i}\cap S_{j}$ commutent avec les restrictions.

$(iii)$ Plaçons-nous sur les ouverts affines. Notons par $k$ le morphisme d'inclusion $S_{i } \times (u_{0},0) \rightarrow S_{i}\times {\bf A}^{1}$ et soit $A_{i}[u,v] \rightarrow A_{i}[u,v]/(u-u_{0},v).A_{i}[u,v]$ le morphisme d'anneaux lui correspondant. 

$(iv)$ Cette assertion se déduit de la précédente en échangeant les rôles de $\calf^{\bullet }$ et $\calg^{\bullet }$.

$(v)$ Exhibons d'abord l'isomorphisme sur les affines $S_{i} \times {\bf G}_{m} \times {\bf A}^1$. Le morphisme d'inclusion $S_{i} \times {\bf G}_{m} \times {\bf A}^1 \rightarrow S_{i} \times {\bf A}^2$ vient du morphisme d'anneaux $A_{i}[u,v] \rightarrow A_{i}[u,u^{-1},v]$. Il suffit donc de vérifier que $B_{i} \otimes_{A_{i}[u,v]} A_{i}[u,u^{-1},v] \cong \sum_{p,q}u^{-p}M_{i,pq}$ est isomorphe au module de Rees obtenu lorsque la deuxième filtration $\calg^{\bullet}$ est triviale. C'est bien le cas puisque ce module n'a pas de terme en $v^{-q}m$ lorsque cette filtration est triviale. Pour conclure on remarque à nouveau que les restrictions commutent aux isomorphismes de recollement.

$(vi)$ C'est le pendant de la précédente en échangeant les deux filtrations.

$(vii)$ Cette assertion se démontre de la même façon que les deux précédentes. On se place sur les ouverts affines $S_{i} \times {\bf G}_{m} \times {\bf G}_m$. L'inclusion $S_{i} \times {\bf G}_{m} \times {\bf G}_{m} \rightarrow S_{i} \times {\bf A}^{2}$ vient du morphisme $A_{i}[u,v] \rightarrow A_{i}[u,u^{-1},v,v^{-1}]$. On a un isomorphisme $B_{i} \otimes_{A_{i}[u,v]} A_{i}[u,u^{-1},v,v^{-1}] \cong M_{i}$ ce qui permet de conclure localement. On peut ensuite recoller, car les restrictions commutent aux isomorphismes de recollement.

\end{preuve}

Sous l'hypothèse $({\bf H})$, on a 

\begin{proposition} Le faisceau réflexif de Rees $\xi(\calv,\calf_{0}^{\bullet },\calf_{1}^{\bullet },\calf_{2}^{\bullet })$ a les propriétés suivantes\\ 
\hspace*{1cm}$\bullet$ Pour tout $s \in S$, $\xi(\calv,\calf_{0}^{\bullet },\calf_{1}^{\bullet },\calf_{2}^{\bullet })\vert_{\{s\}\times {{\bf P}^2}} \cong \xi_{{\bf P}^2}(\calv_{s},\calf_{0}^{\bullet }{}_{s},\calf_{1}^{\bullet }{}_{s},\calf_{2}^{\bullet }{}_{s}).$\\
\hspace*{1cm}$\bullet$ Pour tout $(u_{0}:u_{1}:u_{2}) \in {\bf P}^{2}$ tel que $u_{0}u_{1}u_{2}\neq 0$,  $\xi(\calv,\calf_{0}^{\bullet },\calf_{1}^{\bullet },\calf_{2}^{\bullet })\vert_{S \times  \{(u_{0}:u_{1}:u_{2})\}} \cong \calv$.\\
\hspace*{1cm}$\bullet$ On a aussi $\xi(\calv,\calf_{0}^{\bullet },\calf_{1}^{\bullet },\calf_{2}^{\bullet })\vert_{S \times  \{(1:0:0)\}} \cong Gr_{\calf_{1}}Gr_{\calf_{2}}\calv$,  $\xi(\calv,\calf_{0}^{\bullet },\calf_{1}^{\bullet },\calf_{2}^{\bullet })\vert_{S \times  \{(0:1:0)\}} \cong Gr_{\calf_{0}}Gr_{\calf_{2}}\calv$ et  $\xi(\calv,\calf_{0}^{\bullet },\calf_{1}^{\bullet },\calf_{2}^{\bullet })\vert_{S \times  \{(0:0:1)\}} \cong Gr_{\calf_{0}}Gr_{\calf_{1}}\calv$.\\
Un seul gradué aussi
\end{proposition}

\begin{preuve}
Cette proposition se déduit de la précédente par recollement.

\end{preuve}

\subsubsection{Foncteur de Rees}
Nous introduisons ici quelques notation pour décrire la construction du fibré de Rees relatif parallèlement à ce qui a été fait dans la partie 1 pour les construction de Rees associé aux espaces vectoriels filtrés. Rappelons que nous notions $\Phi_{R}$ le foncteur qui à un espace vectoriel filtré muni de $n$ filtration $(V,F^{\bullet}_{0},...,F_{n-1}^{\bullet })$ associe un faisceau de Rees sur la variété affine ${\bf A}^n$. Nous utilisions la même notation pour le foncteur qui à un espace vectoriel muni de trois filtrations associe le fibré de Rees sur le plan projectif ${\bf P}^2$. 
 
$\Phi_R$ est un foncteur de la catégorie des espaces vectoriels filtrés (resp. bifiltrés, resp. trifiltrés) munie des morphismes strictement compatibles aux filtrations vers la catégorie des faisceaux localement libres sur ${\bf A}^1$ (resp. ${\bf A}^2$, resp. ${\bf P}^2$) équivariants pour une certaine action d'un tore et munie des morphismes équivariants de fibré pour cette action dont le conoyau est sans torsion. 

Nous désignerons par $\Phi_{\calr}$ le foncteur relatif, qui à un fibré vectoriel filtré par des sous-fibrés stricts sur une variété algébrique $S$, $(\calv,\calf^{\bullet}_{0}) \rightarrow S$ (resp. bifiltré par des sous-fibrés stricts $(\calv,\calf^{\bullet}_{0},\calf^{\bullet}_{1}) \rightarrow S$, resp. trifiltré par des sous fibrés stricts $(\calv,\calf^{\bullet}_{0},\calf^{\bullet}_{1},\calf_{2}^{\bullet }) \rightarrow S$) associe le fibré équivariant (resp. faisceau reflexif équivariant) pour l'action d'un tore 
\begin{center}
$ \Phi_{\calr}((\calv,\calf^{\bullet}_{0}) \rightarrow S)=\xi(\calv,\calf^{\bullet}_{0}) \rightarrow S \times {\bf A}^1,$\\
resp. $ \Phi_{\calr}((\calv,\calf^{\bullet}_{0},\calf^{\bullet}_{1}) \rightarrow S)=\xi(\calv,\calf^{\bullet}_{0},\calf^{\bullet}_{1}) \rightarrow S \times {\bf A}^2,$\\
resp. $ \Phi_{\calr}((\calv,\calf^{\bullet}_{0},\calf^{\bullet}_{1},\calf^{\bullet }_{2}) \rightarrow S)=\xi(\calv,\calf^{\bullet}_{0},\calf^{\bullet}_{1},\calf^{\bullet}_{2}) \rightarrow S \times {\bf P}^2.$
\end{center}

Soit $\calc_{nfiltr}(S)$ la catégorie des fibrés vectoriels sur $S$ muni de $n$ filtrations par des sous-fibrés stricts et des morphismes de fibrés vectoriels strictement compatibles aux filtrations.

Soit $\calf ib(S \times {\bf A}^{1}/{\bf G}_{m})$ la catégorie des fibrés ${\bf G}_m$-équivariants sur $S \times {\bf A}^1$ muni de l'action par translation sur le deuxième facteur et des morphismes de faisceaux $G$-équivariants dont le conoyau est sans torsion.

Soit $\calr efl(S \times {\bf A}^{2}/{\bf G}^{2})$ la catégorie des faisceaux réflexifs ${\bf G}^2$-équivariants sur $S \times {\bf A}^2$ muni de l'action par translation sur le deuxième facteur dont les restrictions à $S \times \{(u_{0},v_{0})\}$, $\{s\} \times {\bf A}^2$, $S \times {\bf A}^{1} \times \{ v_{0} \}$, $S \times \{u_{0}\} \times {\bf A}^{1}$ et $S \times {\bf A}^{1} \times \{v_{0}\}$ pour tout $s,u_{0},v_{0}$ sont des fibrés vectoriels et des morphismes de faisceaux $G$-équivariants dont le conoyau est sans torsion.

Soit $\calr efl(S \times {\bf P}^{2}/{\bf T})$ la catégorie des faisceaux réflexifs ${\bf T}$-équivariants sur $S \times {\bf P}^2$ munis de l'action par translation sur le plan projectif dont les restrictions aux trois ouverts standards ${\bf A}^2$ sont dans $\calr efl(S \times {\bf A}^{2}/{\bf G}^{2})$ et des morphismes ${\bf T}$-équivariants dont le noyau est sans torsion.

Alors, la section précédente sur les morphismes de fibrés de Rees montre que l'on a bien défini des foncteurs
\begin{center}
$\Phi_{\calr}:\calc_{1filtr}(S) \rightarrow \calf ib(S \times {\bf A}^{1}/{\bf G}_{m})$\\
$(\calv,\calf^{\bullet }) \mapsto \xi(\calv,\calf^{\bullet })$\\
$( f:(\calv,\calf^{\bullet })) \rightarrow (\calv',\calf^{\bullet }{}'))
\mapsto ( f^{r}:\xi(\calv,\calf^{\bullet }) \rightarrow \xi(\calv',\calf^{\bullet }{}')),$\\
\end{center}

\begin{center}
$\Phi_{\calr}:\calc_{2filtr}(S) \rightarrow \calr efl(S \times {\bf A}^{2}/{\bf G}_{m}^{2})$\\
$(\calv,\calf^{\bullet }_{0},\calf^{\bullet}_{1}) \mapsto \xi(\calv,\calf^{\bullet }_{0},\calf^{\bullet }_{1})$\\
$( f:(\calv,\calf^{\bullet }_{0},\calf^{\bullet}_{1})) \rightarrow (\calv',\calf^{\bullet }_{0}{}',\calf^{\bullet}_{1}{}'))
\mapsto ( f^{r}:\xi(\calv,\calf^{\bullet }_{0},\calf^{\bullet}_{1}) \rightarrow \xi(\calv',\calf^{\bullet }_{0}{}',\calf^{\bullet}_{1}{}'))),$\\
\end{center}

\begin{center}
$\Phi_{\calr}:\calc_{3filtr}(S) \rightarrow \calr efl(S \times {\bf P}^{2}/{\bf T})$\\
$(\calv,\calf^{\bullet }_{0},\calf^{\bullet}_{1},\calf_{2}^{\bullet }) \mapsto \xi(\calv,\calf^{\bullet }_{0},\calf^{\bullet }_{1},\calf^{\bullet }_{2})$\\
$( f:(\calv,\calf^{\bullet }_{0},\calf^{\bullet}_{1},\calf_{2}^{\bullet })) \rightarrow (\calv',\calf^{\bullet }_{0}{}',\calf^{\bullet}_{1}{}',\calf_{2}^{\bullet }{}'))
\mapsto ( f^{r}:\xi(\calv,\calf^{\bullet }_{0},\calf^{\bullet}_{1},\calf_{2}^{\bullet }) \rightarrow \xi(\calv',\calf^{\bullet }_{0}{}',\calf^{\bullet}_{1}{}',\calf_{2}^{\bullet }{}'))).$\\
\end{center}

\subsubsection{Construction relative inverse, foncteur de Rees inverse}
Le but ici est de définir un foncteur inverse du foncteur de Rees relatif qui à un espace vectoriel trifiltré associe un faisceau relatif sur le plan projectif afin d'étendre le dictionnaire ponctuel au cas relatif. Nous allons tout d'abord définir un foncteur inverse relatif ``local'' sur les plans affines pour toute variété algébrique $S$, $$\Phi_{\cali }: \calr efl(S \times {\bf A}^{2}/{\bf G}_{m}^{2}) \rightarrow   \calc_{2filtr}(S).$$

Soit $\xi$ un objet de $ \calr efl(S \times {\bf A}^{2}/{\bf G}_{m}^{2})$. Alors par définition de cette catégorie, $\xi\vert_{S \times \{(1,1)\}}$ est un fibré vectoriel sur $S$ (du moins modulo l'isomorphisme entre $S$ et $S \times \{(1,1)\}$, nous ferons toujours l'abus de langage). Posons 
$$\Phi_{\cali}(\xi):=\xi\vert_{S \times \{(1,1)\}}.$$
L'action permet de munir $\Phi_{\cali}(\xi)$ de deux filtrations.  $\xi$ est cohérent, donc il existe des ouverts affines qui recouvrent $S_{i}=\text{Spec}A_{i}$ et des $A_{i}[u,v]$-modules $M_{i}$ tels que $\xi\vert_{S_{i}\times {\bf A}^2}\cong(M_{i})^{\sim}$. L'action de ${\bf G}_{m}^{2}$ sur le faisceau vient du morphisme d'anneaux $\sigma^{\sharp}: M_{i} \rightarrow M_{i}\otimes_{A_{i}[u,v]}k[x,x^{-1},y,y^{-1}]$. Considérons l'inclusion $j :S_{i}\times {\bf A}^{1} \times \{1\} \rightarrow S_{i}\times {\bf A}^{2}$. L'action de ${\bf G}_{m}$ sur $\xi\vert_{{S_{i}\times {\bf A}^{1} \times \{1\}}}$ qui est un fibré vectoriel par hypothèse vient du morphisme de $A_{i}[u]$-modules  $M_{i}/k[v] \rightarrow M_{i}/k[v]\otimes_{A_{i}[u]}k[x,x^{-1},y,y^{-1}]/k[y,y^{-1}]$. On en déduit une trivialisation par cette action du fibré $\xi\vert{S_{i}\times {\bf A}^{1} \times \{1\}}$ sur $S_{i}\times {\bf G}_{m} \times \{1\}$ 
$$\xi\vert_{{S_{i}\times {\bf G}_{m} \times \{1\}}} \cong k^{*}\xi\vert_{S \times \{(1,1)\}}.$$
Pour tout $i \in I$, l'action à gauche induit une action à droite triviale sur la base $S_{i} \times \{(1,1)\}$ que l'on peut décomposer suivant les caractères
\begin{eqnarray}
\xi\vert_{S_{i} \times \{(1,1)\}}&=&\oplus_{\chi \in X^{*}({\bf G}_{m}(u))}(\xi\vert_{S_{i} \times \{(1,1)\}})^{\chi_u}\\
&=&\oplus_{p \in {\bf Z}}(\xi\vert_{S_{i} \times \{(1,1)\}})^{p_u}
\end{eqnarray}
On a écrit ${\bf G}_{m}(u)$ pour bien indiquer que cette décomposition vient de la ``direction'' $u$ suivant laquelle agit par translation le facteur $x$ dans ${\bf G}_{m}^{2}=\text{Spec}k[x,x^{-1},y,y^{-1}]$. Ceci permet de définir une filtration décroissante et exhaustive par des sous-fibrés stricts de $\Phi_{\cali}(\xi_{i})=\xi\vert_{S_{i} \times \{(1,1)\}}$, $\calf^{\bullet }_{0,i}$ par
$$  \calf^{p}_{0,i}=\oplus_{q \leq p}(\xi\vert_{S_{i} \times \{(1,1)\}})^{p_u}.$$
De manière analogue, en décomposant suivant l'action du deuxième facteur, on peut définir une filtration exhaustive et décroissante par des sous-fibrés stricts de $\Phi_{\cali}(\xi_{i})=\xi\vert_{S_{i} \times \{(1,1)\}}$, $\calf^{\bullet }_{1,i}$,
$$  \calf^{p}_{1,i}=\oplus_{q \leq p}(\xi\vert_{S_{i} \times \{(1,1)\}})^{p_v}.$$
Les filtrations locales par des sous-fibrés stricts se recollent pour former des filtrations par des sous-fibrés stricts de $\Phi_{\cali}(\xi_{})=\xi\vert_{S_{} \times \{(1,1)\}}$.

Le foncteur $\Phi_{\cali}$ est ainsi bien défini au niveau des objets.\\
${}$\\
Considérons un morphisme dans la catégorie $ \calr efl(S \times {\bf A}^{2}/{\bf G}_{m}^{2})$, $g : \xi \rightarrow \xi'$. On en déduit, par restriction à $S \times \{(1,1)\}$, un morphisme de fibrés vectoriels 
$$g^{i}=\Phi_{\cali}(g):=g\vert_{S \times \{(1,1)\}}: \Phi_{\cali}(\xi) \rightarrow \Phi_{\cali}(\xi').$$
Le morphisme $g$ est ${\bf G}_{m}^2$-équivariant donc le morphisme $gî$ respecte les filtrations, c'est un morphisme de fibrés vectoriels filtrés. Reste à montrer qu'il est strict. Comme l'assertion porte sur chacune des filtrations séparement, restreignons (en ne considérant que la première par exemple) le fibré à $S \times {\bf A}^{1} \times \{1\}$.

Le foncteur de Rees relatif inverse $\Phi_{\cali}$ est donc bien défini au niveau des morphismes.\\
${}$\\ 
Justifions son appelation de foncteur inverse.

\begin{proposition}
Pour toute variété algébrique $S$, les foncteurs $\Phi_{\calr}$ et $\Phi_{\cali}$ établissent une équivalence de catégories entre la catégorie $ \calc_{2filtr}(S)$ des fibrés vectoriels sur $S$ munis de deux filtrations exhaustives et décroissantes par des sous-fibrés stricts dont les morphismes sont les morphismes stricts de fibrés vectoriels et la catégorie $ \calr efl(S \times {\bf A}^{2}/{\bf G}_{m}^{2})$ des faisceaux cohérents réflexifs ${\bf G}_{m}^2$-équivariants sur $S \times { \bf G}_{m}^2$ dont les morphismes sont les morphismes équivariants dont le conoyau est sans torsion
$$\xymatrix{
 \calc_{2filtr}(S) \ar@<2pt>[r]^-{\Phi_{\calr}} &   \calr efl(S \times {\bf A}^{2}/{\bf G}_{m}^{2}) \ar@<2pt>[l]^-{\Phi_{\cali}}
}.$$ 
\end{proposition}

\begin{preuve}
Montrons que $\Phi_{\cali}$ établit une correspondance pleinement fidèle et essentiellement surjective. Soit $\calv \in  \calc_{2filtr}(S)$, alors $\Phi_{\calr}(\calv ) \in \calr efl(S \times {\bf A}^{2}/{\bf G}_{m}^{2})$. Par la proposition 31, $ \Phi_{\calr}(\calv ) \vert_{S \times \{(1,1)\}}\cong \calv$. Ainsi $\Phi_{\cali}(\Phi_{\calr}(\calv))\cong \calv$ comme fibrés vectoriels. Les filtrations de $
\Phi_{\cali}(\Phi_{\calr}(\calv))$ coïncident, modulo l'isomorphisme, avec celles de $\calv$.

Montrons que $\Phi_{\cali}$ est essentiellement surjective. Soit $\xi$ et $\xi'$ dans $ \calr efl(S \times {\bf A}^{2}/{\bf G}_{m}^{2})$ et $g,g' \in \text{Hom}_{\calr efl(S \times {\bf A}^{2}/{\bf G}_{m}^{2})}(\xi,\xi')$ telles que $\Phi_{\cali}(g)=\Phi(g')$ dans $\text{Hom}_{ \calc_{2filtr}(S)}(\Phi_{\cali}(\xi),\Phi_{\cali}(\xi'))$. Alors, comme $\Phi_{\calr}(\Phi_{\cali}(g))$ induit un morphisme isomorphe à $g$ entre les objets $\Phi_{\calr}\circ \Phi_{\cali}(\xi)$ et $\Phi_{\calr}\circ \Phi_{\cali}(\xi')$, et que l'isomorphisme est le même pour les deux morphismes identiques $\Phi_{\cali}(g)=\Phi(g')$, il vient $g=g'$, ce qui prouve l'injectivité de la correspondance. La surjectivité est claire par la construction de la section sur le morphisme de Rees, à $f \in \text{Hom}_{ \calc_{2filtr}(S)}(\Phi_{\cali}(\xi),\Phi_{\cali}(\xi'))$ on associe le morphisme $\varphi_{}' \circ \Phi_{\calr}(f) \circ \varphi_{}^{-1} \in  \text{Hom}_{\calr efl(S \times {\bf A}^{2}/{\bf G}_{m}^{2})}(\xi,\xi')$ où $\varphi$ est l'isomorphisme entre $\Phi_{\calr}\circ \Phi_{\cali}(\xi)$ et $\xi$ et $\varphi'$ est l'isomorphisme entre $\Phi_{\calr}\circ \Phi_{\cali}(\xi)$ et $\xi'$.

\end{preuve}

${}$\\
{\bf Remarque :} On déduit bien sûr de cette proposition une équivalence de catégorie pour les objets munis d'une seule filtration. Il suffit pour cela de restreindre la construction précédente faites avec deux filtrations dont l'une est triviale à l'un des deux facteurs, d'où l'équivalence de catégories

$$\xymatrix{
 \calc_{1filtr}(S) \ar@<2pt>[r]^-{\Phi_{\calr}} &   \calf ib(S \times {\bf A}^{1}/{\bf G}_{m}^{}) \ar@<2pt>[l]^-{\Phi_{\cali}}
}.$$

\subsubsection{Dictionnaire relatif}
Définissons quelques catégories qui vont entrer dans le dictionnaire relatif.\\
${}$\\
$\bullet$ Soit $(\calv,\calf^{\bullet }_{0},\calf^{\bullet }_{1},\calf^{\bullet }_{2}) \rightarrow S$ un fibré vectoriel filtré par des sous-fibré stricts et tel que les filtrations soient exhaustives et décroissantes. Comme pour les filtrations d'espaces vectoriels, on dira que les trois filtrations sont simultanéement scindées s'il existe des sous-fibrés $\calv^{p,q,r}$ de $\calv$ tels que pour tout $(p,q,r)$,
$\calf^{p}_{0}= \oplus_{p'\geq p,q,r} \calv^{p',q,r}$, $\calf^{q}_{1}= \oplus_{p,q'\geq q,r} \calv^{p,q',r}$ et $\calf^{r}_{2}= \oplus_{p,q,r'\geq r} \calv^{p,q,r'}$. 

Soit $ \calc_{3filtr,scind}(S)$ la sous-catégorie pleine de $ \calc_{3filtr}(S)$ dont les objets sont simultanéement scindés.

On notera par  $\calf ib_{scind}(S \times {\bf P}^{2}/{\bf T})$ la sous-catégorie pleine de $\calr efl(S \times {\bf P}^{2}/{\bf T})$ dont les objets sont en plus des fibrés vectoriels sommes de fibrés en droites.\\
${}$\\
$\bullet $ Soit $(\calv,\calf^{\bullet }_{0},\calf^{\bullet }_{1},\calf^{\bullet }_{2}) \rightarrow S$ un objet de $ \calc_{3filtr}(S)$. On dit que les filtrations par des fibrés stricts sont opposées si en tout point les filtrations induites sur les fibres ont cette propriété, i.e. si pour tout $s \in S$, $  (\calv_{s},\calf^{\bullet }_{0}{}_{s},\calf^{\bullet }_{1}{}_{s},\calf^{\bullet }_{2}{}_{s})$ est un espace vectoriel muni de trois filtrations opposées.

On notera $ \calc_{3filtr,opp}(S)$ la sous-catégorie pleine de $ \calc_{3filtr}(S)$ dont les objets sont des fibrés vectoriels filtrés dont les filtrations sont opposées.

Soit  $\calr efl_{\mu-semistable/S,\mu=0/S}(S \times {\bf P}^{2}/{\bf T})$ la sous-catégorie pleine de $\calr efl(S \times {\bf P}^{2}/{\bf T})$ dont les objets restreints à $\{s \} \times {\bf P}^2$ pour tout $s \in S$ sont $\mu$-semistables de pente $0$.\\
${}$\\

{\bf Ce dictionnaire ne marche que sur les strates sur lesquelles }$({\bf H})$ {\bf est vérifiée, les faisceaux réflexifs devenant des fibrés}

Nous avons donc montré sous l'hypoth\`ese $({\bf H})$,
\begin{theoreme} Soit $S$ une variété algébrique sur un corps algébriquement clos de caractéristique nulle $k$. La construction des faisceaux de Rees relatifs sur $S \times {\bf P}^2$ établit les équivalences de catégories entre :\\
$\bullet$ La catégorie $ \calc_{3filtr}(S)$ des fibrés vectoriels sur $S$ trifiltrés par des sous-fibrés stricts tels que les filtrations soient exhaustives et décroissantes munie des morphismes de fibrés strictement compatibles aux filtrations et la catégorie $\calf ib(S \times {\bf P}^{2}/{\bf T})$ des fibrés vectoriels ${\bf T}$-équivariants sur $S \times {\bf P}^2$ dont les restrictions ont les propriétés voulues munie des morphismes équivariants de fibrés dont le conoyau est sans torsion :
$$\xymatrix{
 \calc_{3filtr}(S) \ar@<2pt>[r]^-{\Phi_{\calr}} &   \calf ib(S \times {\bf P}^{2}/{\bf T}) \ar@<2pt>[l]^-{\Phi_{\cali}}
}.$$ 
$\bullet$ La catégorie $ \calc_{3filtr,scind}(S)$ des fibrés vectoriels sur $S$ trifiltrés par des sous-fibrés stricts tels que les filtrations soient simultanéement scindées, exhaustives et décroissantes munie des morphismes de fibrés strictement compatibles aux filtrations et la catégorie $\calf ib(S \times {\bf P}^{2}/{\bf T})$ des fibrés ${\bf T}$-équivariants sur $S \times {\bf P}^2$ dont les restrictions ont les propriétés voulues munie des morphismes équivariants de fibrés dont le conoyau est sans torsion : 
$$\xymatrix{
 \calc_{3filtr,scind}(S) \ar@<2pt>[r]^-{\Phi_{\calr}} &   \calf ib_{scind}(S \times {\bf P}^{2}/{\bf T}) \ar@<2pt>[l]^-{\Phi_{\cali}}
}.$$ 
$\bullet$ La catégorie $ \calc_{3filtr,opp}(S)$ des fibrés vectoriels sur $S$ trifiltrés par des sous-fibrés stricts tels que les filtrations soient opposées, exhaustives et décroissantes munie des morphismes de fibrés strictement compatibles aux filtrations et la catégorie $\calf ib(S \times {\bf P}^{2}/{\bf T})$ des fibrés vectoriels ${\bf T}$-équivariants sur $S \times {\bf P}^2$ dont les restrictions ont les propriétés voulues et dont les restrictions à $\{s\} \times {\bf P}^2$ pour tout $s \in S$ sont ${\bf P}^{1}_{0}$-semistables de pente $0$ munie des morphismes équivariants de fibrés dont le conoyau est sans torsion : 
$$\xymatrix{
 \calc_{3filtr,opp}(S) \ar@<2pt>[r]^-{\Phi_{\calr}} &   \calf ib_{{\bf P}^{1}_{0}-semistable/S,\mu=0/S}(S \times {\bf P}^{2}/{\bf T}) \ar@<2pt>[l]^-{\Phi_{\cali}}
}.$$ 
%Si $k={\bf C}$. On a de plus :\\
%$\bullet$ La catégorie $ \calc_{3filtr,opp,{\bf R}}(S)$ des fibrés vectoriels à système local réel sous-jacent sur $S$ trifiltrés par des sous-fibrés stricts tels que les filtrations soient conjuguées, opposées, exhaustives et décroissantes munie des morphismes de fibrés strictement compatibles aux filtrations et la catégorie $\calf ib(S \times {\bf P}^{2}/{\bf T})$ des fibrés vectoriels ${\bf T}^{\tau}$-équivariants sur $S \times {\bf P}^2$ dont les restrictions ont les propriétés voulues et dont les restrictions à $\{s\} \times {\bf P}^2$ pour tout $s \in S$ sont ${\bf P}^{1}_{0}$-semistables de pente $0$ munie des morphismes équivariants de fibrés dont le conoyau est sans torsion : 
%$$\xymatrix{
% \calc_{3filtr,opp,{\bf R}}(S) \ar@<2pt>[r]^-{\Phi_{\calr}} &   \calf ib_{{\bf P}^{1}_{0}-semistable/S,\mu=0/S}(S \times {\bf P}^{2}/{\bf T}^{\tau}) \ar@<2pt>[l]^-{\Phi_{\cali}}
%}.$$ 
\end{theoreme}

%\subsubsection{Foncteurs fibres}

\subsection{Stratification associée à une famille d'espaces vectoriels trifiltrés}

On s' intéresse au cas où l'hypothèse $({\bf H})$ n'est pas forcément satisfaite.

\subsubsection{Stratification platifiante}

Soit $\calf$ un faisceau cohérent sur un schéma projectif $X$ sur un corps algébriquement clos de caractéristique nulle $k$. Alors la caractéristique d'Euler de $\calf$ est
$$\chi(X,\calf)=\sum_{i}(-1)^{i}h^{i}(X,\calf),$$
où $h^{i}(X,\calf)=\text{dim}_{k}H^{i}(X,\calf)$. Si on fixe un fibré en droite ample $\calo(1)$ sur $X$, alors le polyn\^ome de Hilbert de $\calf$, $P(X,\calf)$ est donné par $$m \mapsto \chi(X,\calf \otimes \calo(m)).$$
Supposons, par exemple, que $X$ soit une surface et $\calf $ soit sans torsion et de rang $r$, alors la formule de Riemann-Roch (cf \cite{fri}) donne
$$ \chi(X,\calf)=\dfrac{\text{c}_{1}(\calf).(\text{c}_{1}(\calf)-K_{X})}{2}-\text{c}_{2}(\calf)+r\chi(X,\calo_{X}).$$
Les formules $\text{c}_{1}(\calf \otimes L)=\text{c}_{1}(\calf)+r\text{c}_{1}(L)$ et $\text{c}_{2}(\calf \otimes L)=\text{c}_{2}(\calf)+(r-1)\text{c}_{1}(\calf).\text{c}_{1}(L)+\binom{r}{2} \text{c}_{1}(L)^{2}$ où $L$ est un fibré en droite sur $X$ permettent alors de trouver le polyn\^omes de Hilbert de $\calf$. Nous utiliserons ce résultat dans le cas où $X={\bf P}^{2}$.

Supposons maintenant que nous ayons un faisceau cohérent de $\calo_{X}$-modules $\calf$ sur un schéma $X$ sur $S$ où $f:X \rightarrow S$ est un morphisme de type fini entre schémas n{\oe}theriens. Pour tout $s \in S$ on note la fibre $f^{-1}(s)=\text{Spec}(k(s))\times_{S}X$ de $f$ par $X_{s}$. On note aussi $\calf_{s}$ pour $\calf\vert_{X_s}$. On pense à $\calf$ comme une famille de polyn\^omes paramétrisée par $S$. La notion de ``continuité'' en famille des $\calf_{s}$ est donnée par la définition :

\begin{definition}
Une famille plate de faisceaux cohérents sur les fibres de $f$ est un faisceau de $\calo_{X}$-modules cohérents $\calf$ qui est plat sur $S$ i.e. pour tout $x \in X$, $\calf_{x}$ est plat au dessus de l'anneau local $\calo_{S,f(x)}$. 
\end{definition} 
Nous allons associer une stratification de $S$ à tout faisceau cohérent $\calf $ sur $S$ telle que sur chacune des strates $S_{P}$, strates données par la constance du polyn\^ome des Hilbert des fibres $P(X_{s},\calf_{s})$, la restriction de $\calf$, $\calf\vert_{S_{P}}$ soit une famille plate de faisceaux cohérents. 

Pour cela on utilise un théorème d\^u à Mumford qui permet de platifier tout faisceau cohérent en découpant de façon adéquate la base par des sous-schémas fermés.

\begin{theoremese}(Mumford, cf \cite{huyleh} Theorem 2.1.5.) Soit $f:X \rightarrow S$ un morphisme projectif entre schémas n{\oe}thériens, soit $\calo(1)$ un faisceau inversible sur $X$ qui est très ample par rapport à $S$ et $\calf$ un $\calo_{X}$-module cohérent. Alors l'ensemble $\calp=\{ P(X_{s},\calf_{s}) \vert s \in S \}$ des polyn\^omes de Hilbert des fibres de $\calf $ est fini. De plus il y a un nombre fini de sous-schémas locallement fermés $S_{P}\subset S$, indexés par les polyn\^omes $P \in \calp$ qui satisfont aux propriétés :\\
$\hspace*{2cm}(i)$ Le morphisme naturel $j:\coprod_{P}S_{P} \rightarrow S$ est une bijection.\\
$\hspace*{2cm}(ii)$ Si $g: S'\rightarrow S$ est un morphisme de schémas n{\oe}theriens, alors $g^{*}_{X}\calf$ est plat au dessus de $S'$ si et seulement si $g$ se factorise par $j$ où $g^{*}_{X}:X \times_{S}S' \rightarrow S'$ est le morphisme canonique du produit fibré vers $S'$.
\end{theoremese}

La stratification $S=\coprod_{P \in \calp}S_{P}$ est la stratification platifiante associée au faisceau cohérent $\calf$.

\subsubsection{Application à une famille d'espaces vectoriels filtrés}

Soit $(\calv,\calf^{\bullet }_{0},\calf^{\bullet }_{1},\calf^{\bullet }_{2}) \rightarrow S$ un fibré vectoriel filtré par des sous-fibrés stricts et tel que les filtrations soient exhaustives et décroissantes sur une variété algébrique lisse $S$. La construction de Rees associée à cette famille fournit un faisceau de $\calo_{X \times {\bf P}^2}$-modules cohérent $\xi(\calv,\calf^{\bullet }_{0},\calf^{\bullet }_{1},\calf^{\bullet }_{2})$. Considérons le morphisme de projection sur le deuxième facteur $f:X \times {\bf P}^{2} \rightarrow X$. Le faisceau $f^{*}\calo_{{\bf P}^2}(1)$ est très ample par rapport à $S$. Par le théorème précédent on obtient 
\begin{proposition}
L'ensemble  $\calp=\{ P(X_{s},\xi(\calv,\calf^{\bullet }_{0},\calf^{\bullet }_{1},\calf^{\bullet }_{2})_{s}) \vert s \in S \}$ des polyn\^omes de Hilbert des fibres de $\xi(\calv,\calf^{\bullet }_{0},\calf^{\bullet }_{1},\calf^{\bullet }_{2})$ est fini. De plus il y a un nombre fini de sous-schémas locallement fermés $S_{P}\subset S$, indexés par les polyn\^omes $P \in \calp$ qui satisfont aux propriétés :\\
$\hspace*{1cm}(i)$ Le morphisme naturel $j:\coprod_{P}S_{P} \rightarrow S$ est une bijection.\\
$\hspace*{1cm}(ii)$ Si $g: S'\rightarrow S$ est un morphisme de schémas n{\oe}theriens, alors le tiré en arrière du faisceau de Rees $g^{*}_{X}\calf$ est plat au dessus de $S'$ si et seulement si $g$ se factorise par $j$ où $g^{*}_{X}:X \times_{S}S' \rightarrow S'$ est le morphisme canonique du produit fibré vers $S'$. 
\end{proposition}
Pour relier cette stratification par le polyn\^ome de Hilbert des fibres, déterminés par la première et la deuxième classe de Chern dans $H^{*}({\bf P}^2,k)$ des fibres du faisceau de Rees, aux invariants discrets première et deuxième classe de Chern des fibrés associés aux espaces vectoriels filtrés donnés par les fibres des filtrations, il faut comparer pour tout $s \in S$ d'une part
$$\text{c}_{1}(\xi(\calv,\calf^{\bullet }_{0},\calf^{\bullet }_{1},\calf^{\bullet }_{2})_{s})_{s} \text{ à } \text{c}_{1}(\xi_{{\bf P}^2}(\calv_{s},\calf^{\bullet }_{0}{}_{s},\calf^{\bullet }_{1}{}_{s},\calf^{\bullet }_{2})_{s}{}_{s}),$$
et d'autre part 
$$\text{c}_{2}(\xi(\calv,\calf^{\bullet }_{0},\calf^{\bullet }_{1},\calf^{\bullet }_{2})_{s})_{s} \text{ à } \text{c}_{2}(\xi_{{\bf P}^2}(\calv_{s},\calf^{\bullet }_{0}{}_{s},\calf^{\bullet }_{1}{}_{s},\calf^{\bullet }_{2})_{s}{}_{s}).$$

\subsubsection{Filtrations opposées}
Considérons encore $(\calv,\calf^{\bullet }_{0},\calf^{\bullet }_{1},\calf^{\bullet }_{2}) \rightarrow S$ fibré vectoriel filtré par des sous-fibrés stricts et tel que les filtrations soient exhaustives et décroissantes sur une variété algébrique $S$. Supposons de plus que pour tout $s \in S$ l'espace vectoriel trifiltré défini par ces données $(\calv_{s},\calf^{\bullet }_{0}{}_{s},\calf^{\bullet }_{1}{}_{s},\calf^{\bullet }_{2}{}_{s})$ soit un espace vectoriel muni de filtrations opposées et telles que $\calf^{\bullet }_{1}{}_{s}$ et $\calf^{\bullet }_{2}{}_{s}$ soient positives (tous les termes d'indices négatifs sont égaux à $\calv$) et que $\calf^{\bullet }_{0}{}_{s}$ et $\calf^{\bullet }_{1}{}_{s}$ d'une part et $\calf^{\bullet }_{0}{}_{s}$ et $\calf^{\bullet }_{2}{}_{s}$ d'autre part vérifient l'hypothèse ${\bf (H)}$\footnote{Ces hypothèses seront vérifiées lorsque les données viendront de la théorie de Hodge}. On veut alors relier les strates de la stratification platifiante de $\xi(\calv,\calf^{\bullet }_{0},\calf^{\bullet }_{1},\calf^{\bullet }_{2})_{s}$ aux classes de Chern des fibrés de Rees associés aux espaces vectoriels donnés en chaque point.\\
${}$\\
\textbf{Conjecture :} Soit $(\calv,\calf^{\bullet }_{0},\calf^{\bullet }_{1},\calf^{\bullet }_{2}) \rightarrow S$ un fibré vectoriel filtré par des sous-fibrés stricts tel que les filtrations soient opposées, exhaustives et décroissantes sur une variété algébrique $S$ et telle que de plus que pour tout $s \in S$ l'espace vectoriel trifiltré défini par ces données $(\calv_{s},\calf^{\bullet }_{0}{}_{s},\calf^{\bullet }_{1}{}_{s},\calf^{\bullet }_{2}{}_{s})$ soit un espace vectoriel muni de filtrations opposées et telles que $\calf^{\bullet }_{1}{}_{s}$ et $\calf^{\bullet }_{2}{}_{s}$ soient positives et que $\calf^{\bullet }_{0}{}_{s}$ et $\calf^{\bullet }_{1}{}_{s}$ d'une part et $\calf^{\bullet }_{0}{}_{s}$ et $\calf^{\bullet }_{2}{}_{s}$ d'autre part vérifient l'hypothèse $({\bf H})$, alors la filtration platifiante associée à $\xi(\calv,\calf^{\bullet }_{0},\calf^{\bullet }_{1},\calf^{\bullet }_{2})_{s})$ est au changement d'indice près la stratification associée à la deuxième classe de Chern des fibrés sur les fibres de $f: S \times {\bf P}^{2} \rightarrow {\bf P}^{2}$ associés aux espaces vectoriels filtrés donnés par les fibres $\xi_{{\bf P}^2}(\calv_{s},\calf^{\bullet }_{0}{}_{s},\calf^{\bullet }_{1}{}_{s},\calf^{\bullet }_{2}{}_{s})$.

\newpage
$\,$
\vspace{2cm}
$\,$\\
{\bf Dictionnaire sous l'hypothèse ${\bf (H)}$: Famille d'espaces vectoriels trifiltrés-Faisceaux équivariants.}
${}$\\
${}$\\
\begin{center}
%\begin{tabular}{|c||c|}
\begin{tabular}
{|>{\centering}m{3in}
||>{\centering}m{3in}|}
\hline

Catégorie de familles d'espaces multifiltrés sur $S$ & Catégorie de fibrés équivariants sur $S \times {\bf P}^2$\\  %\hline 

\end{tabular}

\begin{tabular}
{|>{\centering}m{3in}
||>{\centering}m{3in}|}
\hline

 $\calc_{3filtr}(S) $ &  $  {\calf ibr}(S \times {\bf P}^{2}/{\bf T}) $ \\

\end{tabular}

\begin{tabular}
{|>{\centering}m{3in}
||>{\centering}m{3in}|}
\hline

$ \calc_{3filtr,scind}(S) $ &  $  {\calf ib}_{scind}(S \times {\bf P}^{2}/{\bf T}) $ \\

\end{tabular}

\begin{tabular}
{|>{\centering}m{3in}
||>{\centering}m{3in}|}
\hline

$ \calc_{3filtr,opp}(S) $ &  $  {\calf ibr}_{{\bf P}^{1}_{0}-semistable/S,\mu=0/S}(S \times {\bf P}^{2}/{\bf T}) $ \\
\end{tabular}

\begin{tabular}
{|>{\centering}m{3in}
||>{\centering}m{3in}|}
\hline

$ \calc_{3filtr,opp,{\bf R}}(S) $ &  $  {\calf ibr}_{{\bf P}^{1}_{0}-semistable/S,\mu=0/S}(S \times {\bf P}^{2}/{\bf T}^{\tau}) $ \\
\end{tabular}

\begin{tabular}
{>{\centering}m{3in}
>{\centering}m{3in}}
\hline
$\,$ & $\,$ \\
\end{tabular}

%\bottomcaption{Dictionnaire sous l'hypothèse ${\bf (H)}$: Famille d'espaces vectoriels trifiltrés-Faisceaux équivariants.}
\end{center}

\vspace{2cm}
$\,$\\

%On indicera par $1+$ ou $2+$ les catégories pour lesquelles la filtration indexée par $1$ et la filtration indexée par $2$ sont positives.

\vspace{1cm}
$\,$\\
{\bf Dictionnaire sous l'hypothèse} ${\bf (H)}$ {\bf pour} $(\calf^{\bullet}_{0},\calf^{\bullet }_{1})$ {\bf et} $(\calf^{\bullet}_{0},\calf^{\bullet }_{2})$.
\begin{center}
%\begin{tabular}{|c||c|}
\begin{tabular}
{|>{\centering}m{3in}
||>{\centering}m{3in}|}
\hline

Catégorie de familles d'espaces multifiltrés sur $S$ & Catégorie de fibrés équivariants sur $S \times {\bf P}^2$\\  %\hline 

\end{tabular}

\begin{tabular}
{|>{\centering}m{3in}
||>{\centering}m{3in}|}
\hline

$ \calc_{3filtr,opp}(S) $ &  $  {\calr efl}_{{\bf P}^{1}_{0}-semistable/S,\mu=0/S}(S \times {\bf P}^{2}/{\bf T}) $ \\
\end{tabular}

\begin{tabular}
{|>{\centering}m{3in}
||>{\centering}m{3in}|}
\hline

$ \calc_{3filtr,opp,{\bf R}}(S) $ &  $  {\calr efl}_{{\bf P}^{1}_{0}-semistable/S,\mu=0/S}(S \times {\bf P}^{2}/{\bf T}^{\tau}) $ \\
\end{tabular}

\begin{tabular}
{>{\centering}m{3in}
>{\centering}m{3in}}
\hline
$\,$ & $\,$ \\
\end{tabular}

%\bottomcaption{Dictionnaire sous l'hypothèse ${\bf (H)}$ pour $(\calf^{\bullet}_{0},\calf^{\bullet }_{1})$ et $(\calf^{\bullet}_{0},\calf^{\bullet }_{2})$}
\end{center}

\subsection{Transversalité}
\subsubsection{Une filtration}
Revenons à la construction du fibré de Rees sur $S \times {\bf A}^1$ associée à un fibré vectoriel $\calv$ sur $S$ muni d'une filtration $\calf^{\bullet }$ par des sous-fibrés vectoriels stricts de sorte que la filtration soit exhaustive et décroissante. On obtient un fibré $\xi(\calv,\calf^{\bullet})$ ${\bf G}_{m}$-équivariant sur $S \times {\bf A}^1$ pour l'action par translation dans la direction de la droite affine et triviale dans l'autre direction.

Supposons que le fibré vectoriel $\calv$ sur $S$ soit muni d'une connexion intégrable 
$$\nabla : \calv \rightarrow \Omega_{S}\otimes_{\calo_{S}}\calv.$$   
Alors $\calv $ est muni d'une ``action'' du fibré tangent $T(S)$ déduite de la connexion. Elle est donnée par
\begin{center}
$T(S)\otimes_{\calo_S}\calv \rightarrow \calv ,$\\
$ \eta \otimes v \mapsto \nabla(v)(\eta).$
\end{center}

\begin{lemme}\cite{sim1} Soit $\calv$ un fibré vectoriel sur $S$ muni d'une connexion intégrable $\nabla$. Supposons aussi que $\calv$ est muni d'une filtration exhaustive et décroissante par des sous-fibrés stricts. Alors $(\calf^{\bullet},\nabla)$ satisfait la transversalité de Griffith $\nabla(\calf^{p}) \subset \Omega_{S}\otimes_{\calo_{S}}\calf^{p-1}$ si et seulement si l'action de $T(S)$ sur $\calv$ s'étend à une action ${\bf G}_{m}$-invariante du faisceau $T(S \times {\bf A}^{1}/{\bf A}^{1})(-S\times\{0\})$ des champs de vecteurs relatifs qui s'annulent au premier ordre suivant $S \times \{ 0\}$ sur le fibré $\xi(\calv,\calf^{\bullet})$.   

\end{lemme}

\begin{preuve}
Comme la question est locale sur $S$, on peut supposer que $\calv=\oplus_{p}\calv^{p}$ est une somme de fibrés vectoriels triviaux et que la filtration par des sous-fibrés est donnée par $\calf^{p}=\oplus_{q \geq p}\calv^{q}$ (dans le langage de la construction du fibré de Rees relatif on se place sur un ouvert affine $S_{i}$ trivialisant le fibré vectoriel). La connexion $\nabla$ est donnée en terme de la connexion triviale $d$ par $\nabla=d+A$, où $A$ est une $1$-forme à valeur dans $\text{End}(\calv)$. On peut décomposer $A$ en somme directe $A=\oplus_{p,q}A_{p,q}$ où $A_{p,q}$ est une $1$-forme à valeurs dans $\text{Hom}(V^{p},V^{q})$. La transversalité de Griffith est alors équivalente à
$$A_{p,q}=0 \text{ pour } q < p-1.$$  
Comme remarqué précédement, on peut écrire
$$\xi(\calv,\calf^{\bullet })=\oplus_{p}u^{-p}\calo_{{\bf A}^1}\otimes \calv^{p}.$$
L'action de ${\bf G}_m$ induit un isomorphisme $\varphi=\oplus_{p}\varphi_{p}$
$$\varphi:\xi(\calv,\calf^{\bullet })\cong \oplus_{p} \calo_{{\bf A}^1}\otimes \calv^{p},$$
où $\varphi_{p}$ agit sur $\calo_{{\bf A}^1}\otimes \calv^{p}$ par multiplication par $u^p$.

On peut étendre la connexion $\nabla$ à une action ${\bf G}_m$-invariante du champs des vecteurs tangents relatifs dans la direction de $S$. Elle est donnée en terme de l'isomorphisme $\varphi$ par la flèche en pointillés qui rend le diagramme suivant commutatif
$$\xymatrix{
T(S\times {\bf A}^{1}/{\bf A}^{1})\otimes \xi(\calv,\calf^{\bullet }) \ar@{.>}[r]\ar[d]^{id \otimes \varphi} & \xi(\calv,\calf^{\bullet })\\
T(S\times {\bf A}^{1}/{\bf A}^{1})\otimes  \oplus_{p} \calo_{{\bf A}^1}\otimes \calv^{p} \ar[r]^-{id \otimes \nabla} & \oplus_{p} \calo_{{\bf A}^1}\otimes \calv^{p} \ar[u]^{\varphi^{-1}}\\
}.$$
Les éléments du champs tangent sont écrits sous la forme $\sum_{i}u^{i}\eta_{i}$ où les vectreurs $\eta_{i}$ sont tangents à $S$. L'action d'un élément $\sum_{i}u^{i}\eta_{i}$ sur $\xi(\calv,\calf^{\bullet })=\oplus_{p}u^{-p}\calo_{{\bf A}^1}\otimes \calv^{p}$ est donnée par la matrice $\sum_{i}u^{i}A(\eta_{i})$ i.e. que l'image d'un élément $\sum_{p}u^{-p}v_{p}$ par $\sum_{i}u^{i}\eta_{i}$ est $\sum_{i,p,q}u^{i-p}A_{p,q}(\eta_{i})v_{q}$. Lorsqu'on transporte ceci par l'isomorphisme $\varphi$, l'action sur $ \oplus_{p} \calo_{{\bf A}^1}\otimes \calv^{p}$ est donnée par la matrice $\sum_{i,p,q}u^{i-p+q}A_{p,q}(\eta_{i})$. C'est l'action sur un voisinage d'un point de $S \times \{0\}$. L'action s'étend sur $S \times \{0\}$ si et seulement si les termes en puissance de $u$ négatives s'annulent. Ainsi la condition que cette action s'étende à tout champs de vecteur de la forme $\sum_{i\geq 1}u^{i}\eta_{i}$, i.e. les champs de vecteurs relatifs qui s'annulent au premier ordre suivant $S \times \{0\}$, est équivalente à la condition $A_{p,q}=0$ pour $q<p-1$ qui est la condition de transversalité de Griffiths. 
\end{preuve}

%\subsubsection{Plusieurs}

%\subsection{le champs des ${\bf T}^3$-fibré sur ${\bf P}^{2}$}

%\subsubsection{Construction}

%\subsubsection{? est un champs algébrique}

%\subsection{Semi-continuité, faiscaux $U$-cohérents}

%\subsection{Cadre analytique}

\newpage
\section{Fibrés vectoriels sur le plan projectif et structures de Hodge mixtes}

\subsection{Rappels sur les structures de Hodge mixtes}

\subsubsection{Définitions}

\begin{center} 
{\bf Structures de Hodge pures}
\end{center}
\begin{definition}
Une structure de Hodge pure de poids $n$, $(H_{\bf Z},(H_{\bf C},F^{\bullet }))$, est la donnée d'un ${\bf Z}$-module de type fini $H_{\bf Z}$ tel que $H_{\bf C}={\bf C} {\otimes}_{\bf Z} H_{\bf Z}$ soit muni d'une filtration finie d\'ecroissante $F^{\bullet }$ telle que :
 \begin{center}
$H_{\bf C}=F^{0}H_{\bf C}\supset ...\supset F^{p}H_{\bf C} \supset F^{p+1}H_{\bf C} \supset ...\supset F^{n+1}H_{\bf C}=\{ 0 \}$ et\\
$F^{p}H_{\bf C} \oplus \overline{F^{q}}H_{\bf C}=H_{\bf C}$
\end{center} 
pour tous les couples $(p,q)$ tels que $p+q=n+1$, o\`u ${\overline F}^\bullet$ est la filtration compl\`ete d\'eduite de $F^{\bullet }$ par conjugaison complexe par rapport \`a la structure r\'eelle sous-jacente sur $H_{\bf R}={\bf R} {\otimes}_{\bf Z} H_{\bf Z}$.\\
\end{definition}
${}$\\
{\bf Exemple:}
Structure de Hodge de Tate $T \langle -p \rangle$. C'est l'unique structure de Hodge pure de type $(-p,-p)$ et de réseau entier $(2 \pi i)^{p} {\bf Z}$.\\

La définition précédente est équivalente à une décomposition de $H_{\bf C}$ en somme directe $H_{\bf C}=\oplus_{p+q=n} \,H^{p,q}$ telle que $H^{p,q}={\overline H}^{q,p}$. L'équivalence est donnée dans un sens par $\{ H^{p,q} \} \mapsto F^{p}=\oplus_{p'\geq p}H^{p',n-p'}$ et dans l'autre sens par $ F^{\bullet } \mapsto H^{p,q}=F^{p} \cap {\overline F}^{q}$.
\begin{center}
{\bf Structures de Hodge mixtes}
\end{center}

\begin{definition}
Une structure de Hodge mixte sur un espace $H_{\bf Q}$ consiste en les données suivantes :\\
                 \hspace*{.5cm}
{\bf (i)} \ Une filtration croissante $W_{\bullet }$ sur $H_{\bf Q}$ appelée filtration par le poids.\\
                \hspace*{.4cm}
{\bf (ii)} \ Une filtration décroissante appelée filtration de Hodge $F^{\bullet }$ sur $H=H_{\bf Q}\otimes_{\bf Q}{\bf C}$ qui induit avec sa filtration conjuguée par rapport à la structure réelle sous-jacente $H_{\bf R}=H_{\bf Q}\otimes_{\bf Q}{\bf R}$, ${\overline F}^{\bullet }$ une structure de Hodge pure de poids $n$ sur
$Gr^{W}_{n}H=W_{n}\otimes_{\bf Q}{\bf C}/W_{n-1}\otimes_{\bf Q}{\bf C}$.
\end{definition}

Dans le ${\bf (ii)}$, la filtration $F^{\bullet }Gr^{W}_{n}H$ induite par $F^{\bullet }$ sur $Gr^{W}_{n}H$ est donnée par les quotients successifs :
 $$F^{p }Gr^{W}_{n}H=(F^{p} \cap (W_{n}\otimes_{\bf Q}{\bf C})+W_{n-1}\otimes_{\bf Q}{\bf C}))/ W_{n-1}\otimes_{\bf Q}{\bf C}.$$

Pour \^etre plus précis la condition {\bf (ii)} signifie que $F^{\bullet }Gr^{W}_{n}H$ et ${\overline F}^{\bullet }Gr^{W}_{n}H$ sont $n$-opposées sur $Gr^{W}_{n}H$ ce qui signifie, d'après la première partie que $Gr^{p}_{F^{\bullet }Gr^{W}_{n}H}Gr^{q}_{{\overline F}^{\bullet }Gr^{W}_{n}H}H =0$ seulement si $p+q \neq n$ où ${\overline F}^{\bullet}$ est la filtration conjuguée à la filtration de Hodge par rapport à la structure réelle sous-jacente $H_{\bf R}=H_{\bf Q} \otimes_{\bf Q} {\bf R}$.

Nous n'utiliserons par la suite, ni la structure rationnelle $H_{\bf Q}$, ni le réseau entier $H_{\bf Z}$. Il nous suffira de considérer la filtration par le poids sur $H_{\bf C}=H_{\bf Q} \otimes_{\bf Q} {\bf C}$ que nous noterons $W_{\bullet }$ alors que nous la notions $W_{\bf C} \otimes_{\bf Q} {\bf C}$. Précisons les différents contextes dans les quels nous nous placerons.

\begin{definition}${}$\\
${\bf (i)}$ Soit $\crzmhs$ la catégorie des structures de Hodge mixtes comme définies ci-dessus (au sens usuel). Les morphismes sont les morphismes de réseaux entiers tels que les morphismes induits au niveau des filtrations soient strictement compatibles aux filtrations.\\

${\bf (ii)}$ Soit $\crmhs$ la catégorie des structures de Hodge mixtes sans strutures entière et rationnelle sous-jacentes. Les objets sont les ${\bf R}$-espaces vectoriels filtrés $(H_{\bf R}, W_{\bullet })$ (filtration par le poids) tels que $H=H_{\bf R}\otimes_{\bf R}{\bf C}$ soit muni d'une filtration appelée filtration de Hodge de sorte que $F^{\bullet}$ et ${\overline F}^{\bullet}$ soient des filtrations $n$-opposées sur $Gr^{W}_{n}H$. Les morphismes sont les morphismes de ${\bf R}$-espaces vectoriels strictement compatibles à la filtration par le poids qui induisent des morphismes strictement compatibles à toutes les filtrations sur les ${\bf C}$-vectoriels.\\
${\bf (iii)}$ Soit $\cmhs$ la catégorie des structures de Hodge mixtes complexes sans structure réelle, rationnelle ou entière sous-jacente. Les objets sont des ${\bf C}$-espaces vectoriels $H$ munis de trois filtrations $(W_{\bullet},F^{\bullet},{\hat F}^{\bullet })$ opposées. Les morphismes sont les morphismes strictement compatibles aux filtrations.

\end{definition}

On considère les foncteurs oublis suivants entre ces différentes catégories de structures de Hodge:\\
$ \Phi_{z}: \crzmhs \rightarrow \crmhs $ qui est le foncteur oubli de la structure entière et rationnelle sous-jacente à une struture de Hodge mixte au sens usuel.\\
$ \Phi_{r}: \crmhs \rightarrow \cmhs $ qui est le foncteur oubli de la structure réelle sous-jacente.\\

On peut définir un foncteur de $\cmhs$ vers $\crmhs$ de la façon suivante (cf \cite{sim2}) :

Nous avons donc les foncteurs suivants :

$$\xymatrix{  \crzmhs \ar[r]^{\Phi_{z}}& \crmhs  \ar[r]^{\Phi_{r}} &\cmhs \ar@/^/@{.>}[l]^-{H \rightarrow H \oplus {\overline H}} }.$$

Les foncteurs $\Phi_{r}$ et ``$H \rightarrow H \oplus \overline H$'' ne sont pas adjoints. Si l'on part d'un objet de $\crmhs$ dont l'espace vectoriel sous-jacent est de rang $n$, alors l'espace vectoriel sous-jacent de l'image dans cette même catégorie par la composée de ces deux foncteurs est de rang $2n$. 

Par la suite, lorsque rien ne sera précisé, lorsqu'on parlera de structure de Hodge mixte on se placera dans la catégorie $\crmhs$. On pourra toujours penser en terme de structures de Hodge au sens usuel i.e. que ce sont des objets de $\crzmhs$ mais les informations que nous fournira l'étude suivante seront toutes contenues dans les images de ces objets par le foncteur oubli $\Phi_{z}$.    
%$ \Phi_{str}: \cmhs \rightarrow \trif $ qui à une structure de Hodge mixte $H=(H_{\bf C},F^{\bullet},{\hat F}^{\bullet },W_{\bullet})$ associe l'espace vectoriel complexe trifiltré muni de trois filtrations exhaustives et décroissantes $(H_{\bf C},F^{\bullet},{\hat F}^{\bullet },{W_{.}}^{\bullet})$.\\

\begin{definition}
 Soit $H$ une structure de Hodge mixte. Les nombres entiers $$h^{p,q}=dim_{\bf C}Gr_{F^{\bullet }}^{p}Gr_{m}^{W_\bullet }H$$ o\`u $p+q=m$ sont appelés les nombres de Hodge de la structure de Hodge mixte $H$.\\
\end{definition}
Ces nombres généralisent les nombres de Hodge dimension des sous-espaces $H^{p,q}$ où $p+q=n$ pour une structure de Hodge pure de poids $n$. Ils peuvent être non nuls en dehors de la droite d'équation $p+q=n$.\\
${}$\\
{\bf Remarque :} Par la propriété de conjugaison des filtrations $F^{\bullet}$ et ${\overline F}^{\bullet}$ associées à une structure de Hodge mixte dans $\crzmhs$ ou dans $\crmhs$ on voit que les nombres de Hodge des structures de Hodge au sens usuel ou des structures de hodge réelle sont symétriques i.e. pour tout $(p,q)$ on a $h^{p,q}=h^{q,p}$. Une telle propriété n'est pas vérifiée en général pour les structures de Hodge complexes puisque les filtrations de Hodge $F^{\bullet }$ et ${\hat F}^{\bullet }$ sont définies indépendement.\\  

\begin{definition}
Soit $H$ une structure de Hodge mixte donn\'ee on note ${\mathcal E }_{H}$ le sous-ensemble de ${\bf Z} \times {\bf Z} $ formé par les couples $(p,q) \in {\bf Z}\times {\bf Z}$ tels que les nombres de Hodge $h^{p,q}_{H}$ sont non nuls. Cet ensemble est appel\'e le type de la structure de Hodge $H$.\\
\end{definition}

\begin{definition}
La structure de Hodge $T \langle k \rangle$ est l'unique structure de Hodge de rang $1$, de type $(-k, -k)$ et de réseau entier $ (2 \pi i )^{k}{\bf Z}$.\\
\end{definition}
${}$\\
{\bf Remarque :} Soit $H$ une structure de Hodge mixte et $k \in {\bf Z}$, alors les types de $H$ et $H \otimes {T \langle k \rangle}$ sont liés par : ${\mathcal E }_{H \otimes {T \langle k \rangle}}= \{ (p-k,q-k)\in {\bf Z} \times {\bf Z} \vert (p,q) \in {\mathcal E }_{H} \} $.

Pour les filtrations on a les relations pour $p,q,n \in {\bf Z}^{3}$ :

\hspace{1cm}$\left\lbrace \begin{array}{l}
F^{p}(H \otimes {T \langle k \rangle})= F^{p+k}(H),\\   
{\overline F}^{q}(H \otimes {T \langle k \rangle})= {\overline F}^{q+k}(H),\\
W_{n}(H \otimes {T \langle k \rangle})= W_{n+2k}(H),    
\end{array}
\right.$ \hspace{0.5cm} et \hspace{0.5cm} $\left\lbrace \begin{array}{l}
h^{p,q}_{H \otimes {T \langle k \rangle}}=h^{p+k,q+k}_{H }.\\   
\end{array}
\right.$

%\hspace{1cm}$\left\lbrace \begin{array}{l}
%F^{p}(H \otimes {T \langle k \rangle})= F^{p+k}(H),\\   
%{\overline F}^{q}(H \otimes {T \langle k \rangle})= {\overline F}^{q+k}(H),\\
%W_{n}(H \otimes {T \langle k \rangle})= W_{n+2k}(H),    
%\end{array}
%\right.$ \hspace{0.5cm} et \hspace{0.5cm} $\left\lbrace \begin{array}{l}
%t^{p,q}_{H \otimes {T \langle k \rangle}}=t^{p+k,q+k}_{H},\\
%h^{p,q}_{H \otimes {T \langle k \rangle}}=h^{p+k,q+k}_{H }.\\   
%\end{array}
%\right.$

\begin{definition}
${\bf (i)}$ La longueur d'une structure de Hodge mixte est la longueur du plus grand intervalle $[a,b]$ tel que $W_{m}/W_{m-1}\neq0$ pour tout $m \in \{a,b\}$. $a$ et $b$ sont le plus bas et le plus haut poids et la longueur est égale à $b-a$. Une structure de longueur $0$ est une structure pure.\\
${\bf (ii)}$ Le niveau d'une struture de Hodge mixte est la longueur du plus grand intervalle $[c,d]$ tel que $F_{p}/F_{p-1}\neq0$ pour tout $p \in \{c,d\}$. Le niveau est égal à $b-a$. Par symétrie des structures de Hodge. Une structure de niveau $0$ est une tensorisation de structures de Hodge de Tate $T\langle -p \rangle^{k}$ pour $k \in {\bf Z}$.\\ 
\end{definition}

\subsubsection{Structures de Hodge mixtes ${\bf R}$-scindées}

Contrairement au cas des structures de Hodge pures, on ne pas en général trouver de scindement compatible à toutes les filtrations qui composent une struture de Hodge mixte $H$. Ceci est d\^u au fait remarqué dans la partie $1$ que contrairement aux cas des espaces bifiltrés, on ne peut pas en général trouver de scindement compatibles aux trois filtrations composant un espace vectoriel trifiltré.

Rappelons un lemme de Deligne, important dans l'étude des structure de Hodge mixtes. Bien que les structures de Hodge mixte n'admettent pas en général de décomposition en une somme directe de "$(p,q)$-sous-espaces" comme les structures pures, c'est-à-dire de graduation compatible aux trois filtrations, les espaces canoniquement définis ci-dessous ont des propriétés de décomposition trés utiles :
$$ I^{p,q}=(F^{p} \cap W_{p+q}) \cap ({\overline F}^{q} \cap W_{p+q}+\sum_{i \geq 1} \, {\overline F}^{q-i}\cap W_{p+q-i-1}).$$   
En effet,    
\begin{lemme}\cite{del2}\label{lesipq} 

\item{(i)} $I^{p,q}={\overline I}^{q,p} \ \ mod \ W_{p+q-2}. $

\item{(ii)} $W_{m}=\oplus_{p+q \leq m}I^{p,q}$.

\item{(iii)} $F^{p}=\oplus_{i \geq p} \oplus_{q} I^{i,q}$.

\item{(iv)} La projection $W_{m} \rightarrow Gr_{m}^{W}H$ induit un isomorphisme pour $p+q=m$ de $I^{p,q}$ vers le sous-espace de Hodge $(Gr^{W}_{m}H)^{p,q}$.\\ 
\end{lemme}

Quitte à changer les indices, cette décomposition canonique en sous-espaces $I^{p,q}$ donne deux bigraduations canoniques $H_{\bf C}=\oplus_{p,q}I^{p,q}$ associées aux paires de filtrations $(W_{\bullet },F^{\bullet })$ et $(W_{\bullet },{\overline F}^{\bullet })$. La première bigraduation, associée à $(W_{\bullet },F^{\bullet })$ est donnée directement par le lemme précédent, la deuxième associée à $(W_{\bullet },{\overline F}^{\bullet })$ est donnée par $W_{m}=\oplus_{p+q \leq m}{\overline I}^{p,q}$ et ${\overline F}^{q}=\oplus_{i \geq q} \oplus_{p} {\overline I}^{i,p}$. Ces deux bigraduations ne sont pas compatibles à moins que l'on ne soit dans le cadre de la définition suivante :
\begin{definition}
Une structure de Hodge mixte $H=(H_{\bf R},W_{\bullet },F^{\bullet },{\overline F}^{\bullet })$ est dite {\bf R}-scindée si les espaces $I^{p,q}$ définis de façon canonique vérifient la condition : pour tous $p,q$    $I^{p,q}={\overline I}^{q,p}$.
\end{definition}
Dans le cas où la structure de Hodge mixte $H=(H_{\bf R},W_{\bullet },F^{\bullet },{\overline F}^{\bullet })$ est scindée, on a : $F^{p}\cap {\overline F}^{q}=\oplus_{p'\geq p,q' \geq q}I^{p',q'}$. Ce qui signifie que les trois filtrations sont simultanéement scindées.\\
${}$\\
{\bf Remarque :} D'après le lemme \ref{lesipq} p.\pageref{lesipq}, $(i)$, les structures de Hodge mixtes de longueur inférieure à $1$ sont toutes ${\bf R}$-scindées. En particulier les structures de Hodge pures sont ${\bf R}$-scindées.

\subsection{Dictionnaire structures de Hodge-fibrés vectoriels}
Dans cette section nous allons appliquer les constructions établies dans la première partie au cas où les filtrations proviennent d'une structures de Hodge mixte.

\subsubsection{Construction dans le cas où les filtrations proviennent d'une structure de Hodge mixte}
Une structure de Hodge mixte est essentiellement la donnée de trois filtrations $(W_{\bullet }, F^{\bullet },{\overline F}^{\bullet})$ avec certaines relation d'incidence (elles sont opposées). Soit $H=(H_{\bf R},W_{\bullet},F^{\bullet},{\overline F}^{\bullet })$ une structure de Hodge mixte. Nous allons construire un fibré sur ${\bf P}^2$ associé à cette structure de Hodge mixte. Rappelons que dans la construction du fibré de Rees associé à trois filtrations nous partions de trois filtrations exhaustives et décroissantes. Pour ``caractériser'' la structure de Hodge mixte, nous prendrons donc ici non plus la filtration par le poids qui est croissante mais la filtration décroissante associé à la filtration par le poids. Rappelons que l'on définit la filtration décroissante notée $W^{\bullet}$ associée à la filtration croissante $W_{\bullet}$ de $H_{\bf R}$ par,  
$$\text{pour tout } p \in {\bf Z},\,\,\, W^{p}=W_{-p}.$$
Les trois filtrations décroissantes obtenues sont alors opposées. En effet les filtrations de départ forment une structure de Hodge mixte i.e. que $F^{\bullet }$ et ${\overline F}^{\bullet }$ sont $n$-opposée sur $Gr^{W}_{n}$ (dans \cite{del2} on dit que les filtrations de Hodge, sa filtration opposée et la filtration décalée $W[n]^{\bullet}$ sont opposées). On a donc
$$\text{pour tous } (p,q) \text{ tels que } p+q-n\neq 0,\,\,\, Gr_{F}^{p}Gr_{\overline F}^{q}Gr^{W_{\bullet }}_{n}H=0,$$   
et de façon équivalente
$$\text{pour tous } (p,q) \text{ tels que } p+q+n\neq 0,\,\,\, Gr_{F}^{p}Gr_{\overline F}^{q}Gr_{W^{\bullet}}^{n}H=0.$$

Ainsi à une structure de Hodge mixte $H=(H_{\bf R},W_{\bullet},F^{\bullet},{\overline F}^{\bullet })$ on associe un espace vectoriel complexe trifiltré $(H,W^{\bullet},F^{\bullet},{\overline F}^{\bullet })$ dont les filtrations sont opposées. 

\begin{definition}
Le fibré de Rees sur le plan projectif ${\bf P}^2$ associé à la structure de Hodge mixte $H=(H_{\bf R},W_{\bullet},F^{\bullet},{\overline F}^{\bullet })$ est le fibré de Rees construit à partir de l'espace vetoriel trifiltré $(H,W^{\bullet},F^{\bullet},{\overline F}^{\bullet })$.
$$ \xi_{{\bf P}^2}(H):=\xi_{{\bf P}^2}(H,W^{\bullet},F^{\bullet},{\overline F}^{\bullet }).$$  
\end{definition}

De l'étude faite dans la partie $1$ il vient : 
\begin{proposition}${}$\\ \label{gradloc}
${\bf (i)}$ $ \xi_{{\bf P}^2}(H)_{(1:1:1)}\cong H$.\\ 
${\bf (ii)}$ $ \xi_{{\bf P}^2}(H)_{(1:0:0)}\cong \oplus_{p,q}Gr_{F^{\bullet}}^{p}Gr_{{\overline F}^{\bullet}}^{q}H$.\\ 
${\bf (iii)}$ $ \xi_{{\bf P}^2}(H)_{(0:1:0)}\cong \oplus_{p,q}Gr_{F^{\bullet}}^{p}Gr_{W^{\bullet}}^{q}H\cong \oplus_{p,q}Gr_{F^{\bullet}}^{p}Gr_{W_{\bullet}}^{-q}H$.\\ 
${\bf (iv)}$ $ \xi_{{\bf P}^2}(H)_{(0:0:1)}\cong \oplus_{p,q}Gr_{{\overline F}^{\bullet}}^{p}Gr_{W^{\bullet}}^{q}H\cong \oplus_{p,q}Gr_{{\overline F}^{\bullet}}^{p}Gr_{W_{\bullet}}^{-q}H$.\\ 
${\bf (v)}$ $ \xi_{{\bf P}^2}(H)_{(1:1:0)}\cong\oplus_{q}Gr_{{\overline F}^{\bullet}}^{q}H$.\\ 
${\bf (vi)}$ $ \xi_{{\bf P}^2}(H)_{(1:0:1)}\cong \oplus_{p}Gr_{{\overline F}^{\bullet}}^{p}H$.\\ 
${\bf (vii)}$ $ \xi_{{\bf P}^2}(H)_{(0:1:1)}\cong \oplus_{n}Gr_{W^{\bullet}}^{n}H\cong \oplus_{n}Gr_{W_{\bullet}}^{-n}H$.\\ 
Ces isomorphismes donnent une correspondance entre les foncteurs fibres et les foncteurs gradués associés aux trois filtrations.

\end{proposition}

Le fibrés de Rees associé à la structure de Hodge mixte $H$ peut être vu comme une façon de transformer les espaces bigradués les uns en les autres. En effet les fibres canoniques aux trois points origines des cartes affines sont isomorphes aux espaces bigradu\'es associ\'es \`a chacune des paires de filtrations.\\ 
${}$\\
{\bf Exemples :} Donnons quelques exemples :\\
$\bullet$ Le fibré associé à une structure de Hodge de Tate $T \langle -k \rangle$ est le fibré trivial ${\bf T}^{\tau}$-équivariant $$\xi_{{\bf P}^2}^{-2k,k,k}$$ décrit dans la partie $1$.\\
$\bullet$ Le fibré associé à une structure de Hodge pure de poids $n$, $H=\oplus_{p,q,\,\,p+q=n}H^{p,q}$ est le fibré ${\bf T}^{\tau}$-équivariant $$\xi_{{\bf P}^2}(H)=\oplus_{p,q,\,\,p+q=n} \xi_{{\bf P}^2}^{-n,p,q}{}^{\text{dim}_{\bf C}H^{p,q}}.$$

\subsubsection{Dictionnaire structures de Hodge-fibrés sur ${\bf P}^2$ }

Comme conséquence directe des théorèmes \ref{th2} p.\pageref{th2} et \ref{th3} p.\pageref{th3}, il vient : 

\begin{theoreme}\label{equivcatmhs}
La construction du fibré de Rees sur ${\bf P}^2$ associé à un espace vectoriel trifiltré établit les équivalences de catégories entre :\\
${\bullet}$ La catégorie des structures de Hodge mixtes réelles $\crmhs$ et la catégorie des fibrés vectoriels ${\bf T}^{\tau}$-équivariants ${\bf P}^{1}_{0}$-semistables de pente $\mu=0$ munie des morphismes ${\bf T}^\tau$-équivariants dont les singularit\'es des conoyaux sont en codimension $2$, support\'ees au point $(1:0:0)$ : 
$$
\xymatrix{
\{ \crmhs \} \ar@<2pt>[r]^-{\Phi_{R}} &  \{  {\calf ib}_{{\bf P}^{1}_{0}-semistables,\mu=0}({\bf P}^{2}/{\bf T}^{\tau}) \} \ar@<2pt>[l]^-{\Phi_{I}}
}
.$$
${\bullet}$ La catégorie des structures de Hodge mixtes complexe $\cmhs$ et la catégorie des fibrés vectoriels ${\bf T}$-équivariants ${\bf P}^{1}_{0}$-semistables de pente $\mu=0$ munie des morphismes ${\bf T}$-équivariants dont les singularit\'es des conoyaux sont en codimension $2$, support\'ees au point $(1:0:0)$ : 
$$
\xymatrix{
\{ \cmhs \} \ar@<2pt>[r]^-{\Phi_{R}} &  \{  {\calf ib}_{{\bf P}^{1}_{0}-semistables,\mu=0}({\bf P}^{2}/{\bf T}) \} \ar@<2pt>[l]^-{\Phi_{I}}
}
.$$

\end{theoreme}

Ce théorème donne une démonstration géométrique des faits suivants (Théorème $(1.3.16)$, \cite{del2} pour $\crmhs$ et \cite{sim2} pour $\cmhs$) :

\begin{corollaire}\label{corab}${}$\\
$\bullet$ La catégorie des structures de Hodge mixtes réelles $\crmhs$ est abélienne.\\
$\bullet$ La catégorie des structures de Hodge mixtes complexes $\cmhs$ est abélienne.\\ 
$\bullet$ Les foncteurs ``oubli des filtrations'', $Gr_{W}$, $Gr_{F}$, $Gr_{\overline F}$ et $Gr_{W}Gr_{F}\cong Gr_{F}Gr_{W} \cong Gr_{W}Gr_{\overline F} \cong Gr_{\overline F}Gr_{W}$ de $\crmhs$ ou $\cmhs$ dans la catégorie des espaces vectoriels complexes sont exacts. 
\end{corollaire}

\begin{preuve}
Les deux premiers points sont des conséquences directes du fait que les catégories de fibrés vectoriels présentes dans les équivalences sont abéliennes comme établi dans la première partie.\\

Considérons une suite exacte de structures de Hodge mixtes 
$$\xymatrix{ 0 \ar[r]^{i} & A \ar[r] &  H  \ar[r]^{\pi} &  B   \ar[r] & 0}.$$
D'après l'étude faite dans la section ``suites exactes et $\mu$-semistabilité'' da la partie $1$, la suite exacte de faisceaux associée dans ${\calf ib}_{{\bf P}^{1}_{0}-semistables,\mu=0}({\bf P}^{2}/{\bf T}^{\tau}) $ est  
$$\xymatrix{ 0 \ar[r] & \xi_{{\bf P}^2}(A) \ar[r]^{i} & \xi_{{\bf P}^2}(H)  \ar[r]^{\pi^{**}} &  \xi_{{\bf P}^2}(B)   \ar[r] & 0}.$$
Elle provient de la suite exacte de faisceaux
$$\xymatrix{ 0 \ar[r] & \xi_{{\bf P}^2}(A) \ar[r]^{i} & \xi_{{\bf P}^2}(H)  \ar[r]^{\pi} &  \calc oker(i)   \ar[r] & 0},$$
qui fournit la suite exacte
$$\xymatrix{ 0 \ar[r] & \xi_{{\bf P}^2}(A) \ar[r]^{i} & \xi_{{\bf P}^2}(H)  \ar[r]^-{\pi^{**}} & \xi_{{\bf P}^2}(B)   \ar[r] & T \ar[r] & 0},$$
avec $\xi_{{\bf P}^2}(B)=\calc oker(i)^{**}$ et $\xymatrix{0 \ar[r] & \calc oker(i) \ar[r]^{\nu} & \xi_{{\bf P}^2}(B)  \ar[r] &  T   \ar[r] & 0}$.\\

Le support de $T$ est inclu dans le point $(1:0:0)$ donc pour tout $x \in {\bf P}^{2} \backslash \{(1:0:0)\}$ on a la suite exacte courte
$$\xymatrix{ 0 \ar[r] & \xi_{{\bf P}^2}(A)_{x} \ar[r]^{i} & \xi_{{\bf P}^2}(H)_{x}  \ar[r]^{\pi^{**}} &  \xi_{{\bf P}^2}(B)_{x}   \ar[r] & 0},$$
ce qui permet de conclure en choisissant $x$ de façon adéquate et en utilisant la proposition \ref{gradloc} p.\pageref{gradloc}.

\end{preuve}

{\bf Dictionnaire : Structures de Hodge mixtes-Fibrés équivariants.}
${}$\\
\begin{center}
%\begin{tabular}{|c||c|}
\begin{tabular}
{|>{\centering}m{2in}
||>{\centering}m{2in}|}
\hline

Catégorie de structures de Hodge mixtes & Catégorie de fibrés équivariants sur ${\bf P}^2$\\  %\hline 

\end{tabular}

\begin{tabular}
{|>{\centering}m{2in}
||>{\centering}m{2in}|}
\hline

 $\crmhs{}_{,scind} $ &  ${\calf ib}_{{{\bf P}^{1}_{0}}-semistables,\mu=0,scind}({\bf P}^{2}/{\bf T}^{\tau})$  \\

\end{tabular}

\begin{tabular}
{|>{\centering}m{2in}
||>{\centering}m{2in}|}
\hline

 $\crmhs $ &  ${\calf ib}_{{\bf P}^{1}_{0}-semistables,\mu=0}({\bf P}^{2}/{\bf T}^{\tau})$  \\

\end{tabular}

\begin{tabular}
{|>{\centering}m{2in}
||>{\centering}m{2in}|}
\hline

 $\cmhs{}_{,scind} $ &  ${\calf ib}_{{{\bf P}^{1}_{0}}-semistables,\mu=0,scind}({\bf P}^{2}/{\bf T}^{\tau})$  \\

\end{tabular}

\begin{tabular}
{|>{\centering}m{2in}
||>{\centering}m{2in}|}
\hline

$ \cmhs $ &  $  {\calf ib}_{{{\bf P}^{1}_{0}}-semistables,\mu=0}({\bf P}^{2}/{\bf T}) $ \\

\end{tabular}

\begin{tabular}
{>{\centering}m{2in}
>{\centering}m{2in}}
\hline
$\,$ & $\,$ \\
\end{tabular}

%\bottomcaption{Dictionnaire : Structures de Hodge mixtes-Fibrés équivariants.}
\end{center}

\subsubsection{Structure de Hodge mixte ind\'ecomposable}

 Soit $\xi_{{\bf P}^2}(H)$ le fibr\'e de Rees associ\'e \`a une structures de Hodge mixte $H$. On dit qu'un fib\'re \'equivariant $\cale$ sur ${\bf P}^2$ est ind\'ecomposable s'il n' est pas somme directe de deux sous-fibr\'es \'equivariants. On peut d\'ecomposer tout fibr\'e $G$-\'equivariant sur le plan projectif, o\`u $G={\bf T}$ ou ${\bf T}^{\tau}$, sur le plan projectif en somme directe de sous-fibr\'es $G$-\'equivariants
$$
\cale=\oplus_{i}\cale_{i}
,
$$    
o\`u les $\cale_i$ sont des sous-fibr\'es $G$-\'equivariants ind\'ecomposables. 

Notons par $H_{i}$ la structure de Hodge associ\'ee au fibr\'e $G$-\'equivariant $\cale_{i}$, on a alors la d\'ecomposition
$$
H=\oplus_{i}H_{i}.
$$
\begin{definition}Une structure de Hodge mixte est dite ind\'ecomposable si le fibr\'e de Rees \'equivariant qui lui est associ\'e est ind\'ecomposable.
 
La d\'ecomposition $H=\oplus_{i}H_{i}$ est appel\'e d\'ecomposition de la structure de Hodge mixte $H$ en structures de Hodge mixtes ind\'ecomposables.
\end{definition}

Remarquons la relation entre le type de $H$ et les types des termes ind\'ecomposables $H_i$
$$
\cale_{H}=\cup_{i}\cale_{H_i}.$$
\begin{lemme}
Soient $A$ et $B$ deux structures de Hodge mixtes s\'epar\'ees, de d\'ecomposition $A=\oplus_{i}A_{i}$ et $B=\oplus_{j}B_{j}$, alors
$$Ext_{MHS}^{1}(B,A)=\oplus_{i,j}Ext_{MHS}^{1}(B_{j},A_{i}).$$
\end{lemme}

%Nous verrons dans la section ? que l'introduction de la notion de structures de Hodge mixte ind\'ecomposable permet de comparer

\subsection{Invariants discrets de structures de Hodge mixtes, définition du niveau de ${\bf R}$-scindement}

\subsubsection{Classe d'une structure de Hodge mixte dans $K_{0}({\bf P}^{2},{\bf T})$}

Soit $H=(H_{\bf R},W_{\bullet },F^{\bullet },{\overline F}^{\bullet })$ une structure de Hodge mixte. On définit sa classe dans $K_{0}({\bf P}^{2},{\bf T})$ comme étant la classe du fibré équivariant associé
$$[H]=[\xi_{{\bf P}^2}(H)] \in K_{0}({\bf P}^{2},{\bf T}).$$ 
D'après l'étude du $K_{0}$ faite dans la partie $1$, la classe du fibré de Rees associé à trois filtrations est décrite par les entiers dimensions des gradués et bigradués associés aux filtrations. Dans le cas où les filtrations proviennent d'une structure de Hodge mixte ces entiers sont décrits par les nombres de Hodge, les $h^{p,q}$ définis pour tout $(p,q)$ par $h^{p,q}=\text{dim}_{\bf C}Gr_{F}^{p}Gr_{W}^{-p-q}H$ et les nombres entiers suivants
\begin{definition}
On définit les entiers "$s^{p,q}$" par
$$s^{p,q}=\text{dim}_{\bf C}Gr_{F}^{p}Gr_{\overline F}^{q}H.$$
\end{definition} 
Ces entiers vont jouer un rôle important par la suite pour mesurer la complexité d'une structure de Hodge mixte.

D'après la section ``Calcul explicite de $K_{0}({\bf P}^{2},{\bf T})$ et classe des fibré de Rees'' et la proposition \ref{ko} p.\pageref{ko} on a :
\begin{proposition}\label{kor} Soit $H=(H_{\bf R},W_{\bullet},F^{\bullet },{\overline F}^{\bullet })$ un objet de $\crmhs$, alors la classe de $\xi_{{\bf P}^2}(H)$ dans $K_{0}({\bf P}^{2},{\bf T})$ est donnée par
$$[\xi_{{\bf P}^2}(H)]_{K_{0}({\bf P}^{2},{\bf T})}=((P_{A_{0}},P_{A_{1}},P_{A_{2}}),(P_{{\bf G}_{m12}},P_{{\bf G}_{m02}},P_{{\bf G}_{m01}}),P_{{\bf G}_{m}^{2}}),$$
où,\\
$\hspace*{1cm}$ $\bullet$ $P_{A_{0}}=\sum_{p,q}h^{p,q}\,\,\,u_{0}^{p}v_{0}^{q}$,\\
$\hspace*{1cm}$ $\bullet$ $P_{A_{1}}=\sum_{p,q}h^{p,q}\,\,\,u_{1}^{p}v_{1}^{q}$,\\
$\hspace*{1cm}$ $\bullet$ $P_{A_{2}}=\sum_{p,q}s^{p,q}\,\,\,u_{2}^{p}v_{2}^{q}$,\\
$\hspace*{1cm}$ $\bullet$ $P_{{\bf G}_{m12}}=\sum_{n}(\sum_{p,q,p+q=n}h^{p,q})\,\,\,u_{12}^{n}$,\\
$\hspace*{1cm}$ $\bullet$ $P_{{\bf G}_{m02}}=\sum_{p}(\sum_{q}h^{p,q})\,\,\,u_{12}^{p}=\sum_{p}(\sum_{q}s^{p,q})\,\,\,u_{12}^{p}$,\\
$\hspace*{1cm}$ $\bullet$ $P_{{\bf G}_{m01}}=\sum_{q}(\sum_{p}h^{p,q})\,\,\,u_{12}^{q}=\sum_{q}(\sum_{p}s^{p,q})\,\,\,u_{12}^{q}$,\\
$\hspace*{1cm}$ $\bullet$ $P_{{\bf G}_{m}^{2}}=\emph{dim}_{\bf C}V=\sum_{p,q}h^{p,q}=\sum_{p,q}s^{p,q}$.\\ 
\end{proposition}

Soit $H=\oplus_{i} H_{i}$ la d\'ecomposition de la structure de Hodge mixte $H$ en facteurs ind\'ecomposables. On a alors, 
$$
[\xi_{{\bf P}^2}(H)]_{K_{0}({\bf P}^{2},{\bf T})}=\oplus_{i}[\xi_{{\bf P}^2}(H_{i})]_{K_{0}({\bf P}^{2},{\bf T})},$$
et donc 
$$s^{p,q}_{H}=\sum_{i} s^{p,q}_{H_i}.$$

\subsubsection{Définition du niveau de ${\bf R}$-scindement des structures de Hodge mixtes}

Dans cette sous-section, comme partout où l'on parlera de la notion de niveau de ${\bf R}$-scindement d'une structure de Hodge mixte, une structure de Hodge mixte sera une structure de Hodge mixte réelle i.e. un objet de $\crmhs$.
%A une structure de Hodge mixte $H=(H_{\bf Q},W_{\bullet },F^{\bullet },{\overline F}^{\bullet })$ on peut associer les nombres entiers suivants :
%$$\left\lbrace \begin{array}{l}
%         h^{p,q}_{H}=dim_{\bf C}Gr_{{\overline F}^{\bullet }}^{q}Gr_{F^{\bullet }}^{p}Gr_{W_{\bullet }}^{p+q}H_{\bf C} \ \mathrm{pour \ tous}\ p,q \in {\bf Z}, \\
%         s^{p,q}_{H}=dim_{\bf C}Gr_{{\overline F}^{\bullet }}^{q}Gr_{F^{\bullet }}^{p}H_{\bf C}\ \mathrm{ pour\ tous}\ p,q \in {\bf Z}.
         
%\end{array}
%\right.$$ 
Les entiers $h^{p,q}_{H}$ sont les nombres de Hodge classiques associés aux structures de Hodge mixtes. Dans le cas où la structure de Hodge mixte $H$ est ${\bf R}$-scindée, pour tous $p,q \in {\bf Z}$ on a l'égalité : 
$$h^{p,q}_{H}=s^{p,q}_{H},$$
en effet $ h^{p,q}_{H}=\text{dim}_{\bf C}Gr_{{\overline F}^{\bullet }}^{q}Gr_{F^{\bullet }}^{p}Gr_{W_{\bullet }}^{p+q}H_{\bf C}=\text{dim}_{\bf C}I^{p,q}=\text{dim}_{\bf C}Gr_{{\overline F}^{\bullet }}^{q}Gr_{F^{\bullet }}^{p}H_{\bf C}=s^{p,q}_{H}$.

 Ce n'est pas vrai en général comme nous le verrons par la suite. Le calcul des invariants du fibré de Rees associé à une structure de Hodge mixte construit va nous permettra de voir à quel point les entiers $s^{p,q}_{H}$ sont proches des nombres de Hodge $h^{p,q}_{H}$, c'est à dire de voir si la structure de Hodge mixte est proche de la structure de Hodge mixte ${\bf R}$-scindée qui lui est associée \cite{catkapsch}.\\

\begin{proposition}
Soit ${\xi}_{{\bf P}^2}(H_{\bf C},W_{.}^{\bullet },F^{\bullet },{\overline F}^{\bullet })$ le fibré sur ${\bf P}^2$ associé à la structure de Hodge mixte $(H_{\bf Q},W_{\bullet },F^{\bullet},{\overline F}^{\bullet})$, alors :
$$ \emph{ch}({\xi}_{{\bf P}^2}(H_{\bf C},W_{.}^{\bullet },F^{\bullet },{\overline F}^{\bullet }))= \emph{dim}_{\bf C}H_{\bf C}+\frac{1}{2}\sum_{p,q}(p+q)^{2}(s^{p,q}_{H}-h^{p,q}_{H})w^{4}$$
et donc :
$$\emph{c}_{2}( {\xi}_{{\bf P}^2}(H_{\bf C},W_{.}^{\bullet },F^{\bullet },{\overline F}^{\bullet }))=\frac{1}{2}\sum_{p,q}(p+q)^{2}(h^{p,q}_{H}-s^{p,q}_{H})w^{4}.$$

\end{proposition}

\begin{preuve}
Nous pouvons donner deux preuves de cette proposition à l'aide la première partie. Pour la première, on applique la proposition \ref{filtrscin} p.\pageref{filtrscin} dans le cas où l'espace vectoriel trifiltré provient d'une structure de Hodge mixte avec $(V,F^{\bullet }_{0},F^{\bullet }_{1},F^{\bullet }_{2})=(H_{\bf C},W_{.}^{\bullet },F^{\bullet },{\overline F}^{\bullet })$. D'où  $  \text{dim}_{\bf C}Gr_{F_{2}^{\bullet }}^{q}Gr_{F_{1}^{\bullet }}^{p}V= \text{dim}_{\bf C}Gr_{{\overline F}_{}^{\bullet }}^{q}Gr_{{ F}_{}^{\bullet }}^{p}H_{\bf C}=s^{p,q}_{H}$ et $ \text{dim}_{\bf C}Gr_{F_{2}^{\bullet }}^{q}Gr_{F_{1}^{\bullet }}^{p}Gr_{F_{0}^{\bullet }}^{-p-q}V = \text{dim}_{\bf C}Gr_{{\overline F}^{\bullet }}^{q}Gr_{F_{}^{\bullet }}^{p}Gr^{W^{\bullet }_{.}}_{-p-q}H_{\bf C}=h_{H}^{p,q}$ ce qui donne la relation.

La deuxième preuve consiste à déduire le caractère de Chern du fibré de sa classe dans $K_{0}({\bf P}^{2},{\bf T})$ en oubliant l'action.
\end{preuve}
Ceci nous amène à définir le niveau de ${\bf R}$-scindement comme il suit :   
\begin{definition}
Le niveau de ${\bf R}$-scindement d'une structure de Hodge mixte $(H_{\bf Q},W_{\bullet },F^{\bullet},{\overline F}^{\bullet})$ est le nombre entier : $$\alpha (H)=\emph{c}_{2}({\xi}_{{\bf P}^2}(H))=-\emph{ch}_{2}({\xi}_{{\bf P}^2}(H_{\bf C},W_{.}^{\bullet },F^{\bullet },{\overline F}^{\bullet }))=\frac{1}{2}\sum_{p,q}(p+q)^{2}(h^{p,q}_{H}-s^{p,q}_{H}).$$
\end{definition}
Cette définition constitue une généralisation de la notion de structure de Hodge ${\bf R}$-scindée.

\begin{proposition} Soit $H$ une structure de Hodge mixte réelle, alors
$$ \alpha(H)=0    \Longleftrightarrow  \text{ La structure de Hodge mixte } H \text{ est }  {\bf R}\text{-scind\'ee. }$$
\end{proposition}
${}$\\
\begin{preuve}
%L'implication $\Leftarrow$ est \'evidente car les $I^{p,q}$ sindent les trois filtrations simultan\'eement.
Si $(H_{\bf Q},W_{\bullet },F^{\bullet},{\overline F}^{\bullet})$ est une structure de Hodge ${\bf R}$-scindée alors $\alpha (H)=0$. Reciproquement supposons que l'invariant $\alpha $ soit \'egal \`a z\'ero, le fibr\'e de Rees associ\'es aux filtration est alors trivial, ce qui signifie que les scindements sont les m\^emes sur chacunes des cartes affines. Le scindement sur la carte de $(W_{\bullet},F^{\bullet})$ par les $I^{p,q}$ est donc le m\^eme que le scindement par les $S^{p,q}=Gr_{F^\bullet}^{p}Gr_{{\overline F}^{\bullet }}^{q}H_{\bf C}$, scindement sur la carte de $(F^{\bullet },{\overline F}^{\bullet})$. Les $S^{p,q}$ \'etant envoy\'es sur les $S^{q,p}$, on a $I^{p,q}={\overline I}^{q,p}$ pour tout $(p,q)$, ce qui prouve que la structure est ${\bf R}$-scind\'ee. 

La réciproque peut aussi se voir de la façon qui suit : si $\alpha(H)=0$, alors les classes de Chern du fibré sur le plan projectif complexe sont nulles et par la théorie de Donaldson ce fibré est donc trivial et donc la structure de Hodge mixte à laquelle il est associé est scindée. 
\end{preuve}
${}$\\
{\bf Remarques :} $\bullet$ Notons $H'=(H_{\bf C},W_{\bullet },e^{-i.\delta}.F^{\bullet })$ la structure de Hodge mixte ${\bf R}$-scindée associée à la structure de Hodge mixte $H=(H_{\bf C},W_{\bullet },F^{\bullet })$. On vérifie que l'on a bien pour tous $(p,q)$ $h^{p,q}_{H'}=s^{p,q}_{H'}=dim_{\bf C}\,I^{p,q}$ et donc $\alpha(H_{\bf C},W_{\bullet },e^{-i.\delta}.F^{\bullet })=0$.\\
$\bullet$ Pour des raisons de dimension, on a $\sum_{p,q}(h^{p,q}_{H}-s^{p,q}_{H})=0$.\\ 
$\bullet$ De même
\begin{eqnarray*}
\sum_{p,q}(p+q)(h^{p,q}_{H}-s^{p,q}_{H})&=&\sum_{p}p(\sum_{q}\,h^{p,q}_{H}-s^{p,q}_{H})+\sum_{q}q(\sum_{p}\,h^{p,q}_{H}-s^{p,q}_{H})\\
&=&\sum_{p}p(\text{dim}_{\bf C}{\overline F}^{q}-\text{dim}_{\bf C}{\overline F}^{q})+\sum_{q}q(\text{dim}_{\bf C}{ F}^{p}-\text{dim}_{\bf C}{ F}^{p})\\
&=&0.
\end{eqnarray*}\\
$\bullet $ Il para\^{\i}t possible que deux structures de Hodge mixtes $H$ et $H'$ aient les m\^emes nombres de Hodge, des niveaux de ${\bf R}$-scindement $\alpha (H)$ et $\alpha (H')$ égaux mais que leurs nombres $h^{p,q}_H$ et $s^{p,q}_H$ ne soient pas tous égaux.\\
$\bullet$ On peut d\'emontrer directement par des consid\'erations d'alg\`ebre lin\'eaire la deuxi\`eme implication d\'emontr\'ee dans la proposition pr\'ec\'edente. Supposons que $\alpha(H)=0$ pour une structure de Hodge mixte $H$. On veut montrer que dans ce cas pour tous $(p,q)$ $I^{p,q}={\overline I}^{q,p}$. Supposons que ce ne soit pas le cas. Munissons ${\bf Z}^2$ de l'ordre lexicographique. Soit $(p_{0},q_{0})$ le plus grand des éléments $(p,q)$ tels que $I^{p,q}\neq I^{q,p}$ (rappelons que l'on a toujours égalité modulo $W_{p+q-2}$). Soit $I^{p_{0},q_{0}}_{0}$ le sous-espace vectoriel de $I^{p_{0},q_{0}}$ de dimension maximale vérifiant l'égalité $I^{p_{0},q_{0}}_{0}={\overline I}^{q_{0},p_{0}}_{0}$. Il est donné par $I^{p_{0},q_{0}}_{0}=I^{p_{0},q_{0}}\cap {\overline I}^{p_{0},q_{0}}$. On a par hypothèse $I^{p_{0},q_{0}}_{0} \subsetneq I^{p,q}$. Ainsi, comme pour tous $(p,q) \geq (p_{0},q_{0})$ l'égalité $h^{p,q}=s^{p,q}$ est vérifiée, on en déduit l'inégalité $s^{p_{0},q_{0}}_{H}=\text{dim}_{{\bf C}}\,I^{p_{0},q_{0}}_{0} < \text{dim}_{\bf C}\,I^{p_{0},q_{0}}=h^{p_{0},q_{0}}_{H}$ ce qui contredit l'hypothèse car le terme correspondant dans la formule de $\alpha$ contribue strictement positivement. 
 
\subsection{Comportement du niveau de ${\bf R}$-scindement par opérations sur les structures de Hodge mixtes}

Dans cette partie nous allons étudier le comportement de $\alpha$ lors des différentes opérations qui peuvent \^etre faites sur des structures de Hodge mixtes. Ces opération sont héritées de la structure de catégorie abélienne et tensorielle de la catégorie des structures de Hodge mixtes. %Nous allons étudier  Morphismes entre structures de Hodge mixtes, quotients, sommes directes, produits tensoriels et extensions.\\

\begin{theoreme}\label{elemalpha}
Pour $H$ et $H'$ deux structures de Hodge mixtes :\\
\hspace*{1cm}$(i)$ pour tout $k \in {\bf Z}$, $\alpha(H \otimes T \langle k \rangle )=\alpha(H)$.\\ 
\hspace*{1cm}$(ii)$ $\alpha(H^{*})=\alpha(H)$ où $H^{*}=\emph{Hom}_{SHM}(H,T \langle 0 \rangle )$.\\
\hspace*{1cm}$(iii)$ $\alpha ( H \oplus H')=\alpha(H) +\alpha (H')$.\\
 \hspace*{1cm}$(iv)$ $\alpha( H \otimes H')=\emph{dim}(H').\alpha(H)+\emph{dim}(H).\alpha(H')$. 
\end{theoreme}

\begin{preuve}(i)\begin{eqnarray*}\alpha(H \otimes T \langle k \rangle )&=&\frac{1}{2}\sum_{(p,q)\in {\cale}_{H \otimes T \langle k \rangle}}(p+q)^{2}(h^{p,q}_{H \otimes T \langle k \rangle}-s^{p,q}_{H \otimes T \langle k \rangle})\\
&=&\frac{1}{2}\sum_{(p+k,q+k)\in {\cale}_{H}}(p+q)^{2}(h^{p,q}_{H \otimes T \langle k \rangle}-s^{p,q}_{H \otimes T \langle k \rangle})\\
&=&\frac{1}{2}\sum_{(p,q)\in {\cale}_{H}}(p+q)^{2}(h^{p-k,q-k}_{H \otimes T \langle k \rangle}-s^{p-k,q-k}_{H \otimes T \langle k \rangle})\\
&=&\alpha(H),
\end{eqnarray*}
d'après la remarque précédente.\\
(ii) Il suffit d'écrire que ${\cale }_{H^*}=\{(-p,-q) \vert (p,q) \in {\cale}_{ H}$, $s^{p,q}_{H^*}=s^{-p,-q}_{H}$ et $h^{p,q}_{H^*}=h^{-p,-q}_{H}\}$, l'égalité en découle car tous les coefficients $p+q$ sont au carré.\\
(iii) ${\cale }_{H\oplus H'}={\cale}_{H} \cup {\cale}_{ H'}$. Les filtrations de Hodge et la filtration décroissante associée à la filtration par le poids de $H \oplus H'$ se déduisent des filtrations respectives de $H$ et $H'$ par somme directe. Les dimensions des quotients sont donc faciles à calculer, on en déduit : pour $p,q$ entiers, $h^{p,q}_{H \oplus H'}=h^{p,q}_{H}+h^{p,q}_{H'}$ et $s^{p,q}_{H \oplus H'}=s^{p,q}_{H}+s^{p,q}_{H'}$. En découle la formule voulue.\\
(iv) D'après le lemme \ref{sumfibrees} :\\
 ${\xi}_{{\bf P}^2}((H_{\bf C},W^{\bullet }_{.},F^{\bullet},{\overline F}^{\bullet}) \otimes (H_{\bf C} {}',W^{\bullet }_{.} {}',F^{\bullet} {}',{\overline F}^{\bullet} {} ')) \cong {\xi}_{{\bf P}^2}((H_{\bf C} {},W^{\bullet }_{.} {},F^{\bullet} {},{\overline F}^{\bullet} {} )) \otimes {\xi}_{{\bf P}^2}((H_{\bf C} {}',W^{\bullet }_{.} {}',F^{\bullet} {}',{\overline F}^{\bullet} {} ')) $.\\
Donc :

$\text{ch}{\xi}_{{\bf P}^2}((H_{\bf C},W^{\bullet }_{.},F^{\bullet},{\overline F}^{\bullet}) \otimes (H_{\bf C} {}',W^{\bullet }_{.} {}',F^{\bullet} {}',{\overline F}^{\bullet} {} '))=$

\hspace{4cm} $\text{ch}({\xi}_{{\bf P}^2}((H_{\bf C} {},W^{\bullet }_{.} {},F^{\bullet} {},{\overline F}^{\bullet} {} ))).\text{ch}({\xi}_{{\bf P}^2}((H_{\bf C} {}',W^{\bullet }_{.} {}',F^{\bullet} {}',{\overline F}^{\bullet} {} '))) $.\\
Or les trois filtrations constituant une structure de Hodge étant mixte étant opposées, on a : 
$\text{ch}_{1}({\xi}_{{\bf P}^2}((H_{\bf C},W^{\bullet }_{.},F^{\bullet},{\overline F}^{\bullet}))=0$ et de m\^eme $\text{ch}_{1}({\xi}_{{\bf P}^2}((H_{\bf C}{} ',W^{\bullet }_{.}{} ',F^{\bullet} {}',{\overline F}^{\bullet}{} '))=0$. Ainsi :

$\text{c}_{2}({\xi}_{{\bf P}^2}((H_{\bf C},W^{\bullet }_{.},F^{\bullet},{\overline F}^{\bullet}) \otimes (H_{\bf C} {}',W^{\bullet }_{.} {}',F^{\bullet} {}',{\overline F}^{\bullet} {} ')))=$

\hspace{1cm} $\text{dim}_{\bf C}(H_{\bf C}{}').\text{c}_{2}({\xi}_{{\bf P}^2}((H_{\bf C} {},W^{\bullet }_{.} {},F^{\bullet} {},{\overline F}^{\bullet} {} )))+\text{dim}_{\bf C}(H_{\bf C}).\text{c}_{2}({\xi}_{{\bf P}^2}((H_{\bf C} {}',W^{\bullet }_{.} {}',F^{\bullet} {}',{\overline F}^{\bullet} {} '))) $.\\
Qui est l'égalité cherchée.

\end{preuve}
On peut donc, d'apr\`es le $(iii)$ du th\'eor\`eme, calculer $\alpha$ facteurs ind\'ecomposables par facteurs ind\'ecomposables. 
\begin{corollaire} Soit $H$ un structure de Hodge mixte et $H=\oplus_{i}H_{i}$ sa d\'ecomposition en facteurs ind\'ecomposables, on a alors
$$
\alpha(H)=\sum_{i}\,\,\alpha(H_{i}).
$$
\end{corollaire}

Soient $A=(A_{\bf Z},W_{\bullet }^{A},F^{\bullet}_{A},{\overline F}^{\bullet}_{A})$ et $B=(B_{\bf Z},W_{\bullet }^{B},F^{\bullet}_{B},{\overline F}^{\bullet}_{B})$ deux structures de Hodge mixtes. Nous pouvons alors définir  suivant \cite{mor} ou \cite{car} le groupe d'extension de $B$ par $A$ dans la cat\'egorie des structures de Hodge mixtes, noté $\text{Ext}^{1}_{MHS}(B,A)$.
%séparées c'est à dire telles que le plus haut poids de $A$ est strictement inferieur au poids le plus bas de $A$, 
C'est à dire l'ensemble des classes d'équivalence à congruence prés des suites exactes de structures de Hodge mixtes :
%$$ 0 \rightarrow A \rightarrow H \rightarrow B \rightarrow 0$$
$$
\xymatrix{
 0 \ar[r] &  A        \ar[r]^{i}  & H         \ar[r]^{\pi}     &   B    \ar[r]   & 0  
}
$$

Pour une structure de Hodge mixte $H$ remplissant une telle suite exacte, on écrira : $H \in  \text{Ext}^{1}_{MHS}(B,A)$. On a bien s\^ur :
 ${\mathcal E }_{H}={\mathcal E }_{A} \cup {\mathcal E }_{B}$.
Nous noterons pour tous $p,q \in {\bf Z}$ $s^{p,q}_{H}$ (resp. $s^{p,q}_{A}$, $s^{p,q}_{B}$) les entiers suivants $\text{dim}_{\bf C}Gr_{{\overline F}_{H}^{\bullet }}^{q}Gr_{{ F}_{H}^{\bullet }}^{p}H_{\bf C}$ (resp. $\text{dim}_{\bf C}Gr_{{\overline F}_{A}^{\bullet }}^{q}Gr_{F_{A}^{\bullet }}^{p}A_{\bf C}$, $\text{dim}_{\bf C}Gr_{{\overline F}_{B}^{\bullet }}^{q}Gr_{F_{B}^{\bullet }}^{p}B_{\bf C}$). On utilisera le m\^eme type de  notation pour les nombres de Hodge $h^{p,q}_{H}$,$h^{p,q}_{A}$ et $h^{p,q}_{B}$.\\
Pour $H \in \text{Ext}^{1}_{MHS}(B,A)$ nous avons les égalités suivantes pour les nombres de Hodge :
\begin{center}
Pour tous $p,q \in {\bf Z}$, $h^{p,q}_{H}=h^{p,q}_{A}+h^{p,q}_{B}$
\end{center}    
Une telle égalité pour les entiers $s^{p,q}$ n'est pas vraie en général ce qui met en évidence la sur-additivité de l'invariant $\alpha$ :

\begin{theoreme}\label{inegext}
Soient $A$ et $B$ deux structures de Hodge mixtes et $H \in \emph{Ext}^{1}_{MHS}(B,A)$, alors :
$$ \alpha(H) \geq \alpha(A)+\alpha(B)$$

\end{theoreme} 

\begin{preuve}
Considérons la suite exacte de structures de Hodge mixtes qui donne l'extension
$$\xymatrix{ 0 \ar[r] & A \ar[r]^{i} & H \ar[r]^{\pi} &  B  \ar[r] & 0}.$$
On en déduit, comme dans la preuve du corollaire \ref{corab} p.\pageref{corab} la suite exacte de faisceaux sur ${\bf P}^2$
$$\xymatrix{ 0 \ar[r] & \xi_{{\bf P}^2}(A) \ar[r]^{i} & \xi_{{\bf P}^2}(H) \ar[r]^{\pi^{**}} &  \xi_{{\bf P}^2}(B)  \ar[r] & T \ar[r]&  0},$$
où $T$ est de torsion. D 'après le lemme 30, $\text{ch}_{2}(T) \geq 0$. Or $\text{ch}_{2}(T)=-\text{c}_{2}(T)$ car le faisceau $T$ est supporté en codimension $2$ donc $\text{c}_{1}(T)=0$. L'inégalité se déduit donc de la relation
$\text{c}_{2}(\xi_{{\bf P}^2}(A))+\text{c}_{2}(\xi_{{\bf P}^2}(B))=\text{c}_{2}(\xi_{{\bf P}^2}(H))+\text{c}_{2}(T)$. 
\end{preuve}

\begin{corollaire}
Soit $H$ une structure de Hodge mixte, alors :    \,\,\,\,       $\alpha(H) \geq 0$.\\

\end{corollaire}

\begin{preuve}
Toute structure de Hodge mixte peut \^etre décrite comme extension successive de structures de Hodge pures. En effet, si $H$ est une structure de Hodge dont les poids varient entre $0$ et $2n$ par exemple (on peut s'y ramener en tensorisant par une structure de Hodge de Tate, ce qui ne change pas la valeur de $\alpha$ d'après la proposition \ref{elemalpha} p.\pageref{elemalpha}), on a $H \in \text{Ext}_{MHS}^{1}(Gr^{W}_{2n}H,W_{2n-1}H)$. Ainsi par la proposition précédente $\alpha(H) \geq \alpha(W_{2n-1}H) + \alpha(Gr^{W}_{2n}H)=\alpha(W_{2n-1}H)$ car $\alpha(Gr^{W}_{2n}H)=0$, $Gr^{W}_{2n}H$ étant une structure de Hodge pure de poids $2n$. De m\^eme $W_{2n-1}H \in \text{Ext}_{MHS}^{1}(Gr^{W}_{2n-1}H,W_{2n-2}H)$ et donc $\alpha(H) \geq \alpha(W_{2n-1}H) \geq \alpha(W_{2n-2}H)$. En poursuivant ainsi :
$$  \alpha(H) = \alpha(W_{2n}H) \geq \alpha(W_{2n-1}H) \geq \alpha(W_{2n-2}H) \geq ...\geq \alpha(W_{0}H)=0,$$
car $W_{0}H$ est une structure de Hodge pure.
\end{preuve}

{\bf Remarques :} $\bullet$ Avec Drézet et de Le Potier, \cite{drepot}, on voit que les fibrés vectoriels $\xi_{{\bf P}^2}(H^{k}(X,{\bf C}),W^{\bullet }_{.},F^{\bullet },{\overline F}^{\bullet })$ sont ``connus''. En effet, on est dans le cas où $r>1$,$c_1=0$ et $c_2 \leq 0$.\\ 
$\bullet$ Soient $A$ et $B$ deux structures de Hodge mixtes, il existe  $m(A,B) \in {\bf Z}$ ($m=m(h_{A}^{.,.},h_{B}^{.,.},t_{A}^{.,.},t_{B}^{.,.}$)) tel que pour tout $H \in \text{Ext}^{1}_{MHS}(B,A)$ : $\alpha(H) \in [\alpha(A)+\alpha(B),m(A,B)]$. On a m\^eme plus : pour tout $p \in [\alpha(A)+\alpha(B),m(A,B)]$ il existe $H \in \text{Ext}^{1}_{MHS}(B,A)$ tel que $\alpha(H)=p$. On en déduit qu'étant donnés des nombres de Hodge $h^{p,q}$, il existe $m(h^{.,.})$ tel que si $H$ est une structure de Hodge mixte qui a ces nombres de Hodge alors $\alpha(H) \in [0,m(h^{.,.})]$.\\
${}$\\

La suite de cette section est constituée une autre démonstration du théorème sur la sur-additivité de $\alpha$. Cette preuve est une preuve qui utilise des outils d'algèbre linéaire alors que la preuve précédente est géométrique. Elle fournit des expressions de $\alpha$ qui par leurs formes permettront d'étudier les propriétés de semi-continuité de cet invariant.

Considérons la suite exacte de structures de Hodge mixtes donnée par la classe d'extension $H \in \text{Ext}^{1}_{MHS}(B,A)$. D'après le travail de la section $1.6.2$ de la partie $1$, on peut écrire le diagramme

%$$
%\xymatrix{
%    {} & A          \ar[d]_{\alpha} & H         \ar[r]^{\pi} \ar[d]_{\nu}     &   B    \ar[r] \ar[d]_{\beta}  & 0  \\
%          {}  & A'         \ar[r]^{i'} &  H'         \ar[r]^{\pi'}  &      B     &   {}
%}
%$$

%$$
%\xymatrix{
%{}  &    0   \ar[d] &   0   \ar[d] & 0  \ar[d] &  {}\\
% 0  \ar[r] & A        \ar[r]  \ar[d]_{\alpha} & H         \ar[r]^{\pi} \ar[d]_{\nu}     &   B    \ar[r] \ar[d]_{\beta}  & 0  \\
%          {}  & A'   \ar[d]      &  H'    \ar[d]      &      B  \ar[d]   &   {}\\
%{}  & A'   \ar[d]      &  H'    \ar[d]      &      B  \ar[d]   &   {}\\
%{}  & 0         &  0          &      0   &   {}
%}
%$$

\scalebox{0.70}[1]{
$$
\xymatrix{
{}  &    0   \ar[d] &   0   \ar[d] & 0  \ar[d] &  {}\\
 0  \ar[r] &  {\xi}_{\widetilde{{\bf P}^2}}(A_{\bf C},{W^{A}_{. }}^{\bullet },F_{A}^{\bullet },{\overline F}_{A}^{\bullet },Triv^{\bullet })       \ar[r]  \ar[d] & {\xi}_{\widetilde{{\bf P}^2}}(H_{\bf C},{W^{H}_{.}}^{\bullet },F_{H}^{\bullet },{\overline F}_{H}^{\bullet },Triv^{\bullet })         \ar[r] \ar[d]     &   {\xi}_{\widetilde{{\bf P}^2}}(B_{\bf C},{W^{B}_{.}}^{\bullet },F_{B}^{\bullet },{\overline F}_{B}^{\bullet },Triv^{\bullet })    \ar[r] \ar[d]  & 0  \\
          {}  &  e^{*}{\xi}_{{{\bf P}^2}}(A_{\bf C},{W^{A}_{.}}^{\bullet },F_{A}^{\bullet },{\overline F}_{A}^{\bullet })  \ar[d]      &   e^{*}{\xi}_{{{\bf P}^2}}(H_{\bf C},{W^{H}_{.}}^{\bullet },F_{H}^{\bullet },{\overline F}_{H}^{\bullet })    \ar[d]      &  e^{*}{\xi}_{{{\bf P}^2}}(B_{\bf C},{W^{B}_{.}}^{\bullet },F_{B}^{\bullet },{\overline F}_{B}^{\bullet })      \ar[d]   &   {}\\
{}  &  {\calf}_{A}   \ar[d]      &   {\calf}_{H}    \ar[d]      &       {\calf}_{B}  \ar[d]   &   {}\\
{}  & 0         &  0          &      0   &   {}
}
$$
}

%$$
%\xymatrix{
   
%0 \ar[r]   &   A    \ar[r] \ar[d]    &   H    \ar[r] \ar[d]  &   B   \ar[r]  \ar[d]
%&     0    \\

 %0 \ar[r] &   {\xi}_{\widetilde{{\bf P}^2}}(A_{\bf C},W^{A}_{\bullet },F_{A}^{\bullet },{\overline F}_{A}^{\bullet },Triv^{\bullet })        \ar[r]  &   {\xi}_{\widetilde{{\bf P}^2}}(H_{\bf C},W^{H}_{\bullet },F_{H}^{\bullet },{\overline F}_{H}^{\bullet },Triv^{\bullet })       \ar[r]      &   {\xi}_{\widetilde{{\bf P}^2}}(B_{\bf C},W^{B}_{\bullet },F_{B}^{\bullet },{\overline F}_{B}^{\bullet },Triv^{\bullet })    \ar[r]   & 0  

%}
%$$

%           {} & e^{*}{\xi}_{{{\bf P}^2}}(A_{\bf C},W^{A}_{\bullet },F_{A}^{\bullet },{\overline F}_{A}^{\bullet })         \ar[d]  &  e^{*}{\xi}_{{{\bf P}^2}}(H_{\bf C},W^{H}_{\bullet },F_{H}^{\bullet },{\overline F}_{H}^{\bullet })         \ar[d]  &      e^{*}{\xi}_{{{\bf P}^2}}(B_{\bf C},W^{B}_{\bullet },F_{B}^{\bullet },{\overline F}_{B}^{\bullet })  \ar[d]   &   {}\\

%{}  & {\calf}_{A}  \ar[d]  &  {\calf}_{H}  \ar[d]    & {\calf}_{B}  \ar[d]  &   {}\\

% {}    &  0   \ar[d]     &   0   \ar[d]  &  0  \ar[d] &   {}   

%\begin{center}
%$0 \rightarrow {\xi}_{\widetilde{{\bf P}^2}}(V,F_{0}^{\bullet },F_{1}^{\bullet },F_{2}^{\bullet },Triv^{\bullet }) \rightarrow e^{*}{\xi}_{{\bf P}^2}(V,F_{0}^{\bullet },F_{1}^{\bullet },F_{2}^{\bullet }) \rightarrow \calf \rightarrow 0$.
%\end{center}

où, d'après la section précédente, les colonnes sont des suites exactes associées à chacunes des structures de Hodge mixtes $A$, $H$ et $B$. La première ligne est une suite exacte car le foncteur ${\xi}_{\widetilde{{\bf P}^2}}(.,.,.,.,Triv^{\bullet })$ est un foncteur exact de la catégorie des structures de Hodge mixtes vers la catégorie des fibrés vectoriels sur ${\widetilde{{\bf P}^2}}$. On obtient donc la relation entre les niveaux de ${\bf R}$-scindement des structures de Hodge mixtes :
$$  \alpha (H)= \alpha (A) + \alpha (B)+ \frac{1}{2}\sum_{p,q}\,(p+q)^{2}(s^{p,q}_{A}+s^{p,q}_{B}-s^{p,q}_{H}).$$
Cette relation est importante dans cette autre démonstration du fait que l'invariant $\alpha$ est sur-additif par extension :\\
${}$\\
{\bf Notation :} pour une structure de Hodge mixte $\bullet $, on définit
 $$f^{p,q}_{\bullet}=\text{dim}_{\bf C}F^{p}_{\bullet}\cap{\overline F}^{q}_{\bullet}.$$

La démonstration du théorème utilise le lemme suivant ainsi que son corollaire :\\
  
\begin{lemme}
Soient $V_1$ et $V_2$ deux espaces vectoriels, $V=V_{1}\oplus V_{2}$ et $\pi$ la projection de $V$ sur $V_2$, $i$ l'injection de $V_1$ dans $V$. Soient $W_1$, $W_{1}'$ (resp. $W_2$, $W_{2}'$) des sous-espaces vectoriels de $V_1$ (resp. $ V_{2}$). Si $W$ et $W'$ sont des sous-espaces vectoriels tels que $i(W_{1})=W \cap V_{1},i(W_{1}')=W' \cap V_{1}$ et $\pi (W)=W_{2},\pi (W')=W_{2}'$, alors :

$  \emph{dim}_{\bf C}(W_{1} \cap W_{1}')+ \emph{dim}_{\bf C}(W_{2} \cap W_{2}')-\min(\emph{dim}_{\bf C}(W_{2} \cap W_{2}'),\emph{dim}_{\bf C}(V_{1}))$
\begin{eqnarray*}
 \hspace{4cm}  &   \leq   &   \emph{dim}_{\bf C}( W\cap W')\\
\hspace{4cm} & \leq & \emph{dim}_{\bf C}(W_{1} \cap W_{1}')+\emph{dim}_{\bf C}(W_{2} \cap W_{2}'). 
\end{eqnarray*}

\end{lemme}

\begin{corollaire}
Soient $A$ et $B$ deux structures de Hodge mixtes et $H \in \emph{Ext}^{1}_{MHS}(B,A)$, alors :
$$ \forall (p,q) \in {\bf Z} \times {\bf Z} \,\,\,\, f^{p,q}_{H}-f^{p,q}_{A}-f^{p,q}_{B} \leq 0$$
 
\end{corollaire}

\begin{preuve}(du corollaire)
D'après la construction de $H \in \text{Ext}^{1}_{MHS}(B,A)$ (c'est à dire une classe de suite exacte  $0 \longrightarrow A \stackrel{i}{\longrightarrow} H \stackrel{ \pi}{\longrightarrow} B \longrightarrow 0 $) suivant \cite{mor},\cite{car}, pour tout $(p,q) \in {\bf Z} \times {\bf Z}$, on est exactement dans le cadre du lemme précédent avec : $H_{\bf C}=A_{\bf C}\oplus B_{\bf C}$, $W_1=F^{p}_{A}$, $W_{1}'={\overline F}^{q}_{A}$, $W_1=F^{p}_{B}$, $W_{1}'={\overline F}^{q}_{B}$, $W=F^{p}_{H}$, $W'={\overline F}^{q}_{H}$. Par construction de la filtration de Hodge et sa filtration opposée sur l'extension, on a
$i(F^{p}_{A_{\bf C}})=A_{\bf C} \cap F^{p}_{H_{\bf C}}$, $i({\overline F}^{q}_{A_{\bf C}})=A_{\bf C} \cap {\overline F}^{q}_{H_{\bf C}}$ et $\pi(F^{p}_{H_{\bf C}})=F^{p}_{B_{\bf C}}$, $\pi({\overline F}^{q}_{H_{\bf C}})= {\overline F}^{q}_{B_{\bf C}}$ par stricte compatibilité des morphismes de structures de Hodge.
\end{preuve}

\begin{preuve}(du lemme)
Ecrivons sous forme matricielle les coordonnées des sous espaces vectoriels $W$ et $W'$ dans $V=V_{1}\oplus V_{2}$ i.e. les représentations matricielles des points $W$ et $W'$ dans les grassmanniennes $G(V,\text{dim}_{\bf C}W)$ et $G(V,\text{dim}_{\bf C}W')$. On note $M_{i}Sev$ la matrice représentant $Sev \subset V_{i}$ dans $G(V_{i},\text{dim}_{\bf C}Sev)$ pour $i \in \{\Box ,1,2\}$. La base de $V$ prise pour la représentation matricielle est la réunion des bases de $V_{1}$ et $V_{2}$. Alors :\\
$MW=\left(
    \begin{array}{ccc}
    M_{1}W_{1}   &   \vert    &      0   \\
    A_{1}     &    \vert       &     M_{2}W_{2}

\end{array}
\right)$    et $MW'=\left(
    \begin{array}{ccc}
    M_{1}W_{1}'   &   \vert    &      0   \\
    A_{1}'     &    \vert       &     M_{2}W_{2}'

\end{array}
\right)$.\\
où $A_{1}$ et $A_{2}$ sont des matrices quelconques de dimension $\text{dim}_{\bf C}W_{2} \times \text{dim}_{\bf C}V_{1} $ et $\text{dim}_{\bf C}W_{2}' \times \text{dim}_{\bf C}V_{1}$ respectivement (on vérifie que ce sont bien des sous-espaces vectoriels tels que $i(Sev_{1})=Sev \cap V_{1}$ et $\pi(Sev)=Sev_{2}$). Comme pour deux sous-espaces vectoriels on a l'égalité sur les dimensions : $\text{dim}(Sev_{1}\cap Sev_{2})+\text{dim}(Vect(Sev_{1},Sev_{2}))=\text{dim}(Sev_{1})+\text{dim}(Sev_{2})$, conna\^{\i}tre la dimension de $Vect(W,W')$ nous donne la dimension cherchée. Il s'agit donc de trouver le maximum des dimensions des matrices extraites de déterminant non nul de la matrice : 

\hspace{2cm}$MW=\left(
    \begin{array}{ccc}
    M_{1}W_{1}   &   \vert    &      0   \\
    A_{1}     &    \vert       &     M_{2}W_{2}\\
    --------&   \vert  & --------\\
    M_{1}W_{1}'   &   \vert    &      0   \\
    A_{1}'     &    \vert       &     M_{2}W_{2}'

\end{array}
\right)$\\
Une telle matrice contient des lignes et colonnes obtenues à partir des matrices extraites de déterminant non nul de dimension maximale pour les sous-espaces engendrés $Vect( W_{1},  W_{1}')$ et $Vect( W_{2},  W_{2}')$ i.e. les lignes et colonnes complétées à partir de telle matrices extraites dans :

\hspace{2cm} $MW'=\left(
    \begin{array}{c}
    M_{1}W_{1} \\
    M_{1}W_{1}' 
\end{array}
\right)$ et $MW'=\left(
    \begin{array}{c}
    M_{2}W_{2} \\
    M_{2}W_{2}' 
\end{array}
\right)$.\\
On peut donc trouver une matrice extraite de déterminant non nul de dimension au moins égale à $ \text{dim}\,Vect( W_{1},  W_{1}')+\text{dim}\,Vect( W_{2},  W_{2}')$. On peut ajouter à cela au plus $\text{min}(\text{dim}\,( W_{2} \cap  W_{2}'), \text{dim} \,V_{1})$ lignes et colonnes pour obtenir une matrice de déterminant non nul. D'où l'inégalité :
\begin{center}
$ \text{dim}(Vect(W_{1},  W_{1}'))+\text{dim}(Vect(W_{2},  W_{2}'))   \leq  \text{dim}(Vect(W_{},  W_{}')) $\\
$  \leq  \text{dim}(Vect(W_{1},  W_{1}'))+\text{dim}(Vect(W_{2},  W_{2}'))+ \text{min}(\text{dim}\,( W_{2} \cap  W_{2}'), \text{dim} \,V_{1})$
\end{center}
En découle l'inégalité voulue en passant aux codimensions. 

%$B=\left(
%    \begin{array}{cccccccccc}
%     1  & 0                  & ... & ...              & 0                & \vert   & 0 & ... & 0 \\
%     \lambda_{1} & 1    &   -1 & 0             & 0                & \vert   & ... & ... & ...           \\
%     ...              & 0                 & 1  & -1               &...               & \vert    & ... & ... & ...           \\
%      ...             &...                 &...  &...               &...               & \vert     & ... &
%c_{ij}=log( \frac{\theta(Q_{j}-p_{i}-\frac{1}{2} (1+ \tau ))}{\theta(Q_{j}-p_{i+1}-\frac{1}{2}
%(1+\tau ))})-log( \frac{\theta(P_{j}-p_{i}-\frac{1}{2} (1+ \tau ))}{\theta(P_{j}-p_{i+1}-\frac{1}{2} (1+ \tau ))}) & ...   \\
%      ...             &...                 &...  &...               &...               & \vert      & ... & ... & ...          \\
%       \lambda_{k-1}              & ...                & 0   & 1  & -1  & \vert         & ... & ... & ...
%\end{array}
%\right)$\\

\end{preuve}

\begin{preuve}(du théorème)
 Comme ${\mathcal E }_{H}={\mathcal E }_{A} \cup {\mathcal E }_{B}$ est un sous-ensemble fini de ${\bf Z} \times {\bf Z}$, on peut trouver $k \in {\bf Z}$ et $N \in {\bf Z}$ tels que ${\mathcal E }_{H \otimes {T \langle k \rangle }}={\mathcal E }_{A \otimes {T \langle k \rangle }} \cup {\mathcal E }_{B\otimes {T \langle k \rangle }} \subset [0,N] \times [0,N]$. D'après le théorème \ref{elemalpha} $(i)$, l'égalité n'est pas modifiée par tensorisation des structures de Hodge $A$,$B$ et $H$ par une structure de Hodge de Tate. Tensoriser les structures de Hodge dans une extension par une structure de Hodge de Tate ${T \langle k \rangle }$ induit un isomorphisme entre $\text{Ext}^{1}_{MHS}(B,A)$ et $\text{Ext}^{1}_{MHS}(B\otimes {T \langle k \rangle },A\otimes {T \langle k \rangle })$. Quitte à remplacer les structures de Hodge $A,B$ et $H$ par leurs tensorisées par ${T \langle k \rangle }$, on peut ainsi supposer que  ${\mathcal E }_{H } \subset [0,N] \times [0,N]$. La proposition se ramène donc à montrer que : 
$$ \alpha(H) - \alpha(A)-\alpha(B) =\frac{1}{2} \sum_{(p,q) \in {[0,N] \times [0,N]}} (p+q)^{2}(s^{p,q}_{A}+s^{p,q}_{B}-s^{p,q}_{H}) \geq 0$$
Pour une structure de Hodge mixte $\bullet$, on a les suites exactes suivantes :

$$
\xymatrix{
{}   &   0  \ar[d]    &    0   \ar[d]   &   {}    &   {}\\
{}   &   F^{p+1}_{\bullet}\cap{\overline F}^{q+1}_{\bullet}  \ar[d]    &    F^{p}_{\bullet}\cap{\overline F}^{q+1}_{\bullet}   \ar[d]   &   {}    &   {}\\
  {}   &   F^{p+1}_{\bullet}\cap{\overline F}^{q}_{\bullet}  \ar[d]    &    F^{p}_{\bullet}\cap{\overline F}^{q}_{\bullet}   \ar[d]   &   {}    &   {}\\
 0  \ar[r]  &   F^{p+1}_{\bullet}Gr_{{\overline {F_{\bullet}}}}^{q}   \ar[r] \ar[d] &  F^{p}_{\bullet}Gr_{{\overline {F_{\bullet}}}}^{q}   \ar[r] \ar[d]   &  Gr_{{ {F_{\bullet}}}}^{p}Gr_{{\overline {F_{\bullet}}}}^{q}  \ar[r]  &   0\\
{}    &   0    &    0    &  {}   &   {}
}
$$
Et donc la relation sur les dimensions :$$ s^{p,q}_{\bullet}=f^{p,q}_{\bullet}-f^{p+1,q}_{\bullet}-f^{p,q+1}_{\bullet}+f^{p+1,q+1}_{\bullet}.$$
Posons $\alpha(\bullet)=\alpha^{+}(\bullet)-\alpha^{-}(\alpha)$ où $\alpha^{+}(\bullet)$ est constitué des termes de la somme avec les nombres de Hodge et $\alpha^{-}(\bullet)$ est constitué des termes de la somme qui comportent les entiers $s^{p,q}_{\bullet }$. Comme pour tout $(p,q) \in {\bf Z}$ $h^{p,q}_{H}=h^{p,q}_{A}+h^{p,q}_{B}$, $\alpha^{+}$ est additif. Montrer que $\alpha$ est sur-additif revient donc à montrer que $\alpha^{-}$ est sous-additif. On a :
$$\alpha^{-}( {\bullet})= \frac{1}{2} \sum_{(p,q) \in {[0,N] \times [0,N]}} (p+q)^{2}\,s^{p,q}_{\bullet }=\frac{1}{2}\sum_{(p,q) \in {[0,N] \times [0,N]}} (p+q)^{2}(f^{p,q}_{\bullet}-f^{p+1,q}_{\bullet}-f^{p,q+1}_{\bullet}+f^{p+1,q+1}_{\bullet}).$$ 
On effectue un changement d'indice sur les sommes, comme pour tout $p$, $f^{p,N+1}_{\bullet}=0$, pour tout $q$, $f^{N+1,q}_{\bullet}=0$ :
\begin{center}
$\alpha^{-}( {\bullet})= \frac{1}{2} \sum_{(p,q) \in {[1,N] \times [1,N]}} \Bigl[ (p+q)^{2}-(p+q-1)^{2}-(p+q-1)^{2}+(p+q-2)^{2}              \Bigr] f^{p,q}_{\bullet}+\frac{1}{2} \sum_{q \in [0,N]}\Bigl[ (0+q)^{2}-(0+q-1)^{2}            \Bigr] f^{0,q}_{\bullet}$\\

$+\frac{1}{2} \sum_{p \in [1,N]}\Bigl[ (p+0)^{2}-(p-1+0)^{2}            \Bigr] f^{p,0}_{\bullet} $.
\end{center}
$\alpha^{-}( {\bullet})= \frac{1}{2} \sum_{(p,q) \in {[1,N] \times [1,N]}} ( 2pq+4) f^{p,q}_{\bullet}+\frac{1}{2} \sum_{q \in [0,N]} ( 2q-1) f^{0,q}_{\bullet}+\frac{1}{2} \sum_{p \in [1,N]}(2p-1) f^{p,0}_{\bullet} $,\\
que nous écrirons, afin que les coefficients des $f^{p,q}_{\bullet}$ sous les $\sum$ soient tous positifs, sous la forme :\\
$\alpha^{-}( {\bullet})= \frac{1}{2} \sum_{(p,q) \in {[1,N] \times [1,N]}} ( 2pq+4) f^{p,q}_{\bullet}+\frac{1}{2} \sum_{q \in [1,N]} ( 2q-1) f^{0,q}_{\bullet}+\frac{1}{2} \sum_{p \in [1,N]}(2p-1) f^{p,0}_{\bullet}- \frac{1}{2} f^{0,0}_{\bullet}$\\
Ainsi :\\
$\alpha( H)-\alpha(A)-\alpha(B)=\alpha^{-}( H)-\alpha^{-}(A)-\alpha^{-}(B)= \frac{1}{2} \sum_{(p,q) \in {[1,N] \times [1,N]}} ( 2pq+4) \Bigl[ f^{p,q}_{H}-f^{p,q}_{A}-f^{p,q}_{B} \Big]+\frac{1}{2}\sum_{q \in [1,N]} ( 2q-1) \Bigl[ f^{0,q}_{H}-f^{0,q}_{A}-f^{0,q}_{B} \Big]+\frac{1}{2} \sum_{p \in [1,N]}(2p-1) \Bigl[ f^{p,0}_{H}-f^{p,0}_{A}-f^{p,0}_{B} \Big]- \frac{1}{2} \Bigl[ f^{0,0}_{H}-f^{0,0}_{A}-f^{0,0}_{B} \Big]$.\\
$f^{0,0}_{H}-f^{0,0}_{A}-f^{0,0}_{B}=\text{dim}_{\bf C}H-\text{dim}_{\bf C}A-\text{dim}_{\bf C}B=0$, donc tous les coefficients des $f^{p,q}_{\bullet}$ dans les $\sum$ étant positifs, d'après le lemme précédent, pour tous $(p,q) \in [0,N] \times [0,N]$, $f^{p,q}_{H}-f^{p,q}_{A}-f^{p,q}_{B} \leq 0$. D'où le résultat.

%$$ 0 \rightarrow F^{p+1}_{\bullet}Gr_{{\overline {F_{\bullet}}}}^{q} \rightarrow F^{p}_{\bullet}Gr_{{\overline {F_{\bullet}}}}^{q} \rightarrow Gr_{{ {F_{\bullet}}}}^{p}Gr_{{\overline {F_{\bullet}}}}^{q} \rightarrow 0$$ 
% $$0 \rightarrow F^{p+1}_{\bullet}Gr_{{\overline {F_{\bullet}}}}^{q} \rightarrow F^{p}_{\bullet}Gr_{{\overline {F_{\bullet}}}}^{q} \rightarrow Gr_{{ {F_{\bullet}}}}^{p}Gr_{{\overline {F_{\bullet}}}}^{q} \rightarrow 0$$

\end{preuve} 
${}$\\
De la formule explicite de $\alpha$ dans la preuve du th\'eor\`eme  nous tirons le lemme suivant :

\begin{lemme}\label{alphaplus} Si l'on considère une famille de structures de Hodge mixtes paramétrées par une base $S$, alors en tout point $s \in S$, $\alpha_{s}=\alpha^{+}_{s}-\alpha^{-}_{s}$ où $\alpha^{+}_{s}$ est constitué des nombres de Hodge en $s$ et  
 $$\alpha^{-}_{s}= \frac{1}{2} \sum_{(p,q) \in {[1,N] \times [1,N]}} ( 2pq+4) f^{p,q}_{s}+\frac{1}{2} \sum_{q \in [1,N]} ( 2q-1) f^{0,q}_{s}+\frac{1}{2} \sum_{p \in [1,N]}(2p-1) f^{p,0}_{s}- \frac{1}{2} f^{0,0}_{s}$$.
\end{lemme}

\subsection{Calculs de niveaux de ${\bf R}$-scindemment $\alpha$ de la cohomologie}

D'après Deligne \cite{del2}, \cite{del3}, les groupes de cohomologie des variétés algébriques (schémas de type fini sur ${\bf C}$) sont munis de structures de Hodge mixtes. Dans cette section, on donne des exemples de calcul de l'invariant $\alpha$ pour des structures de Hodge mixtes qui viennent de la cohomologie de variétés algébriques. Pour ce faire nous explicitons quelques constructions de ces structures de Hodge mixtes. Précisons d'abord certains cas pour lesquels l'invariant $\alpha $ est nul. 

$\alpha$ est nul pour toutes les structures de Hodge mixtes ${\bf R}$-scindées, ce qui comprend :\\
\hspace*{1.2cm}$\bullet$ Toutes les structures de Hodge pures.\\
\hspace*{1.2cm}$\bullet$ Plus généralement, toutes les structures de Hodge mixtes de longueur inférieure ou égale à $2$, comme par exemple :\\
\hspace*{2cm}$-$ Les structures de Hodge mixtes sur la cohomologie des variétés projectives à poids (\cite{dol}).\\   
\hspace*{2cm}$-$ La cohomologie des variétés à singularit\'es logarithmiques cf \cite{usu}.

Soit $X$ un schéma de type fini sur ${\bf C}$, pour tout $k \in {\bf Z}$, on peut lui associer l'entier
$$\alpha_{k}(X)=\alpha((H^{k}(X,{\bf C}),W^{\bullet},F^{\bullet},{\overline F}^{\bullet})).$$

\begin{proposition}
$(i)$ Soit $X$ une variété algébrique, $X'$ sa normalisée et $\overline X$ sa complétée, alors pour tout $k \in {\bf Z}$ :
$$\alpha_{k}(X) \geq \alpha_{k}(X')   \text{ et }  \alpha_{k}(X) \geq \alpha_{k}({\overline X}).$$
$(ii)$ Plus généralement, si $f:X \rightarrow Y$ est un morphisme qui induit un morphisme injectif (resp. surjectif) sur la cohomologie $f^{*}:H^{k}(Y,{\bf C}) \rightarrow H^{k}(X,{\bf C})$, alors, pour tout $k \in {\bf Z}$ :
$$\alpha_{k}(X) \geq \alpha_{k}(Y)   (\text{resp. }  \alpha_{k}(Y) \geq \alpha_{k}({\overline X})).$$ 
\end{proposition}

\begin{preuve}
$(i)$ se déduit de $(ii)$ à l'aide de \cite{del3}, proposition (8.2.6) : l'application de la normalisée $X'$ vers $X$ induit un morphisme surjectif en cohomologie, l'application de $X$ vers sa complétée $\overline X$ induit un morphisme injectif  en cohomologie. La catégorie des structures de Hodge mixtes étant abélienne, on peut écrire des suites exactes à partir des surjections et injections de structures de Hodge mixtes et donc pour prouver $(ii)$, on écrit les extensions correspondantes à ces flèches injectives ou surjectives et on utilise le théorème \ref{inegext} p.\pageref{inegext}.

\end{preuve}

\begin{proposition}
Soient $X$ et $Y$ deux variétés algébriques, alors :
$$\alpha_{k}(X \times Y)=\sum_{i=0}^{i=k}\,\,\emph{dim}_{\bf C}(H^{k-i}(Y,{\bf C}))\alpha_{i}(X)+\emph{dim}_{\bf C}(H^{i}(X,{\bf C}))\alpha_{k-i}(Y).$$
\end{proposition}

\begin{preuve}
  D'après \cite{del3}, proposition $(8.2.10)$, les isomorphismes de Künneth sur la cohomologie $H^{*}(X \times Y,{\bf C}) \cong H^{*}(X,{\bf C}) \otimes H^{*}(Y,{\bf C})$ sont des isomorphismes de structures de Hodge mixtes. On applique $(iii)$ du théorème \ref{elemalpha} p.\pageref{elemalpha} à la somme donnée par la formule de Künneth, puis le $(iv)$ à chacun des facteurs.  

\end{preuve}

\subsubsection{D\'etails des filtrations et calcul de $\alpha$ pour les courbes alg\'ebriques complexes de genre $0$ ou $1$}

Soit $X$ une courbe alg\'ebrique sur ${\bf C}$, $X'$ la normalis\'ee de $X$ et $r: X' \rightarrow X$ le morphisme qui s'en d\'eduit.
Soit ${\overline X}'$ la courbe projective non singuli\`ere dont $X'$ est un ouvert dense et $\overline X$ la courbe d\'eduite de
${\overline X}'$ en contractant chacun des $r^{-1}(s)$ pour $s \in X$.
Notons $\overline r$ le morphisme ${\overline X}' \rightarrow {\overline X}$, et $j$,$j'$ les inclusions respectives de $X$ dans
$\overline X$ et  $X'$ dans ${\overline X}'$.
Soit $S$ l'ensemble fini ${\overline X}-X$. D'où le carré cartésien suivant :
$$
\xymatrix{
     X' \ar@{^{(}->}[r]^{j'} \ar[d]_r  & {\overline X}' \ar[d]^{\overline r} \\
     X \ar@{^{(}->}[r]_{j} & {\overline X}
   }
$$

\begin{proposition}\cite{del3}
 La cohomologie de $X$ est donn\'ee par donnée par :
\begin{center}
$H^{1}(X,{\bf C})={\bf H}^{1}({\overline X},[ {\mathcal O}_{\overline X} \stackrel{d}{ \rightarrow } r_{*} {\Omega }_{{\overline
X}'}^{1}(logS)])$
\end{center}
\end{proposition}
de plus\\
\begin{proposition}\cite{del3}
La filtration par le poids est donnée par :\\
${}$\\
$W_{1}(H^{1}(X,{\bf C}))=Im(H^{1}({\overline X},{\bf C}) \rightarrow H^{1}(X,{\bf C})$\\
$W_{0}(H^{1}(X,{\bf C}))=Ker(H^{1}({\overline X},{\bf C}) \rightarrow H^{1}({\overline X}',{\bf C}).$\\

De plus la suite spectrale d\'efinie par la filtration b\^ete de $[ {\mathcal O}_{\overline X} \stackrel{d}{\rightarrow} r_{*}
{\Omega}_{{\overline X}'}^{1}(logS)]$ d\'eg\'en\`ere  en $E_{1}$ et aboutit \`a la filtration de Hodge de $H^{\bullet }(X,{\bf C})$.\\
\end{proposition}
Calculons l'hypercohomologie du complexe $[ {\mathcal O}_{\overline X} \stackrel{d}{\rightarrow} r_{*}{\Omega}_{{\overline
X}'}^{1}(logS)]$. Pour un complexe de faisceaux $({\mathcal K^{\bullet }},d)$ sur une vari\'et\'e $X$ muni d'un recouvrement acyclique ${\bf \calu}=(\calu_{i})_{i \in I}$, l'hypercohomologie ${\bf H}^{\bullet }(X,{\mathcal K}^{\bullet })$ du complexe est d\'efinie comme \'etant la cohomologie du complexe total $ C^{(p,q)}=C^{p}({\bf \calu},\delta ,d )$ o\`u $d $ est l'op\'erateur du complexe et $\delta $ celui de $\check{C}$ech associ\'e au recouvrement. ${\bf H}^{\bullet }(X,{\mathcal K}^{\bullet })$ est l'aboutissement des suites spectrales d\'efinies par les filtrations \'evidentes de $C^{(p,q)}$ (cf \cite{gh}).
\\
\\
Le complexe double associ\'e \`a $[ {\mathcal O}_{\overline X} \stackrel{d}{\rightarrow} r_{*}{\Omega}_{{\overline
X}'}^{1}(logS)]$ est donc ici, relativement à ${\bf \calu}$ :
\\
${\bf (*)}$

%\begin{equation*}
%\begin{CD}
%{C^{0}({\mathcal U} \coprod {\mathcal V},{\mathcal O}_{\overline X})} @>{\delta}>>           {C^{1}({\mathcal {UV}},{\mathcal
%O}_{\overline X})}\\
%@VdVV                                                                                         @VVdV\\
%C^{0}({\mathcal U} \coprod {\mathcal V},r_{*}{\Omega}_{{\overline X}'}^{1}(logS)) @>{\delta}>>
%C^{1}({\mathcal {UV}},r_{*}{\Omega}_{{\overline X}'}^{1}(logS))\\
%\end{CD}
%\end{equation*}
$$
\xymatrix{
{C^{0}({\coprod_{i} \calu_{i}} ,{\mathcal O}_{\overline X})} \ar[r]^{\delta} \ar[d]_{d}     &    {C^{1}({\coprod_{i,j} \calu_{i} \cap \calu_{j}},{\mathcal O}_{\overline X})}   \ar[d]_{d}\\
C^{0}({\coprod_{i} \calu_{i}} ,r_{*}{\Omega}_{{\overline X}'}^{1}(logS)) \ar[r]_{\delta}  & C^{1}({\coprod_{i,j} \calu_{i} \cap \calu_{j}},r_{*}{\Omega}_{{\overline X}'}^{1}(logS))
}
$$

Pour voir dans $H^{1}(X,{\bf C})$ les différents sous-espaces associés à la filtration par le poids $W_{0} ,W_{1} ,W_{2} $, on utilise la proposition pr\'ec\'edente qui se traduit par :
\begin{center}
$W_{1}H^{1}({X},{\bf C})={\bf H}^{1}({\overline X},[ {\mathcal O}_{\overline X} \stackrel{d}{\rightarrow} r_{*}
{\Omega}_{{\overline X} '}^{1}])$
\end{center}
et comme annonc\'e : 
\begin{center}
$W_{2}H^{1}(X,{\bf C})={\bf H}^{1}({\overline X},[{\mathcal O}_{\overline X} \stackrel{d}{\rightarrow} r_{*} {\Omega}_{{\overline
X}'}^{1}(logS)])$.
\end{center}
On peut donc "d\'emonter" l'\'etude de l'hypercohomologie comme il suit :
\\
\\
$\bf (**)$
\scalebox{0.8}[0.8]{
$$
\xymatrix {
   {C^{0}({\coprod_{i} \calu_{i}} ,{\mathcal O}_{\overline X})} \ar[rr] \ar[dd] \ar[dr] && {C^{1}({\coprod_{i,j} \calu_{i} \cap \calu_{j}},{\mathcal O}_{\overline X})} \ar[dr] \ar[dd] |!\hole \\
   & {C^{0}({\coprod_{i} \calu_{i}} ,{\mathcal O}_{\overline X})} \ar[rr] \ar[dd] && {C^{1}({\coprod_{i,j} \calu_{i} \cap \calu_{j}},{\mathcal O}_{\overline X})} \ar[dd] \\
   C^{0}({\coprod_{i} \calu_{i}} ,r_{*}{\Omega}_{{\overline X}'}^{1}) \ar[rr] |!\hole \ar[dr] && C^{1}({\coprod_{i,j} \calu_{i} \cap \calu_{j}},r_{*}{\Omega}_{{\overline X}'}^{1}) \ar[rd] \\
   &  C^{0}({\coprod_{i} \calu_{i}} ,r_{*}{\Omega}_{{\overline X}'}^{1}(logS)) \ar[rr] && C^{1}({\coprod_{i,j} \calu_{i} \cap \calu_{j}},r_{*}{\Omega}_{{\overline X}'}^{1}(logS)) \\
 }  
$$
}
${}$\\
où le diagramme au deuxième plan fournit la cohomologie de la complétée de la courbe $\overline X$ qui est isomorphe à $W_{1}H^{1}({ X},{\bf C})$.

%\begin{equation*}
%\begin{CD}
%{C^{0}({\mathcal U} \coprod {\mathcal V},{\mathcal O}_{\overline X})} @>{\delta}>>           {C^{1}({\mathcal {UV}},{\mathcal
%O}_{\overline X})}\\
%@VdVV                                                                                         @VVdV\\
%C^{0}({\mathcal {U}} \coprod {\mathcal {V}},r_{*}{\Omega}_{{\overline X}'}^{1}) @>{\delta}>>
%C^{1}({\mathcal {UV}},r_{*}{\Omega}_{{\overline X}'}^{1})\\
%\end{CD}
%\end{equation*}
%qui fournit la cohomologie venant du fait que $X$ soit singuli\`ere. Et le premier diagramme donn\'e.
%\\
%Effectuons le calcul du ${\bf H}^1$ de $\bf (**)$:
%\\
%\\

\subsubsection{Exemple de ${\bf P}^1$ avec quatre points marqués}

Explicitons la décomposition précédente sur le premier groupe de la cohomologie ${ H}^1(X,{\bf C})$ de la courbe $X$ obtenue à partir de ${\bf P}^1$ avec quatre points distincts distingués $\{p_{1},p_{2},P_{1},Q_{1} \}$(dans l'ordre) dont on a enlevé les deux premiers $p_{1}$ et $p_{2}$ et recollé les deux autres $P_{1}$,et $Q_{1}$. (Pour une justification du recollement, voir \cite{ser}, chapitre IV., partie 3.).  

Pour effectuer le calcul, utilisons le recouvrement standard de  ${\bf P}^1$ par les cartes affines de coordonnées $\calu ={\bf P}^1-\{\infty \}$ munie de la coordonnée $u$ et $\calv ={\bf P}^1-\{0\}$ munie de la coordonnée $v$. On notera $\calu \calv$ l'ouvert intersection des deux cartes.

Nous commencerons par calculer la cohomologie de courbe complétée $\overline X$. Soit  $f \in C^{1}({\mathcal {UV}},{\mathcal O}_{\overline X})$ (on choisira $u=v^{-1}$ comme coordonn\'ee d'écriture sur ${\mathcal {UV}}$). Alors $f$ et $df$ peuvent s'\'ecrirent :
\begin{center}
$f(u)= \sum_{- \infty}^{\infty} a_{n}u^{n}$ et $df(u)= \sum_{- \infty}^{\infty} na_{n}u^{n-1}$
\end{center}
Prenons pour $g \in C^{0}({\mathcal U} \coprod {\mathcal V},r_{*}{\Omega}_{{\overline X}'}^{1})$, $g_{0}(u)$ \'egal \`a la partie de
$dg(u)$ \`a exposants positifs en $u$ et $g_{1}(v)$ la partie de $dg(u)$ correspondant aux exposants n\'egatifs en $u$. Alors, on a, pour tout $u \in \calu \calv$ : 
\begin{center}
$\delta (g_{0},g_{1})(u)=dg(u)$
\end{center}
Donc tout \'el\'ement de $C^{1}({\mathcal {UV}},{\mathcal O}_{\overline X})$ voit son
image par $d$ dans $C^{1}({\mathcal {UV}},r_{*}{\Omega}_{{\overline X}'}^{1})$ annul\'ee par l'image par $\delta$ d'un \'el\'ement
de $C^{0}({\mathcal U} \coprod {\mathcal V},r_{*}{\Omega}_{{\overline X}'}^{1})$, donc l'hypercohomologie ${\bf H}^{1}({\overline X},[ {\mathcal O}_{\overline X} \stackrel{d}{\rightarrow} r_{*}{\Omega}_{{\overline X} '}^{1}])$ c'est à dire la cohomologie du complexe
$$
\xymatrix{
{C^{0}({\mathcal U} \coprod {\mathcal V},{\mathcal O}_{\overline X})} \ar[r]^-{(\delta,d)} &{C^{1}({\mathcal {UV}},{\mathcal O}_{\overline X})} \oplus C^{0}({\mathcal U} \coprod {\mathcal V},r_{*}{\Omega}_{{\overline X}'}^{1})  \ar[r]^-{d+\delta} & C^{1}({\mathcal {UV}},r_{*}{\Omega}_{{\overline X}'}^{1})
}
$$
%\begin{equation*}
%\begin{CD}
%{C^{0}({\mathcal U} \coprod {\mathcal V},{\mathcal O}_{\overline X})} @>{\delta}>>           {C^{1}({\mathcal {UV}},{\mathcal
%O}_{\overline X})}\\
%@VdVV                                                                                         @VVdV\\
%C^{0}({\mathcal U} \coprod {\mathcal V},r_{*}{\Omega}_{{\overline X}'}^{1}(logS)) @>{\delta}>>
%C^{1}({\mathcal {UV}},r_{*}{\Omega}_{{\overline X}'}^{1}(logS))\\
%\end{CD}
%\end{equation*}
est d\'etermin\'e par l'image par $\delta$ de $C^{0}({\mathcal U} \coprod {\mathcal V},{\mathcal O}_{\overline X})$ dans $C^{1}({\mathcal {UV}},{\mathcal
O}_{\overline X})$.\\

Nous avons donc à chercher quelles sont les images des  \'el\'ements $(h_{0},h_{1}) \in C^{0}({\mathcal U} \coprod {\mathcal V},{\mathcal
O}_{\overline X})$. Comme $\overline X$ est donnée par ${\bf P}^1$ avec deux  points $P_{1}$ et $Q_{1}$ identifi\'es, on doit donc avoir $f(P_{0})=f(Q_{0})$ et comme plus haut, on peut montrer qu'il existe un \'element $(h_{0},h_{1}) \in C^{0}({\mathcal U} \coprod {\mathcal V},{\mathcal O}_{{\bf P}^1})$ tel que $\delta (h_{0},h_{1})=f$. Cet \'el\'ement n'est cependant pas forc\'ement dans  $C^{0}({\mathcal U} \coprod {\mathcal V},{\mathcal
O}_{\overline X})$ mais l'on a :
\\
\begin{center}
$f(P_{1})=h_{0}(P_{1})-h_{1}(P_{1})=h_{0}(Q_{1})-h_{1}(Q_{1})=f(Q_{1})$
donc il existe $\lambda $ tel que
$h_{0}(Q_{1})=h_{0}(P_{1})+\lambda $ et $h_{1}(Q_{1})=h_{1}(P_{1})+\lambda $
\end{center}
Posons $h_{0}'=h_{0}- \lambda (\frac{u-{u_{P_1}}}{{u_{Q_1}}-{u_{P_1}}})$ et  $h_{1}'=h_{1}- \lambda
(\frac{v-v_{P_1}}{v_{Q_1}-v_{P_1}})$. Alors $h_{0}'(Q_{1})=h_{0}'(P_{1})$ et  $h_{1}'(Q_{1})=h_{1}'(P_{1})$ donc ces \'el\'ements sont bien dans $C^{0}({\mathcal U} \coprod {\mathcal V},{\mathcal O}_{\overline X})$ et $h_{0}(u)-h_{1}(u)=f(u)- \lambda (\frac{u-u_{P_1}}{u_{Q_1}-u_{P_1}}-\frac{v-v_{P_1}}{v_{Q_1}-v_{P_1}})$ et donc le coker de $\delta$ est engendr\'e par l'\'el\'ement $(\frac{u-u_{P_1}}{u_{Q_1}-u_{P_1}}-\frac{v-v_{P_1}}{v_{Q_1}-v_{P_1}})$.

Cet \'element est l'image d'un g\'en\'erateur de $H^{1}({\overline X},{\bf C})$, il provient de la boucle dans $\overline X$
form\'ee par l'identification de deux points dans ${\bf P}^1$ (de premier groupe de cohomologie trivial).

Afin de déterminer $\alpha_{1}(X)=\alpha(H^{1}({\overline X},{\bf C}))$, calculons $F^{1}{\bf H}^{1}({\overline X},[ {\mathcal O}_{\overline X} \stackrel{d}{ \rightarrow } r_{*} {\Omega }_{{\overline X}'}^{1}(logS)])$. C'est la partie de la cohomologie qui vient du noyau de la flèche $$\delta \,:\,C^{0}({\calu \coprod \calv} ,r_{*}{\Omega}_{{\overline X}'}^{1}(logS)) \rightarrow C^{1}({\calu \cap \calv},r_{*}{\Omega}_{{\overline X}'}^{1}(logS)).$$ 
Il suffit donc de prendre une forme qui n'est pas exacte $\omega \in C^{0}({\calu \cap \calv} ,r_{*}{\Omega}_{{\overline X}'}^{1}(logS))$ et prendre $(\omega,\omega) \in 
C^{0}({\calu \coprod \calv} ,r_{*}{\Omega}_{{\overline X}'}^{1}(logS))$. D'où le lemme,

\begin{lemme}

On peut donc prendre pour g\'en\'erateur de $F^{1}H^{1}(X,{\bf C})$ la forme $$\omega=(\frac{1}{u-p_{1}}-\frac{1}{u-p_{2}})du.$$
\end{lemme}
${}$\\

\begin{center}
{\bf Calcul effectif du niveau de ${\bf R}$-scindement de la courbe $X$}
\end{center}
${}$\\

Pour calculer le niveau de ${\bf R}$-scindement associé à la structure de Hodge mixte du premier groupe de cohmologie de la courbe il faut déterminer en termes de dimension les positions relatives de la filtration de Hodge et sa filtration conjuguée, c'est à dire les dimensions des intersections deux à deux des sous-espaces vectoriels qui composent les deux filtrations, de Hodge et conjuguée. Comme elle est de niveau $1$, il suffit de trouver la dimension  $$s^{1,1}_{H^{1}(X,{\bf C})}=\text{dim}_{\bf C}\,F^{1}H^{1}(X,{\bf C}) \cap {\overline F}^{1}H^{1}(X,{\bf C})$$
qui est $0$ ou $1$. La filtration de Hodge induit une filtration sur le dual du premier groupe de cohomologie par :
\begin{center}
$H^{1}_{DR}(X, {\bf C}) \cong H^{1}_{Betti}(X, {\bf C}) \cong (H_{1}(X, {\bf C}))^{*}=(H_{1}(X, {\bf Z}) \otimes {\bf R})^{*}\otimes {\bf C}$
\end{center}
Choisissons des générateurs de $H_{1}(X, {\bf Z})$, ${\gamma}_0$ et ${\gamma}_1$. Ces générateurs serviront de base de $H_{1}(X, {\bf R})$ puis de base invariante par conjugaison complexe de $H_{1}(X, {\bf C})$. Pour ${\gamma}_0$ prenons le
lacet qui serait homologue \`a $0$ dans $X \cup p_0$ et pour ${\gamma}_1$ un lacet d\'ecrivant la boucle formée par l'identification des deux points. Le resultat ne d\'epend pas du choix des g\'en\'erateurs du $H_1$ puisqu'un choix différent revient à multiplier la matrice des périodes par une matrice à coefficients rationnels et préserve la colinéarité complexe. Int\'egrons donc $\omega $ sur les cycles :

\begin{eqnarray*}
<w,{\gamma}_0>&=&\int_{{\gamma}_0}w=2\pi i \text{ (c'est le r\'esidu de } \omega  \text{ en } p_1),\\
<w,{\gamma}_1>&=&\int_{{\gamma}_1}w\\
{}&=& (\int_{0}^{P_1}-\int_{0}^{Q_1}) (\frac{1}{u-p_{1}}-\frac{1}{u-p_{2}})du\\
{}&=&\lbrack \log(\frac{u-p_1}{u-p_2}) {\rbrack}^{Q_1}_{P_1}\\
{}&=&\log(\frac{Q_{1}-p_1}{Q_{1}-p_2})-\log(\frac{P_{1}-p_1}{P_{1}-p_2})
\end{eqnarray*}
L'int\'egrale entre les points $P_{1}$ et $Q_{1}$ sur la courbe $X'$ d\'esingularis\'ee de $X$ repr\'esente en cohomologie le cycle
cr\'e\'e par l'identification de ces points dans $X'$, boucle non nul-homologue dans $X$.\\

Notons $(Q_{1},P_{1},p_{1},p_{2}):=( \frac{Q_{1}-p_{1}}{Q_{1}-p_{2}})/ (\frac{P_{1}-p_{1}}{P_{1}-p_{2}})$ le
birapport des quatre points $Q_{1},P_{1},p_{1},p_{2}$ (lorsque au moins trois d'entre eux sont diff\'erents). Si l'un d'eux est
l'$\infty$ le birapport peut \^etre d\'efini en passant \`a la limite par $(\infty
,P_{1},p_{1},p_{2}):=\frac{P_{1}-p_{2}}{P_{1}-p_{1}}$

Avec cette notation on a $<w,{\gamma}_1>=\log(Q_{1},P_{1},p_{1},p_{2})$. Supposons qu'il existe $\lambda
\in {\bf C}$ tel que le vecteur de ${\bf C}^2$ $(<w,{\gamma}_0>,<w,{\gamma}_1>)$ soit proportionnel (complexe) au vecteur obtenu par conjugaison :
\begin{center}
$( 2\pi i,\log(\frac{Q_{1}-p_1}{Q_{1}-p_2})-\log(\frac{P_{1}-p_1}{P_{1}-p_2}))= \lambda \overline{( 2\pi
i,\log(\frac{Q_{1}-p_1}{Q_{1}-p_2})-\log(\frac{P_{1}-p_1}{P_{1}-p_2}))}$
\end{center}
Comme l'action de $PGL(1)$  sur ${\bf P}^1$ est transitive sur trois points, on peut supposer sans changer la condition de
colin\'earit\'e que l'on a $p_{1}=0, p_{2}=1$ et $P_{1}=\infty$ (on garde la notation $Q_{1}$ pour son image par $\tau$), d'o\`u $\,\,\log(\frac{Q_{1}-p_1}{Q_{1}-p_2})-\log(\frac{P_{1}-p_1}{P_{1}-p_2})=\log(\frac{Q_{1}}{Q_{1}-1})$.
Or la colin\'earit\'e impose que $\lambda =-1$ et donc que $\log(\frac{Q_{1}}{Q_{1}-1})$ soit un imaginaire pur c'est \`a dire que $Q_1$ soit sur la droite $\Re(u)=\frac{1}{2}$. Donc :
$$
\left\lbrace \begin{array}{l}
        s^{1,1}_{H^{1}(X,{\bf C})}=\text{dim}_{\bf C}F^{1}H^{1}(X,{\bf C}) \cap {\overline F}^{1}H^{1}(X,{\bf C})=1 \text{ si } Q_{1} \in \Re(u)=\frac{1}{2},\\
        s^{1,1}_{H^{1}(X,{\bf C})}=0 \text{ autrement, et alors } s^{1,0}_{H^{1}(X,{\bf C})}=1 \text{ et } s^{1,0}_{H^{1}(X,{\bf C})}=1.
\end{array}
\right.
$$
Comme
\begin{eqnarray*}
\alpha_{1}(X_{Q_1})&=&\frac{1}{2}((0+0)^{2}(h^{0,0}_{H^{1}(X_{Q_1},{\bf C})}-s^{0,0}_{H^{1}(X_{Q_1},{\bf C})})+(1+0)^{2}(h^{1,0}_{H^{1}(X_{Q_1},{\bf C})}-s^{1,0}_{H^{1}(X_{Q_1},{\bf C})}))\\
{}&=&\frac{1}{2}((0+1)^{2}(h^{0,1}_{H^{1}(X_{Q_1},{\bf C})}-s^{0,1}_{H^{1}(X_{Q_1},{\bf C})})+(1+1)^{2}(h^{1,1}_{H^{1}(X_{Q_1},{\bf C})}-s^{1,1}_{H^{1}(X_{Q_1},{\bf C})}))
\end{eqnarray*}
${}$\\
\begin{proposition}
Ainsi le niveau de ${\bf R}$-scindement du premier groupe de cohomologie de $X=X_{Q_1}$ est donné par :

$$
\left\lbrace \begin{array}{l}
         \alpha_{1}(X)=\alpha(H^{1}(X,{\bf C}))=0 \,\, si \,\, Q_{1} \in {\mathcal M}_{0,4}- {\{\Re=\frac{1}{2} \}}\\
         \alpha_{1}(X)=\alpha(H^{1}(X,{\bf C}))=1 \,\, si Q_{1} \in \{\Re=\frac{1}{2}\}
\end{array}
\right.
$$
\end{proposition}
\subsubsection{G\'en\'eralisation sur ${\bf P}^1$}

Soit $X={\bf P}^{1}$ avec $m$ points enlev\'es $p_{1},...,p_{m}$ et $2n$ points identifi\'es $P_{1},Q_{1},...,P_{n},Q_{n}$deux à deux (on identifie $P_{k}$ à $Q_{k}$ pour tout $k \in [1,n]$). $H^{1}(S,{\bf C})$ est de rang $m+n-1$, il est engendr\'e par $m+n-1$ \'el\'ements dont $n$ proviennent des boucles form\'ees par l'identification des points et $m-1$ repr\'esentent les lacets autour cr\'e\'es par la non-complétude de $X$. Le calcul de l'hypercohomologie est similaire au précédent, on recouvre ${\overline X}'={\bf P}^{1}$ par les deux ouverts standards $\calu $ $\calv $ munis des coordonnées $u$ et $v$. Comme générateurs de $F^{1}H^{1}(X,{\bf C})$  on peut prendre les \'el\'ements $ \omega_{i}=( \frac{1}{u-p_{i}} -\frac{1}{u-p_{i+1}}) du $ avec $i \in [1,m-1]$. ''Filtrons'' l'homologie par la filtration de Hodge :  
D\'esignons par $\gamma_j$ $ j \in [1,m-1]$ les g\'en\'erateurs de $ H_{1}(X,{\bf
C})$ qui sont homologues à zéro dans $X \cup p_j$ et par $\beta_j$, $j \in [1,n]$ les lacets \'el\'ements du $ H_{1}(X,{\bf C})$
repr\'esentants les boucles donn\'ees par $P_{j}=Q_j$ respectivement et ne passant par aucun autre point $P_{k}=Q_k$ pour $k \in
[1,n]$.

Les $\beta_j$ ne sont pas d\'efinis canoniquement mais modulo les $\gamma_j$. Changer la base du premier groupe d'homologie ne modifie pas les conditions de colinéarité dans la matrice des périodes car cela revient à la multiplier par une matrice à coefficients rationnels.\\ 
Filtrons ${\bf C}^{m+n-1}$ par $F^1$, on obtient alors $m-1$ vecteurs :\\
\begin{center}
$(<\alpha_{1}, \omega_{i}>,...,<\alpha_{m-1}, \omega_{i}>,<\beta_{1}, \omega_{i}>,...,<\beta_{n}, \omega_{i}>)$ pour $i \in
[1,m-1]$
\end{center}
et on regarde s'ils sont colin\'eaires complexes \`a leur conjugu\'e. On obtient une matrice $A=(B,C) \in Mat_{(m-1)\times (m+n-1)}({\bf C}) $ telle que les coefficients de $B=(b_{i,j})_{i \in [1,m-1], j \in
[1,m+-1]}$ sont donn\'es par :
$$
\left\lbrace \begin{array}{l}
        b_{i,j}= 2\pi \sqrt{-1}  \,\, si \,\, i=j\\
        b_{i,j}= -2\pi \sqrt{-1}  \,\, si \,\, i=j+1\\
        b_{i,j}=0 \,\, sinon
\end{array}
\right.
$$
où les coefficients sont donnés par les résidus des formes qui $\omega_i$. Les coefficients de $C=(c_{i,j})_{i  \in [1,m-1],j \in [m,n+m-1]}$ sont :
\begin{center}
$c_{i,j}=\log(p_{i},p_{i+1},P_{j},Q_{j})$
\end{center}
La matrice des périodes $A=(B,C)$ est donc :

$$A=\left(
    \begin{array}{ccccccccc}
     2 \sqrt{-1} \pi  & 0                  & ... & ...              & 0                & \vert   & ... & ... & ... \\
     -2 \sqrt{-1} \pi & 2 \sqrt{-1} \pi    &   0 & ...              & 0                & \vert   & ... & ... & ...           \\
     ...              &...                 &...  &...               &...               & \vert    & ... & ... & ...           \\
      ...             &...                 &...  &...               &...               & \vert     & ... &
c_{ij}=(\log(p_{i},p_{i+1},P_{j},Q_{j})) & ...            \\
      ...             &...                 &...  &...               &...               & \vert      & ... & ... & ...          \\
       0              & ...                & 0   & -2 \sqrt{-1} \pi & 2 \sqrt{-1} \pi  & \vert         & ... & ... & ...
\end{array}
\right)$$
Les vecteurs lignes de la matrice $A$ engendrent l'image de $F^{1}H^{1}(X,{\bf C})$ sur $H_{1}(X,{\bf C})$ et donc calculer $s^{1,1}_{H^{1}(X,{\bf C})}=dim_{\bf C}\,F^{1}H^{1}(X,{\bf C}) \cap {\overline F}^{1}H^{1}(X,{\bf C})$ revient à calculer le rang de la matrice $\left(
\begin{array}{c}
A\\
{\overline A}
\end{array}
\right) \in  Mat_{(2m-2)\times (m+n-1)}({\bf C}) $. D'aprés la forme de $B$ ceci revient à voir si chaque vecteur colonne de $A$ est opposé à son conjugué, or :
$$c_{i,j}=-{\overline c}_{i,j} \Leftrightarrow \log(p_{i},p_{i+1},P_{j},Q_{j})=-{\overline {\log(p_{i},p_{i+1},P_{j},Q_{j})}} \Leftrightarrow \vert c_{i,j} \vert =1$$
Finalement notons ${\mathcal M}_{0,m+2n}$ l'espace des modules des courbes stables de genre $0$ sur ${\bf C}$ avec $m+2n$ points marqués $(p_{1},...,p_{m},P_{1},Q_{1},...,P_{n},Q_{n})$. Il est naturellement isomorphe à $({\bf P}^{1}-\{0,1,\infty\})^{m+n-1} - \Delta$ où $\Delta $ est la grande diagonale définie par l'ensemble des $n-3$-uplets dec points de $({\bf P}^{1}-\{0,1,\infty\})$ tels que deux au moins co\"{\i}ncident. Ainsi :\\
\begin{proposition} Le niveau de ${\bf R}$-scindement de $X$ est donné par
$$
\alpha_{1}(X)=(m-1)-\emph{card} \{ i \in [1,m-1]\, \vert \, \forall j \in [m,m+n-1]\,\,, \vert \log (p_{i},p_{i+1},P_{j},Q_{j}) \vert=1 \}.
$$
\end{proposition}
${}$\\
{\bf Remarque :} Pour tout $i,j \in [1,m-1] \times [m,m+n-1]$, il existe $\tau \in PGL(1)$ tel que $\tau(p_{i},p_{i+1},P_{j},Q_{j})=(0,1,\infty,\tau(Q_{j}))$ et donc le vecteur ligne $i$ est colinéaire à son conjugué ssi pour tout $j \in  [m,m+n-1]$ il exite un élément $ \tau \in PGL(1)$  tel que $\tau(Q_{j})$ appartienne à la droite de partie réèlle $\frac{1}{2}$.\\ 
${}$\\

${}$\\
{\bf Remarque :}(Sur les éléments de ${\bf H}^{1}({\overline X},[{\mathcal O}_{\overline X} \stackrel{d}{\rightarrow} r_{*} {\Omega}_{{\overline X}'}^{1}(\log S)])$ qui ne sont pas dans\\
 $F^{1}{\bf H}^{1}({\overline X},[{\mathcal O}_{\overline X} \stackrel{d}{\rightarrow} r_{*} {\Omega}_{{\overline X}'}^{1}(\log S)])$).\\
D'après le travail fait sur ${\bf P}^{1}$ avec deux points identifiés, la partie du $H^{1}(X,{\bf C})$ ne provenant pas des formes logaritmiques est donnée par :
\begin{center}
$<(\frac{u-u_{P_1}}{u_{Q_1}-u_{P_1}})...(\frac{u-u_{P_n}}{u_{Q_1}-u_{P_n}})-
(\frac{v-v_{P_1}}{v_{Q_1}-v_{P_1}})...(\frac{v-v_{P_n}}{v_{Q_1}-v_{P_n}}),....,(\frac{u-u_{P_1}}{u_{Q_n}-u_{P_1}})...
(\frac{u-u_{P_n}}{u_{Q_n}-u_{P_n}})-
(\frac{v-v_{P_1}}{v_{Q_n}-v_{P_1}})...(\frac{v-v_{P_n}}{v_{Q_n}-v_{P_n}})>$.
\end{center}
Ce sont les \'el\'ements de $H^{1}({ X},{\bf C})$ qui proviennent de $H^{1}({\overline X},{\bf C})$, cr\'e\'es par les boucles
sur ${\bf P}^1$ venant de l'identification des points $P_i$ et $Q_i$. Le reste de la cohomologie est donn\'e par la non
compl\`etude de $X$ i.e. par les \'el\`ements de $H^{1}({ X'},{\bf C})$, exhibés plus haut.
\\

%-----------------------------------------------------------------------
%-------------------------------------------------------------------------

\subsubsection{Les courbes de genre $1$}

Soit $X$ une courbe alg\'ebrique de type finie sur ${\bf C}$ et de genre arithmétique $1$ avec $m$ points enlev\'es $p_{1},...,p_{m}$ et $2n$ points identifi\'es deux à deux $P_{1},Q_{1},...,P_{n},Q_{n}$. $\overline X'$ est isomorphe au quotient de ${\bf C}$ muni de la coordonn\'ee $z$ par le r\'eseau ${\Lambda}_{\bf Z}= {\bf Z}+{\bf Z} \tau$ où $\tau \in {\bf C}$ et $Im(\tau )>0$. Cherchons les g\'en\'erateurs de $F^{1}H^{1}(X,{\bf C})$. D'apr\`es \cite{del3} ce sont les \'el\'ements qui proviennent de $ H^{1}(X',{\bf C})$. \\
\begin{proposition}\cite{mum} La fonction
$\Psi(z)= \sum_{i=1}^{i=k-1} \,\lambda_{i} \, \frac{d}{dz} \log( \theta(z-a_{i}))+cste$ avec $\sum_{i=1}^{i=k-1} \lambda_i=1$ est p\'eriodique pour ${\Lambda}_{\bf Z}$ avec des p\^oles simples aux points $a_i+ \frac{1}{2} (1+ \tau )$ et r\'esidus $\lambda_{i}$, o\`u $\theta $ est la fonction theta sur la courbe elliptique donn\'ee par ${\Lambda}_{\bf Z}$, $\theta(z)= \sum_{n \in {\bf Z}} \exp( \pi \sqrt{-1}
n^{2} + 2 \pi \sqrt{-1} n \tau )$.\\
\end{proposition}

On peut donc prendre comme formes g\'en\'eratrices de $F^{1}H^{1}(S,{\bf C})$ :
$$
\left\lbrace \begin{array}{l}
\omega_{0}=dz \text{ et}\\
\omega_{i}=d(\log( \theta(z-p_{i}-\frac{1}{2} (1+ \tau )))-\log( \theta(z-p_{i+1}-\frac{1}{2} (1+ \tau )))) \text{ o\`u}  i \in [1,m-1]
\end{array}
\right.
$$

Soient $\alpha_{0},\alpha_{1},...,\alpha_{m},\beta_{1},...,\beta_{n}$ des g\'en\'erateurs de $H_{1}(X, {\bf C})$ tels que $\alpha_{j}$ pour $ j \in [1,m]$ soit homologue \`a z\'ero dans $ X \cup p_{j}$ et $\beta_{j}$ est un \'el\'ement representant en homologie la boucle dans $X'$ obtenue par le recollement $ P_{j}=Q_{j}$ (les éléments $\beta_{j}$ sont d\'efinis modulo les $\alpha_{j}$).\\
Fitrons ${\bf C}^{m+n+1}$ par $F^{1}$, cela donne $m$ vecteurs :\\
$(<H_{1},\omega_{i}>)=(<\alpha_{0},\omega_{i}>,<\alpha_{1},\omega_{i}>,...,<\alpha_{m},\omega_{i}>,<\beta_{1},\omega_{i}>,...,<\beta_{n},\omega_{i}>)$ pour $i \in [1,m]$. Les $<\alpha_{j},\omega_{i}>$ sont donn\'es par les r\'esidus des g\'en\'erateurs de $F^{1}H^{1}$ et les $<\beta_{j},\omega_{i}>$ sont donn\'es par int\'egration le long des boucles $\beta_{j}$ i.e. sur $X'$, par l'int\'egration de la forme $\omega_{i}$ de $P_j$ \`a $Q_j$ :\\

\begin{center}
$<\beta_{j},\omega_{i}>= \int_{P_j}^{Q_j} \, \omega_{i}=\int_{P_j}^{Q_j} \,d(\log( \theta(z-p_{i}-\frac{1}{2} (1+ \tau )))-\log(
\theta(z-p_{i+1}-\frac{1}{2} (1+ \tau )))$\\
i.e. $ <\beta_{j},\omega_{i}>=\log( \frac{\theta(Q_{j}-p_{i}-\frac{1}{2} (1+ \tau ))}{\theta(Q_{j}-p_{i+1}-\frac{1}{2}
(1+\tau ))})-\log( \frac{\theta(P_{j}-p_{i}-\frac{1}{2} (1+ \tau ))}{\theta(P_{j}-p_{i+1}-\frac{1}{2} (1+ \tau ))})$
\end{center}
d'où la matrice $A=(B,C)$ o\`u $B$ est une matrice de dimension $m \times m+1$ et $C$ une matrice $m \times n$.\\
$$A=\left(
    \begin{array}{cccccccccc}
     1  & 0                  & ... & ...              & 0                & \vert   & 0 & ... & 0 \\
     \lambda_{1} & 1    &   -1 & 0             & 0                & \vert   & ... & ... & ...           \\
     ...              & 0                 & 1  & -1               &...               & \vert    & ... & ... & ...           \\
      ...             &...                 &...  &...               &...               & \vert     & ... &
c_{ij}=\log( \frac{\theta(Q_{j}-p_{i}-\frac{1}{2} (1+ \tau ))}{\theta(Q_{j}-p_{i+1}-\frac{1}{2}
(1+\tau ))})-\log( \frac{\theta(P_{j}-p_{i}-\frac{1}{2} (1+ \tau ))}{\theta(P_{j}-p_{i+1}-\frac{1}{2} (1+ \tau ))}) & ...   \\
      ...             &...                 &...  &...               &...               & \vert      & ... & ... & ...          \\
       \lambda_{m-1}              & ...                & 0   & 1  & -1  & \vert         & ... & ... & ...
\end{array}
\right)$$
o\`u $i \in [1,m-1]$ et $j \in [1,n]$ et les $\lambda $ sont dans ${\bf Z}$.\\

Pour calculer le niveau de ${\bf R}$-scindement, il faut donc déterminer si $\log( \frac{\theta(Q_{j}-p_{i}-\frac{1}{2} (1+ \tau
))}{\theta(Q_{j}-p_{i+1}-\frac{1}{2} (1+\tau ))})-\log( \frac{\theta(P_{j}-p_{i}-\frac{1}{2} (1+ \tau
))}{\theta(P_{j}-p_{i+1}-\frac{1}{2} (1+ \tau ))})= \overline{\log( \frac{\theta(Q_{j}-p_{i}-\frac{1}{2} (1+ \tau
))}{\theta(Q_{j}-p_{i+1}-\frac{1}{2} (1+\tau ))})-\log( \frac{\theta(P_{j}-p_{i}-\frac{1}{2} (1+ \tau
))}{\theta(P_{j}-p_{i+1}-\frac{1}{2} (1+ \tau ))})}$ suivant la configuration des points $p_{i},p_{i+1},Q_{j},Q_{j+1}$. D'où :

\begin{proposition} Le niveau de ${\bf R}$-scindement de $X$ est donné par
$
\alpha_{1}(X)=(m-1)-\emph{card} \{ i \in [1,m-1]\, \vert \, \forall j \in [m,m+n-1]\,\,, \vert \,\,$\\
\hspace*{0.5cm}$\log( \frac{\theta(Q_{j}-p_{i}-\frac{1}{2} (1+ \tau
))}{\theta(Q_{j}-p_{i+1}-\frac{1}{2} (1+\tau ))})-\log( \frac{\theta(P_{j}-p_{i}-\frac{1}{2} (1+ \tau
))}{\theta(P_{j}-p_{i+1}-\frac{1}{2} (1+ \tau ))})=$\\
\hspace*{4cm}$\overline{\log( \frac{\theta(Q_{j}-p_{i}-\frac{1}{2} (1+ \tau
))}{\theta(Q_{j}-p_{i+1}-\frac{1}{2} (1+\tau ))})-\log( \frac{\theta(P_{j}-p_{i}-\frac{1}{2} (1+ \tau
))}{\theta(P_{j}-p_{i+1}-\frac{1}{2} (1+ \tau ))})} \}.$

\end{proposition}
${}$\\

\subsection{Traduction des extensions de structures de Hodge mixtes en termes de fibrés de Rees}

\subsubsection{$\text{Ext}^{1}$ des structures de Hodge mixtes}
 
Nous allons étudier les extensions de structures de Hodge mixtes (cf Annexe A pour les définitions et propriétés classiques des extensions de structures de Hodge) dans le language des fibrés de Rees équivariants. 

On se limite dans cette section aux structures de Hodge mixtes de rang $2$ obtenues comme extensions de structures de Hodge de rang $1$ dans la catégorie $\crmhs$.

Considérons une structure de Hodge de rang deux $H$ qui est une extension de deux structures de Hodge de Tate $A=T \langle - p \rangle $ et $B=T \langle - q \rangle $ telles que $p <q$. Vient la suite exacte de structures de Hodge mixtes
$$ 0 \rightarrow A \rightarrow H \rightarrow B \rightarrow 0.$$
Cette suite exacte se traduit par le théorème \ref{equivcatmhs} en la suite exacte dans la catégorie\\ $ {\calf ib}_{{\bf P}^{1}_{0}-semistables,\mu=0}({\bf P}^{2}/{\bf T}^{\tau})$

$$  0 \rightarrow \xi_{A} \rightarrow \xi_{H} \rightarrow \xi_{B} \rightarrow 0.$$
Cette suite exacte dans la catégorie des faisceaux se transforme en couple de suite exactes
$$\left\lbrace \begin{array}{l}
0 \rightarrow \xi_{A} \rightarrow \xi_{H} \rightarrow \text{Coker} \rightarrow 0,\\    
{}\\
0 \rightarrow \text{Coker} \rightarrow \xi_{B} \rightarrow \calo_{T} \rightarrow 0,
\end{array}
\right.$$
où $T$ est un sous-schéma de codimension $2$ supporté au point $(0:0:1)$. Soit $\cali_{T}$ le faisceau d'idéaux qui correspond à $T$, on a alors la suite exacte 
$$  0 \rightarrow \xi_{B}\otimes\cali_{T} \rightarrow \xi_{B} \rightarrow \calo_{T} \rightarrow 0.$$
Ceci permet d'avoir à considérer les extensions (équivariantes) de faisceaux de rang $1$ dans la catégorie des faisceaux cohérents de la forme
$$  0 \rightarrow \xi_{A} \rightarrow \xi_{H} \rightarrow \xi_{B} \otimes \cali_{T} \rightarrow 0.$$

\begin{center}
{\bf Etude $\text{Ext}^{1}(L'\otimes \cali_{T},L)$}
\end{center}

Ici $L$ et $L'$ sont deux fibrés en droite sur une surface lisse projective $X$ et $T$ un schéma de codimension $2$.

L'étude de ces extensions est faite dans le chapitre $2$ sur les faisceaux cohérents de \cite{fri}. Pour calculer les extensions globales, on passe par une suite spectrale dont le deuxième terme est
$$E^{p,q}_{2}=H^{p}(X,\cale xt^{q}(L'\otimes \cali_{T},L)) \Rightarrow \text{Ext}^{p+q}( L'\otimes \cali_{T},L)).$$
\begin{lemme}(Lemme 7, p.36, \cite{fri}) Soit $R$ un anneau local régulier et $I=(f,g)R$ un idéal de $R$ engendré par deux éléments premiers entre eux, alors :\\
$(i)$ $\emph{Hom}_{R}(I,R) \cong R $ et l'isomorphisme est engendré par l'application naturelle de restriction $\emph{Hom}_{R}(R,R) \rightarrow \emph{Hom}_{R}(I,R)$,\\
$(ii)$ $\emph{Hom}_{R}(I,I)\cong R $ et l'isomorphisme est encore induit par l'application de restriction,\\
$(iii)$ $\emph{Ext}^{1}_{R}(I,R) \cong \emph{Ext}^{2}_{R}(R/I,R) \cong R/I$,\\
$(iv)$ $\emph{Ext}^{k}_{R}(I,R)=0$ pour $k \geq 2$.
\end{lemme} 

La suite spectrale se ramène par le lemme à une suite exacte longue
$$ 0 \rightarrow H^{1}(L'{}^{-1}\otimes L) \rightarrow \text{Ext}^{1}(L'\otimes \cali_{T},L) \rightarrow
H^{0}(\cale xt^{1}(L'\otimes \cali_{T},L)) \rightarrow H^{2}(L'{}^{-1}\otimes L)\rightarrow ...$$
Dans le cas d'une surface $\cale xt^{1}(L'\otimes \cali_{T},L)=\calo_{T}$.   
 
Reste à voir si l'extension est juste un faisceau cohérent ou bien un faisceau localement libre... 

\begin{center}
{\bf Extensions équivariantes}
\end{center}

Pour calculer les $\text{Ext}^{p}_{G}(\calf,\calg)$ dans la catégorie abélienne des faisceaux équivariants on utilise la suite spectrale 
$$ E_{2}^{p,q}=H^{p}(G,\text{Ext}^{q}(\calf, \calg)) \Rightarrow \text{Ext}^{p+q}_{G}(\calf,\calg).$$ 
Lorsque $G$ est un groupe réductif, il n'a pas de cohomologie supérieure et alors les groupes d' extensions équivariantes sont donnés par les extensions invariantes par l' action du groupe, i.e.
$$\text{Ext}^{n}_{G}(\calf,\calg)\cong \text{Ext}^{n}(\calf,\calg)^{G}.$$
Nous utilisons ce résultat dans le cas de l'action du groupe réductif $G$ qui est le tore qui agit sur le plan projectif complexe.\\

\subsubsection{$\text{Ext}^{p}$ des structures de Hodge mixtes, $p>1$}
Les extensions supérieures de structures de Hodge mixtes sont étudiées dans \cite{carhai} à l'aide de resolutions projectives relatives.\\

Nous donnons une démonstration géométrique du fait suivant

\begin{proposition}\cite{carhai}
Soient $A $ et $B$ deux structures de Hodge mixtes, alors, pour tout $p>1$,
$$\text{Ext}^{p}_{MHS}(B,A)=0.$$
\end{proposition}

\begin{preuve}
Le lemme précédent $(iii)$ montre que tous les $\cale xt^p$ locaux sont nuls pour $p>1$, il en va donc de même pour les $\text{Ext}^p$ globaux. On conclut par
$$ \text{Ext}^{p}_{MHS}(B,A)\cong \text{Ext}^{p}_{G}(\xi_{B},\xi_{A})\cong \text{Ext}^{p}(\xi_{B},\xi_{A})^{G}=0.$$

\end{preuve}

%---------------------------------------------------------------------------------------------------------------

%-----------------------------------------------------------------------------------------------------------------

%--------------------------------------------------------------------------------------------------------------

%--------------------------------------------------------------------------------------------------------------

\newpage

\section{Stratifications et variations de structures de Hodge, perspectives}

\subsection{Variations de structures de Hodge mixtes}

\subsubsection{Rappel : variations de structures de Hodge}
 Soit $\calv_{\bf C}$ un système local complexe sur une variété $S$. Ici $S$ désignera un schéma lisse projectif sur ${\bf C}$ ou respectivement une variété complexe lisse. Alors le fibré vectoriel algébrique (resp. holomorphe) $\calv:= \calo_{X} \otimes_{\bf C} \calv_{\bf C}$ associé à ce système est muni d'une connexion intégrable $\nabla$ dont l'espace des sections horizontales est $\calv_{\bf C}$. 

La définition qui suit peut-être donnée en remplaçant ${\bf R}$ par tout sous-anneau $A$ de ${\bf R}$, mais nous n'avons pas besoin de ce contexte plus large car toutes les constructions de l'étude précédente ne concernent que les structures sur ${\bf R}$ et ${\bf C}$.

\begin{definition} Une variation polarisée de structures de Hodge de poids $n$ sur $S$ consiste en la donnée d'un système local de ${\bf R}$-vectoriels de rang fini $\calv_{\bf R}$ ainsi que \\
\hspace*{1cm}$(i)$ d'une filtration décroissante $...\, \calf^{p+1} \subset \calf^{p}...$ du fibré vectoriel $\calv$ par des sous-fibrés vectoriels,\\
\hspace*{1cm}$(ii)$ une forme bilinéaire plate $Q: \calv_{\bf R} \otimes \calv_{\bf R} \rightarrow {\bf R}(-n)$ tel que :\\
\hspace*{2cm}$(a)$ la transversalité de Griffith soit satisfaite, i.e.
$$\nabla(\calf^{p}) \subset \Omega_{S}^{1}\otimes_{\calo_{S}} \calf^{p-1},$$
\hspace*{2cm}$(b)$ pour tout $s \in S$, le ${\bf R}$-module $\calv_{\bf R}{}_{s}$ muni de la filtration $\calf^{\bullet }_{s}$ et de la forme bilinéaire $Q_{s}$ sont la donnée d'une structure de Hodge polarisée de poids $n$.
\end{definition}    
  La situation typique dans laquelle arrivent les variations de structures de Hodge polarisées est lorsque est donné un morphisme lisse et projectif $f:X \rightarrow S$ entre deux variétés quasi-projectives sur ${\bf C}$, le système local étant donné par $\calv_{\bf R}:=R^{n}f_{*}{\bf R}$, cf \cite{gri1}.  

\subsubsection{Définition des variations de structures de Hodge mixtes}
 
Dans la situation géométrique décrite ci-dessus, si l'on enlève la condition de lissité et/ou de propreté du morphisme $f$, alors sur un ouvert dense, la famille paramétrée par $f$ donne lieu à une famille de structures de Hodge mixtes. Pour tout $s$ dans l'ouvert dense, les groupes de cohomologie de la fibre $f^{-1}(s)$ sont équipés de structures de Hodge mixtes. Ceci mène à la notion de variation graduellement polarisée de structures de Hodge mixtes.

Comme pour les variations de structures de Hodge, on se limite ici à définir les variations de structures de Hodge mixtes à système local sous-jacent réel, mais la définition peut être donnée pour tout système local de $A$-module avec $A$ sous-anneau de ${\bf R}$.
\begin{definition}\cite{stezuc}
 Une variation graduellement polarisée de structures de Hodge mixtes sur une variété $S$ consiste en la donnée d'un système local de ${\bf R}$-vectoriels de rang fini $\calv_{\bf R}$ ainsi que \\
\hspace*{1cm}$(i)$ d'une filtration croissante $...\calw^{p} \subset \calw^{p+1}...$ du système local $\calv_{\bf R}$,\\
\hspace*{1cm}$(ii)$ d'une filtration décroissante $...\,\calf^{p+1} \subset \calf^{p}...$ du fibré vectoriel $\calv=\calv_{\bf R} \otimes_{\calo_{S}} \calo_{S}$ par des sous-fibrés vectoriels,\\
\hspace*{1cm}$(iii)$  une suite de forme bilinéaires plates,

$$Q_{k}:Gr^{\calw}_{k}(\calv_{\bf R})\otimes Gr^{\calw}_{k}(\calv_{\bf R}) \rightarrow {\bf R}(-k),$$ tels que :\\
\hspace*{2cm}$(a)$ la transversalité de Griffith soit satisfaite, 
$$\nabla(\calf^{p}) \subset \Omega_{S}^{1}\otimes_{\calo_{S}} \calf^{p-1},$$
\hspace*{2cm}$(b)$ la donnée de $Gr^{\calw}_{k}$, de la filtration induite par $\calf^{\bullet}$ sur $\calo_{S}\otimes Gr^{\calw}_{k}$ et de la polarisation $Q_{k}$ soit celle d'une variation de structures de Hodge de poids $k$ .
\end{definition}

\subsection{Etude des stratifications associées aux variations de structures de Hodge mixtes : notions préliminaires}
Les premières sections sont préparatoires, on introduit d'abord les fibrés universels sur les variétés de drapeaux, on rappelle quelques notions d'amplitude qui seront nécessaires pour appliquer les théorèmes sur les dégénérescences de fibrés vectoriels.  

\subsubsection{Fibrés universels, variétés des drapeaux généralisés, amplitude }

\begin{center}
{\bf Variétés des drapeaux}
\end{center}

Soit $V$ un espace vectoriel complexe de dimension finie muni d'une filtration $F^{\bullet }$ exhaustive et décroissante :$$  \{ 0 \} =F^{n+1} \subset F^{n} \subset ...\subset F^{0}=V.$$

Nous construisons ici une variété qui paramétrise toutes les filtrations de cette forme sur $V$. Notons $f^{p}=dim_{\bf C}F^{p}V$. Il s'agit donc de trouver une variété dont les points paramétrisent l'ensemble des filtrations exhaustives et décroissantes de $V$, $F^{\bullet}$ telles que  $f^{p}=dim_{\bf C}F^{p}V$.\\

\begin{proposition}
Il existe une variété algébrique lisse compacte ${\bf Flag}(f^{0},f^{1},...,f^{n},V)$ qui paramétrise toutes les filtrations exhaustives et décroissantes, on l'appelle variété des drapeaux associée à $(V,f^{0},f^{1},...,f^{n})$.\\ 
\end{proposition}  
${}$
{\bf Remarque :}
Dans l'écriture de ${\bf Flag}(f^{0},f^{1},...,f^{n},V)$, on peut omettre $f^{0}=\text{dim}_{\bf C}V$ qui est implicite comme dimension de $V$.\\
${}$\\
{\bf Remarque :}
La variété des drapeaux décrite ici est une généralisation de ce que l'on entend par variété des drapeaux habituellement où les drapeaux sont pris comme étant les suites de sous-espaces vectoriels $  \{ 0 \} =F_{0} \subset F^{1} \subset ...\subset F_{n}=V$ avec $dim_{\bf C}F_{i}=i$, i.e. des drapeaux complets.\\

Avant de donner la preuve de la proposition, pour fixer les notations, nous allons rappeler la construction des grassmaniennes (\cite{adagri}, \cite{hir} par exemple).\\

Définissons la grassmanienne des sous-espaces vectoriels complexes de dimension $k$ de l'espace vectoriel complexe $V$ de dimension $m$, $grass(k,m)$ (que l'on pourra aussi noter $grass(k,V)$ suivant le contexte, s'il est évident ou non que l'on regarde des sous-espaces vectoriels de $V$).\\
Soient $(z_{1},...,z_{m})$ les coordonnées complexes standard()s sur $V={\bf C}^{m}$. Et soit $L_{k}$ le sous-espace vectoriel défini par l'équation $z_{k+1}=...=z_{m}=0$. Le groupe $GL(m,{\bf C})$ agit transitivement sur les espaces vectoriels de dimension $k$. Le sous-groupe $GL(k,m-k,{\bf C})$ des matrices inversibles qui laisse $L_k$ invariant est constitué des éléments de la forme suivante :
$$  A \in GL(k,m-k,{\bf C}) \text{\,si et seulement si \,}
A=\left(
\begin{array}{cc}
           A'    &     B       \\
           0         &   A''
\end{array}
\right)
$$
où $A' \in GL(k,{\bf C})$, $A'' \in GL(m-k,{\bf C})$ et $B$ est une matrice $k \times m-k$ quelconque. La grassmanienne $grass(k,V)$ peut donc \^etre identifiée à $GL(m,{\bf C})/GL(k,m-k,{\bf C})$ et est donc une variété algébrique lisse sur ${\bf C}$. Elle ainsi naturellement munie de l'action du groupe algébrique d'automorphismes $GL(m,{\bf C})$. $GL(k,m-k,{\bf C})$ est de codimension $k\times(m-k)$ dans $GL(m,{\bf C})$ donc $grass(k,V)$ est de dimension $k\times(m-k)$.\\

Munissons $V={\bf C}^m$ du produit scalaire hermitien standard. Le groupe unitaire standard associé $U(m,{\bf C})$ agit également de façon transitive sur les sous-espaces vectoriels de dimension $k$. Le sous-groupe du groupe unitaire qui fixe $z_{k+1}=...=z_{m}=0$ sera noté $U(k,m-k,{\bf C})$. Il est naturellement isomorphe à $U(k,{\bf C}) \times U(m-k,{\bf C})$. On a donc : 
$$grass(k,m)= U(m,{\bf C}) /(U(k,{\bf C}) \times U(m-k,{\bf C})).$$

\begin{preuve}
Nous allons construire ${\bf Flag}(f^{0},f^{1},...,f^{n},V)$ comme sous-variété de $ \prod_{p=1}^{n} grass(f^{p},V)$. Notons, pour $x \in grass(k,m)$, par $<x>$ le sous-espace vectoriel de $V$ engendré par $x$. Ensemblistement, ${\bf Flag}(f^{0},f^{1},...,f^{n},V)$ est l'ensemble des $n$-uplets $(x_{1},...,x_{n}) \in  \prod_{p=1}^{n} grass(f^{p},V)$ tels que $<x_{i+1}> \subset <x_{i}>$ pour $i \in [1,n-1]$. Rappelons que d'après la description des grassmaniennes, si l'on note par $m$ la dimension de $V$ :
$$ \prod_{p=1}^{n} grass(f^{p},V) = \prod_{p=1}^{n}  U(f^{0},{\bf C}) /(U(f^{p},{\bf C}) \times U(f^{0}-f^{p},{\bf C}))$$
Considérons maintenant la suite de sous-espaces vectoriels emboités de $V$ $\{ L_{f^{i}} \}_{i \in [0,f^{n+1}]} $ :
$$   \{ 0 \}=L_{f^{n+1}} \subset L_{f^{n}} \subset ... \subset L_{f^{i}} \subset L_{f^{i-1}} \subset ...\subset L_{f^{1}} \subset L_{f^{0}}=V$$
où $L_{f^{i}}$ est défini par l'équation $z_{f^{i}+1}=z_{f^{i}+2}=...=z_{f^{0}}=0$ dans $V={\bf C}^{f^{0}}$ muni des coordonnées $(z_{1},...,z_{f^{0}})$. C'est un drapeau de la forme définie par $(V,f^{0},f^{1},...,f^{n})$.
Le groupe $U(f^{0},{\bf C})$ agit transitivement sur l'ensemble des drapeaux définis de la sorte. Notons $U(f^{n},f^{n-1}-f^{n},f^{n-2}-f^{n-1},...,f^{1}-f^{2},f^{0}-f^{1},{\bf C})$ le sous-groupe qui laisse le drapeau standard fixe. Dans la base choisie, les éléments de $U(f^{n},f^{n-1}-f^{n},f^{n-2}-f^{n-1},...,f^{1}-f^{2},f^{0}-f^{1},{\bf C})$ sont de la forme :  
$$\left(
\begin{array}{cccccc}
           A_{n}    &     *   &     *   &  ...  &  * &    *    \\
          0    &     A_{n-1}  &     *    &    ...   &   *   &   *   \\
           0    &  0           &         A_{n-2}  &   ... &  *     &   *  \\
          ...    &  ...    &  ... &   ... &... & ...\\       
          0 & 0 & 0 & ... & A_{1} & *\\
          0 & 0 & 0 & ... & 0 & A_{0}
\end{array}
\right)
$$
où $A_{i} \in U(f^{i}-f^{i+1},,{\bf C})$. Comme c'est un sous-groupe du groupe unitaire, on a aussi la relation $*=0$. Donc $U(f^{n},f^{n-1}-f^{n},f^{n-2}-f^{n-1},...,f^{1}-f^{2},f^{0}-f^{1},{\bf C})\cong U(f^{n},{\bf C}) \times U(f^{n-1}-f^{n},{\bf C})\times ... \times U(f^{0}-f^{1},{\bf C})$ où le morphisme injectif $i: U(f^{i}-f^{i+1},{\bf C}) \rightarrow 
U(f^{0},{\bf C})$ est donné par $A \mapsto Diag(I_{f^{n}},...,I_{f^{i+1}-f^{i+2}},A,I_{f^{i-1}-f^{i}},...,I_{f^{0}-f^{1}})$. Ainsi,
$${\bf Flag}(f^{0},f^{1},...,f^{n},V)=\frac{U(f^{0},{\bf C})}{U(f^{n},{\bf C}) \times U(f^{n-1}-f^{n},{\bf C})\times ... \times U(f^{0}-f^{1},{\bf C})}$$
est donc une variété algébrique lisse compacte et est munie de l'action transitive du groupe algébrique $U(f^{0},{\bf C})$.\\

\end{preuve}
 
On notera l'immersion dans le produit des grassmaniennes par :
$$i_{\bf F}: {\bf Flag}(f^{0},f^{1},...,f^{n},V) \hookrightarrow \prod_{p=1}^{n} grass(f^{p},V) .$$
{\bf Remarque :} Si l'on ne regarde que les drapeaux au sens habituel du terme (cf remarque suivant la proposition précédente), alors 
$U(f^{n},{\bf C}) \times U(f^{n-1}-f^{n},{\bf C})\times ... \times U(f^{0}-f^{1},{\bf C})$ est le sous-groupe des matrices unitaires diagonales ${\bf T}^{f^{0}}={\bf \Delta}(f^{0},{\bf C}) \cap U(f^{0},{\bf C})$ i.e. un tore complexe de dimension $f^0$.\\

\begin{center}
{\bf Fibrés universels}
\end{center}
Pour tout $k$ les grassmaniennes $grass(k, f^{0})$ sont munies de fibré universels $\calu_{k}$, sous-fibrés du fibré trivial $grass(k, f^{0}) \times V$ :

%$$
%\xymatrix{
% 0 \ar[r] & A          \ar[r]^{i} \ar[d]_{id_{A}} & H         \ar[r]^{\pi} \ar[d]_{(r,\pi)} \ar[dl]_{r}     &   B    \ar[r] \ar[d]_{id_{B}}  & 0  \\
%          0 \ar[r] & A         \ar[r]^{(1,0)} &  L         \ar[r]^{(0,1)}  &      B  \ar[r]   &0
%}
%$$

$$
\xymatrix{
\calu_{k} \lhook\rightarrow \ar[dr] &     grass(k,f^{0}) \times V \ar[d]\\
                      &  grass(k,f^{0})
}
$$   
La fibre au dessus de $x \in grass(k,f^{0})$ est le sous-espace vectoriel que représente $<x>\subset V$ de la fibre du fibré trivial en $x$.\\

Notons, pour $q \in [1,n]$, par $\pi_{q}$ la projection :
$$ \pi_{q}: \prod_{p=1}^{n} grass(f^{p},V) \rightarrow  grass(f^{q},V)$$
et $\pi_{q}^{*}\calu_{k}$ le ''pull-back'' du fibré universel le $q$-ième terme du produit des grassmaniennes.\\

\begin{proposition}  
Le fibré trivial de rang $f^{0}$ sur $ {\bf Flag}(f^{0},f^{1},...,f^{n},V) \times V$ est muni d'une filtration $\calf^{\bullet}_{univ}(f^{0},f^{1},...,f^{n},V)$ exhaustive et décroissante de sous-fibrés ''universels'' tel que pour tout $x \in 
 {\bf Flag}(f^{0},f^{1},...,f^{n},V)$ la filtration induite sur $V$, $(\calf^{\bullet}_{univ}(f^{0},f^{1},...,f^{n},V))_{x}$, par les fibres en $x$ des sous-fibrés soit la filtration paramétrée par le point $x$ dans la variété des drapeaux.\\
\end{proposition}

\begin{preuve}
Considérons pour tout $q \in [1,n]$ les morphismes :
$$\pi_{q}\circ i_{\bf F}: {\bf Flag}(f^{0},f^{1},...,f^{n},V) \hookrightarrow \prod_{p=1}^{n} grass(f^{p},V) \rightarrow  grass(f^{q},V)$$
Et les pull-back de $\calu_{q}$ par ces morphismes :
$$
\xymatrix{
  ( i_{\bf F}^{*}\circ \pi_{q}^{*})\,\calu_{q} \ar[d]\\
  {\bf Flag}(f^{0},f^{1},...,f^{n},V) 
}
$$
Ce sont des sous-fibrés du fibré trivial $ {\bf Flag}(f^{0},f^{1},...,f^{n},V) \times V$ tels que la fibre au dessus d'un point $x \in {\bf Flag}(f^{0},f^{1},...,f^{n},V) $ soit le sous-espace $<\pi_{q}\circ i_{\bf F}(x)> \subset V$ paramétrisé par $\pi_{q}\circ i_{\bf F}(x)$ dans la grassmanienne $grass(f^{p},V)$. Posons pour tout $q \in [1,n]$ :
$$\calf^{q}_{univ}(f^{0},f^{1},...,f^{n},V)= (i_{\bf F}^{*}\circ \pi_{q}^{*})\,\calu_{q}$$
et on a la filtration de $ {\bf Flag}(f^{0},f^{1},...,f^{n},V) \times V$ voulue.
\end{preuve}

\begin{center}
{\bf Amplitude}
\end{center}

Soit $X$ une variété projective. Rappelons qu'un fibré vetcoriel $\mathcal L$ sur $X$ est ample s'il exite une entier positif $m$ tel que ${\mathcal L}^{\otimes m}$ soit la restriction à $X$ du fibré en droite ${\mathcal O}(1)$ pour un plongement de $X$ dans un espace projectif. 
Notons par $E^{*}$ le fibré dual de $E$, $\text{Hom}(E,\calo_{X})$.\\

\begin{lemme} Les fibrés $((i_{\bf F}^{*}\circ \pi_{q}^{*})\,\calu_{q})^{*}$ sur ${\bf Flag}(f^{0},f^{1},...,f^{n},V)$ sont amples.\\
\end{lemme}

\begin{preuve}
cf \cite{gri2}
\end{preuve}

\subsubsection{Lieu de dégénérescence des morphismes de fibrés vectoriels}

\begin{center}
{\bf Introduction}
\end{center}
 Considérons deux familles $\cale$ et $\calf$ de sous-espaces vectoriels d'un m\^eme espace vectoriel $V$ paramétrées par une variété $X$. En tout point $x \in X$ on a deux espaces vectoriels $\cale_{x}$ et $\calf_{x}$. Le but de cette section est d'étudier les variations de $\text{dim}_{\bf C}\,\cale_{x}\cap \calf_{x}$ sur $X$dans le cas où les deux familles sont données par des sous-fibrés vectoriels d'un m\^eme fibré vectoriel.

Soient $\phi : \cale \rightarrow \calf $ un morphisme de fibrés vectoriels complexes de rang $m$ et $n$ sur une variété complexe lisse. Il existe un recouvrement de $X$ par des ouverts sur lesquels les fibrés $\cale$ et $\calf$ sont triviaux. Soit $U$ un tel ouvert. Sur $U$, $\phi $ est donnée par une matrice $m \times n$ de fonctions à valeurs complexes. On peut donc définir le rang de $\phi$ sur $U$. Ce rang ne depend pas des ouverts de trivialisations choisis. On peut donc poser :
$$D_{r}(\phi)=\{ x \in X \vert \text{rang}(\phi(x)) \leq r \}.$$     
C'est une sous-variété de $X$. On obtient donc une suite de sous-variétés de $X$ :
$$ \emptyset \subset D_{0}(\phi) \subset D_{1}(\phi) \subset ... \subset D_{\text{min}(m,n)}(\phi)=X.$$
Si $\phi$ est suffisament générique, pour tout $r$ tel que $0 \leq r \leq \text{min}(m,n)$, $D_{r}(\phi)$ a pour codimension $(m-r).(n-r)$.
Considérons maintenant deux sous-fibrés $\calg$ et $\calh$ d'un m\^eme fibré $\calv$. En tout point $x \in X$ on a deux sous-espaces vectoriels $\calg_{x}$ et $\calh_{x}$ de $\calv_{x}$. Nous voulons calculer $\text{rang}(\phi (x))$ en fonction de $x$.
Notons $i : \calg \rightarrow \calv $ l'inclusion de $\calg$ dans $\calv$ et $\pi : \calv \rightarrow \calv/\calh$ la projection. Considérons le morphisme de fibrés composition des deux morphismes :
$$ \pi \circ i : \calg \longrightarrow \calv \longrightarrow \calv/\calh.$$
En tout point $x \in X$ on a un morphisme d'espaces vectoriels :
$$ (\pi \circ i)_{x} : \calg_{x} \longrightarrow \calv_{x} \longrightarrow \calv_{x}/\calh_{x}.$$ 

On a $\text{rang}((\pi \circ i)_{x})=\text{dim}_{\bf C}{\calg_{x}}-\text{dim}_{\bf C}\calg_{x}\cap \calh_{x}$, i.e. $\text{rang}((\pi \circ i)_{x})=\text{rang}(\calg)-\text{dim}_{\bf C}\calg_{x}\cap \calh_{x}$. On est dans la situation décrite plus haut avec $\cale=\calg$, $\calf=\calv/\calh$ et $\phi = \pi \circ i$. On peut donc définir les lieux :
\begin{eqnarray*}
D_{r}(\calg,\calh)&=&\{ x \in X \vert \text{dim}_{\bf C}(\calg_{x}\cap \calh_{x}) \geq r \}\\
 &=&\{x \in X \vert \text{rang}((\pi \circ i)_{x}) \leq \text{rang}(\calg-r \}\\
 &=&\{x \in X \vert \text{rang}(\phi)_{x}) \leq \text{rang}(\calg)-r \}\\
 &=&D_{\text{rang}_{\bf C}{\calg}-r}(\phi).
\end{eqnarray*}
Pour tout $r \in [\text{max}(0,\text{rang}(\calg)+\text{rang}(\calh)-\text{rang}(\calv)), \text{min}(\text{rang}(\calg),\text{rang}(\calh))]$. 
D'où la suite de sous-variétés de $X$ :
$$ X = D_{\text{max}(0,\text{rang}(\calg)+\text{rang}(\calh)-\text{rang}(\calv))}(\phi) \supset  ... \supset D_{\text{min}(\text{rang}(\calg),\text{rang}(\calh))}(\phi) \supset \emptyset.$$

\begin{center}
{\bf Dégénérescence sur le produit de variétés des drapeaux}
\end{center}
Nous allons appliquer le travail de la section précédente aux fibrés tautologiques sur les variétés des drapeaux ${\bf Flag}(f^{0},f^{1},...,f^{n},V)$. Pour $p \in [0,n]$, désignons par $\calf^{p}$ le fibré $i_{\bf F}^{*}\circ \pi_{p}^{*}\,\calu_{p}$ sur ${\bf Flag}(f^{0},f^{1},...,f^{n},V)$. Notons par $\pi_{i}$ pour $i \in \{1,2 \}$ les projections respectives de ${\bf Flag}(f^{0},f^{1},...,f^{n},V) \times {\bf Flag}(f^{0},f^{1},...,f^{n},V)$ sur le premier et deuxième facteur. On dispose donc, pour tous $(p,q) \in [0,n]^2$ de fibrés "tautologiques" $\pi_{1}^{*}(\calf^{p})$ et $\pi_{2}^{*}(\calf^{q})$, sous-fibrés vectoriels du fibré trivial $\calv ={\bf Flag}(f^{0},f^{1},...,f^{n},V) \times {\bf Flag}(f^{0},f^{1},...,f^{n},V) \times V$. Au dessus d'un point $(x,y)$, la fibre de $\pi_{1}^{*}(\calf^{p})$ est le $p$-ième sous-espace vectoriel de $V$ du drapeaux paramétré par $x \in  {\bf Flag}(f^{0},f^{1},...,f^{n},V)$. On a les morphismes suivants : 

$$
\xymatrix{
\calv  & \pi_{2}^{*}(\calf^{q}) \ar@{^{(}->}[l] \ar[r]^{\pi_{2}^{*}} \ar[d] & \calf^{q} \ar[d]\\
\pi_{1}^{*}(\calf^{p}) \ar@{^{(}->}[u] \ar[d]^{\pi_{1}^{*}} \ar[r] &   {\bf Flag}(f^{0},f^{1},...,f^{n},V) \times {\bf Flag}(f^{0},f^{1},...,f^{n},V) \ar[r]^-{\pi_{2}} \ar[d]^{\pi_{1}} &  {\bf Flag}(f^{0},f^{1},...,f^{n},V)\\
\calf^{p} \ar[r] &  {\bf Flag}(f^{0},f^{1},...,f^{n},V) &   
}
.$$
Considérons le morphisme de fibrés vectoriels :
$$\phi^{p,q} : \pi_{1}^{*}(\calf^{p}) \longrightarrow \calv \longrightarrow \calv /\pi_{2}^{*}(\calf^{q}),$$
et les applications pour tous $(p,q) \in [0,n]^{2}$ :
$$
\xymatrix{
\varphi^{p,q} :& {\bf Flag}(f^{0},f^{1},...,f^{n},V) \times {\bf Flag}(f^{0},f^{1},...,f^{n},V) \ar[r] & {\bf Z}\\
 & (x,y) \ar@{|->}[r] & \text{rang} \phi^{p,q}_{(x,y)}
}  
.$$
En tout point $(x,y)$, $\varphi^{p,q}_{(x,y)}=\text{dim}_{\bf C} \calf^{p}_{x} \cap \calf^{q}_{y}$.\\

\begin{lemme}\label{flagsci}
Pour tous $(p,q) \in [0,n]^2$ $\varphi^{p,q}$ est semi-continue supérieurement sur le produit des variétés de drapeaux ${\bf Flag}(f^{0},f^{1},...,f^{n},V) \times {\bf Flag}(f^{0},f^{1},...,f^{n},V)$ et prend toutes les valeurs de $[\text{max}(0,f^{p}+f^{q}-f^{0}),\text{min}(f^{p},f^{q})]$.\\
\end{lemme}

\begin{preuve}
Pour la première assertion il suffit de prouver que l'application du produit des Grassmaniennes $grass(f^{p},V) \times grass(f^{q},V)$ dans ${\bf Z}$ donnée par :
$$ (G,H) \longmapsto \text{dim}_{\bf C}\, G \cap H $$
est semi-continue supérieurement.\\ 

Donnons la démonstration suivant \cite{catdelkap}. Soient $G'$ et $H'$ proches de $G$ et $H$. Il existe un automorphisme proche de l'identité $g \in GL(V)$ tel que $g.G'=G$. On peut supposer que $G'=G$. Soit $I_1$ un supplémentaire de $G \cap H$ dans $G$. Soit $I_2$ un supplémentaire de $G+H$ dans $V$. On a ainsi $G=G \cap H \oplus I_1$ et $V=G + H \oplus I_2$. $I_1$ et $I_2$ sont en somme directe. Posons $I=I_{1}\oplus I_{2}$, on a alors $V=G \oplus I$ et $G=G \cap H \oplus G \cap I$. Lorsque $H'$ est assez proche de $H$, $H'$ peut \^etre vu comme le graphe d'une application $h': H \mapsto I$. L'application qui a $H'$ associe l'application $h'$ est continue. Soit $x$ un élément de $H'$. $x$ s'écrit $(u,h'(u))$. $x $ est dans $G$ si et si seulement si $u \in G$ et $h'(u) \in I$ sont dans $G$ ce qui implique que $u \in G \cap H$. En découle l'inégalité sur les dimensions.

\end{preuve}

On peut donc définir pour tous $(p,q) \in [0,n]^2$ conformément à la section précédente les sous-ensembles analytiques du produit des variétés des drapeaux :
$$ D_{r}(\varphi^{p,q})=\{ (x,y) \in {\bf Flag}(f^{0},f^{1},...,f^{n},V) \times {\bf Flag}(f^{0},f^{1},...,f^{n},V) \vert \text{rang}\, \varphi^{p,q}_{(x,y)} \leq r \}$$
pour obtenir une suite de sous-variétés analytiques :
$$  \emptyset \subset ...\subset D_{r}(\varphi^{p,q}) \subset  D_{r+1}(\varphi^{p,q}) \subset ...\subset {\bf Flag}(f^{0},f^{1},...,f^{n},V) \times {\bf Flag}(f^{0},f^{1},...,f^{n},V).$$
Pour des raisons de dimensions si $r < \text{max}(0,f^{p}+f^{q}-f^{0})$ alors $ D_{r}(\varphi^{p,q})=\emptyset$ et si $r > \text{min}(f^{p},f^{q})$ alors $ D_{r}(\varphi^{p,q})=X$. 

Les variétés des drapeaux sont des variétés projectives. On peut en effet plonger la variété des drapeaux ${\bf Flag}(f^{0},f^{1},...,f^{n},V)$ dans ${\bf P}^{N}$ pour $N$ assez grand par l'application de Veronese itérée, de même pour le produit. Par G.A.G.A. les sous-ensembles analytiques $ D_{r}(\varphi^{p,q})$ sont donc des sous-ensembles algébriques. 

Posons pour tous $(p,q) \in [0,n]^2$ et tous $r \in {\bf Z}$ $\cals_{r}(\varphi^{p,q})= D_{r}(\varphi^{p,q})- D_{r-1}(\varphi^{p,q})$. C'est le lieu sur ${\bf Flag}(f^{0},f^{1},...,f^{n},V) \times {\bf Flag}(f^{0},f^{1},...,f^{n},V)$ des points tels que le rang de $\varphi$ soit exactement égal à $r$.\\

\begin{proposition}La ${\bf Z}$-décomposition définie par les strates algébriques $ \cals_{r}(\varphi^{p,q})$ définit une stratification de Whitney du produit des variétés de drapeaux 
$${\bf Flag}(f^{0},f^{1},...,f^{n},V) \times {\bf Flag}(f^{0},f^{1},...,f^{n},V).$$

\end{proposition}

%\textit{Exemple :}Celui de Strasbourg.

%\textit{Exemple :}Les drapeaux complets.

%\textit{Exemple :}Les drapeaux "incomplets"
%\subsubsection{Stratifications sur des quotients}
%\subsubsection{Stratifications sur des sous-variétés}

\subsubsection{Espaces classifiants de structures de Hodge mixtes}
\begin{center}
{\bf Cadre}
\end{center}
Suivant le contexte, on se placera dans les catégories suivantes :
\begin{itemize}\parindent=1.5cm
\item {(i)} $\crzmhs$ la catégorie des structures de Hodge mixtes au sens usuel.
\item {(ii)} $\crmhs$ la catégorie des structures de Hodge mixtes sans strutures entières sous-jacentes.
\item {(iii)} $\cmhs$ la catégorie des structures de Hodge mixtes complexes sans structures réelles sous-jacentes.
\item {(iv)} $\trif$ la catégorie des espaces vectoriels munis de trois filtrations exhaustives et décroissantes .

\par\end{itemize}
${}$\\
On a les foncteurs oublis suivants entre ces différentes catégories :\\

$ \Phi_{z}: \crzmhs \rightarrow \crmhs $ qui est le foncteur oubli de la structure entière sous-jacente à une struture de Hodge mixte.

$ \Phi_{r}: \crmhs \rightarrow \cmhs $ qui est le foncteur oubli de la structure réèlle sous-jacente.

$ \Phi_{str}: \cmhs \rightarrow \trif $ qui à une structure de Hodge mixte $H=(H_{\bf C},F^{\bullet},{\hat F}^{\bullet },W_{\bullet})$ associe l'espace vectoriel complexe trifiltré muni de trois filtrations exhaustives et décroissantes $(H_{\bf C},F^{\bullet},{\hat F}^{\bullet },{W_{.}}^{\bullet})$.
\\
Pour exhiber les constuctions d'espaces classifiants, nous aurons besoin d'équiper les morphismes entre structures de Hodge mixtes. Ceci est possible d'après le théorème de Deligne que nous avons démontrer géométriquement dans la partie 3 qui affirme que la catégorie des structures de Hodge mixtes $\crmhs $ est abélienne. De plus elle est fermée par sommes directes, produits tensoriels et dualisation.\\

D'après la décomposition de Deligne \ref{lesipq}, la donnée d'une structure de Hodge mixte sur $V=V_{\bf Q}\otimes {\bf C}$ induit une strucure de Hodge mixte sur $gl(V)$ de bigraduation associée :
$$ gl(V)^{r,s}=\{ \alpha \in gl(V) \vert \forall p,q \, \alpha : I^{p,q} \rightarrow I^{p+r,q+s} \}.$$

Soit $\cals$ une polarisation graduée de la structure de Hodge mixte $(V,W_{\bullet },F^{\bullet })$ et $G_{\bf C}=\{ g \in GL(V)^{ W_{\bullet }} \vert Gr(g) \in Aut_{\bf C}(\cals) \}$.
 Alors $(W_{\bullet },F^{\bullet })$ définit une structure de hodge mixte sur $\mathfrak{g}=Lie(G_{\bf C})$ par la bigraduation :
$$ {\mathfrak{g}}^{r,s}=gl(V)^{r,s}\cap Lie(G_{\bf C}).$$
\begin{definition} 
Soit $W_{\bullet }$ une filtration croissante d'un espace vectoriel complexe de dimension finie $V$. Alors un endomorphisme semi-simple $Y \in gl(V)$ gradue() $W_{\bullet }$ si pour tout $k$ $W_{k}=W_{k-1}\oplus E_{k}(Y)$ où $E_{k}(Y)$ est le sous-espace propre de $Y$ associé à la valeur propre $k$.\\
\end{definition}
Notons $\caly (W_{\bullet})$ l'ensemble des graduations de $W_{\bullet}$.\\
\\
Pour décrire $\caly (W_{\bullet})$, considérons $Lie_{-1}$ l'idéal nilpotent de $gl(V)$ défini par :
$$ Lie_{-1}=\{ g \in gl(V) \vert \forall k \,\alpha(W_{k})\subset W_{k-1} \}$$
Il vient alors :\\
\begin{theoremese}\cite{catkapsch}\label{explie}Le groupe de Lie unipotent $\exp (Lie_{-1})$ agit de façon simple et transitive sur $\caly(W_{\bullet })$.\\
\end{theoremese}

\begin{center}
{\bf Rappel des classifiants de structures de Hodge pures polarisées}
\end{center}
Pour fixer les notations et préparer le travail pour les structures de Hodge mixtes, rappelons la construction des espaces classifiants des structures de Hodge pures polarisées d'après Griffiths (\cite{gri1}, \cite{gri2} et \cite{grisch}).\\
Soit $V$ un espace vectoriel de dimension finie et muni d'une structure rationnelle $V_{\bf Q}$ et d'une forme bilinéaire non-dégénérée :
$$ Q: V_{\bf Q} \otimes V_{\bf Q} \rightarrow {\bf Q}$$
de la m\^eme parité que $k$, pour $k$ un entier fixé. Donnons-nous une partition de $dim \,V$ par une somme d'entiers naturels $\{ h^{p,k-p} \}$ symétriques par rapport à la diagonales i.e. $h^{p,q}=h^{q,p}$. On peut alors construire l'espace classifiant de toutes les structures de Hodge pures de poids $k$ sur $V$, polarisées par $Q$ et satifaisant $dim_{\bf C}H^{p,k-p}=h^{p,k-p}$. 
$$ \cald = \cald (V,Q,\{ h^{p,k-p} \})$$
On peut équiper cet ensemble d'une structure de variété complexe de la façon suivante. Posons $f^{p}= \sum_{r \geq p}h^{r,k-r}$ et soit $\check \calf$ la variété des drapeaux consistant en toutes les filtrations décroissantes $F^{\bullet}$ de $V$ telles que $dim_{\bf C}F^{p}=f^{p}$.\\
\begin{lemme}
$\check \calf$ est une variété algébrique lisse.\\
\end{lemme}
\begin{preuve}A été faite plus haut, dans l'étude des variétés des drapeaux.

\end{preuve}

Soit ${\check \cald} \subset \check \calf$ le sous-ensemble des filtrations qui vérifient la première des relations bilinéaires de Hodge-Riemann $Q(F^{p},F^{k-p+1})=0$. Le lemme qui suit est une conséquence du lemme précédent :\\

\begin{lemme}
$\check \cald$ est une variété algébrique lisse.\\
\end{lemme}
\begin{preuve}
Le groupe de Lie complexe
$G_{\bf C}=\text{Aut}_{\bf C}(Q) $ agit transitivement sur $\check \cald$ et donc $\check \cald$ est une sous-variété lisse de $\check \calf$. 
\end{preuve} 

L'espace ainsi définit est classifiant dans le sens suivant : soit $\calv \rightarrow S$ une variation de structures de Hodge pures polarisées.  Le choix d'un point de base $s_{0} \in S$ détermine une représentation de  monodromie :
$$ \rho : {\pi}_{1}(S,s_{0}) \rightarrow \text{Aut} ( {\calv_{s_{0}}}).$$
Notons $\Gamma $ le sous-groupe $\rho({\pi}_{1}(S,s_{0}))$, la représentation de la monodromie peut s'écrire $ \rho : {\pi}_{1}(S,s_{0}) \rightarrow \Gamma $. Le choix du point de base nous donne donc une application :
$$\Phi: S \rightarrow \cald /\Gamma.$$

\begin{center}
{\bf Classifiants dans $\crzmhs \,$,$\, \crmhs$ }
\end{center}
%()expliquer le choix d'une structure de Hodge de reference

Nous avons définit suivant \cite{stezuc} la notion de variation de structure de Hodge mixte graduellement polarisée $\calv \rightarrow S$. On peut associer à une telle donnée un espace classifiant $\calm$ et un morphisme de la base de la variation vers cet espace classifiant (modulo le choix d'un point de base). Cet espace classifiant se construit de la façon suivante (\cite{usu}, \cite{pea}) :

Soit $V$ un espace vectoriel complexe muni d'une structure rationnelle sous-jacente $V_{\bf Q}$ et d'une filtration croissante $W_{\bullet}$ définie aussi sur ${\bf Q}$, la filtration par le poids. Donnons-nous aussi un ensemble $\cals $ de polarisations sur les espaces gradués associés à la filtration par le poids, c'est à dire une suite de formes biliéaires de $(-1)^{k}$-paires $${\cals }_{k}: Gr_{k}^{W} \otimes Gr_{k}^{W} \rightarrow {\bf C}$$ 
et une partition de $\text{dim}_{\bf C}V$ par des entiers $\{h^{p,q}\}$.

 On peut alors construire l'espace classifiant $\calm $ des structures de Hodge mixtes graduellement polarisées par $\cals $ et vérifiant la condition de dimension 
$$\text{dim}_{\bf C}Gr_{F}^{p}Gr_{p+q}^{W}V=\text{dim}_{\bf C}I^{p,q}_{(F^{\bullet },W_{\bullet })}=h^{p,q}.$$ 

Commençons, de façon similaire au cas des structures de Hodge pures, par construire la variété des drapeaux $\check \calf$ de toutes les filtrations décroissantes $F^{\bullet }$ de $V$ telles que $\text{dim}_{\bf C}F^{p}=f^{p} $ où $f^{p}=\sum_{r \geq p,s}h^{r,s}$. C'est une variété algébrique lisse d'après l'étude des variétés des drapeaux.

Définissons ensuite $\check \calf(W_{\bullet})$ comme étant la sous-variété de $\check \calf$ correspondant aux filtrations $F^{\bullet}$ de $V$ qui vérifient en plus la propriété 
$$\text{dim}_{\bf C}F^{p}Gr_{k}^{W_{\bullet }}=f^{p}_{k},$$ 
où $f^{p}_{k}=\sum_{r \geq p}h^{r,k-r}$. Le lemme suivant prouve que $\check \calf(W_{\bullet})$ est lisse, sa démonstration est similaire à celle du lemme précédent :
\begin{lemme}
Le groupe de Lie complexe $GL(V)^{W_{\bullet}}=\{ g \in GL(V) \vert \forall k g : W_{k} \rightarrow W_{k} \}$ agit transitivement sur $\check \calf(W_{\bullet})$. $\check \calf(W_{\bullet})$ est donc lisse.\\ 

\end{lemme}

Soit ensuite $\check \calm$ l'ensemble des filtrations $F^{\bullet} \in \check \calf(W_{\bullet})$ qui satifont à la première des relations bilinéaires de Hodge-Riemann $\forall k \,\cals_{k}(F^{p}Gr_{k}^{W_{\bullet}},F^{k-p+1}Gr_{k}^{W_{\bullet }})=0$. Alors :\\
\begin{lemme}
Le groupe de Lie complexe $G_{\bf C}=\{ g \in GL(V)^{W_{\bullet}} \vert Gr(g) \in Aut_{\bf C}({\cals}) \}$ agit transitivement sur $\check \calm$ qui est donc lisse.\\
\end{lemme}

\begin{preuve}
Soit $\caly (W_{\bullet })$ l'ensemble des graduations de $V$ par $W_{\bullet }$, 
$$\caly (W_{\bullet })=\{ Y \vert Y:Gr^{W_{\bullet }}V \cong {\rightarrow } V \}.$$ 
Soit 
$$\check \cald =\bigoplus_{k} {\check \cald}(Gr_{k}^{W_{\bullet}},{\cals_{k}},h^{p,k-p}),$$ 
et, 
$$\check \calx = \check \cald \times \caly (W_{\bullet }).$$ 
Soit $\pi :\check \calx \rightarrow \check \calm $ la projection qui envoie un point $ (\oplus_{k}F^{\bullet}_{k},Y) \in {\check \calx}$ sur la filtration $F^{\bullet } \in {\check \calm}$ déterminée par la filtration $\oplus_{k}F^{\bullet}_{k}$ sur $Gr^{W_{\bullet }}$ et l'isomorphisme $Y: Gr^{W_{\bullet }}V \cong{\rightarrow } V$ c'est à dire : $F^{p}=\oplus_{k}Y(F^{p}_{k})$. L'application $\pi $ est surjective et le diagramme suivant est commutatif :

$$
\xymatrix{
{\check \calx } \times G_{\bf C}  \ar[d] \ar[r]^{(\pi,id)} & {\check \calm} \times G_{\bf C} \ar[d]\\
{\check \calx }  \ar[r]^{\pi} & {\check \calm}, 
}
$$
et l'application $\pi :\check \calx \rightarrow \check \calm $ est donc $G_{\bf C}$-équivariante. Ainsi si l'on montre que $G_{\bf C}$ agit transitivement sur $\check \calx$, on aura aussi prouvé que l'action de $G_{\bf C}$ est transitive sur $\check \calm$. Soit un point $P \in {\check \calx}$, $P=(\oplus_{k}F^{\bullet}_{k},Y)$. Nous allons montrer la transitivité dans les deux directions du produit $\check \calx = \check \cald \times \caly(W_{\bullet })$. D'une part le sous-groupe de $G_{\bf C}$, $G_{\bf C}^{Y}=\coprod_{k} G_{\bf C}(Gr_{k}^{W_{\bullet}},{\cals_{k}},h^{p,k-p}) $ agit transitivement sur $\check \cald$ tout en maintenant le scindement $Y$ constant. D'autre part, d'après le théorème \ref{explie} p.\pageref{explie} le groupe de Lie unipotent $\exp(Lie_{-1})$ agit de façon simple et transitive sur $\caly(W_{\bullet })$ et laisse $\check \cald $ fixe. D'où la transitivité de l'action sur $\check \calx$ et donc sur $\check \calm$. 

\end{preuve} 

Pour toute variation graduellement polarisée de structures de Hodge mixtes au dessus d'une variété $S$ nous avons donné la construction d'un espace classifiant $\calm$. Cet espace est classifiant dans le sens où toute variation de structures de Hodge mixtes associée au choix d'un point de base $s_{0}\in S$ donne lieu à un morphisme $$\Phi: S \rightarrow \calm /\Gamma,$$ où $\Gamma$ est l'image du groupe fondammental de $S$ par la représentation de monodromie $\rho$ associée au système local sous-jacent à la variation de structures de Hodge mixtes, $\rho : \pi_{1}(S,s_{0}) \rightarrow \text{Aut}(\calv_{s_{0}})$. L'application $\Phi$ admet un relèvement au revêtement universel $\tilde S$ de $S$, qui mène au diagramme commutatif :
$$\xymatrix{ {\tilde S} \ar[r]^{\tilde \Phi} \ar[d] &  \calm \ar[d]\\
S \ar[r]^{\Phi} & \calm /\Gamma.}$$

Au dessus de $\calm$ il y a un fibré vectoriel universel $\calv$ filtré par des sous-fibrés stricts $\calf^{\bullet }$ tels que le tiré en arrière de ces données par le morphisme $\Phi$ permet de retrouver la variation de structure de Hodge mixtes associée à $\Phi$.

\begin{center}
{\bf Strate des structures de Hodge mixtes ${\bf R}$-scindées}
\end{center}

Nous donnons ici une remarque (cf \cite{pea}, p.17) sur les structures de Hodge qui sont ${\bf R}$-scindées. Par les mêmes méthodes que pour montrer que l'action de $G_{\bf C}$ est transitive, on montre que le groupe de Lie réel 
$$G_{\bf R}=\{ g \in GL(V_{\bf R})^{W} \vert Gr(g) \in \text{Aut}_{\bf R}(\cals) \}$$
agit transitivement sur la variété $\calc^{\infty}$ $\calm_{\bf R} \subset \calm$ qui paramétrise toutes les filtrations de Hodge qui donne lieu à une structure de Hodge mixte qui est ${\bf R}$-scindée. En particulier l'action de ${ G}_{\bf R}$ préserve donc la sous-variété $\calm_{\bf R}$.

 On peut généraliser cette remarque aux autres strates que la strate $\alpha =0$, $\calm_{\bf R}$. Puisque l'action de $G_{\bf R}$ ne modifie pas les positions relatives de la filtration de Hodge et de sa conjuguée :
\begin{proposition}
L'action de $G_{\bf R}$ préserve les strates sur lesquelles l'invariant $\alpha$ est constant, i.e. préserve le niveau de ${\bf R}$-scindement.  
\end{proposition}

\subsection{Etude des stratifications associées aux variations de structures de Hodge mixtes, exemple des variations de structures de Hodge mixtes de Tate}

Considérons une variation de structures de Hodge mixtes $(\calv,\calw_{\bullet },\calf^{\bullet})$ sur une variété $S$. En tout point $x \in S$, le niveau de ${\bf R}$-scindement de la structure de Hodge mixte au point $x$ est donné par
$$\alpha(H_{x})=\frac{1}{2}\sum_{p,q}\,(p+q)^{2}(h^{p,q}_{x}-s^{p,q}_{x}).$$

Les nombres de Hodge sont constants sur la base. En effet ils sont donnés par $h^{p,q}_{x}=\text{dim}_{\bf C}\,(Gr^{\calf}_{p}Gr_{\calw}^{p+q}\calv)_{x}$ et correspondent par définition d'une variation de structures de Hodge mixtes aux nombres de Hodge d'une variation de structures de Hodge pures de poids $p+q$. Les sommes $\sum_{p'\geq p}\, h^{p',p+q-p'}$ sont constantes pour tout $p$ car elles donnent les rangs des fibrés induits par les $\calf^{p}$ sur $\calw_{p+q}/\calw_{p+q-1}$. 

Pour tout $x$, $h^{p,q}_{x}=h^{p,q}$.

En terme des fibrés de Hodge, sous-fibrés stricts de $\calv$, la formule devient
$$\alpha(H_{x})=\sum_{p,q}\,(p+q)^{2}(h^{p,q}-\text{dim}_{\bf C}Gr_{\calf_x}^{p}Gr_{{\overline \calf}_x}^{q}).$$
Dans la partie ``comportement du niveau de ${\bf R}$-scindement par opérations sur les structures de Hodge mixtes'' on a décomposé $\alpha$ en une partie avec les nombres de Hodge $h^{p,q}$ et une partie avec les $s^{p,q}$, $\alpha(\bullet)=\alpha^{+}(\bullet)-\alpha^{-}(\bullet)$. Avec
$$\alpha^{-}_{s}= \frac{1}{2} \sum_{(p,q) \in {[1,N] \times [1,N]}} ( 2pq+4) f^{p,q}_{s}+\frac{1}{2} \sum_{q \in [1,N]} ( 2q-1) f^{0,q}_{s}+\frac{1}{2} \sum_{p \in [1,N]}(2p-1) f^{p,0}_{s}- \frac{1}{2} f^{0,0}_{s}.$$
$\alpha^{+}(\bullet)$ est constant.

La filtration conjuguée à la filtration de Hodge varie de façon anti-holomorphe sur la base de la variation de structures de Hodge mixtes. On ne peut donc utiliser directement l' étude des morphismes de fibré sur S de la forme
$$ \calf^{p} \rightarrow  \calv \rightarrow \calv / {\overline \calf}^{q},$$ 
car ces morphismes ne sont pas holomorphes mais $\calc^{\infty}$. Pour nous placer dans un cadre holomorphe, nous allons "déplier" la variation. 

Considérons la filtration de Hodge universelle $ \calf^{\bullet }$ au dessus du classifiant $\calm$ et les deux familles de filtrations $ \calf^{\bullet }_{1}$ et $ \calf^{\bullet }_{2}$ au dessus du produit $\calm \times \calm$ données comme tirées en arrière de $ \calf^{\bullet }$ sur $\calm$ par les deux projections sur le produit. $\calm$ est naturellement équipé d'une conjugaison complexe car chaque point de $\calm$ représente une filtration de Hodge dont la conjuguée est aussi dans le classifiant. Soit $z \in \calm$, on a alors,
$$ \text{dim}_{\bf C}( \calf^{p }_{1} \cap \calf^{q}_{2})_{(z,{\overline z})}=\text{dim}_{\bf C}( \calf^{p }_{} \cap {\overline \calf}^{q }_{})_{z}.$$

Ainsi nous voulons étudier les sauts du niveau de ${\bf R}$-scindement des structures de Hodge mixtes données par la variation en fonction des sauts de la quantité $\alpha$ donnée en fonction des $f^{p,q}$ et des $h^{p,q}$ sur le produit des classifiants par le morphisme $$\calm \rightarrow \calm \times \calm,\,\,\,z\mapsto (z,{\overline z}).$$ 

On veut donc étudier la stratification associée à l'invariant $\alpha$ par l'intermédaire de la stratification du classifiant par le morphisme $$\phi : S \rightarrow \calm \rightarrow \calm \times \calm .$$

Pour définir la stratification de $\calm \times \calm$ nous avons le résultat :

\begin{proposition}\label{alphasci}
L'invariant $alpha$ sur le produit du classifiant $\calm \times \calm$ associ\'e \`a la variation de structures de Hodge mixtes $\calv \rightarrow S$ est semi-continu supérieurement sur la base de la variation à valeurs entières et définit donc une stratification 
$$\calm \times \calm=\coprod_{i \in I}(\calm \times \calm)_{i}.$$
\end{proposition}
\begin{preuve}
C'est une application directe de la formule de $\alpha^{-}$ et du lemme \ref{flagsci} p.\pageref{flagsci}.
\end{preuve}
Par l'étude précédente de dégénérescence sur la variété des drapeaux nous pouvons donner un expression explicite des strates pour $\alpha $. La semi-continuité supérieure des $\varphi^{p,q}$ sur $\calm \times \calm$, donnés ici par les rangs des morphismes de fibrés
$$ p_{1}^{*}\calf^{p} \rightarrow  \calv \rightarrow \calv / p_{2}^{*}{\calf}^{q},$$
est équivalente à la semi-continuité inférieure des $f^{p,q}$. 

Etudions pour tout $n $ l'ouvert $\{ s \in \calm \times \calm \vert \alpha(s) \geq \} $. Rappelons que l'on note
$D_{r_{p,q}}=\{ s \in \calm \times \calm \vert f^{p,q} \leq r \}.$
Alors, on rappelle que $\alpha^{+}(s)=\alpha^{+}(s_{0})$ est constant sur la base,
\begin{eqnarray*}
{} & \alpha(s) &  \geq n\\
\Leftrightarrow & \alpha^{-}(s) &  \leq  \alpha^{+}(s_{0})-n\\
\Leftrightarrow & \sum_{(p,q) \in {[1,N] \times [1,N]}} ( 2pq+4) f^{p,q}_{s} & \leq  \alpha^{+}(s_{0})-\frac{1}{2} \sum_{q \in [1,N]} ( 2q-1) f^{0,q}_{s}\\
{} & {} &- \sum_{p \in [1,N]}(2p-1) f^{p,0}_{s}+ f^{0,0}_{s}-n,
\end{eqnarray*}
Notons par $C$ tout ce qui est constant dans le menbre de droite, il vient,
\begin{eqnarray*}
\alpha(s)   \geq  C-n  \Leftrightarrow & {}\\
{}   &     s \in \cup_{ \{ r_{p,q},p,q\, \in {\bf Z}^{3} \vert \sum_{(p,q) \in {[1,N] \times [1,N]}} ( 2pq+4) r^{p,q}_{s} \leq C-n \} } \cap_{p,q}\,\,D_{r_{p,q}}
\end{eqnarray*}
Par intersection avec l'image de $\calm $ dans le produit on en déduit une stratification qui donne le niveau de ${\bf R}$-scindement des structures de Hodge mixtes, d' après la proposition 

\begin{proposition}\label{alpharesci}
Le niveau de scindement $\alpha$ associ\'e \`a une variation de structures de Hodge mixtes sur $S$ est semi-continu inférieurement sur le revêtement universel de $S$.
\end{proposition}

\begin{preuve}
Montrons d'abord qu'il est semi-continu inférieurement sur $\calm$. Comme $\alpha=\alpha^{+}-\alpha^{-}$ et comme $\alpha^{+}$ est constant, il faut montrer que $-\alpha^{-}$ est semi-continu inférieurement. $-\alpha^{-}$ est somme à coefficients positifs des fonctions $s \mapsto -f^{p,q}$ d'après le lemme \ref{alphaplus} p. \pageref{alphaplus} et de la fonction $s \mapsto f^{0,0}$ qui elle est constante comme rang de l'espace vectoriel fibre. Ces fonctions sont semi-continues inférieurement d'après le lemme \ref{flagsci} p.\pageref{flagsci} ce qui permet de conclure (on peut en effet toujours supposer que l'une des deux filtration est fixe dans le lemme bien qu'elles soient liées car conjuguées, l'automorphisme $g$ pour passer de l'un à l'autre ne varie plus holomorphiquement mais de façon $\calc^{\infty}$ mais seule la continuité est nécessaire pour conclure ici).

Le revêtement universel de $S$ est alors stratifié par $\alpha$ d'après l'existence d'un morphisme de ce revêtement vers $\calm$ associé à la variation de structures de Hodge mixtes :
$${\widetilde S}=\coprod_{i \in I}\,{\widetilde S}_{i}.$$
\end{preuve}

\begin{center}
{\bf Strates $\alpha=cste$ et hypothèse $({\bf H})$}
\end{center}
${}$\\

Considérons une variation de structures de Hodge mixtes $\calv \rightarrow S$ sur une variété $S$ que l'on supposera simplement connexe pour simplifier. On en déduit une stratification de $S$ par le niveau de ${\bf R}$-scindement
$${ S}=\coprod_{i \in I}\,{ S}_{i}.$$

\begin{theoreme}
Soit $X$ une sous-variété de $S$ incluse dans l'une des strates $S_{i}$ de la décomposition associée à la variation de structures de Hodge mixtes sur $S$ alors l'hypothèse $({\bf H})$ est vérifiée sur $X$ pour la variation de structures de Hodge mixtes associée sur $X$.

Réciproquement si $({\bf H})$ est vérifiée pour la restriction de la variation à une sous-variété $X$ alors $X$ est incluse dans l'une des composantes de la stratification associée à la variation de structures de Hodge mixtes.
\end{theoreme}

\begin{preuve}
L'invariant $\alpha$ est somme des $x \mapsto -f^{p,q}_{x}$ qui sont des fonctions semicontinues supérieurement d'après le lemme $40$ et de la fonction $x \rightarrow f^{0,0}_{x}$ qui est constante comme rang du système local sous-jacent à la variation. La somme étant constante, ces fonctions sont constantes donc d'après la partie II $({\bf H})$ est vérifiée.

La réciproque est immédiate puisque les $t^{p,q}$ sont constants si $({\bf H})$ est vérifiée, en effet cette hypothèse est équivalente à la constance des fonctions $x \mapsto \text{dim}_{\bf C} Gr_{\calf^{\bullet}}^{p}{}_{s}Gr_{{\overline \calf}^{\bullet}}^{q}{}_{s}$ qui ne sont autres que les fonctions $x \mapsto t^{p,q}_{x}$.
\end{preuve}

\subsection{Perspectives}
\subsubsection{Géométriser les variations de structures de Hodge mixtes}
Comme précisé ci-dessus, on ne peut appliquer directement les constructions de Rees faites dans la partie 2 aux variations de structures de Hodge, même par strates. En effet, le principal obstacle est le fait que la filtration conjuguée à la filtration de Hodge varie de façon anti-holomorphe. 

Considérons une variation de structures de Hodge mixtes sur une base $S$. Cette variation telle que nous voulons l'étudier ici est la donnée d'un ${\bf R}$-système local $\calv_{\bf R}$ filtré par une filtration décroissante de sous-systèmes locaux $\calw^{\bullet}_{\bf R}$ (comme pour le cas ponctuel, pour se ramener à la situation de trois filtrations décroissantes, on associe à la filtration croissante $\calw_{\bullet}$ par des sous-systèmes locaux du système local $\calv_{\bf R}$ une filtration décroissante $\calw^{\bullet }$ en posant, pour tout entier $k$, $\calw^{k}=\calw_{-k}$), d'une filtration décroissante $\calf^{\bullet }$ du fibré vectoriel $\calv:=\calo_{S}\otimes \calv_{\bf R}$ tels que la connexion canonique associée au système local, $\nabla$, satisfasse la tranversalité de Griffith pour la filtration $\calf^{\bullet }$ et qu'en tout point $s \in S$, $(\calv_{s},\calw^{\bullet}_{s}, \calf^{\bullet}_{s}, {\overline \calf}^{\bullet}_{s})$ définissent un triplet (ordonné) de trois filtrations opposées sur un espace vectoriel. La filtration ${\overline \calf}^{\bullet}_{s}$ se déduit de la filtration de Hodge de la fibre au point $s$ par conjugaison par rapport à la structure réelle sous-jacente. Cette donnée nous fournit un fibré vectoriel plat muni de trois filtrations par des sous-espaces vectoriels stricts $(\calv, \calw^{\bullet}, \calf^{\bullet }, {\overline \calf}^{\bullet })$, le problème pour une utilisation directe de la partie 2 est que les sous-fibrés définis par la filtration conjuguée à la filtration de Hodge ne sont pas holomorphes mais anti-holomorphes. On ne peut donc pas utiliser directement le dictionnaire par strates donné dans l'étude des constructions de Rees relatives.

Pour palier ce problème, nous allons ``déplier'' les données qui correspondent à une variation de structures de Hodge mixtes sur $S$, qui contiennent de l'holomorphie et de l'antiholomorphie, en données holomorphes sur le produit de $S$ par elle-même, de façon à ``rendre holomorphe'' la filtration conjuguée à la filtration de Hodge et retrouver la situation initiale en considérant les tirés en arrière des objets définis sur le produit $S\times S$ par le morphisme $(i,{\overline i})$ défini par :
$$(i,{\overline i}) :S \rightarrow S  \times S ,\,\,\, (i,{\overline i})(s)=(s,{\overline s}).$$
Nous voulons voir la variété originale $S$ comme "l'antidiagonale" du produit $S \times S$. 
  
Considérons la situation directement sur le produit du classifiant par lui-même. Nous voudrions étudier les variations de structures de Hodge à partir des objets qu'elles permettent de définir sur $\calm \times \calm$.
Ces données sont, en notant $V=\calv_{s_0}$ la fibre choisie pour référence :\\
\hspace*{1.3cm}$(i)$ le système local $\calv:=\calm \times \calm \times V$,\\
\hspace*{1.3cm}$(ii)$ la filtration décroissante associée à la filtration par le poids $\calw^{\bullet}$,\\
\hspace*{1.3cm}$(iii)$ la filtration $\calf^{\bullet}_{1}:=p_{1}^{*}\calf^{\bullet}$ par des sous-fibrés vectoriels stricts, et,\\
\hspace*{1.3cm}$(iv)$ la filtration $\calf^{\bullet}_{2}:=p_{2}^{*}\calf^{\bullet}$ par des sous-fibrés vectoriels stricts,\\
où $p_{1}$ et $p_{2}$ sont les projections sur le premier et deuxième facteur de $\calm \times \calm$ vers $\calm$,\\
\hspace*{1.3cm}$(v)$ deux connexions intégrables $\nabla_{i}$ où $i \in \{1,2\}$ telles que $\calf^{\bullet }_{i}$ vérifie la transversalité de Griffith pour $\nabla_{i}$.

Les deux connexions pouvant se déduire l'une de l'autre, on peut remplacer les données $(i),(ii),(iii),(iv),(v)$ par $(i),(ii),(iii),(iv),(v)'$ où $(v)'$ ne porte que sur $\nabla_{1}$ par exemple.    

${}$\\

Alors $(\calv,\calw^{\bullet},\calf^{\bullet }_{1},\calf^{\bullet }_{2}) \rightarrow \calm \times \calm$ est un fibré vectoriel filtré par des sous-fibrés stricts et tel que les filtrations soient exhaustives et décroissantes et tels que de plus pour tout $x \in \calm \times \calm$, la donnée de $(V,\calw^{\bullet}_{x},\calf^{\bullet }_{1}{}_{x},\calf^{\bullet }_{2}{}_{x})$ est celle d'un espace vectoriel muni de trois filtrations opposées. 

On peut faire la construction de Rees relative exhibée plus haut pour obtenir un faisceau réflexif équivariant pour l'action du groupe ${\bf T}$ sur $\calm \times \calm \times {\bf P}^2$ qui provient de l'action de ce groupe sur ${\bf P}^2$, $$\xi(\calv,\calw^{\bullet},\calf^{\bullet}_{1},\calf^{\bullet }_{2}) \rightarrow \calm \times \calm \times {\bf P}^2.$$ 

Rappelons l'hypothèse ${\bf (H)}$. Elle est satisfaite si pour tout $(p,q)$, les faisceaux cohérents $Gr^{p}_{\calf_{1}} Gr^{q}_{\calf_{2}}\calv$ sont des faisceaux localement libres sur $\calm \times \calm$, cela revient au même de demander que les fonctions sur cette variété $x \mapsto \text{dim}_{\bf C} Gr^{p}_{\calf_{1}{}_x} Gr^{q}_{\calf_{2}{}_x}\calv_{x}$ soient constantes. On sait par la partie 2 qu'il existe une stratification de $\calm \times \calm=\coprod_{i \in I}(\calm \times \calm)_{i}$ telle que sur chaque strate les restrictions sont des fibrés vectoriels et tellles que l'on peut décrire le fibré vectoriel trifiltré en fonction du fibré équivariant associé.
 
Supposons que l'on ait un ensemble de propriétés $(P)$ à ajouter aux données\\
$(i),(ii),(iii),(iv),(v)$ tel que l'on puisse établir, une équivalence entre variation de structures de Hodge mixtes sur $\calm $ et objets holomorphes (ou algébriques suivant le contexte) donnés par $(i),(ii),(iii),(iv),(v)$ et $(P)$ sur le produit $\calm \times \calm$, alors au moins sur les strates sur lesquelles l'hypothèse ${\bf (H)}$ est vérifiée, on pourrait traduire la donnée d'une variation de structures de Hodge mixtes sur $S$ en termes de fibrés sur $\calm \times \calm \times {\bf P}^{2}$, équivariant pour l'action habituelle du tore sur le produit héritée de l'action sur le plan projectif, ${\bf P}^{1}$-semistable par rapport à $\calm \times \calm$ plus une condition $(P)'$ associée à $(P)$. Les construction de Rees sur les différents $\calm \times \calm \times {\bf A}^{1}$ dans cette description feraient appara\^{\i}tre la propriété de transversalité de Griffiths (pour les filtrations issues de Hodge et conjuguée) et la platitude de la filtration par le poids comme décrit dans le lemme qui conclut la partie 2, intégrons ces propriétés dans $(P)'$, on aurait alors la description géométrique des variations de structures de Hodge mixtes suivante, en notant $VSHMP$ les variations de structures de Hodge mixtes graduellement polarisées,  

$$\xymatrix{ VSHMP \text{ sur } \calm \ar@{.>}[r] \ar@{.>}[rd] & \text{Sys. locaux }+\,(i),(ii),(iii),(iv),(v)+\,(P)   \ar@{.>}[d] \ar@{.>}[l]\\ & \calf ib_{{\bf P}^{1}-semistables/\calm \times \calm,\mu=0/\calm \times \calm}(\calm \times \calm \times {\bf P}^{2}/{\bf T})+\,(P)' \ar@{.>}[lu] \ar@{.>}[u] }$$

\subsubsection{Structures de Hodge limites}
La possibilité de géométriser les notions de structures de Hodge mixtes limites dépend évidemment directement du dictionnaire escompté ci-dessus. Modulo ces équivalences définir des objets géométriques associés aux structures de Hodge mixtes limites reviendrait à appréhender les limites d'objets dans les espaces de modules ou champs des objets de 
$$\calf ib_{{\bf P}^{1}-semistables/\calm \times \calm,\mu=0/\calm \times \calm}(\calm \times \calm \times {\bf P}^{2}/{\bf T})+\,(P)'.$$   

\subsubsection{Rigidité de variations de structures de Hodge mixtes}
Rappelons le théorème suivant, d\^u \`a Deligne,
\begin{theoremese} Une variation admissible de structures de Hodge mixtes sur une variété k\"alherienne compactifiable est complétement déteminée par sa représentation de monodromie et sa valeur en un point, i.e. que, un système local filtré $(\calv,\calw^{\bullet})$ étant donné sur une variété du type considéré $S$ et en un point $s \in S$ une filtration décroissante $F^{\bullet }_s$ telle que $(\calv_{s},\calw^{\bullet}_{s},F^{\bullet }_{s})$ déterminent une structure de Hodge mixte, alors il existe au plus une variation de structures de Hodge mixtes admissible qui prolonge ces données. 
\end{theoremese}
On propose ici une idée pour trouver des exemples, à partir des données du théorème, i.e. une représentation du groupe fondamental et d'une structure de Hodge mixtes en un point, de telles données qui ne se prolongent pas en variation admissible de structures de Hodge mixtes. 

Considérons donc un système local filtré sur une variété $S$, $(\calv,\calw^{\bullet})$ et une filtration qui donne une structure de Hodge en un point. On peut alors en déduire une stratification du rev\^etement universel de $S$, $\tilde S$ par le niveau de ${\bf R}$-scindement. Une telle stratification doit donner des informations géométriques (homologie) sur $\tilde S$, théoriques puisqu'elles correspondent à une représentation du groupe fondamental qui ne provient pas forcément d'une variation de structures de Hodge mixtes. Si ces informations sont contraires à ce que l'on connait de la variété, alors c'est qu'il n'y a aucune variation de structures de Hodge mixtes correspondant aux données initiales.

Si l'on a une description des variations de structures de Hodge mixtes sur un classifiant $\calm$ en termes de fibrés trifiltrés sur $\calm \times \calm$ comme explicité au-dessus, alors la géométrie des strates est bien connues dans de nombreux cas car elle est donnée par la dégénérescence de morphismes de fibrés vectoriels. Par exemple sur $\calm \times \calm$, les fibrés $p_{1}^{*}\calf^{p}$ et $p_{2}^{*}\calf^{q}$ sont des sous-fibrés des fibrés universels, donc le fibré $\text{Hom}(p_{1}^{*}\calf^{p},\calv p_{2}^{*}\calf^{q})$ est ample et donc suivant \cite{fullaz}, si la dimension attendue du lieu où le morphisme est de rang inférieur ou égal à $r$, $2.\text{dim}(\calm)-(f^{p}-r)(f^{q}-r)$ est au moins $1$, alors ce lieu est connexe.
 
On a des résultats plus fort lorsqu'on peut équiper les fibrés dont on étudie les dégénérescences de formes bilinéaires (cf \cite{deb}), qui pourraient provenir ici des polarisations associées aux variations de structures de Hodge mixtes.

%\subsection{dessins}

%\begin{figure}[ht] 

%\begin{center}

%\scalebox{1}{\input{try.pstex_t}}

%\end{center}

%\caption {Exemple saut de $\alpha$}

%\end{figure}

%blablablabla  gGGGjhklhj

%\begin{figure}[ht] 

%\begin{center}

%\scalebox{1}{\input{essi.pstex_t}}

%\end{center}

%\caption {Exemple saut de $\alpha$}

%\end{figure}

%\begin{figure}[ht] 

%\begin{center}

%\scalebox{1}{\input{tryagain.pstex_t}}

%\end{center}

%\caption {Exemple saut de $\alpha$}

%\end{figure}

%\begin{center}
%\begin{tabular}{|c||c|}
%\hline
%$X_{\Rr}$ non séparante &$X_{\Rr}$  séparante\\
%\hline \hline

%$Y_{\Rr}$ est maximal et &$Y_{\Rr}$ est séparant et\\
% $\widetilde{X_{\Rr}}$ est non séparante &$\widetilde{X_{\Rr}}$ est 
%séparante\\

%\hline

%$Y_{\Rr}$ est  non séparant et & $Y_{\Rr}$ est non séparant
%et\\$\widetilde{X_{\Rr}}$ est un nid et $p - n \equiv k^2 \pm 2 \;
%\mbox{mod} \; 16.$ & $\widetilde{X_{\Rr}}$ est non séparante\\

%\hline
%\end{tabular}
%\end{center}

\clearpage
\section*{Annexe  : Structures twistorielles mixtes.}
\addcontentsline{toc}{section}{Annexe 1 : Structures twistorielles mixtes}
\setcounter{section}{5}
\setcounter{theoreme}{0}

%{\bf ANNEXE A : Structures twistorielles mixtes.}\\
$${}$$
Dans cet appendice on rappelle la notion de structures twistorielles mixtes d'après \cite{sim2}. Les entiers $s^{p,q}$ introduits dans la section d'application à la théorie de Hodge permettent de donner une description explicite des structures twistorielles mixtes en fonction des dimensions d'intersection des sous-espaces vectoriels que définissent les filtrations qui composent une structure twistorielle équivariante.

On fixe une droite projective ${\bf P}^1$ avec ses points $0,1,\infty$ et le fibré en droite standard ${\calo}_{{\bf P}^1}$.
\begin{definition}Une structure twistorielle est un fibré vectoriel $\cale$ sur la droite projective ${\bf P}^{1}=\emph{Proj}\,{\bf C}[u,v]$. L'espace vectoriel sous-jacent est $\cale_{1}$, la fibre en $1$ du fibré.\\
\end{definition}

Elle est dite pure de poids $\omega$ si $\cale$ est semistable de pente $\omega$. D'après Grothendieck, \cite{gro1}, tout fibré vectoriel sur ${\bf P}^{1}$ est somme directe de fibrés en droite. Cela est équivalent à dire que $\cale=\oplus_{i} \,\calo_{{\bf P}^{1}}(a_{i})^{i}$.\\
\begin{definition}
Une structure twistorielle mixte est la donnée d'une filtration croissante $\calw_{\bullet}\cale$ par des sous-fibrés stricts d'une structure twistorielle $\cale$ telle que pour tout $i$ les fibrés quotients $\calw_{i}\cale / \calw_{i-1}\cale$ soit une structure twistorielle pure de poids $i$. La filtration $\calw_{\bullet }\cale$ est appelée filtration par le poids. Elle est dite pure si le gradué associé à la filtration par le poids est non trivial en un seul degré. \\  
\end{definition}  

Rappelons comment l'on peut fabriquer un fibré vectoriel sur ${\bf P}^1$ à partir de deux filtrations. Le travail qui a été fait sur ${\bf P}^2$ est une généralisation de cette construction. 

Soit $(V,F^{\bullet},{\hat F}^{\bullet })$ un espace vectoriel muni de deux filtrations décroissantes et exhaustives. On peut alors faire la construction du fibré de Rees de la droite affine pour chacune des filtrations. Pour $(V,F^{\bullet})$ par exemple, on construit le faisceau localement libre $\xi(V,F^{\bullet})$ comme précisé dans la partie $1$. Si l'on note $j:{\bf G}_{m} \rightarrow {\bf A}^1$ l'inclusion, ce faisceau est le sous-faisceau de $j_{*}(V \otimes_{\bf C} \calo_{{\bf G}_m}$ engendré par les sections de la forme $u^{-p}v_{p}$ pour $v_{p} \in F^{p}$ (avec ${\bf A}^{1}=\text{Spec}k[u]$). Nous avons vu, partie $1$, section 1.5.2, que c'est un faisceau ${\bf G}_{m}$-équivariant pour l'action qui recouvre l'action de ${\bf G}_m$ par translation sur la base ${\bf A}^1$ (cf plus bas). Comme décrit plus haut, ou suivant \cite{sim1}, la construction $(V,F^{\bullet}) \mapsto \xi(V,F^{\bullet})$ admet une construction inverse qui à partir d'un faisceau ${\bf G}_{m}$-équivariant sur la droite affine donne un espace vectoriel filtré. Si l'on part d'un fibré vectoriel obtenu à partir d'un espace vectoriel filtré par le constrcution de Rees, la construction inverse donne la même filtration sur le même espace vectoriel, la fibre en $1$.

Revenons au deux filtrations. Par la construction décrite, nous obtenons deux faisceaux localement libres ${\bf G}_{m}$-équivariants sur la droite affine. Nous allons recoller ces deux faisceaux par l'intermédaire de l'automorphisme de ${\bf G}_{m}$, $u \mapsto u^{-1}$. Ceci est possible car les restrictions à ${\bf G}_m$ des faisceaux de Rees associés à chacune des filtrations sont trivialisés par l'action. On en déduit deux isomorphismes $\xi(V,F^{\bullet})\cong {\bf G}_{m}\times V$ et $\xi(V,{\hat F}^{\bullet})\cong {\bf G}_{m}\times V$. Le recollement se fait donc par l'automorphismes de ${\bf G}_{m} \times V$ donné par $(t,v) \mapsto (t^{-1},v)$. Le fibré de Rees obtenu est le fibré ${\bf G}_m$ équivariant
$$\xi(V,F^{\bullet},{\hat F}^{\bullet }).$$

On a alors le résultat suivant : $F^{\bullet}$ et ${\hat F}^{\bullet}$ sont $n$-opposées si et seulement si $\xi(V,F^{\bullet},{\hat F}^{\bullet }) $ est une somme directe de ${\calo}_{{\bf P}^1}(n)$. 

Ceci se montre facilement en utilisant un scindement commun $V=\oplus_{p,q}V^{p,q}$ aux deux filtrations $F^{\bullet}$ et ${\hat F}^{\bullet}$ ; ce qui est toujours possible. Alors $(V,F^{\bullet},{\hat F}^{\bullet})$ est somme directe d'objets munis de deux filtrations triviales décalées $\oplus_{p,q}(V^{p,q},Dec^{p}Triv^{\bullet},Dec^{q}Triv^{\bullet})$. Or,
\begin{eqnarray*}
\xi(V^{p,q},Dec^{p}Triv^{\bullet},Dec^{q}Triv^{\bullet})&=&V^{p,q}\otimes_{\bf C}\calo_{{\bf P}^1}(p.O+q.\infty)\\
{}&\cong & V^{p,q} \otimes_{\bf C}\calo_{{\bf P}^1}(n),
\end{eqnarray*}    
donne l'expression de $\xi(V,F^{\bullet},{\hat F}^{\bullet})$ en prenant la somme directe sur les couples $(p,q)$, ce qui permet de conclure l'équivalence annoncée.\\
${}$\\   
Rappelons que le groupe ${\bf G}_m$ agit sur ${\bf P}^1$ par translation
\begin{eqnarray*}
\mu:& {\bf G}_{m} \times {\bf P}^{1} \rightarrow  & {\bf P}^{1}\\
{} & (t,u) \mapsto &  t.u
\end{eqnarray*}
Une structure twistorielle mixte $(\cale,\calw^{\bullet})$ est dite ${\bf G}_{m}$-équivariante s'il existe un isomorphisme de faisceaux cohérents filtrés sur ${\bf G}_{m} \times {\bf P}^{1}$
$$\rho : \mu^{*}(\cale,\calw^{\bullet})\cong p_{2}^{*}(\cale,\calw^{\bullet}),$$
où $p_{2}$ est la projection sur le deuxième facteur, tel que les deux morphismes que l'on peut en déduire sur ${\bf G}_{m} \times {\bf G}_{m} \times {\bf P}^{1}$ 
$$(\mu \circ \mu)^{*}(\cale,\calw^{\bullet}) \cong  p_{3}^{*}(\cale,\calw^{\bullet}),$$
où $p_{3}$ est la projection sur le troisième facteur.

\begin{proposition}\cite{sim2}\label{equivcar} La catégorie des structures twistorielles mixtes ${\bf G}_{m}$-équivariantes est équivalente à la catégorie des structures de Hodge mixtes complexes.
\end{proposition}
Donnons une idée de la preuve. Comme nous l'avons au dessus et dans l'étude des constructions de Rees sur la droite affine, un espace vectoriel filtré est la même chose qu'un fibré ${\bf G}_m$-équivariant sur la droite affine ${\bf A}^1$. D'un tel fibré équivariant on peut donc déduire deux filtrations sur la fibre au dessus du point $1$. La filtration par le poids est données par la filtration de cette fibre au dessus de $1$ par $\calw^{\bullet}$. Dans l'autre sens, on fait la construction du fibré de Rees associé à deux filtrations ($F^{\bullet}$ et ${\hat F}^{\bullet}$). Le fibré obtenu est filtré par des sous-fibrés stricts associés à la filtration par le poids car à chaque niveau de cette filtration, en  partant du plus bas, on a un sous-fibré strict. C'est bien une structure twistorielle mixte car chaque composante du gradué associé à la filtration par le poids du fibré de Rees, notée $\calw_{\bullet}$, est isomorphe à la construction du fibré de Rees associé au filtrations $F^{\bullet}$ et ${\hat F}^{\bullet}$ induite sur gradué associé par la filtration par le poids de la structures de Hodge complexe. Or les filtrations $F^{\bullet}$ et ${\hat F}^{\bullet}$ sont par hypothèse opposées sur chaque composante du gradué, ainsi chaque composante du gradué par la filtration $\calw_{\bullet}$ est pur i.e. de façon équivalente pour un fibré vectoriel sur ${\bf P}^1$, semistable de pente l'indice de la composante de la graduation.

Soit $H=(H_{\bf R},W_{\bullet}, F^{\bullet},{\overline F}^{\bullet})$ une structure de Hodge mixte réelle\\
(resp. $H=(H_{\bf C},W_{\bullet}, F^{\bullet},{\hat F}^{\bullet})$ une structure de Hodge mixte complexe). 
Les entiers $s^{p,q}=\text{dim}_{\bf C}Gr_{F}^{p}Gr_{\overline F}^{q}$ (resp. $s^{p,q}=\text{dim}_{\bf C}Gr_{F}^{p}Gr_{\hat F}^{q}$) permettent d'expliciter la forme scindée (qui existe toujours d'après le résultat de Grothendieck) de la structure twistorielle mixte associée par la proposition \ref{equivcar} p.\pageref{equivcar} à la structure de Hodge mixte réelle (resp. complexe). En effet d'après la décomposition en somme directe explicitée pour l'équivalence entre espce vectoriel muni de filtrations opposées d'une part et fibrés de Rees sur ${\bf P}^1$ semistable, il vient, avec $H_{\bf C}=H_{\bf R}\otimes_{\bf R} {\bf C}$,
\begin{eqnarray*}
\xi(H_{\bf C},F^{\bullet},{\overline F}^{\bullet})&=&\oplus_{(p,q)\in \cale_{H}}\xi(V^{p,q},DecTriv^{p},DecTriv^{q})\\
{}&\cong&\oplus_{(p,q)\in \cale_{H}}\calo(p+q)^{\text{dim}_{\bf C}Gr_{F}^{p}Gr_{\overline F}^{q}},
\end{eqnarray*} 
et ainsi
$$\xi(H_{\bf C},F^{\bullet},{\overline F}^{\bullet })\cong\oplus_{(p,q) \in \cale_{H}}\calo(p+q)^{s^{p,q}},$$ 
(resp. $\xi(H_{\bf C},F^{\bullet},{\hat F}^{\bullet })\cong\oplus_{(p,q) \in \cale_{H}}\calo(p+q)^{s^{p,q}}$).
Le fibré de Rees est filtré par des sous-fibrés stricts associés à la filtration par le poids. La construction étant fonctorielle, pour tout $n \in {\bf Z}$, on a la suite exacte
$$0 \rightarrow \calw_{n-1}\xi(H_{\bf C},F^{\bullet},{\overline F}^{\bullet }) \rightarrow \calw_{n}\xi(H_{\bf C},F^{\bullet},{\overline F}^{\bullet }) \rightarrow Gr_{\calw}^{n}\xi(H_{\bf C},F^{\bullet},{\overline F}^{\bullet }) \rightarrow 0.$$
Or $Gr_{\calw}^{n}\xi(H_{\bf C},F^{\bullet},{\overline F}^{\bullet })\cong \xi(Gr_{W}^{n}H_{\bf C},F^{\bullet},{\overline F}^{\bullet })$. Les filtrations de Hodge et conjuguée sont opposées sur le gradué associé à la filtration par le poids, donc 
\begin{eqnarray*}
Gr_{\calw}^{p}\xi(H_{\bf C},F^{\bullet},{\overline F}^{\bullet })&\cong &\oplus_{(p,q) \in \cale_{H},p+q=n}\calo(p+q)^{\text{dim}_{\bf C}Gr_{F}^{p}Gr_{W}^{n}}\\
{}& \cong &\oplus_{(p,q) \in \cale_{H},p+q=n}\calo(p+q)^{h^{p,q}}
\end{eqnarray*}

Donc $\xi(H_{\bf C},F^{\bullet},{\overline F}^{\bullet })$ est construit par extensions successives à partir de fibrés de la forme $\oplus_{(p,q) \in \cale_{H},p+q=n}\calo(p+q)^{h^{p,q}}$ 
En calculant d'une part le rang et d'autre part la première classe de Chern de la structure twistorielle mixte associée à une strcuture de Hodge mixte par l'expression explicite en fonction des entiers $s^{p,q}$ ou des entiers de Hodge $h^{p,q}$, on retrouve les relations (qui se déduisent aussi de relations sur les dimensions) $\sum_{(p,q) \in \cale_{H}}(h^{p,q}-s^{p,q})=0$ et $\sum_{(p,q) \in \cale_{H}}(p+q)(h^{p,q}-s^{p,q})=0$, ce qui justifie l'emploi d'une expression quadratique pour définir l'invariant $\alpha$, 
$$\alpha(H)=\sum_{(p,q) \in \cale_{H}}(p+q)^{2}(h^{p,q}- s^{p,q}).$$  

Le fibré $\xi(H_{\bf C},F^{\bullet},{\overline F}^{\bullet })$ associé à la structure de Hodge mixte réelle $H=(H_{\bf R},W_{\bullet}, F^{\bullet},{\overline F}^{\bullet}) \in \crmhs$ ne conserve que les données de la structure de Hodge mixte complexe associée par le morphisme oubli de la structure réelle $\crmhs \rightarrow \cmhs$. Dans \cite{sim2}, on introduit une involution antiholomorphe $\tau$ de ${\bf P}^1$ et une notion de faisceau $\tau$-équivariant. Cette construction dans le cas de ${\bf P}^1$ a été suivie pour la partie ``Structures réelles'' du dictionnaire entre espaces vectoriels trifiltrés et fibrés sur le plan projectif. Ici l'information de la conjugaison entre les filtrations qui permettent de former le fibré de Rees est codée par la $\tau$-équivariance du fibré de Rees associé à ces filtrations. Un fibré est dit ${\bf G}_{m}^{\tau}$-équivariant s'il est à la fois ${\bf G}_m$ et $\tau$-équivariant.

\begin{proposition}\cite{sim2}\label{equivcar} La catégorie des structures twistorielles mixtes ${\bf G}_{m}^{\tau}$-équivariantes est équivalente à la catégorie des structures de Hodge mixtes réelles.
\end{proposition}

\end{document}